\documentclass[11pt,PhD]{muthesis}

\usepackage{verbatim}
\usepackage{graphicx}
\usepackage{url} % typeset URL's reasonably
\usepackage[round]{natbib}
\usepackage{amsmath}
\usepackage{amssymb}
\usepackage{listings} % for printing code listings nicely
\usepackage{textcomp}
\usepackage{tikz}
\usepackage[T1]{fontenc}
\usepackage{xspace}
\usepackage{amssymb}
\usepackage{amsmath}
\usepackage{amsthm}
\usepackage{xcolor}
\usepackage{comment}
\usepackage{ wasysym} 
\usepackage{hyperref}
\usepackage{array}
\usepackage{listofitems}
\newcommand{\myStep}[2]{{\bf Step #1} --- #2\\}
\newenvironment{dedication}
  {\clearpage           % we want a new page
   \thispagestyle{empty}% no header and footer
   \vspace*{\stretch{1}}% some space at the top 
   \itshape             % the text is in italics
   \raggedleft          % flush to the right margin
  }
  {\par % end the paragraph
   \vspace{\stretch{3}} % space at bottom is three times that at the top
   \clearpage           % finish off the page
  }
\subject{Mathematics}
\campus{Palmerston North}
\begin{document}

%% Uncomment the following lines to leave out list of figures and tables until final printing
%\figurespagefalse
\tablespagefalse

%% Uncomment to make it singlespaced to save paper for draft printing
%\singlespace

%% Latex will put the current year on the title page by default. If you wish to change it, modify and uncomment this line
%\year=2010

\title{Globally resonant homoclinic tangencies}
\author{Sishu Shankar Muni\\
Under the supervision of Dr. David Simpson and Dist. Prof. Robert McLachlan
\vspace{3cm}
\newline
Thesis examined by:  Prof. Alan Champneys, University of Bristol, UK\\
           \hspace*{4cm}          Assoc. Prof. Bruce Van Brunt, Massey University, NZ\\
           \hspace*{4cm}       Assoc. Prof. Vivien Kirk, University of Auckland, NZ \\}
\beforeabstract
\prefacesection{Abstract}
%The abstract should be no more than 350 words long. It should present an overview of your thesis in accessible language.

%Remember that the purpose of an abstract is to help the potential reader decide whether or not they want to read your thesis. 
%You should therefore mention the important features: motivation, statement of the problem, your approach, the quality of your results, and what conclusions you drew.

The attractors of a dynamical system govern its typical long-term behaviour. The presence of many attractors is significant as it means the behaviour is heavily dependent on the initial conditions. To understand how large numbers of attractors can coexist, in this thesis we study the occurrence of infinitely many stable single-round periodic solutions associated with homoclinic connections in two-dimensional maps.  We show this phenomenon has a relatively high codimension requiring a homoclinic tangency and `global resonance', as has been described previously in the area-preserving setting.  However, unlike in that setting, local resonant terms also play an important role. To determine how the phenomenon may manifest in bifurcation diagrams, we also study perturbations of a globally resonant homoclinic tangency. We find there exist sequences of saddle-node and period-doubling bifurcations. Interestingly, in different directions of parameter space, the bifurcation values scale differently resulting in a complicated shape for the stability region for each periodic solution. In degenerate directions the bifurcation values scale substantially slower as illustrated in an abstract piecewise-smooth $C^{1}$ map.
\begin{dedication}
In dedication to my family
\end{dedication}
\prefacesection{Acknowledgements}
%\iffalse
To begin with I would esteem it as a privilege to thank my supervisors Dr.\,David John Warwick Simpson and Dr.\,(Dist. Prof.)\,Robert Ian McLachlan for their invaluable advice and continuous supervision throughout my PhD studies. Their doors were always open for clearing my queries. In fact, the more I stayed with them the less I met difficulties in accomplishing my doctoral research. I am really beholden to them. 

I would like to thank Dr.\,David for introducing me to the world of homoclinic tangencies. I admire him for his esteemed comments/ suggestions which helped me to learn many new things. I got so many words like ``great person'', ``great scientist'', and so on reserved for him to acknowledge his contribution for preparing this thesis. His comments/ suggestions list were always full for me from which I learnt many new things and helped me improve the quality of the thesis. He actually showed me how even roadblocks could be turned into stepping stones.

Further, I would like to thank Prof.\,Robert, for putting up with me for such a long time and giving me nice suggestions and holding meaningful discussions which really helped my studies go the extra mile. I would like to thank my supervisors for respecting and supporting my thoughts, choices, to express my results more clearly.

I also thank Massey University for making my journey through my Ph.D. a smooth sailing one by funding my research through Massey University doctoral scholarship. It was nice to be a Ph.D. scholar and gaining a lot of student information. I also want to thank SFS, NZMS, ANZIAM, SIAM, SMB for funding my research travel for various conferences I attended and got feedbacks which improved my presentation skills.

Now I would like to take it a pleasure to mention here all my family members, who are staying continents apart. First, thanks a bunch to my mother Dr.\,(Prof.)\,Anita Kumari Panda, R.C.M. Science College, Khallikote, Odisha, India. Apart from doing a mother’s job towards her children, she being herself a Ph.D. was guiding me all along my studies in India as well as here in New Zealand. It was purely a mother’s love that defied the definition of location on the globe. She has nurtured me, taken care of me. She always wanted me to be independent and to stand on my own foot. I am thankful that she always wanted to hear how my meeting with my supervisors went. Now a few lines of gratitude to my father, Shri Bijaya Kumar Muni. Maybe, by God’s grace I born as a man, but it is my father who, by profession a banker, turned me into a human being. Always cool and caring as he is, whatever I am today, are all because of him only. He just not financed me, inspired me and incited me, he even insulated all the predicaments that I might encounter during my studies. He taught me to be independent and always respect time. His motivation towards his work and family ignites a spark in me to give my best in whatever I do. Further, I would be doing a great injustice if I fail to add another name to the list of persons before whom I bestow all my gratitudes. She is Smt. Kumudini Muni, my grandmother. Of course, she herself is not a highly literate, but does possess the required capabilities to turn any other into a highly educated. So I feel proud to be one of her grandchildren. I would like to thank her for her humourous interactions, backing me up throughout and asking about my cooking to be sent through teleportation. My grandfather, the late Master Devaraj Muni had also his share of contribution towards whatsoever achievements I got today. Hats off to that departed soul. I admire his humour (especially humorous interactions with my grandmother), showing us the culture of \textit{Utkala}, character, morning walks, braveness and will always be remembered. Finally comes my life’s hero, i.e. my elder brother late Sitansu Muni. Being a brilliant par excellence at his studies which he completed in leading institutes in India, left no stone unturned to make me more than what he was. Such was his brilliance that the tasks appearing for others as skating on thin ice, were always a cake walk for him. His funda in every subject was a crystal clear. He was just not my brother, he was my ‘role model’. He monitored me for every letter I scribbled and guided me at every problem in mathematics I unriddled. His personality, character improved me a lot as a person. His brilliance and interest in mathematics and learning a lot has ignited spark in me to pursue science. His teachings on arithmetic progression, trigonometry, JAVA during my high school days motivated me a lot to pursue mathematics. I know you will support me in spirit however difficult the situation is. The warm memories of eating fritters, walks, fights, your spiritual devotion is etched in me. As he was the eldest of all our brothers and sisters, and the most brilliant also, we respected and loved him very much. But perhaps, he was loved more by God than us, he preferred God’s lap thereby leaving a huge gap in us. Today I reap his seeds and weep for missing his deeds. I would like to thank my maternal grandparents who have nourished, brought up and supported me all the way.

I would like to thank my colleagues I.A. Shepelev, A. Provata, V. Dos Santos, K. Rajagopal with whom it was interesting to discuss, learn  and collaborate on other topics of research on nonlinear dynamics. I hope to learn, participate with other researchers in future. I also would like to thank many researchers whom I communicated with and have given me many suggestions at conferences to improve my research work.  

I would like to thank the administrative staffs of the School of Fundamental Sciences for their support in my confirmation event and formal procedures during my Ph.D. research. 
I would like to thank Dr.\,Anil Kaushik and his family members who have welcomed and helped me in my endeavour. But for their help, my journey to New Zealand would have never been a smooth sailing. I am lucky to have many other friends in India as well as here in New Zealand who played their roles hush-hush in my entire endeavour. I am much thankful to my friends, office mates with whom I was able to interact freely and enjoyed many tramps, places, cooking and amazing badminton, squash matches.

SALUTE TO ALL.
%\fi
\afterpreface

% These include the actual text

\chapter{Introduction\label{cha:intro}}
\graphicspath{ {./images/} }
\newtheorem{theorem}{Theorem}[section]
\newtheorem{corollary}[theorem]{Corollary}
\newtheorem{lemma}[theorem]{Lemma}
\newtheorem{prop}[theorem]{Proposition}
\newcommand*{\blue}{\textcolor{black}}
\theoremstyle{definition}
\newtheorem{definition}{Definition}[section]
This thesis provides new results for \textit{homoclinic tangencies}, defined below, for smooth dynamical systems. These objects have an important role in the general theory of dynamical systems
as they are one of the simplest mechanisms for the onset of chaotic dynamics
and are the starting point for transverse homoclinic connections which were a central feature of Poincar{\'e}’s celebrated analysis of the three-body problem in his book \citep{Po92}.

In this chapter we first introduce maps --- discrete-time dynamical systems --- in which the results are formulated. We then review some fundamental concepts of nonlinear dynamical systems in the context of maps. \blue{Much of the fundamental theory on dynamical systems can also be found in standard textbooks, such as \cite{Kuz04,ArPl92,Meiss07}.} Next we introduce the concept of a homoclinic tangency and summarise results related to those obtained in later chapters. Finally \S \ref{sec:Overview} outlines the content of the remaining chapters.

\section{Maps}
\label{sec:maps}
Let $f : \mathbb{R}^n \mapsto \mathbb{R}^n$ be $C^{k}$ for some $k \geq 1$.
Since the domain and range of $f$ are the same we can refer it as a \textit{map} on $\mathbb{R}^n$.
Given $\rm \mathbf{x} \in \mathbb{R}^n$, we can iterate it to obtain the sequence $ \mathbf{x}, f(\mathbf{x}), f^2( \mathbf{x}), \ldots $ known as the \textit{forward orbit} of $\rm \mathbf{x}$ under $f$. Here and throughout the thesis we write $f^{k}$ to denote the composition of $f$ with itself $k$ times.

Maps are often used as mathematical models of physical phenomena. In this context $f^k(\rm \mathbf{x})$ represents the state of the system at the $k^{\rm{th}}$ time step, starting from the initial state $\mathbf{x}$. Maps also arise as Poincar{\'e} maps, or return maps, of systems of ordinary differential equations. For these reasons we study the dynamics of maps with the ultimate goal of better understanding the underlying physical phenomena from which they are derived.

A set $\Omega \in \mathbb{R}^n$ is said to be {\em invariant} under $f$ if $f(\Omega) \subseteq \Omega$.  The simplest invariant sets of a map are fixed points so we consider these first.
\begin{definition}
A point $\mathbf{x^{*}} \in \mathbb{R}^{n}$ is a \textit{fixed point} of $f$ if $f(\mathbf{x^{*}}) = \mathbf{x^{*}}$. 
\label{def:FP}
\end{definition}
\subsection{The stability of fixed points}
In this section we explain what it means for a fixed point $\mathbf{x^*}$ of a map $f$ to be stable and connect stability to the eigenvalues of ${\rm D} f(\mathbf{x^*})$.

\begin{definition}
A fixed point $\mathbf{x^{*}}$ is \textit{stable} or \textit{(Lyapunov stable)} if for every neighborhood $U$ of $\mathbf{x^{*}}$ there exists a neighborhood $V$ of $\mathbf{x^{*}}$ such that for all $\mathbf{x} \in V$, $f^{k}(\mathbf{x}) \in U$ for all $k \ge 0$, otherwise $\mathbf{x^{*}}$ is \textit{unstable}. 
\label{def:LSFP}
\end{definition}

\begin{definition}
A fixed point $\mathbf{x^{*}}$ is \textit{asymptotically stable} if $\mathbf{x^*}$ is stable and there exists a neighborhood $U$ of $\mathbf{x^*}$ such that for all $\mathbf{x} \in U$ we have $f^{k}(\mathbf{x}) \rightarrow \mathbf{x^{*}}$ as $k \rightarrow \infty$.
\label{def:ASFP}
\end{definition}

The linear approximation $g(\mathbf{x}) = {\rm D} f(\mathbf{x^*})(\mathbf{x} - \mathbf{x^*}) + {\rm \mathbf{x^{*}}}$ can provide a useful approximation to the dynamics of $f$ near $\mathbf{x^*}$, where ${\rm D} f$ is the $n \times n$ Jacobian matrix of first derivatives.  If ${\rm D} f(\mathbf{x^{*}}) v = \lambda v$, where $v \in \mathbb{R}^n$ and $\lambda \in \mathbb{R}$, then $g^k(\mathbf{x^*} + v) = \mathbf{x^*} + \lambda^k v$ for all $k \geq 0$.  So if $|\lambda| < 1$ this orbit converges to $\mathbf{x^*}$ while if $|\lambda| > 1$ this orbit diverges \blue{from $\mathbf{x^*}$}.  It follows that if all eigenvalues of $Df({\rm \mathbf{x^*}})$ have modulus less than $1$ then all orbits of $g$ converge to ${\rm \mathbf{x^*}}$.  The following theorem tells us this is true for $f$ for all initial points sufficiently close to $\rm \mathbf{x^*}$.
\begin{theorem}
Let $\rm \mathbf{x^*}$ be a fixed point of $f$. If all eigenvalues of $\textnormal{D}f(\mathbf{x^{*}})$ have modulus less than $1$, then $\rm \mathbf{x^*}$ is asymptotically stable.
If at least one eigenvalue of $\textnormal{D}f(\mathbf{x^{*}})$ \blue{has modulus greater than $1$,} then $\mathbf{x^{*}}$ is unstable.
\label{TheoFP}
\end{theorem}
Theorem \ref{TheoFP} can be viewed as a consequence of the following version of the Hartman-Grobman theorem \citep{Hart60}. \blue{Before we state this theorem we first remind the reader of three very standard definitions.}
\blue{
\begin{definition}
 A fixed point $x^*$ is said to be {\em hyperbolic} if no eigenvalues of ${\rm D}f(x^*)$ have unit modulus.
\label{def:hyperbolic}
\end{definition}
\begin{definition}
A function $h : \mathbb{R}^n \to \mathbb{R}^n$ is said to be a \textit{homeomorphism} if $h$ is one-to-one, onto, continuous, and has a continuous inverse.
\label{def:homeomorphism}
\end{definition}
\begin{definition}
A map $f$ is said to be \textit{conjugate} to another map $g$ on a set $S$ if there exists a homeomorphism $h$ such that $g = h^{-1} \circ f \circ h$ on $S$.
\label{def:conjugacy}
\end{definition}}
\begin{theorem}
{\rm [Hartman-Grobman theorem]} Let $f$ be $C^{1}$ and let $\mathbf{x^*}$ be a hyperbolic fixed point. Then there exists a neighborhood of $\mathbf{x^*}$ within which $f$ is conjugate to its linearisation. 
\label{Th:HartGrob}
\end{theorem}
In this thesis we also use the following terminology when the fixed point is hyperbolic. 
\begin{definition}
Suppose $\rm \mathbf{x^*}$ is a hyperbolic fixed point of $f$.  If all eigenvalues of $\textnormal{D}f(\rm \mathbf{x^*})$ have modulus less than $1$ then $\rm \mathbf{x^*}$ is {\em attracting}, if all eigenvalues of $\textnormal{D}f(\rm \mathbf{x^*})$ have modulus greater than $1$ then $\rm \mathbf{x^*}$ is {\em repelling}, otherwise $\rm \mathbf{x^*}$ is a {\em saddle}.
\label{def:AttrRepSaddle}
\end{definition}

\subsection{The stability of fixed points for two-dimensional maps}
In the case of two-dimensional maps (which will be the focus of later chapters), the eigenvalues of ${\rm D} f(\mathbf{x^*})$ are determined from its trace and determinant.  Algebraically it is easier to work with the trace and determinant than the eigenvalues because the formula for the eigenvalues in terms of the entries of ${\rm D} f(\mathbf{x^*})$ involves a square-root. Suppose $n=2$ and let $\tau = \rm{trace}(\textnormal{D}$$f$$(\mathbf{x^{*}}))$ and $\delta = \rm{det}(\textnormal{D}$$f$$(\mathbf{x^{*}}))$. The eigenvalues $\lambda$ are the solutions to the characteristic equation
$$
\lambda^2 - \tau \lambda + \delta = 0.
$$

The following result is a simple consequence of this equation and Theorem \ref{TheoFP}. Below we provide a proof for completeness.
\begin{theorem}
The eigenvalues of ${\rm D}f(\mathbf{x^*})$ both have modulus less than $1$ if and only if
\begin{equation}
    |\tau| - 1 < \delta < 1.
    \label{stab_eqn}
\end{equation}
\end{theorem}

\begin{proof}
Let $\lambda_{1},\lambda_{2}$ denote the eigenvalues of ${\rm D}f(\mathbf{x}^{*})$.  First suppose $\delta \leq \frac{\tau^{2}}{4}$ in which case the eigenvalues are real-valued. Then ${\rm max}(|\lambda_{1}|,|\lambda_{2}|) = \frac{|\tau| + \sqrt{\tau^2 - 4\delta}}{2}$. Observe
\begin{equation}
\begin{aligned}
& \   \frac{|\tau| + \sqrt{\tau^2 - 4\delta}}{2} < 1\\
& \iff \sqrt{\tau^2 - 4\delta} < 2 - |\tau|\\
& \iff \tau^2 - 4 \delta < (2-|\tau|)^2 \quad {\rm and}\quad 2 - |\tau| > 0\\
& \iff |\tau| - 1 < \delta \quad {\rm and} \quad \delta < 1
\end{aligned}
    \label{eq:RealCase}
\end{equation}
Second suppose $\delta > \frac{\tau^2}{4}$ in which case the eigenvalues are complex and we automatically have $|\tau| - 1 < \delta$. Here
\begin{equation}
    \label{eq:Imgcase}
    \begin{aligned}
    {\rm max}(|\lambda_{1}|,|\lambda_{2}|)  &= \left|\frac{\tau + i\sqrt{4\delta - \tau^2}}{2}\right|\\
        &=\sqrt{\left|\frac{\tau^2}{4} + \frac{4\delta - \tau^2}{4}\right|} = \delta \\
    \end{aligned}
\end{equation}
\end{proof}
The condition \eqref{stab_eqn} can be represented as three inequalities:
\begin{align*}
\delta &> -\tau - 1,\\
\delta &> \tau  -1,\\
\delta &< 1.    
\end{align*}

\begin{figure}[h!]
\begin{center}
\includegraphics{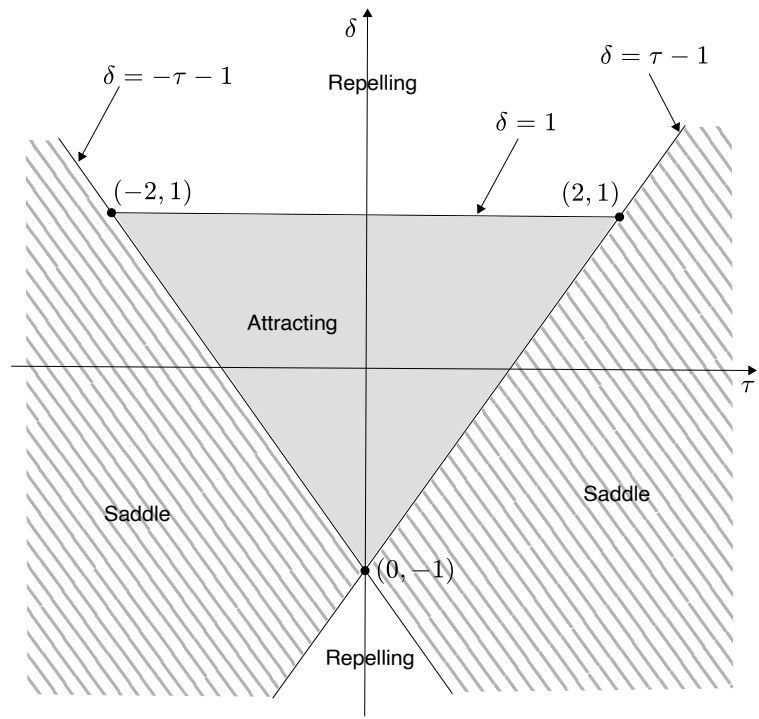}
\caption[Stability triangle in two-dimensional maps.]{ The stability of a fixed point of $f$ in terms of the trace $\tau$ and determinant $\delta$ of the Jacobian matrix $D f$ evaluated at \blue{the fixed point}.}
\label{fig:stability_triangle}
\end{center}
\end{figure}

If we visualize $\delta$ as a function of $\tau$ see Fig.~\ref{fig:stability_triangle}, the above three inequalities give us the stability region as the interior of the triangle bounded by the lines $\delta = -\tau - 1$, $\delta = \tau  -1$ and $\delta = 1$. Further these lines also divide the $(\tau,\delta)$-plane in regions where $\rm \mathbf{x^*}$ is attracting, repelling, and a saddle.

\subsection{Attractors and basins of attraction}
So far we have explained what it means for a fixed point (the simplest type of invariant set) to be attracting.  We now generalise this notion to arbitrary invariant sets of an $n$-dimensional map $f$.  This is the concept of a `topological attractor'.  An important alternate measure-theoretic viewpoint is that of a `Milnor attractor' \citep{Mi85}.
\begin{definition}
A compact set $\Omega \subset \mathbb{R}^{n}$ is said to be a \textit{trapping region} for $f$ if $f(\Omega) \subset {\rm int}(\Omega)$ (where {\rm int} denotes interior).
\label{def:Trapping}
\end{definition}
\begin{definition}
A set $\Lambda \subset \mathbb{R}^{n}$ is said to be an \textit{attracting set} of $f$ if there exists a trapping region $\Omega$ such that $\Lambda = \cap_{k \geq 1} f^{k}(\Omega)$.
\label{def:AttrSet}
\end{definition}
An `attractor' is essentially an attracting set that cannot be decomposed into multiple attracting sets.  There are many non-equivalent ways this is commonly done \citep{Meiss07} but this will not be an issue for the results of this thesis.  For definiteness we provide the following definition.
\begin{definition}
A set $\mathcal{A}$ is an \textit{attractor} if it is an attracting set and contains a dense orbit.
\label{def:Attractor}
\end{definition}
The presence of attractors is common in dynamical systems. The set of points whose forward orbits converge to an attractor is defined as the basin of attraction of that attractor. It is interesting to study which initial conditions converge to which attractor. To differentiate usually colour coded plots are generated. The set of initial conditions which converges to one attractor is marked with one colour and the set of initial conditions which converges to another attractor is marked with another colour and so on. Different basins of attraction can be highly intermingled due to nonlinearity in the system. Basins of attraction are usually bounded by the stable manifold of a saddle invariant set, and these boundaries can be smooth curves or fractals. 
\section{Local bifurcations}
Maps often contain parameters (constants).  Naturally we wish to understand how the dynamics of a map is different for different values of a parameter.  Roughly speaking as the value of a parameter is varied the dynamics only changes in a fundamental way at special values of the parameters --- these are termed {\em bifurcations}. \blue{Formally this can be defined in terms of the existence of a conjugacy (see Definition \ref{def:conjugacy}). Suppose $f_\mu$ is a family of maps.  A value $\mu_0$ is a bifurcation value of this family if for every neighbourhood $\mathcal{N}$ of $\mu_0$ there exist $\mu_1, \mu_2 \in \mathcal{N}$ for which 
$f_{\mu_1}$ and $f_{\mu_2}$ are not conjugate.}

Arguably the simplest types of bifurcations are those where a fixed point loses stability (or more generally where the number of stable eigenvalues \blue{of the Jacobian matrix of the map evaluated at the fixed point} changes).  These are examples of {\em local} bifurcations because they only concern the dynamics of the map in an arbitrarily small region of phase space. There are three different ways by which a fixed point can lose stability in a generic (codimension-one) fashion, that is when an eigenvalue attains the value $1$, the value $-1$, or the value ${\rm e}^{i \theta}$ for some $0 < \theta < \pi$.

These three scenarios correspond to saddle-node, period-doubling, and Neimark-Sacker bifurcations, respectively.  For each bifurcation to occur generically, the eigenvalue has to pass through its critical value linearly with respect to the parameter change (this is known as the {\em transversality condition}).  Also the map needs to have appropriate nonlinear terms (this leads to a {\em non-degeneracy condition}).  Formal statements of these are quite involved, see \cite{Kuz04}. Here we simply describe the effect that these bifurcations have.

\begin{figure}[t!]
\begin{center}
\includegraphics[width=0.5\textwidth]{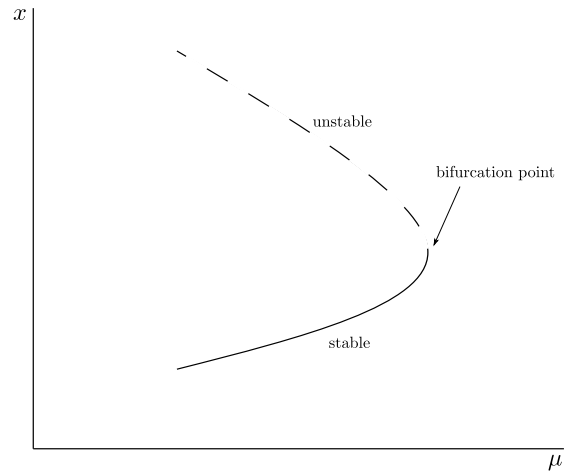}
\end{center}
\caption[Saddle-node bifurcation diagram.]{Saddle-node bifurcation in the case of $x \mapsto f(x)$ when a parameter $\mu$ is varied. A stable and unstable fixed point approach, collide and annihilate as parameter $\mu$ is varied.}
\label{fig:SNbifIntr}
\end{figure}

At a saddle-node bifurcation a fixed point has an eigenvalue $1$.  Under parameter variation two fixed points collide and annihilate at the bifurcation.  Fig.~\ref{fig:SNbifIntr} shows a typical bifurcation diagram.  The number of stable eigenvalues associated with the two fixed points differs by one.

\begin{figure}[ht!]
\begin{center}
\includegraphics[width=0.5\textwidth]{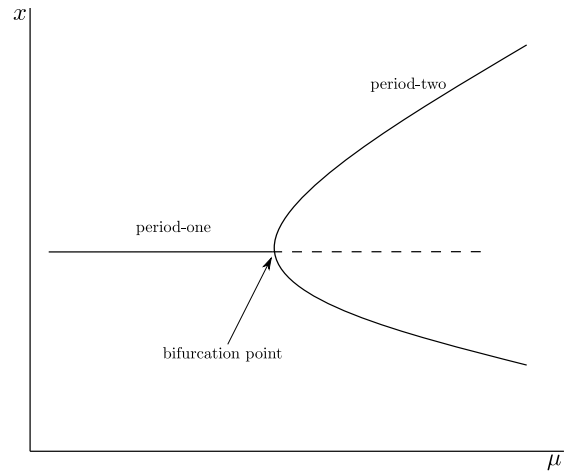}
\end{center}
\caption[Period-doubling bifurcation diagram.]{ Period-doubling bifurcation in the case of $x \mapsto f(x)$ when a parameter $\mu$ is varied. As $\mu$ is varied, the fixed point changes type and gives rise to a period-two solution.}
\label{fig:PDbifIntr}
\end{figure}

At a period-doubling bifurcation a fixed point has an eigenvalue $-1$.  Under parameter variation a period-2 solution is created at the bifurcation (periodic solutions are described in more detail in \S \ref{sec:persoln}).  Fig.~\ref{fig:PDbifIntr} shows a typical bifurcation diagram.  The stability of the period-2 solution matches that of the fixed point on the other side of the \blue{bifurcation, and so depending on which side of the bifurcation the period-$2$ solution emerges, the bifurcation can be classified as either supercritical or subcritical. To explain this, let $\lambda$ be the eigenvalue that has value $-1$ at the bifurcation and varies smoothly with respect to the parameter that we are varying. Then on one side of the bifurcation we have $\lambda < -1$ and the other side of the bifurcation we have $\lambda > -1$. If the period-$2$ solution coexists with the fixed point when $\lambda > -1$ then the period-doubling bifurcation is said to be \textit{subcritical}.  If the period-$2$ solution coexists with the fixed point when $\lambda < -1$ then the period-doubling bifurcation is said to be \textit{supercritical}}. 
%When the period-2 solution is more stable (as in Fig.~\ref{fig:PDbifIntr}) the period-doubling bifurcation is said to be {\em supercritical}, otherwise it is {\em subcritical}.

Finally, at a Neimark-Sacker bifurcation a fixed point has a complex conjugate pair of eigenvalues ${\rm e}^{\pm {\rm i} \theta}$ with $0 < \theta < \pi$. If $\theta \ne \frac{\pi}{2},\frac{\pi}{3},\frac{2\pi}{3},\frac{\pi}{4},\frac{3\pi}{4}$ then under parameter variation an invariant circle is created.  The disallowed fractional values of $\theta$ are `strong' resonances, see \cite{Kuz04}.

In the case of two-dimensional maps these bifurcations can be matched to lines in Fig.~\ref{fig:stability_triangle}. The straight line $\delta =1$ represents the Neimark-Sacker bifurcation. The straight line $\delta = \tau  - 1$ represents the saddle-node bifurcation. The straight line  $\delta = -\tau - 1$ represents the period-doubling bifurcation.
\section{Periodic solutions}
\label{sec:persoln}
In this section we discuss periodic solutions. These are orbits of periodic points where periodic points are fixed points of some iterate of the map.  Consequently the concepts described above for fixed points (particularly stability and bifurcations) extend easily to periodic points.
\begin{definition}
Let $\mathbf{x} \in \mathbb{R}^n$. If there exists a minimal number $k$ such that $f^{k}(\mathbf{x}) = \mathbf{x}$, then $\mathbf{x}$ is a periodic point of period $k$ of the map $f$.
\label{def:periodicpoint}
\end{definition}
\blue{For such $k$,} the collection of points
\begin{equation}
     \{ \mathbf{x}, f(\mathbf{x}), f^{2}(\mathbf{x}), \ldots, f^{k-1}(\mathbf{x})\},
    \label{eq:periodicsoln}
\end{equation} 
is termed a {\em periodic solution}. Observe the periodic point $\rm{\mathbf{x}}$ in Definition \ref{def:periodicpoint} is a fixed point of $f^k$. This observation allows us to reduce the problem of the stability of a periodic solution to that of a fixed point. The stability of the periodic solution \eqref{eq:periodicsoln} determined by the eigenvalues of $\textnormal{D}f^{k}(\mathbf{x})$ due to Theorem \ref{TheoFP}. Further, periodic solutions can be classified as attracting, saddles, and repelling, explained for fixed points in \S \ref{sec:maps} (note, Fig.~\ref{fig:stability_triangle} only applies to $2$d maps).

We observe that $\textnormal{D}f^{k}(\mathbf{x})$ can be written as a product of Jacobian matrices evaluated along each  point of the periodic solution. From chain rule of differentiation
\begin{equation}
    \textnormal{D}f^{k}(\mathbf{x}) = \textnormal{D}f(f^{k-1}(\mathbf{x})) \ldots \textnormal{D}f(f(\mathbf{x})) \textnormal{D}f(\mathbf{x}).
    \label{eq:Df_notation}
\end{equation}
Equation \eqref{eq:Df_notation} gives us a practical way to evaluate ${\rm D} f^k(\mathbf{x})$ from the derivative of $f$ and knowledge of the points of the periodic solution.  Further, from \eqref{eq:Df_notation} we can see that the eigenvalues are independent of the choice of the point of the periodic solution that is used to form \eqref{eq:Df_notation}.  Specifically if we replace $\mathbf{x}$ with $f^i(\mathbf{x})$ then ${\rm D} f^k \left( f^i(\mathbf{x}) \right)$ is a product of the same matrices in \eqref{eq:Df_notation} but multiplied together under a different cyclic permutation.  The cyclic permutation does not change the eigenvalues of the product, thus these eigenvalues are a well-defined feature of the periodic solution \eqref{eq:periodicsoln}.
\section{Non-invertible maps}
 \label{twodimnoninv}
 \begin{definition}
 A map $\textbf{x} \mapsto f(\textbf{x})$ is \textit{invertible} when its inverse $f^{-1}(\textbf{x})$ exists and is unique for each point in the domain. Otherwise it is said to be \textit{non-invertible}.
 \label{def:noninvmaps}
 \end{definition}
 We begin with a two-dimensional example. Consider the map $(x',y') = T(x,y)$ where
 $$T(x,y) =  \begin{bmatrix} ax+y\\ x^2 + b \end{bmatrix},$$
 where $a,b \in \mathbb{R}$ are the parameters of the map. With $a=0.5$ and $b=-2$, for example, the points $(-2,2)$ and $(2,0)$ both map under $T$ to the point $(1,2)$.  So certainly $T$ is not invertible on $\mathbb{R}^2$, at least for these values of the parameters. 
 More generally by solving for $x$ and $y$ in terms of $x'$ and $y'$ we obtain two solutions:\begin{align*}  x &= \sqrt{y' - b},\\ y&=x' - a\sqrt{y'-b}, \end{align*} and \begin{align*} x &= -\sqrt{y' - b},\\ y&=x' + a\sqrt{y' - b},\end{align*}
 assuming $y' > b$.  Thus any point $(x',y') \in \mathbb{R}^2$ with $y' > b$ has two {\em preimages} under $T$.  If instead $y' < b$ the point has no preimages.  The line $y' = b$, where the number of preimages changes, is an example of a {\em critical curve}, referred to \blue{as} $LC$ (\blue{for \textit{Ligne Critique} in French}) following \cite{CM96}. 
% Here we consider two dimensional non-invertible maps $T:\mathbb{R}^2 \rightarrow \mathbb{R}^2$, for $X = (x,y)$, $T(X): (x,y) \rightarrow [f,g]$, where $x,y$ are real variables, $f,g$ are single valued functions either continuously differentiable or piecewise differentiable. 
 %\textcolor{blue}{Formal definition of critical set $LC$, $LC_{-1}$}. 
 
\blue{More generally the critical curves of a map divide its phase space into distinct regions $R_{i}$ for $i=1,2,\ldots,n$ in each of which the map has a constant number of preimages, say $k_i$.  The map can then be classified as $Z_{k_{1}}$-$Z_{k_{2}}$-$Z_{k_{3}}$-$\cdots$-$Z_{k_{n}}$ as determined by the types of regions $Z_{k_{i}}$ that appear \cite{CM96}.} For example, \blue{the map $T$ defined above} is of type $Z_0$-$Z_2$ because points have zero preimages above the line $y'=b$ and two preimages below this line. The action of such a map is sometimes referred to as ``fold and pleat''.

  Maps may also exhibit another kind of complexity related to the presence of one or several cusp points on a critical curve $LC$. We encounter this phenomenon later in Chapter \ref{cha:ExampleInfinte}.  
\section{Stable and unstable manifolds}
\label{sec:stunmnfld}
Stable and unstable manifolds are perhaps the simplest global invariant structures of a map and often play an important role in influencing its global dynamics. 

Let $\mathbf{x^{*}}$ be a fixed point of a map $f$. The stable manifold of $\mathbf{x^{*}}$ is defined as the set
$$W^{s}(\mathbf{x^{*}}) = \{\mathbf{x} \,|\, f^{n}(\mathbf{x}) \rightarrow \mathbf{x^{*}} ~\text{as}~ n \rightarrow \infty \}\backslash \{\mathbf{x^{*}}\}.$$
Similarly the unstable manifold of $\mathbf{x^{*}}$ is 
$$W^{u}(\mathbf{x^{*}}) = \{\mathbf{x} \,|\, f^{-n}(\mathbf{x}) \rightarrow \mathbf{x^{*}} \text{as}~n \rightarrow \infty ~\textrm{for some sequence of preimages}\}\backslash \{\mathbf{x^{*}}\}.$$
The definition of the unstable manifold is more complicated than that of the stable manifold because we are not assuming $f$ is a homeomorphism. For $\mathbf{x} \in \mathbb{R}^n$ to belong to the unstable manifold of ${\rm \mathbf{x^*}}$ we require there to exist a sequence of preimages converging \blue{under backwards iteration of $f$} to ${\rm \mathbf{x^*}}$. Also notice $W^s(\mathbf{x}^*)$ and $W^u(\mathbf{x}^*)$ are forward invariant, that is if $\mathbf{x} \in W^s(\mathbf{x}^*)$ then $f(\mathbf{x}) \in W^s(\mathbf{x}^*)$, and similarly for $W^u(\mathbf{x}^*)$.
\subsection{Homoclinic orbits of maps}
Now suppose $\mathbf{x}^*$ is hyperbolic, that is it has no eigenvalues with unit modulus \blue{(see Definition \ref{def:hyperbolic})}.  Let $a$ be the number of eigenvalues of ${\rm D} f(\mathbf{x}^*)$ that have modulus less than $1$.  The stable manifold theorem (see \citep{Hart60}) tells us that $W^s(\mathbf{x}^*)$ is $a$-dimensional and emanates from $\mathbf{x}^*$ tangent to the stable subspace of \blue{the linearised system of $f$.  This subspace is given by} the span of the real parts of \blue{the eigenvectors and  generalised} eigenvectors corresponding to eigenvalues with modulus less than $1$.  Moreover, $W^u(\mathbf{x}^*)$ is $(n-a)$-dimensional and emanates from $\mathbf{x}^*$ tangent to the analogous unstable subspace.  

If $f$ is invertible then any unstable manifold cannot have any self-intersections. This is not true for non-invertible maps, as will be seen in Chapter \ref{cha:ExampleInfinte}.  Also, for non-invertible maps the numerical computation of stable manifolds is particularly challenging (as will be discussed in Chapter \ref{cha:ExampleInfinte}) because to `grow' the stable manifold outwards we need to iterate backwards under $f$. 

If $W^s(\mathbf{x}^*)$ and $W^u(\mathbf{x}^*)$ are both non-empty then $\mathbf{x}^*$ must be a saddle.
\begin{figure}[t!]
\begin{center}
\includegraphics[width=0.8\textwidth]{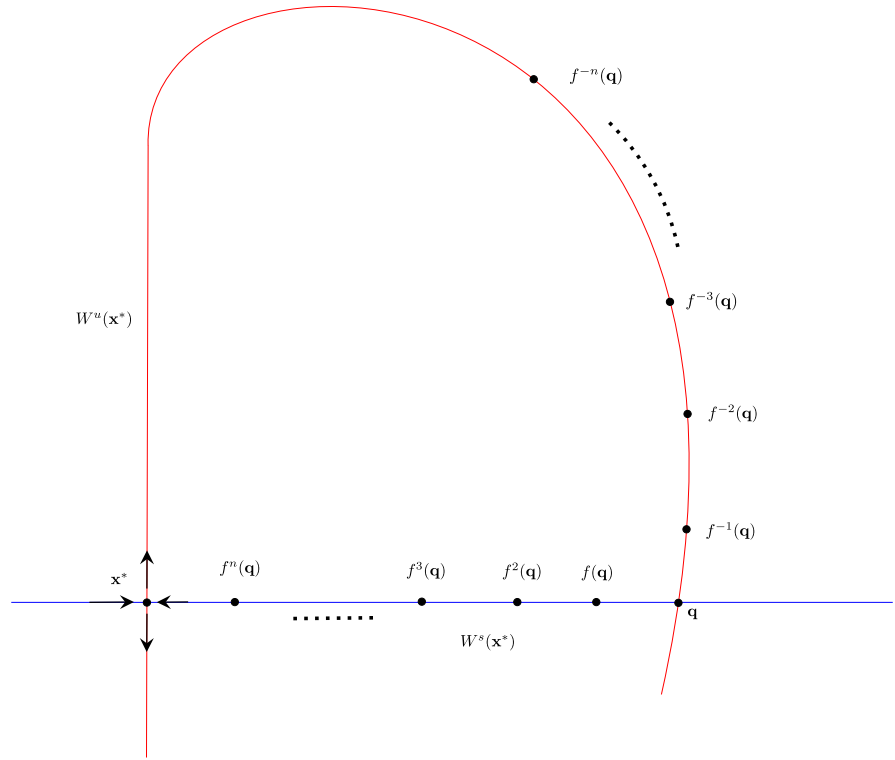}
\caption[Homoclinic points and homoclinic orbit of a saddle fixed point.]{Homoclinic point $\mathbf{q}$ and homoclinic orbit associated with a saddle fixed point $\mathbf{x^{*}}$. The stable manifold (blue) and the unstable manifold (red) of the saddle fixed point $\mathbf{x^{*}}$ intersect transversally at $\mathbf{q}$.}
\label{homoclinic_theo}
\end{center}
\end{figure}
\newtheorem{defn2}[theorem]{Definition}
\begin{defn2}
Let $\mathbf{x^{*}}$ be a saddle fixed point of $f$. If $$\mathbf{q}\in W^{s}(\mathbf{x^{*}}) \cap W^{u}(\mathbf{x^{*}}),$$ then $\mathbf{q}$ is called a {\em homoclinic point} and the orbit of $ \mathbf{q}$ is called a {\em homoclinic orbit}.
\end{defn2}
\newtheorem{defn3}[theorem]{Definition}

Figure \ref{homoclinic_theo} sketches a two-dimensional example.  Shown in blue and red, respectively, are parts of the stable and unstable manifolds of a fixed point $\mathbf{x^{*}}$ as they emanate from $\mathbf{x^{*}}$.  They have been followed outwards far enough that one intersection point $\mathbf{q}$ is visible.  By definition this is a homoclinic point and its orbit is homoclinic to $\mathbf{x^{*}}$.  Moreover every point of the homoclinic orbit is itself a homoclinic point.  Thus if we were to grow the stable and unstable manifolds outwards indefinitely we would see that they pass through every point of the homoclinic orbit.  The manifolds consequently have an extremely complicated geometry, known as a homoclinic tangle.  This structure was first described by Poincar{\'e}. Quoting from his book \citep{Po92} on page $389$:

\textit{Que l'on cherche {\`a} se repr{\'e}senter la figure form{\'e}e par ces deux courbes et leurs intersections en nombre infini dont chacune correspond {\`a} une solution doublement asymptotique, ces intersections forment une sorte de treillis, de tissu, de r{\'e}seau {\`a} mailles infiniment serr{\'e}es; chacune des deux courbes ne doit jamais se recouper elle-m$\hat{e}$me, mais elle doit se replier sur elle-m$\hat{e}$me d'une mani{\`e}re tr{\`e}s complexe pour venir recouper une intinit{\'e} de fois toutes les mailles du r{\'e}seau.
On sera frapp{\'e} de la complexit{\'e} de cette figure, que je ne cherche m$\hat{e}$me pas {\`a} tracer.}

\noindent Here is a translation given by \citep{And94}, \citep{Tour12}:

\textit{When one tries to imagine the figure formed by these two curves and
their infinitely many intersections each corresponding to a doubly asymptotic
solution, these intersections form a kind of lattice, web or network with
infinitely tight loops; neither of the two curves must ever intersect itself, but
it must bend in such a complex fashion that it intersects all the loops of the
network infinitely many times. One is struck by the complexity of this figure which I am not even attempting to draw.}

\begin{figure}[t!]
\begin{center}
\includegraphics[width=0.7\textwidth]{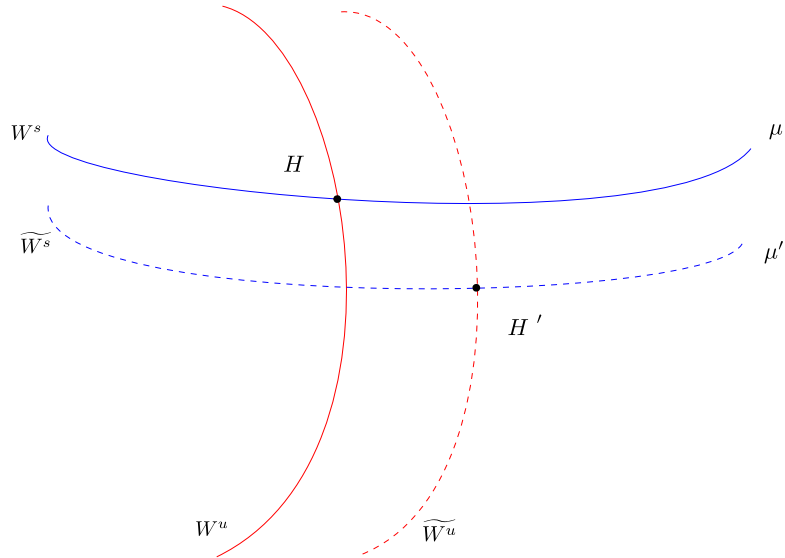}
\end{center}
\caption[Structurally stable transversal intersection.]{
 Structurally stable homoclinic orbits $($ due to the transversal intersection of $W^{s}$ and $W^{u}).$ The stable manifold $W^{s}$ and the unstable manifold $W^{u}$ at parameter $\mu$ intersect transversally at $H$. The transversal intersection is structurally stable as after the parameter $\mu$ is changed to $\mu'$, the stable manifold \blue{$\widetilde{W^{s}}$} intersects transversally with the unstable manifold \blue{$\widetilde{W^{u}}$} at $H'$.}
\label{fig:myFigure1}
\end{figure}
Generically, the intersection of the stable and the unstable manifold is transverse in which case the intersection persists as the map is perturbed (say by adjusting the value of a parameter of the map), see Fig.~\ref{fig:myFigure1}.  We say the intersection, and therefore also the corresponding homoclinic orbit, is {\em structurally stable}. As an example we identify a  homoclinic orbit in the case of a \blue{planar} map given in \cite{PaTa93}. Define a linear map $\phi : \mathbb{R}^2 \rightarrow \mathbb{R}^2$ by 
\begin{equation}
\phi(x,y) = \begin{bmatrix} 2 & 0\\ 0 & \frac{1}{2} \end{bmatrix} \begin{bmatrix} x\\ y \end{bmatrix}.
\label{eq:phi}
\end{equation}
The origin is a saddle-type fixed point. The $x$-axis is the unstable manifold of the origin and the $y$-axis is the stable manifold of the origin. 

\begin{figure}[t!]
\begin{center}
\includegraphics[width=0.8\textwidth]{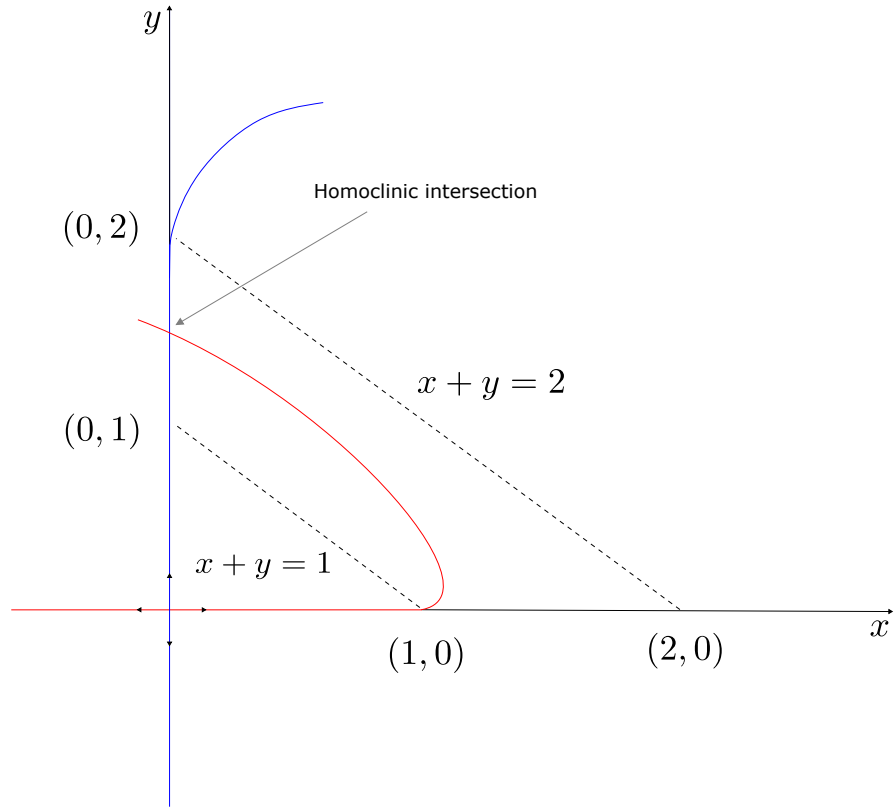}
\caption[Existence of a homoclinic orbit in shear map.]{A homoclinic intersection of the perturbed linear map $\psi \circ \phi$. The blue and red curves denote the stable and unstable manifold of the origin.}
\label{fig:myFigure2}
\end{center}
\end{figure}

We next consider the composition $\psi \circ \phi$, where $\psi(x,y) = (x- h(x+y),y+h(x+y))$ and $h$ is some continuous function with $h(z)=0$ for all $z \leq 1$ and $h(2) > 2$. To understand the effect of $\psi$, write $x' = x- h(x+y)$, $y'=y+h(x+y)$, then $x'+y' = x+y$, that is if points $(x,y)$ lie on the line $l_{c} (x+y=c)$, then the points $(x',y')$ also lie on the line $l_{c}$ but the coordinates are such that they get shifted from $(x,y)$ by the vector $(-h(c),h(c))$. Under the composition $\psi \circ \phi$, the point $(\frac{1}{2},0)$ gets mapped to the point $(1,0)$ and the point $(1,0)$ gets mapped to $(2-h(2),h(2))$ which lies to the left of the $y$-axis. Thus the line segment from $(\frac{1}{2},0)$ to $(1,0)$, which is part of the unstable manifold of $(0,0)$ for the map $\psi \circ \phi$, maps to the curve connecting $(1,0)$ to a point left of the $y$-axis.  It follows that this curve must intersect the line segment from $(0,1)$ to $(0,2)$, which is part of the stable manifold of $(0,0)$, see Fig.~\ref{fig:myFigure2}. Thus $\psi \circ \phi$ has a homoclinic orbit.

%See \eqref{eq:phi}.

\begin{figure}[t!]
\begin{center}
\includegraphics[width=0.8\textwidth]{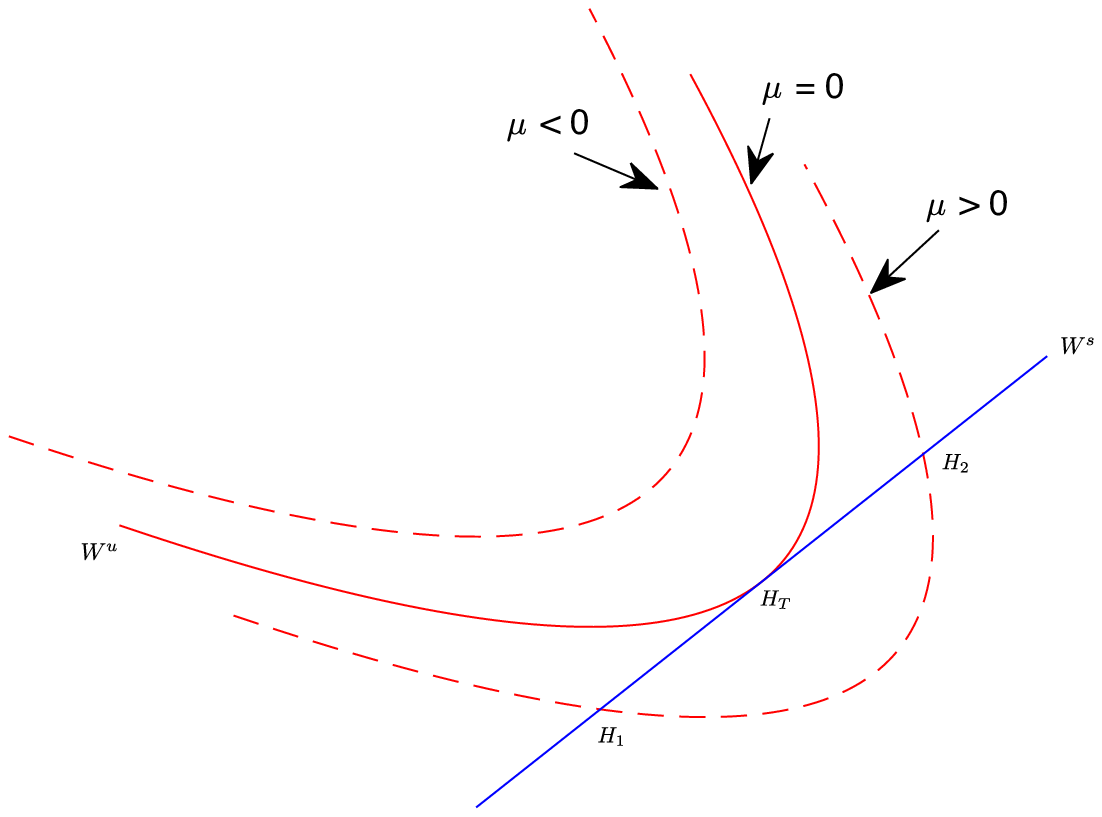}
\caption[Structural instability of homoclinic tangencies.]{Structural instability of homoclinic tangencies. The unstable manifold (solid red curve) has a tangential intersection with the stable manifold (blue curve) when $\mu = 0$ at the point $H_{T}$. The tangency is destroyed with the variation of the parameter $\mu$. For $\mu<0$, there is no tangency. For $\mu>0$, there is a transversal intersection of $W_{u}$ and $W_{s}$ at $H_{1}$ and $H_{2}$. 
\label{fig:myFigure3}
}
\end{center}
\end{figure}

\section{Homoclinic tangencies}
When the stable and unstable manifolds of a fixed point first intersect in response to parameter variation, the intersection is usually tangential.  Such an intersection is called a {\em homoclinic tangency}. This intersection is not structurally stable. As shown in Fig.~\ref{fig:myFigure3} for a two-dimensional map, the homoclinic points $H_1$ and $H_2$ collide and annihilate when $\mu=0$. Consequently homoclinic tangencies are usually codimension-one phenomena, that is they occur on a $(p-1)$-dimensional subset of $p$-dimensional parameter space.  However this is complicated by the fact that if a map has one homoclinic tangency it typically has other homoclinic tangencies on a dense subset of parameter space \citep{PaTa93}.  Intuitively this occurs because transverse homoclinic intersections imply a homoclinic tangle (mentioned in \S \ref{sec:stunmnfld}) so the manifolds wind around each other infinitely many times producing infinitely many chances, so to speak, to intersect tangentially. Also \cite{Ne74} showed that \blue{planar maps have} infinitely many attractors for a dense set of parameter values near a generic homoclinic tangency.

\begin{figure}[htbp!]
\begin{center}
\setlength{\unitlength}{1cm}
\hspace*{-7cm}
\begin{picture}(5.1,22)
\put(0,0){\includegraphics[width=10cm]{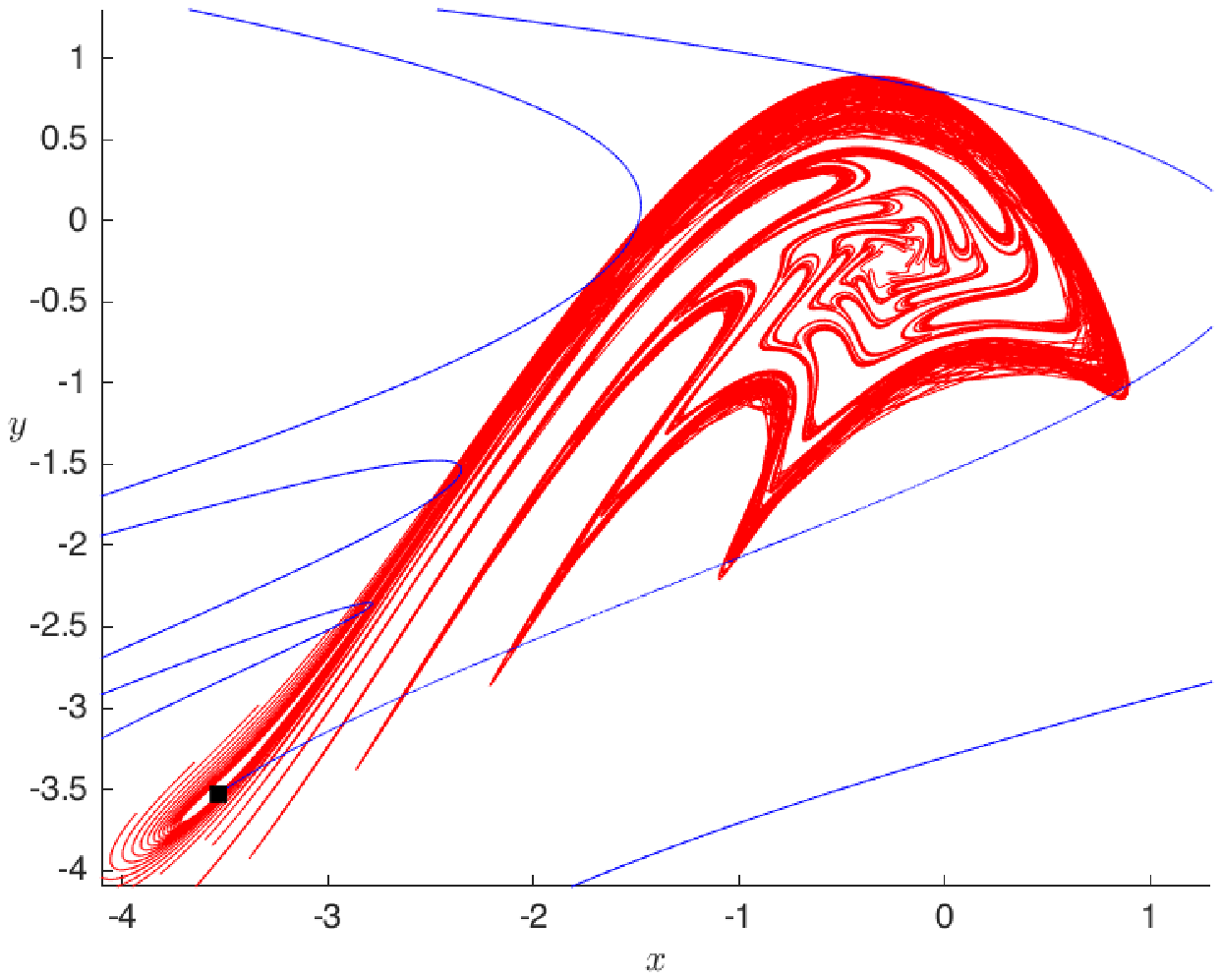}}
\put(0,7.2){\includegraphics[width=10cm]{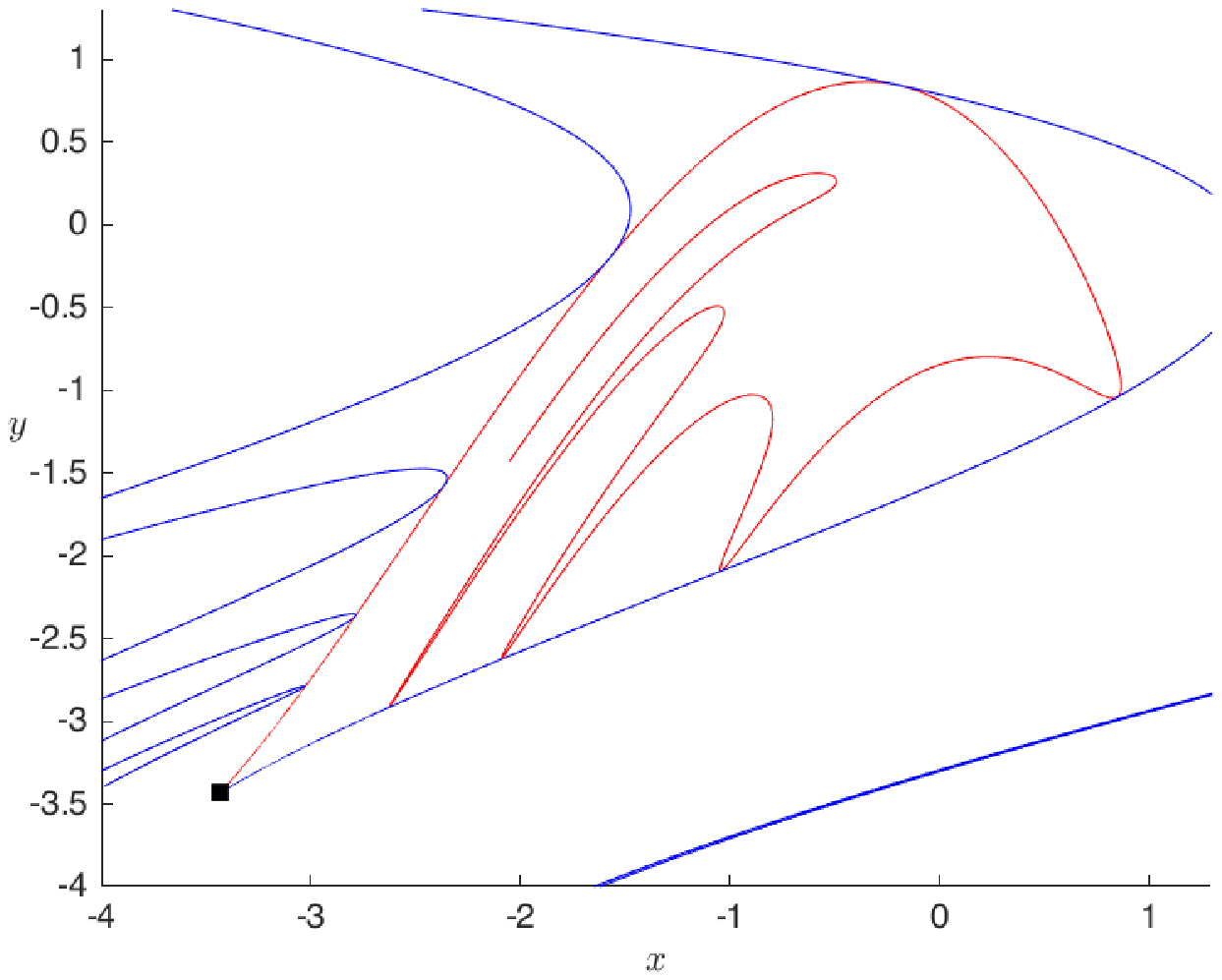}}
\put(0,14.4){\includegraphics[width=10cm]{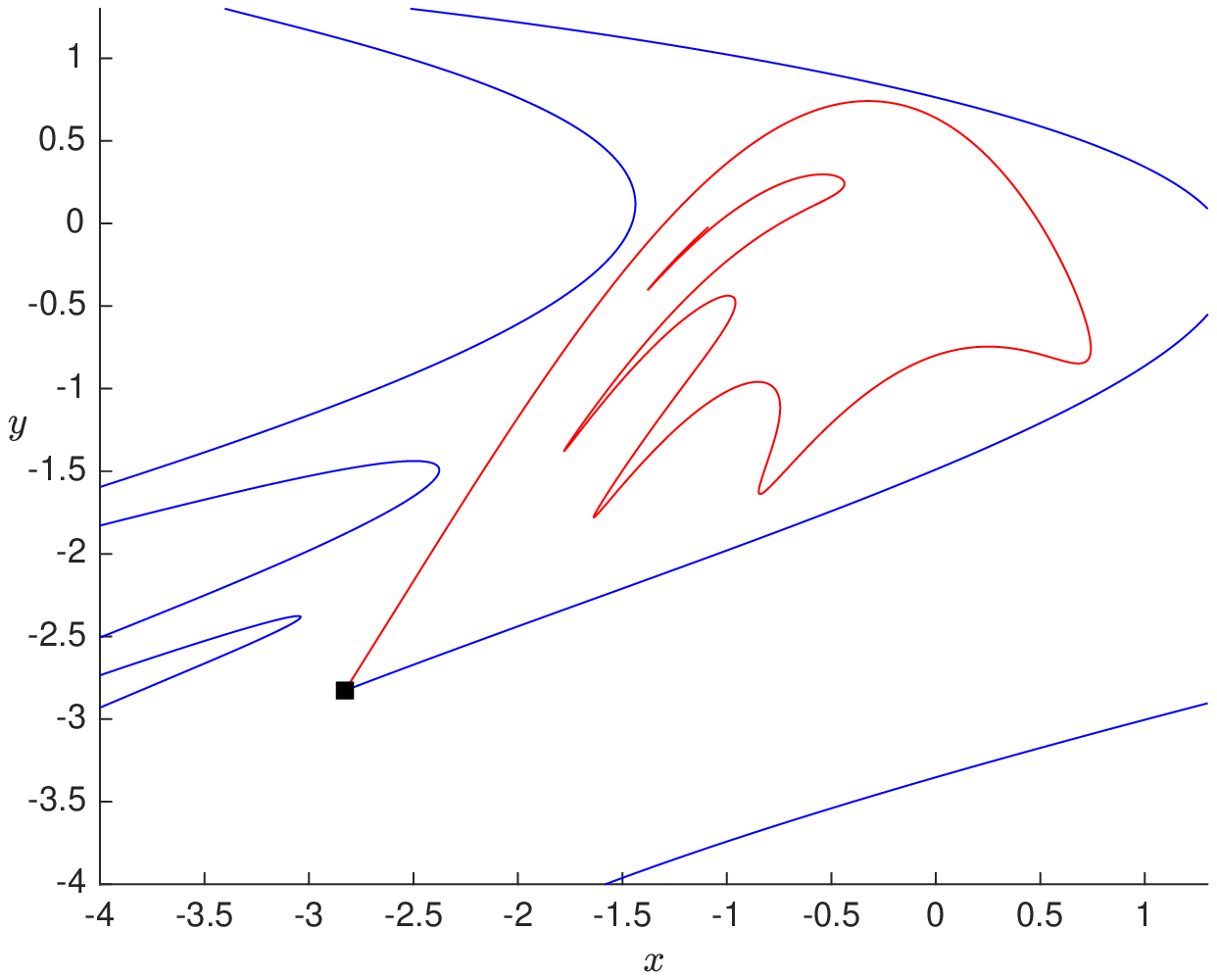}}
\put(0,7){\small \bf c)}
\put(0,14.2){\small \bf b)}
\put(0,21.4){\small \bf a)}
\end{picture}
\caption[Homoclinic tangle in the generalised H{\'e}non map. ]{Homoclinic tangle in the \blue{GHM}. The parameters are set as $S = -0.125, \beta = 0.8, \alpha = -0.4$. In the left, $R = 0.06$, we see no intersections between the stable manifold (blue) and unstable manifold (red). In the middle, $R = 0.08048175$ and we see a homoclinic tangency. In the right, $R = 0.08094$, the stable and unstable manifolds develop transversal intersections.}
	\label{fig:GHM_Tangle}
\end{center}
\end{figure}

Here we illustrate \blue{homoclinic tangles with} the generalised H\'enon map \blue{(GHM)}. \blue{The GHM generalises the well known H{\'e}non map \blue{and is given by}
\begin{align*}
x' &= y,\\
y' &= \alpha - \beta x - y^2 + Rxy + Sy^3,   
\end{align*}
where $S, R, \alpha, \beta$ are parameters of the map. The \blue{GHM} was derived from a general analysis of  homoclinic tangencies by composing a local map with a global map, similar to what is done here in later chapters, and taking a certain limit to justify the removal of higher-order terms \cite{GoKu05}. } We fix parameters $S = -0.125, \beta = 0.8, \alpha = 0.4$ and vary $R$. For $R = 0.06$, we observe no intersection between the stable and unstable manifold of a saddle fixed point (small black square), see Fig.~\ref{fig:GHM_Tangle}(a). For approximately $R = 0.08048175$, the stable and unstable manifolds intersect tangentially, see Fig.~\ref{fig:GHM_Tangle}(b).  For larger values of $R$ they have transversal intersections, see Fig.~\ref{fig:GHM_Tangle}(c). Now the manifolds form a homoclinic tangle.  Its complexity is striking and we can certainly symphasize with Poincar\'e for not attempting to draw it! \blue{The GHM is analysed more in \S \ref{sec:GHManalyse}.}

\section{Single-round and multi-round periodic solutions}
\begin{figure}[htbp!]
\begin{center}
\setlength{\unitlength}{1cm}
\hspace*{-4cm}
\begin{picture}(10.1,4.1)
\put(0,0){\includegraphics[width=0.45\textwidth]{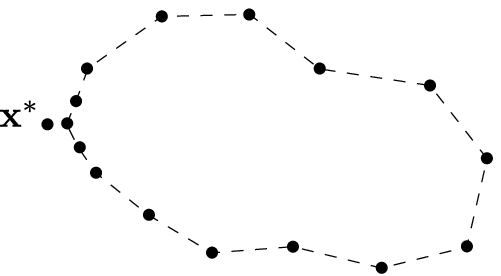}}
\put(7.1,0){\includegraphics[width=0.45\textwidth]{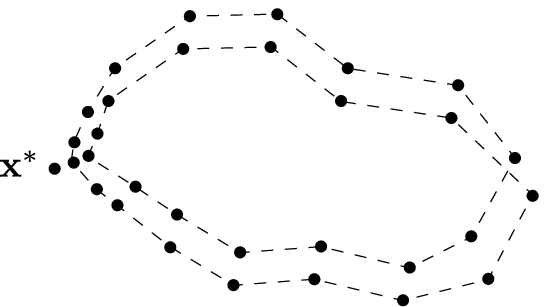}}
\put(0,0){\small \bf a)}
\put(7.1,0){\small \bf b)}

\end{picture}
\caption[A sketch of a single-round and multi-round periodic solution.]{ A single-round periodic solution of period $15$ is shown on the left. By single-round we mean it comes close to the fixed point ${\bf x}^*$ once. A double-round periodic solution of period $30$ is shown on the right.}
	\label{fig:SingMulRound}
\end{center}
\end{figure}
Here we introduce the notion of single-round and multi-round periodic solutions.  A precise definition is provided in \S \ref{sec:sec2p2}.  Essentially a single-round solution is a periodic solution which \blue{comes near the fixed point} once before returning to its starting point. A single-round periodic solution of period $15$ is shown in Fig.~\ref{fig:SingMulRound}(a). In this thesis we will focus on the single-round periodic solutions. 

A \textit{multi-round} periodic solution of $f$ is a periodic solution which \blue{comes near the fixed point} several times before returning to its starting point. A double-round periodic solution of period $30$ is shown in Fig.~\ref{fig:SingMulRound}(b).

\iffalse
\begin{figure}[!htbp]
 %\hspace{0.8cm}
	\begin{tabular}{c c}
		%\includegraphics[width=0.55\textwidth]{case1.eps}&
		\includegraphics[width=0.45\textwidth]{SingleR15.eps}&
		\includegraphics[width=0.45\textwidth]{DoubleR30.eps}&
	\end{tabular}
	\caption[A sketch of a Single-round and multi-round periodic solution.]{ A single-round periodic solution of period $15$ which shadows neighborhood of $\mathbf{x^{*}}$ once before returning to its starting point in shown on left. A double-round periodic solution of period $30$ which shadows neighborhood of $\mathbf{x^{*}}$ twice before returning to its starting point is shown on right.}
	\label{fig:SingMulRound}
\end{figure}
\fi
\iffalse
\begin{figure}[h!]
\begin{center}
\includegraphics[width=0.5\textwidth]{SingRoundSketch.eps}
\caption[A sketch of a Single-round periodic solution.]{ A single-round periodic solution of period $8$ which shadows neighborhood of $\mathbf{x^{*}}$ once before returning to its starting point.}
\label{fig:SingleRoundSketch}
\end{center}
\end{figure}
\begin{figure}[h!]
\begin{center}
\includegraphics[width=0.5\textwidth]{MulRoundSketch.eps}
\caption[Double-round periodic solution.]{A double-round periodic solution of period $16$ which shadows neighborhood of $\mathbf{x^{*}}$ twice before returning to its starting point.}
\label{fig:MultiRoundSketch}
\end{center}
\end{figure}
\fi
\section{Degenerate homoclinic tangencies}
\subsection{Globally resonant homoclinic tangencies}
At a generic homoclinic tangency the number of single-round periodic solutions that are stable must be finite (in typical examples this number is rarely more than one) \citep{GaSi72}. One of the main aims of this thesis is to explain what degeneracies permit this number to instead be infinite.
Consider a homoclinic tangency associated with a fixed point of a two-dimensional map.
Then the eigenvalues associated with the fixed point, call them $\lambda$ and $\sigma$, satisfy \blue{$0 < |\lambda| < 1 < |\sigma|$}. We will show (see Chapter \ref{cha:GRHT}) that one of the required conditions (degeneracies) is that $|\lambda \sigma| = 1$.

Note that this is a local condition. In contrast another required condition is global and more difficult to describe. It can be interpreted as saying that the reinjection mechanism converts a displacement in the stable direction to an equal displacement in the unstable direction (see Chapter \ref{cha:GRHT}). Following \cite{GoGo09} we refer to this condition as {\em global resonance}.

When combined the conditions described above form a codimension-three phenomenon. This is
because (i) a homoclinic tangency, (ii) $|\lambda \sigma| = 1$, and (iii) global resonance
are all codimension-one and independent.
As will be seen, in the case $\lambda \sigma = 1$ an extra codimension-one condition (on resonant terms in the `local map') is also required to have
infinitely many stable single-round periodic solutions.
In Chapter \ref{cha:GRHT} we further identify sufficient conditions for the periodic solutions to be asymptotically stable.

A second aim of this thesis is to unfold a globally resonant homoclinic tangency.
In Chapter \ref{cha:UnfGRHT} we introduce several parameters, but focus on one-parameter families $f_\varepsilon$ that exhibit a globally resonant homoclinic tangency when $\varepsilon = 0$. For each $K \in \mathbb{Z}$ there exists an open interval of values of $\varepsilon$, containing $0$, for which $f_{\varepsilon}$ has $K$ asymptotically stable single-round periodic solutions. In particular in Chapter 4 we show that typically the width of this interval is asymptotically proportional to $|\lambda|^{2K}$ as $K \to \infty$.

\subsection{The area-preserving case}
The results in this thesis apply to families of $C^\infty$ maps.
If we instead restrict our attention to area-preserving maps the phenomenon is codimension-two because $|\lambda \sigma| = 1$ is automatic and the extra condition on the resonance terms of the local map turns out to be automatic similarly. This scenario was considered in \cite{GoSh05} where the periodic solutions are elliptic instead of being asymptotically stable. The phenomenon was unfolded in \cite{GoGo09}.
Their results show that $K$ elliptic single-round periodic solutions coexist
for parameter values in an interval of width asymptotically proportional to $|\lambda|^{2 K}$, matching our result. A more recent summary of these results can be found in \cite{DeGo15}.

\subsection{The piecewise-linear case}
\label{subsec:PWL}
Continuous piecewise-linear maps exhibit homoclinic tangencies in a codimension-one fashion except at any point of intersection one manifold forms a kink (instead of a quadratic tangency) while the manifold is locally linear. In a generic unfolding only finitely many single-round periodic solutions can be stable \citep{Si16b}. In this setting the analogy of a globally resonant homoclinic tangency is codimension-three and was first described in general in \cite{Si14}. Here the parameter width for $K$ is asymptotically proportional to $|\lambda|^K$ \citep{Si14b},
which, interestingly, differs from our result. In this setting the global resonance condition implies that branches of the stable and unstable manifolds not only intersect but are in fact coincident (i.e.~intersect at all points).

\subsection{Other degenerate homoclinic tangencies}

Various other degenerate types of homoclinic tangencies have been considered and unfolded.
For example branches of the stable and unstable manifolds of a fixed point of a smooth two-dimensional map can be coincident in a codimension-two phenomenon. In \cite{HiLa95} it was shown that in a generic two-parameter unfolding, each single-round periodic solution exists between two lines of saddle-node bifurcations, and period-doubling bifurcations can occur nearby. Consequently in one-parameter bifurcation diagrams one typically sees a sequence of `isolas' bounded by the saddle-node bifurcations. This was extended further in \cite{ARC99} where the authors unfolded a codimension-three Shilnikov-Hopf bifurcation at which branches intersect transversely. Their results explain how one-parameter bifurcation diagrams change from isolas to a single wiggly curve
by asymptotic calculations of surfaces of Hopf bifurcations and homoclinic tangencies in three-dimensional parameter space. Also bifurcations associated with cubic and other higher-order tangencies are described in \cite{Da91,GoKu05}.

\section{Thesis overview}
\label{sec:Overview}
After some background of homoclinic theory of maps in Chapter \ref{cha:intro}, we proceed to state the backbone of the linearisation of maps --- the Sternberg's linearisation theorem in Chapter \ref{ch:StLin}. Chapters \ref{cha:GRHT}--\ref{cha:OpenProblems} contain the original research of this thesis. Chapter \ref{cha:GRHT} begins by setting up the mechanism for two-dimensional maps to have infinite coexistence of single-round periodic solutions. We will prove that such infinite coexistence is codimension-four in the orientation-preserving case, while it is codimension-three in the orientation-reversing case. Further we see how the number of periodic solutions change when we vary the parameters near the codimension points. This unfolding scenario is shown in Chapter \ref{cha:UnfGRHT}. We show that the bifurcations scale differently in different directions of the parameter space. Finally we showcase the infinite coexistence of single-round periodic solutions in a piecewise-smooth $C^{1}$ map for both the codimension-three and codimension-four scenarios in Chapter \ref{cha:ExampleInfinte}. We also give a theoretical framework behind the numerics used. In Chapter \ref{cha:OpenProblems}, we give future directions of this research. Preliminary research to find globally resonant homoclinic tangencies in the \blue{GHM} is discussed and other open problems are discussed. \blue{Then Section \ref{sec:conclusion} provides conclusions and Section \ref{sec:significance} discusses the broader significance of the results in the thesis.}

\chapter{Linearisation of maps\label{ch:StLin}}
\graphicspath{ {./images/} }
The standard approach for analysing any sort of homoclinic tangency is to first change coordinates to simplify the form of the map as much as possible in a neighbourhood of the saddle point. Already the Hartman-Grobman theorem (Theorem \ref{Th:HartGrob}) tells us that in a neighbourhood of a hyperbolic saddle the map is conjugate to its linearisation. However, the conjugacy (i.e. the coordinate change that takes the map to its linearisation) may not be differentiable. Sternberg's linearisation theorem \blue{\citep{St58}} tells us that the conjugacy is differentiable if the eigenvalues associated with the saddle satisfy a certain non-resonance condition (explained below). If the saddle is resonant, there exists a differentiable (actually $C^\infty$) coordinate change that transforms the map to something close to the linearisation --- specifically the linearisation plus resonant terms that cannot be removed. 

For the globally resonant homoclinic tangencies studied in later chapters, the saddle is resonant so we need to understand which terms are resonant and cannot be removed. This is the purpose of the present chapter. \blue{From now on $\mathcal{O}$ refers to the big-O notation discussed in \cite{De81}.}
\blue{\begin{definition}
 A function $f(x)$ on $\mathbb{R}^n$ is said to be  $\mathcal{O}(g(x))$ if $\limsup\limits_{x \rightarrow 0} \frac{f(x)}{g(x)}$ is finite.
 \label{defn:bigo}
\end{definition}
Often we will have $g(x) = \| x \|^m$ -- the $m^{\rm th}$ power of the norm of $x$, in which case we might abbreviate $\mathcal{O}(g(x))$ to $\mathcal{O}(m)$.} 

\section{Sternberg's linearisation}

Let $f$ be a smooth map on $\mathbb{R}^n$ for which the origin is a hyperbolic fixed point. Let $y = \phi(x)$ represent a change of coordinates, defined in a neighbourhood of the origin.  If the transformed map $\phi \circ f \circ \phi^{-1}$ is linear, we say that $\phi$ \textit{linearises} $f$. \blue{The theorem below gives conditions on the eigenvalues under which $f$ can be linearised.}
\begin{theorem}
Let $f$ be a $C^{\infty}$ map  defined in some neighborhood of the origin in $\mathbb{R}^n$ for which the origin is a fixed point. Let $\lambda_{1},\lambda_{2},\ldots,\lambda_{n}$ be the eigenvalues of ${\rm D}f(0)$, counting algebraic multiplicity. If \begin{equation}\lambda_{i} \neq \lambda_{1}^{m_{1}} \ldots \lambda_{n}^{m_{n}} \label{eq:ResoIneq}\end{equation} for all $i\in \{ 1,2,\ldots, n \}$ and for all non-negative integers $m_{1},m_{2},\ldots,m_{n}$ with $\sum \limits_{i=1}^{n} m_{i}> 1$, then there exists a neighborhood $\mathcal{N}$ of the origin on which a $C^{\infty}$ change of coordinates linearises $f$.
\label{Th:StLinear}
\end{theorem}
We now rewrite the condition of non-resonance \eqref{eq:ResoIneq} in a simpler form in the case of two-dimensional maps. Let $\lambda, \sigma$ denote the eigenvalues of the Jacobian matrix evaluated at the hyperbolic fixed point at the origin. We assume that $\lambda, \sigma \neq 0$. In the two-dimensional case, the non-resonant condition in Theorem \ref{Th:StLinear} can be written as $\lambda \neq \lambda^{m_{1}} \sigma^{m_{2}}, \blue{\sigma \neq \lambda^{n_{1}} \sigma^{n_{2}}}$ for all $m_{1},m_{2},\blue{n_{1},n_{2}} \geq 0$. It then follows that $\lambda^{m_{1}-1} \sigma^{m_{2}} \neq 1$. Let $p = m_{1}-1, q=m_{2}$, then $p+q = m_{1} + m_{2} - 1 \geq 1$ as $m_{1}+m_{2} >1$. Similarly \blue{$\lambda^{n_{1}} \sigma^{n_{2}-1} \neq 1$ as $\sigma \neq 0$. So here it is convenient to let $p=n_{1}, q = n_{2}-1$ giving again $p+q = n_{1}+n_{2}-1 \geq 1$}. This motivates the following definition of the resonance condition in the case of two-dimensional maps.

Let $f$ be a two-dimensional $C^{\infty}$ map defined in some neighborhood of the origin and for which the origin is a fixed point. Let $\lambda, \sigma$ denote the eigenvalues of ${\rm D}f(0)$. 
\begin{definition}
 The eigenvalues $\lambda,\sigma$ are said to be \textit{resonant} if there exist integers $p \geq -1$ and $q \geq -1$ with $p+q \geq 1$ for which $\lambda^{p} \sigma^{q}=1$.
\label{def:Reso}
\end{definition}
The following result is Theorem \ref{Th:StLinear} in the two-dimensional setting.
\begin{theorem}
 If the eigenvalues $\lambda, \sigma$ are non-resonant, there exists a neighborhood $\mathcal{N}$ of the origin on which a change of coordinates  linearises $f$. 
\label{Th:gencase}
\end{theorem}
In this chapter we show how the change of variables can be constructed and how the resonances arise. We will not prove that the change of coordinates is $C^{\infty}$. We start by showing how quadratic terms can be eliminated. 

Consider a two-dimensional map of the form
\begin{equation}
\label{eq:ResoMapEXthm}
\begin{aligned}
    x' &= \lambda x + a_{1} x^2 + a_{2} x y + a_{3} y^{2} + \mathcal{O}(3),\\
    y' &= \sigma y + b_{1} x^{2} + b_{2} x y + b_{3} y^{2} + \mathcal{O}(3),
\end{aligned}
\end{equation}
where we assume $0 < \lambda < 1 < \sigma$. Below we construct a coordinate change that transforms \eqref{eq:ResoMapEXthm} to
\begin{equation}
    \label{eq:linearCoordTransf}
    \begin{aligned}
    u' &= \lambda u + \mathcal{O}(3),\\
    v' &= \sigma v + \mathcal{O}(3).
    \end{aligned}
\end{equation}

Consider a near identity coordinate change of the form
\begin{equation}
    \label{eq:uvCoordChange}
    \begin{aligned}
    u &= x + c_{1}x^{2} + c_{2}xy + c_{3}y^{2},\\
    v &= y + d_{1}x^{2} + d_{2}xy + d_{3}y^{2}.
    \end{aligned}
\end{equation}
Inverting \eqref{eq:uvCoordChange}, we get $x,y$ in terms of $u,v$ as follows 
\begin{equation}
 \label{eq:InvertSeruv}
    \begin{aligned}
    x &= u - c_{1} u^{2} - c_{2} u v - c_{3} v^{2} + \mathcal{O}(3),\\
    y &= v - d_{1} u^{2} - d_{2} u v - d_{3} v^{2} + \mathcal{O}(3).\\
    \end{aligned}
\end{equation}
Rewriting the map \eqref{eq:ResoMapEXthm} in $(u,v)$-coordinates by substituting the expressions of $x,y$ in terms of $u,v$ as in \eqref{eq:InvertSeruv}, we get
\begin{equation}
  \label{eq:UVmapquad}
    \begin{aligned}
    u' &= x' + c_{1} (x')^{2} + c_{2} x' y'  + c_{3} (y')^{2}\\
       &= \lambda x + a_{1} x^{2} + a_{2} xy + a_{3} y^{2} + c_{1} \lambda^{2} x^{2} + c_{2} \lambda \sigma x y + c_{3} \sigma^{2} y^{2} + \mathcal{O}(3)\\
       &= \lambda(u-c_{1}u^{2}-c_{2}uv-c_{3}v^{2}) + a_{1}u^{2} + a_{2}uv + a_{3}v^{2} + c_{1}\lambda^{2}u^{2} + c_{2}\lambda\sigma uv + c_{3}\sigma^{2}v^{2} + \mathcal{O}(3)\\
       &= \lambda u + u^{2} (-\lambda c_{1} + a_{1} + c_{1}\lambda^{2}) + uv(-\lambda c_{2} + a_{2} + c_{2}\lambda\sigma )+ v^{2}(-\lambda c_{3} + a_{3} + c_{3}\sigma^{2}) + \mathcal{O}(3).
    \end{aligned}
\end{equation}
Observing the expression of $u'$ above in \eqref{eq:UVmapquad}, we try to choose the values of $c_1$, $c_2$, and $c_3$ such that the quadratic terms in the expression of $u'$ are eliminated. The $u^{2}$-term can be eliminated by setting $c_{1}(\lambda^{2} -\lambda) + a_{1} =0$ that is $c_{1} = \frac{a_{1}}{\lambda(1-\lambda)}$. The $uv$-term can be eliminated by choosing $c_{2} = \frac{a_{2}}{\lambda (1- \sigma)}$. Similarly, the $v^{2}$-term can also be eliminated by choosing $c_{3} = \frac{a_{3}}{\lambda - \sigma^{2}}$. We see that all quadratic terms $u^{2}, uv$, and $v^2$ can be eliminated via a coordinate change because we have assumed $0 < \lambda < 1 < \sigma$. If instead $\lambda = 0$, $1$, or $\sigma^2$, or $\sigma = 1$, then the eigenvalues would be resonant and at least one of the quadratic terms could not be eliminated.

Next for the $v$ variable, we obtain similarly 
\begin{equation}
    \begin{aligned}
    v' &= y' + d_{1}(x')^{2} + d_{2} x' y' + d_{3} (y')^{2}\\
       &= \sigma y + b_{1} x^{2} + b_{2} x y + b_{3} y^{2} + d_{1} \lambda^{2} x^{2} + d_{2} \lambda \sigma x y + d_{3} \sigma^{2} y^{2} + \mathcal{O}(3)\\
       &= \sigma v + u^{2}(-\sigma d_{1} + b_{1} + \lambda^{2}d_{1}) + u v(-\sigma d_{2} + b_{2} + \lambda \sigma d_{2}) + + v^{2} (-\sigma d_{3} + b_{3} + d_{3} \sigma^{2}) +\mathcal{O}(3).
    \end{aligned}
    \label{eq:Vmapquad}
\end{equation}
Observing the expression of $v'$ above in \eqref{eq:Vmapquad}, we try to choose the values of $d_1$, $d_2$, and $d_3$ such that the quadratic terms in the expression of $v'$ are eliminated. We can eliminate the $u^2$-term by choosing $d_{1} = \frac{b_{1}}{\sigma - \lambda^{2}}$. The $uv$-term can be eliminated by choosing $d_{2} = \frac{b_{2}}{\sigma(1-\lambda)}$. The $v^2$-term can be eliminated by choosing $d_{3} = \frac{b_{3}}{\sigma(1-\sigma)}$. 

\section{Resonant terms}
In this section, we assume that $\lambda$ and $\sigma$ are resonant.  In this case there exists a $C^\infty$ coordinate change to a map equal to the linearisation plus some resonant terms that cannot be eliminated. Our goal here is to explain exactly which terms are the resonant terms.

Let $n \ge 2$ be an integer.  Consider a map of the form
\begin{equation}
    \label{eq:xyresoterms}
    \begin{aligned}
    x' &= \lambda x + P(x,y) + \sum_{i=0}^{n} a_{i} x^{n-i}y^i + \mathcal{O}(n+1),\\
    y' &= \sigma y + Q(x,y) + \sum_{i=0}^{n} b_{i} x^{n-i} y^{i} + \mathcal{O}(n+1),
    \end{aligned}
\end{equation}
where $P(x,y)$ and $Q(x,y)$ represent the resonant terms of degree two to degree $n-1$. The two series contain all terms of order $n$.  We now perform a coordinate change from which we can see exactly which of these order $n$ terms can be eliminated based on the values of $\lambda$ and $\sigma$.

We perform a near identity coordinate transformation of the form
\begin{equation}
    \label{eq:gencoordtransf}
    \begin{aligned}
    u &= x + \sum_{i=0}^{n} c_{i} x^{n-i} y^{i},\\
    v &= y + \sum_{i=0}^{n} d_{i} x^{n-i} y^{i}.\\
    \end{aligned}
\end{equation}
The inverse of \eqref{eq:gencoordtransf} is 
\begin{equation}
    \label{eq:InvGenxy}
    \begin{aligned}
    x &= u + \sum_{i=0}^{n} -c_{i} u^{n-i} v^{i}+ \mathcal{O}(n+1),\\
    y &= v + \sum_{i=0}^{n} -d_{i} u^{n-i} v^{i}+ \mathcal{O}(n+1).
    \end{aligned}
\end{equation}
The $u$-component of the map in new coordinates is then
\begin{equation}
    \label{eq:GenUVCCoord}
    \begin{aligned}
    u' &= x' + \sum_{i=0}^{n} c_{i} (x')^{n-i} (y')^{i}\\
       &= \lambda x + P(x,y) + \sum_{i=0}^{n} a_{i} x^{n-i} y^{i} + \sum_{i=0}^{n} c_{i} \lambda^{n-i} x^{n-i} \sigma^{i} y^{i} + \mathcal{O}(n+1)\\
       &= \lambda u +  P(u,v)+ \sum_{i=0}^{n} (-\lambda c_{i} + a_{i} + c_{i}\lambda^{n-i}\sigma^i) u^{n-i} v^{i} + \mathcal{O}(n+1).
    \end{aligned}
\end{equation}
Observe that the term $u^{n-i} v^{i}$ can be eliminated by choosing $c_{i} = \frac{a_{i}}{\lambda - \lambda^{n-i}\sigma^{i}}$ for $i=0,\ldots,n$, assuming the denominator in this expression is non-zero. The denominator vanishes if $\lambda^{n-i-1}\sigma^i = 1$ which is an example of the eigenvalues being resonant, as in Definition \ref{def:Reso}.
Similarly, for the $v$-coordinate we have
\begin{equation}
    \label{eq:GenVCoord}
    \begin{aligned}
    v' &= y' + \sum_{i=0}^{n} d_{i} (x')^{n-i} (y')^{i},\\
       &= \sigma y + Q(x,y) + \sum_{i=0}^{n} b_{i} x^{n-i} y^{i} + \sum_{i=0}^{n}d_{i} \lambda^{n-i} x^{n-i} \sigma^{i} y^{i} + \mathcal{O}(n+1)\\
       &= \sigma v + Q(u,v) + \sum_{i=0}^{n} (-\sigma d_{i} + b_{i} + \sigma^{i} \lambda^{n-i} d_{i}) u^{n-i} v^{i} + \mathcal{O}(n+1). 
    \end{aligned}
\end{equation}
Similarly we observe that the term $u^{n-i} v^{i}$ can be eliminated by choosing $d_{i} = \frac{b_{i}}{\sigma (1-\sigma^{i-1}\lambda^{n-i} )}$ for $i=0,1,\ldots,n$. The condition $\sigma^{i-1}\lambda^{n-i} = 1$ represents the resonance condition for this term. \blue{The following theorem provides a simpler form of the map $f$ after the coordinate change has been applied.}
\begin{theorem}
\label{Th:LinearLS}
Let $f$ be a $C^{\infty}$ map on $\mathbb{R}^{2}$ with a saddle fixed point whose eigenvalues are $\lambda, \sigma$. If $0<|\lambda|<1<|\sigma|$ then the map $f$ can be brought to the map $T_{0}$ of the form \begin{equation}
T_{0}(x,y) = \begin{bmatrix}
\lambda x \left( 1 + x y F(x,y) \right) \\
\sigma y \left( 1 + x y G(x,y) \right)
\end{bmatrix}
\label{eq:T00}
\end{equation} under coordinate change, where $F,G$ are $C^{\infty}$ functions.
\end{theorem}
Observe the terms in \eqref{eq:GenUVCCoord} and \eqref{eq:GenVCoord} with $i=0$ and with $i=n$ can always be eliminated by choosing $c_{i} = \frac{a_{i}}{\lambda - \lambda^{n-i}\sigma^{i}}$ and $d_{i} = \frac{b_{i}}{\sigma (1-\sigma^{i-1}\lambda^{n-i} )}$ respectively and can be written as \eqref{eq:T00}.

In this thesis we will be particularly interested in the resonant cases $\lambda \sigma = 1$ and $\lambda\sigma=-1$. Two-dimensional maps can be \blue{brought} into the following normal forms by $C^{\infty}$ changes of variable. 
We have already determined the resonant \blue{(non-removable)} terms above. \cite{St58} proved that the required change of variables to eliminate all non-resonant terms does exist and is $C^{\infty}$.

\blue{The following theorems provide a simpler form of the map $f$ after coordinate change in the orientation-preserving and orientation-reversing cases respectively.}
\begin{theorem}
The map 
\begin{equation}
    \label{eq:T0eqn}
    \begin{aligned}
    x' &= \lambda x \Big(1 + xyF(x,y)\Big),\\
    y' &= \frac{1}{\lambda} y\Big(1 + xyG(x,y)\Big),
    \end{aligned}
\end{equation}
can be transformed to 
\begin{equation}
    \label{eq:OriPresf}
    \begin{aligned}
    x' &= \lambda x \Big(1 + xyA(xy)\Big),\\
    y' &= \frac{1}{\lambda} y \Big(1 + xyB(xy)\Big),
    \end{aligned}
\end{equation}
where $F,G,A$ and $B$ are $C^{\infty}$ functions.
\label{Th:OriPres}
\end{theorem}

\begin{theorem}
The map 
\begin{equation}
    \label{eq:T0eqn11}
    \begin{aligned}
    x' &= \lambda x \Big(1 + xyF(x,y)\Big),\\
    y' &= -\frac{1}{\lambda} y\Big(1 + xyG(x,y)\Big)
    \end{aligned}
\end{equation}
can be transformed to 
\begin{equation}
    \label{eq:OriRevf}
    \begin{aligned}
    x' &= \lambda x \Big(1 + x^2y^2A(x^2y^2)\Big),\\
    y' &= -\frac{1}{\lambda} y \Big(1 + x^2y^2B(x^2y^2)\Big),
    \end{aligned}
\end{equation}
where $F,G,A,$ and $B$ are $C^{\infty}$ functions.
\label{Th:OriRev}
\end{theorem}
Theorems \ref{Th:OriPres} and \ref{Th:OriRev} follow from the results of \cite{St58}. We will see in Chapters \ref{cha:GRHT} and \ref{cha:UnfGRHT} that it is extremely important that the \blue{orientation}-preserving case (\ref{Th:OriPres}) contains cubic resonant terms while the \blue{orientation}-reversing case (\ref{Th:OriRev}) does not. This will mean that the bifurcation due to a homoclinic tangency studied here is codimension-four in the first case but codimension-three in the second.

\blue{The presence of cubic resonant terms in \eqref{eq:OriPresf} but not \eqref{eq:OriRevf} may appear to be a contradiction for the following reason.  When one constructs a map, call it $g$, as the second iterate of a map in the orientation-reversing form \eqref{eq:T0eqn11}, $g$ is orientation-preserving, so why can some of its cubic terms not be removed, when all cubic terms of \eqref{eq:T0eqn11} can be removed?  The answer is that all cubic resonant terms vanish when forming the second iterate.  This is explained by the following result.}
\blue{
\begin{lemma}
The second iterate of \eqref{eq:T0eqn11} is given by 
\begin{equation}
    \label{eq:seconditerate}
    \begin{aligned}
    x'' &= \lambda^{2} x (1 + \mathcal{O}(x^2y)),  \\
    y'' &= \frac{1}{\lambda^{2}}y(1 + \mathcal{O}(xy^2)).
    \end{aligned}
\end{equation}
\end{lemma}}
\blue{
\begin{proof}
Considering the second iterate of the map in \eqref{eq:T0eqn11}, we get
\begin{equation}
\label{eq:interfirst}
\begin{aligned}
x'' &= \lambda (\lambda x(1+x y F(x,y)))  \left( 1 - xy (1+ xyF(x,y))(1+xyG(x,y))F) \right), \\
y'' &= -\frac{1}{\lambda} (-\frac{1}{\lambda}y(1+x y G(x,y)))  \left( 1 - xy (1+ xyF(x,y))(1+xyG(x,y))G) \right). \\
\end{aligned}
\end{equation}
A further expansion of terms in brackets leads to
 \begin{equation}
 \label{eq:intersecond}
\begin{aligned}
x'' &= \lambda^{2} x (1 + \lambda^{2}xyF -\lambda^{2}xyF + \mathcal{O}(x^2y)),  \\
y'' &= \frac{1}{\lambda^{2}}y(1 + \frac{1}{\lambda^{2}} x yG - \frac{1}{\lambda^{2}} x y G + \mathcal{O}(xy^2) ).
\end{aligned}
\end{equation}
Observe that in \eqref{eq:intersecond}, the $x^2 y$ term in the first component and the $xy^2$ term in the second component vanish, and the map becomes \eqref{eq:seconditerate}.
\end{proof}
}

\chapter{Globally Resonant Homoclinic Tangencies\label{cha:GRHT}}
\graphicspath{ {./images/} }
\def\cN{\mathcal{N}}
\def\cO{\mathcal{O}}
\def\ee{\varepsilon}
\newtheorem{conjecture}[theorem]{Conjecture}
\newtheorem{proposition}[theorem]{Proposition}
\theoremstyle{definition}
\newtheorem{example}[definition]{Example}
\theoremstyle{remark}
\newtheorem{remark}{Remark}
Here we study the mechanism behind the coexistence of an infinite number of asymptotically stable single-round periodic solutions. This phenomenon occurs when a homoclinic tangency occurs simultaneously with a number of other special conditions that we shall discover and prove in this chapter. This special tangency is known as globally resonant homoclinic tangency (GRHT). The results in this chapter were published in \cite{MuMcSi21}. 
\section{Main results}
\label{sec:setup}
\subsection{Local coordinates}
\label{sec:sec2p1}
\begin{figure}[t!]
\begin{center}
\includegraphics{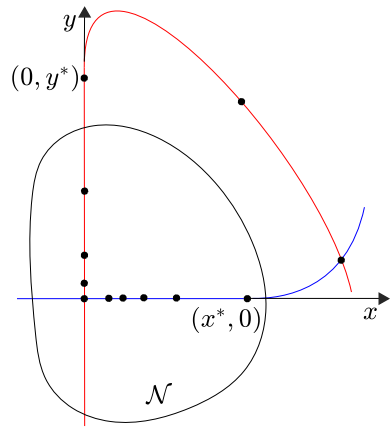}
\caption[Homoclinic orbit and local analysis problem setup]{ A homoclinic orbit is indicated with black dots in the case $0<\lambda<1$ and $\sigma > 1$.}
\label{fig:transv_intn_article}
\end{center}
\end{figure}

Let $f$ be a $C^\infty$ map on $\mathbb{R}^2$.
Suppose $f$ has a homoclinic orbit to a saddle fixed point
with eigenvalues $\lambda, \sigma \in \mathbb{R}$ where
\begin{equation}
0 < |\lambda| < 1 < |\sigma|.
\label{eq:eigenvalueAssumption}
\end{equation}
By applying an affine coordinate change
we can assume the fixed point lies at the origin, $(x,y) = (0,0)$,
about which the stable manifold is tangent to the $x$-axis
and the unstable manifold is tangent to the $y$-axis.
Then by Theorem \ref{Th:LinearLS}
there exists a locally valid $C^\infty$ coordinate change under which $f$ is transformed to
\begin{equation}
T_0(x,y) = \begin{bmatrix}
\lambda x \left( 1 + x y F(x,y) \right) \\
\sigma y \left( 1 + x y G(x,y) \right)
\end{bmatrix},
\label{eq:T01}
\end{equation}
where $F$ and $G$ are $C^\infty$.
In these new coordinates let $\cN$ be a bounded convex neighbourhood of the origin for which
\begin{equation}
f(x,y) = T_0(x,y), \quad \text{for all}~(x,y) \in \cN,
\label{eq:N}
\end{equation}
see Fig.~\ref{fig:transv_intn_article}. 
%Fig B below

It is possible to make $F$ and $G$ identically zero if $\lambda$ and $\sigma$ are  non-resonant, Theorem \ref{Th:gencase}.

In $\cN$ the local stable and unstable manifolds of the origin coincide with the $x$ and $y$-axes, respectively.
Let $(x^*,0)$ and $(0,y^*)$ be some points of the homoclinic orbit.
Without loss of generality we assume $x^* > 0$ and $y^* > 0$.
We further assume $x^*$ and $y^*$ are chosen sufficiently small so that
\begin{equation}
(x^*,0), (\lambda x^*,0), \left( 0, \tfrac{y^*}{\sigma} \right), \left( 0, \tfrac{y^*}{\sigma^2} \right) \in \cN.
\label{eq:fourPoints}
\end{equation}

Note that if $\lambda,\sigma > 0$ then the conditions that $(\lambda x^{*},0)$ and $\left(0,\frac{y^{*}}{\sigma^{2}}\right)$ belong to $\cN$ are redundant. We do not require $(0,y^{*})$ to belong to $\cN$ because it can be interpreted as the starting point of the excursion so does not need to map under $T_{0}$. Equation $\eqref{eq:fourPoints}$ ensures that the forward orbit of $(x^*,0)$ and a backward orbit of $(0,y^*)$ (these are both part of the homoclinic orbit) converge to the origin in $\cN$. By assumption $f^m(0,y^*) = (x^*,0)$ for some $m \ge 1$.
We let $T_1 = f^m$ and expand $T_{1}$ in a Taylor series centred at $(x,y) = (0,y^*)$:
\begin{equation}
T_1(x,y) = \begin{bmatrix}
x^* + c_1 x + c_2 (y - y^*) + \cO \left( \left( x, y-y^* \right)^2 \right) \\
d_1 x + d_2 (y - y^*) + d_3 x^2 + d_4 x (y - y^*) + d_5 (y - y^*)^2 + \cO \left( \left( x, y-y^* \right)^3 \right)
\end{bmatrix},
\label{eq:T11}
\end{equation}
where here we have explicitly written the terms that will be important below.

The value of $m$ and the values of the coefficients $c_i, d_i \in \mathbb{R}$
depend on our choice of $x^*$ and $y^*$.
It is a simple exercise to show that $d_1 = d_1(x^*,y^*) = \frac{\xi y^*}{x^*}$, for some constant $\xi$, as shown in the next theorem. \blue{We prove this for the orientation-preserving case below ; the orientation-reversing case can be proved similarly.}
\begin{theorem}
The quantity $\frac{d_{1}x^{*}}{y^{*}}$ is independent of the choice of homoclinic points $(0,y^*)$ and $(x^*,0)$.
\end{theorem}
\begin{proof}
Consider the map $T_{0}$ and $T_{1}$ as given in \eqref{eq:T01} and \eqref{eq:T11} respectively. \blue{Let $T_{1,2}$ represent the second component of  the map $T_{1}$}. We have $d_{1} = \frac{\partial{T_{1,2}}}{\partial{x}}(0,y^*)$. 
\begin{figure}\hspace{3cm}
\includegraphics{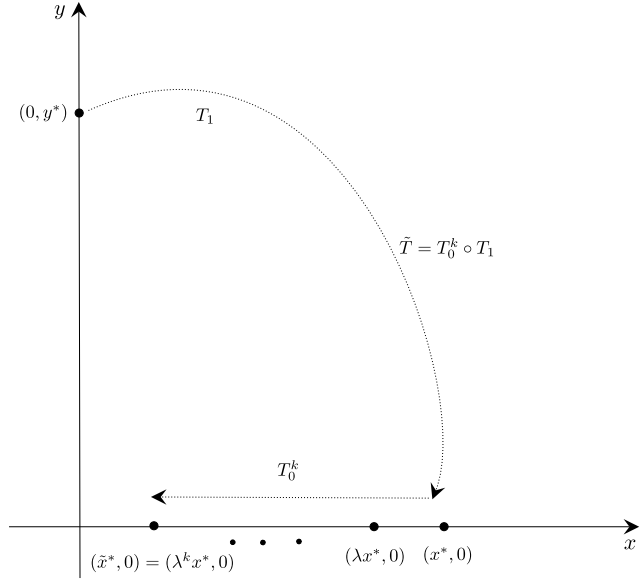}
\caption[A sketch to accompany the argument in the text that changing the homoclinic point on the $x$-axis leaves the quantity $\frac{|d_{1}|x^{*}}{y^{*}}$  unchanged.]{A sketch to accompany the argument in the text that changing the homoclinic point on the $x$-axis \blue{from $(x^{*},0)$ to $(\tilde{x}^{*},0)$} leaves the quantity $\frac{|d_{1}|x^{*}}{y^{*}}$ unchanged. \blue{The excursion from the homoclinic point $(0,y^{*})$ to the homoclinic point $(\tilde{x}^{*},0)$ determined by the map $\tilde{T}$ is shown in black dots.}}
\label{fig:thirdcondnY}
\end{figure}
Instead of using $y^*$, let us use another point say $\tilde{y}^* = \frac{y^*}{\sigma^k}$ as shown in Fig.~\ref{fig:thirdcondnY}. So, the map $T_{1}$ changes and we represent the new map by $\tilde{T}_{1}:$
\begin{equation}
    \begin{aligned}
    \tilde{T}_{1}(x,y) &= T_{1} \circ T_{0}^k,\\
                       &= \begin{bmatrix} x^* + \blue{c_{1} \lambda^{k} x} + c_{2}(\sigma^k y - y^*) + \mathcal{O}(2)\\ d_{1}\lambda^k x + \mathcal{O}(2)\end{bmatrix}.
    \end{aligned}
    \label{eq:NewMapT1}
\end{equation}

Thus the new value of $d_1$ is $\tilde{d}_1 =\frac{\partial{\tilde{T}}_{1,2}}{\partial{x}} = d_{1} \lambda^k$. Therefore the new value of the quantity of interest is $\frac{\tilde{d}_1 x^*}{\tilde{y}^*} = \frac{d_1 \lambda^k \sigma^k x^*}{y^*}$. As $\lambda \sigma  =1$, this simplifies to the original quantity $\frac{d_1 x^*}{y^*}$. This shows that the quantity is independent of the choice of $y^*$. By similarly altering our choice for the point $x^*$, or by altering both $x^*$ and $y^*$, we arrive at the same conclusion.
\begin{figure}\hspace{3cm}
\includegraphics{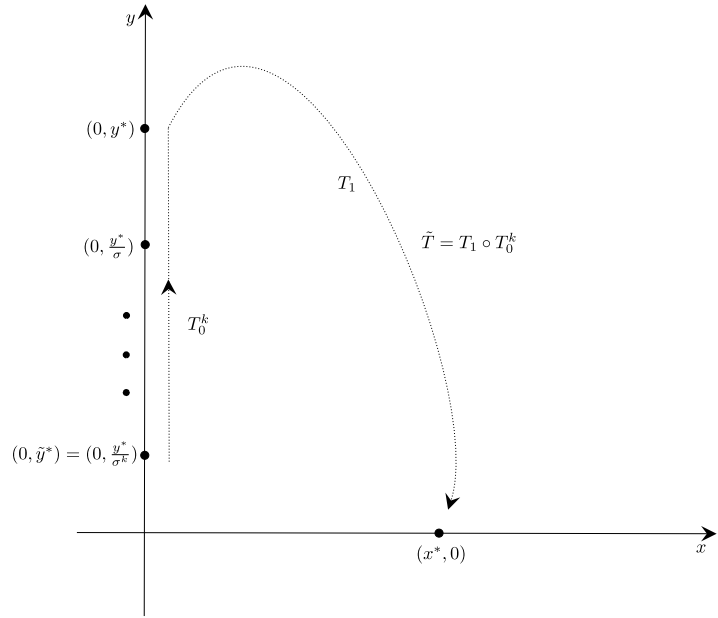}
\caption[Change of homoclinic point on the $y$-axis leaves the quantity $\frac{|d_{1}|x^{*}}{y^{*}}$ unchanged.]{Change of homoclinic point on the $y$-axis \blue{from $(0,y^{*})$ to $(0,\tilde{y}^{*})$} leaves the quantity $\frac{|d_{1}|x^{*}}{y^{*}}$ unchanged. \blue{The excursion from $(0,\tilde{y}^{*})$ to $(x^*,0)$ determined by the map $\tilde{T}$ is shown in black dots}.}
\label{fig:thirdcondnX}
\end{figure}
That is, $\frac{d_1 x^*}{y^*}$ is invariant with respect to a change in our choice of $x^*$ and $y^*$.
\end{proof}

This invariant (analogous to a quantity denoted $\tau$ in  \cite{DeGo15,GoSh87}) will be important below. From \eqref{eq:T11} we have det$\Big({\rm D}T_{1}\big(0,y^{*}\big)\Big) = c_{1}d_{2} - c_{2}d_{1}$. Thus if $c_{1}d_{2} - c_{2}d_{1} \neq 0$ then $f$ is locally invertible along the homoclinic orbit.
The stable and unstable manifolds of the origin intersect tangentially at $(0,y^*)$ if and only if $d_2 = 0$.
From basic homoclinic theory \cite{PaTa93,Ro04}
it is known that a tangential intersection is necessary to have stable single-round periodic solutions
near the homoclinic orbit.
For completeness we prove this as part of Theorem \ref{th:necessaryConditions} below.
\subsection{Three necessary conditions for infinite coexistence}
\label{sec:sec2p2}
It can be expected that points of single-round periodic solutions in $\cN$
converge to the $x$ and $y$-axes as the period tends to infinity.
This motivates our introduction of the set
\begin{equation}
\cN_\eta = \left\{ (x,y) \in \cN \,\big|\, |x y| < \eta \right\},
\label{eq:Nalt}
\end{equation}
where $\eta > 0$.
Below
we control the resonance terms in \eqref{eq:T01} by
choosing the value of $\eta$ to be sufficiently small. We first provide an $\eta$-dependent definition of single-round periodic solutions (SR abbreviates single-round).
\begin{figure}\hspace{3cm}
\includegraphics{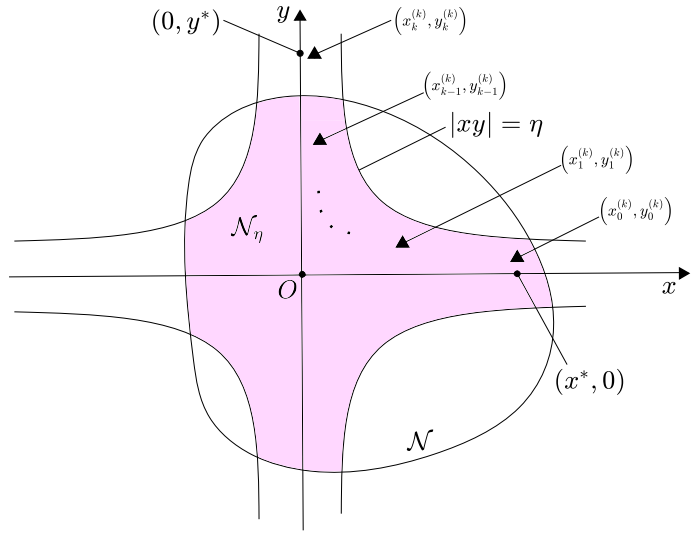}
\caption[Analytical definition of single-round periodic solution.]{Selected points of an {\em ${\rm SR}_k$}-solution (single-round periodic solution satisfying Definition \ref{df:SRkSolution}) in the case $\lambda > 0$ and $\sigma > 0$.
The region $\cN_\eta$ is shaded.}
\label{fig:SRP_defn}
\end{figure}
%.................................................................................
\begin{definition}
An {\em ${\rm SR}_k$-solution}
is a period-$(k+m)$ solution of $f$
involving $k$ consecutive points in $\cN_\eta$.
\label{df:SRkSolution}
\end{definition}

Below we consider ${\rm SR}_k$-solutions that approach the homoclinic orbit as $k \rightarrow \infty$ (see Theorem \ref{th:necessaryConditions} ). We denote points of an ${\rm SR}_k$-solution
by $\left( x^{(k)}_j, y^{(k)}_j \right)$, for $j = 0,\ldots,k+m-1$,
with $\left( x^{(k)}_j, y^{(k)}_j \right) \in \cN_\eta$ for all $j = 0,\ldots,k-1$
as shown in Fig.~\ref{fig:SRP_defn}.
The point $\left(x_{k}^{(k)},y_{k}^{(k)}\right)$ is a fixed point of $T_{0}^{k} \circ T_{1}$, thus the eigenvalues of ${\rm D} \left( T_0^k \circ T_1 \right) \left( x^{(k)}_k, y^{(k)}_k \right)$
determine the stability of the ${\rm SR}_k$-solution.
Let $\tau_k$ and $\delta_k$ respectively denote the trace and determinant of this matrix.
If the point $(\tau_k,\delta_k)$ lies in the interior of the triangle shown in Fig.~\ref{fig:stability_triangle}
then the ${\rm SR}_k$-solution is asymptotically stable.
If it lies outside the closure of this triangle then the ${\rm SR}_k$-solution is unstable.

%.................................................................................
\begin{theorem}
Suppose \eqref{eq:N} and \eqref{eq:fourPoints} are satisfied, with $x^*, y^* > 0$,
and $c_1 d_2 - c_2 d_1 \ne 0$. There exists $\eta > 0$ such that
if $f$ has a stable ${\rm SR}_k$-solution for infinitely many values of $k \ge 1$
and points $\left( x^{(k)}_0, y^{(k)}_0 \right)$ of these solutions converge to $(x^*,0)$ as $k \to \infty$, then 
\begin{equation}
d_2 = 0.
\label{eq:tangencyCondition}
\end{equation}
If also $d_5 \ne 0$, then
\begin{equation}
|\lambda \sigma| = 1,
\label{eq:det1Condition}
\end{equation}
and
\begin{equation}
|d_1| = \frac{y^*}{x^*},
\label{eq:globalResonanceCondition}
\end{equation}
with $d_1 = \frac{y^*}{x^*}$ in the case $\lambda \sigma = 1$.
\label{th:necessaryConditions}
\end{theorem}

\blue{Theorem \ref{th:necessaryConditions} is proved in $\S\ref{sec:necessaryConditions}$}. Here we provide some technical remarks. Equations \eqref{eq:tangencyCondition}--\eqref{eq:globalResonanceCondition}
are independent scalar conditions, thus together represent a codimension-three scenario.
The condition $d_2 = 0$ is equivalent to $(0,y^*)$ being a point of homoclinic tangency.
With also $d_5 \ne 0$ the tangency is quadratic.
The condition $|\lambda \sigma| = 1$ is equivalent to $f$ being area-preversing at the origin.

\begin{figure}[h!]
\begin{center}
\includegraphics{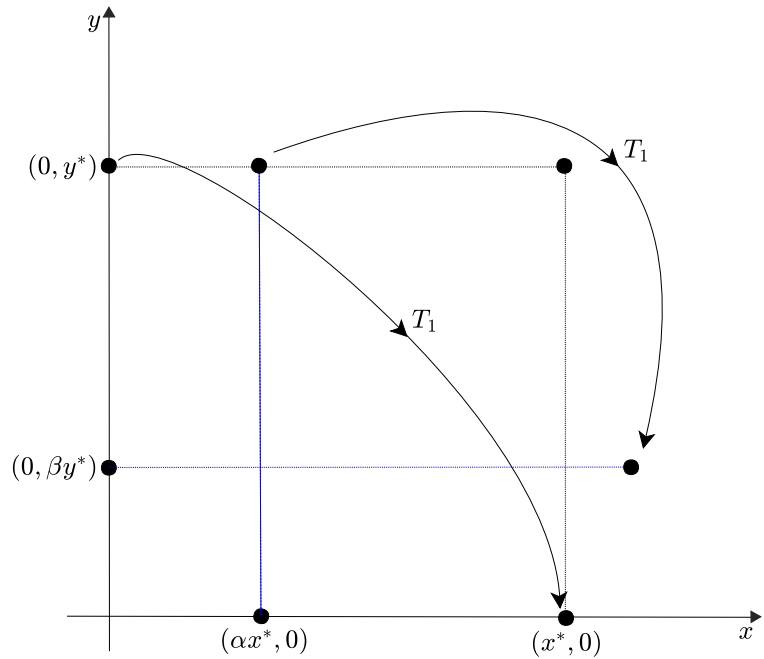}
\caption[A geometrical interpretation of the global resonance condition.]{ A geometrical interpretation of the global resonance condition \eqref{eq:globalResonanceCondition}. }
\label{fig:third_condn}
\end{center}
\end{figure}

The condition $|d_1| = \frac{y^*}{x^*}$ is a global condition,
termed {\em global resonance} in the area-preserving setting \cite{GoGo09}. This condition is well-defined because the value of $\frac{d_1 x^*}{y^*}$ is independent of our choice of $x^*$ and $y^*$.
To give geometric meaning to~\eqref{eq:globalResonanceCondition},
consider the perturbed point $(\alpha x^*, y^*)$, where $\alpha \in \mathbb{R}$ is small.
Under $T_1$ this point maps to
\begin{equation}
\Big( x^* \left( 1 + \cO(\alpha) \right), y^* \left( \beta + \cO \left( \alpha^2 \right) \right) \Big),
\nonumber
\end{equation}
where $\beta = \frac{d_1 x^* \alpha}{y^*}$. Notice we are writing the $x$ and $y$-components of these points as multiples of $x^*$ and $y^*$ so that these values provide a relative sense of scale.
Then condition~\eqref{eq:globalResonanceCondition} is equivalent to $|\beta| = |\alpha|$. Therefore~\eqref{eq:globalResonanceCondition} implies that when a point that is perturbed from $(0,y^*)$ in the $x$-direction by an amount $\alpha x^*$
is mapped under $T_1$, the result is a point that is perturbed from $(x^*,0)$ in the $y$-direction by an amount $\pm \alpha y^*$, to leading order as shown in Fig. \ref{fig:third_condn}. 
\subsection{The orientation-preserving case}
\label{sec:sec2p3}
Theorem~\ref{th:necessaryConditions} tells us that the infinite coexistence requires $|\lambda \sigma| = 1$.  Here we consider the case $\lambda \sigma = 1$ for which $f$ is orientation-preserving at the origin.
In this case $\lambda$ and $\sigma$ are resonant and by Theorem \ref{Th:OriPres}, $T_0$ can be reduced to the form
\begin{equation}
T_0(x,y) = \begin{bmatrix}
\lambda x \left( 1 + x y A(x y) \right) \\
\frac{1}{\lambda} y \left( 1 + x y B(x y) \right)
\end{bmatrix},
\label{eq:T02}
\end{equation}
where $A$ and $B$ are scalar $C^\infty$ functions.
Let
\begin{align}
a_1 &= A(0), & b_1 &= B(0).
\label{eq:a1b1}
\end{align}
The following theorem tells us that here infinite coexistence is only possible if $a_1 + b_1 = 0$.  This condition is satisfied automatically for area-preserving maps because 
\begin{equation}
    \begin{aligned}
    1 &= {\rm det} \left( {\rm D}f(x,y) \right)\\
       &= {\rm det}\left( {\rm D} T_{0}(x,y) \right) {\rm for} ~ {\rm all} ~ (x,y) \in \mathcal{N}\\
       &= 1 + 2(a_{1}+b_{1})xy + \mathcal{O}(x^2y^2).
    \end{aligned}
    \label{eq:FourthCondn}
\end{equation}

%.................................................................................
\begin{theorem}
Suppose \eqref{eq:N} and \eqref{eq:fourPoints} are satisfied with \eqref{eq:T02} and $x^*, y^* > 0$.  Further suppose $\lambda \sigma = 1$, $d_1 = \frac{y^*}{x^*}$, $d_2 = 0$, and $d_5 \ne 0$.
There exists $\eta > 0$ such that if $f$ 
has a stable ${\rm SR}_k$-solution for infinitely many values of $k \ge 1$
and points $\left( x^{(k)}_0, y^{(k)}_0 \right)$ of these solutions converge to $(x^*,0)$ as $k \to \infty$,
then $a_1 + b_1 = 0$.
\label{th:a1b1}
\end{theorem}
\blue{Theorem \ref{th:a1b1} is proved in \S \ref{sec:positiveCase}}. Together Theorems~\ref{th:necessaryConditions} and~\ref{th:a1b1} provide four necessary scalar conditions for infinite coexistence.  The next result provides sufficient conditions \blue{for the existence of infinitely many stable and saddle} single-round periodic solutions. 
%..................................................................................
\begin{theorem}
Suppose \eqref{eq:N} and \eqref{eq:fourPoints} are satisfied with \eqref{eq:T02} and $x^*, y^* > 0$. Further suppose $\lambda \sigma = 1$, $d_1 = \frac{y^*}{x^*}$, $d_2 = 0$, $a_1 + b_1 = 0$,
$d_5 \ne 0$ and $\Delta > 0$, where
\begin{equation}
\Delta = \left( 1 - \frac{c_2 y^*}{x^*} - \frac{d_4 y^*}{d_1} \right)^2 - 4 d_5 \left( d_3 x^{*^2} + c_1 d_1 x^* \right),
\label{eq:discriminant}
\end{equation}
and
\begin{equation}
-1 < \frac{c_2 y^*}{x^*} < 1 - \frac{\sqrt{\Delta}}{2}.
\label{eq:stabilityCondition}
\end{equation}
Then there exists $k_{\rm min} \in \mathbb{Z}$ such that for all $k \ge k_{\rm min}$
$f$ has an asymptotically stable ${\rm SR}_k$-solution and a saddle ${\rm SR}_k$-solution.
\label{th:sufficientConditionsPositiveCase}
\end{theorem}

From the proof of Theorem~\ref{th:sufficientConditionsPositiveCase} (given below in \S\ref{sec:positiveCase}) it can be seen that the condition $\Delta > 0$ ensures the existence of two ${\rm SR}_k$-solutions for large values of $k$, while~\eqref{eq:stabilityCondition} ensures that one of these is asymptotically stable and the other is a saddle.
%...............................................................................................

\subsection{The orientation-reversing case}
\label{sec:sec2p4}
We now consider the case $\lambda \sigma = -1$, in which $f$ is orientation-reversing at the origin. In this case $\lambda$ and $\sigma$ are again resonant but $T_0$ can be reduced further than in the $\lambda \sigma = 1$ case.
By Theorem \ref{Th:OriRev} we may assume
\begin{equation}
T_0(x,y) = \begin{bmatrix} \lambda x \left( 1 + \cO \left( x^2 y^2 \right) \right) \\
-\tfrac{1}{\lambda} y \left( 1 + \cO \left( x^2 y^2 \right) \right) \end{bmatrix}.
\label{eq:T0rev}
\end{equation}
Analogous to Theorem~\ref{th:sufficientConditionsPositiveCase}, the following result provides sufficient conditions for infinite coexistence .  We see that again it is only required that $\Delta > 0$ and \eqref{eq:stabilityCondition} is satisfied. However, here the periodic solutions only exist for either even values of $k$, or odd values of $k$, as determined by the sign of $d_1$.
%.................................................................................
\begin{theorem}
Suppose \eqref{eq:N} and \eqref{eq:fourPoints} are satisfied with \eqref{eq:T0rev} and $x^*, y^* > 0$. Further suppose $\lambda \sigma = -1$, $|d_1| = \frac{y^*}{x^*}$, $d_2 = 0$,
$d_5 \ne 0$, $\Delta > 0$ (where $\Delta$ is given by \eqref{eq:discriminant})
and \eqref{eq:stabilityCondition} is satisfied.
Then there exists $k_{\rm min} \in \mathbb{Z}$ such that
for all even $k \ge k_{\rm min}$ in the case $d_1 = \frac{y^*}{x^*}$,
and for all odd $k \ge k_{\rm min}$ in the case $d_1 = -\frac{y^*}{x^*}$,
$f$ has an asymptotically stable ${\rm SR}_k$-solution
and a saddle ${\rm SR}_k$-solution.
\label{th:sufficientConditionsNegativeCase}
\end{theorem}
\blue{Theorem \ref{th:sufficientConditionsNegativeCase} is proved in \S \ref{sec:negativeCase}.}
%..................................................................................
%===============================================================================

%===============================================================================
\section{Derivation and proof of necessary conditions for infinite coexistence}
\label{sec:necessaryConditions}
In this section we work towards a proof of Theorem \ref{th:necessaryConditions}.  Since an ${\rm SR}_k$-solution has period $k + m$, its stability depends on the eigenvalues of ${\rm D} f^{k+m}$ evaluated at any point of the solution.  Below we work with the point $\left( x^{(k)}_k, y^{(k)}_k \right)$.  Since $\left( x^{(k)}_k, y^{(k)}_k \right) = T_0^k \left( x^{(k)}_0, y^{(k)}_0 \right)$, and $\left( x^{(k)}_0, y^{(k)}_0 \right) = T_1 \left( x^{(k)}_k, y^{(k)}_k \right)$ (see Fig.~\ref{fig:SRP_defn}), we have
\begin{equation}
{\rm D} f^{k+m} \left( x^{(k)}_k, y^{(k)}_k \right) = {\rm D} T_0^k \left( x^{(k)}_0, y^{(k)}_0 \right)
{\rm D} T_1 \left( x^{(k)}_k, y^{(k)}_k \right).
\label{eq:DfkPlusm}
\end{equation}

To obtain information on the eigenvalues of \eqref{eq:DfkPlusm} we first construct bounds on the values of the points of an ${\rm SR}_k$-solution (Lemmas \ref{le:y0Bound} and \ref{le:xjyjBound}).  We then estimate the entries of ${\rm D} T_0^k \left( x^{(k)}_0, y^{(k)}_0 \right)$ (Lemmas \ref{le:MjBoundsl} and \ref{le:PjBoundsl}).  Next we estimate the contribution of the resonant terms in $T_0^k$ (Lemma \ref{le:T0k}).  These finally enable us to prove Theorem \ref{th:necessaryConditions}.  Essentially we show that if the conditions \eqref{eq:tangencyCondition}--\eqref{eq:globalResonanceCondition} are not all met then the trace of \eqref{eq:DfkPlusm}, denoted $\tau_k$ above, diverges as $k \to \infty$.  This implies that the eigenvalues of  \eqref{eq:DfkPlusm} cannot lie inside the shaded region of Fig.~\ref{fig:stability_triangle} for more than finitely many values of $k$.

%\setcounter{equation}{0}
%General aspects.

We first observe that the resonant terms $F$ and $G$ of \eqref{eq:T01} are continuous and $\cN$ is bounded, so there exists $R > 0$ such that
\begin{equation}
|F(x,y)| \le R, \quad |G(x,y)| \le R, \quad \text{for all}~ (x,y) \in \cN.
\label{eq:FGBound}
\end{equation}
Let
\begin{equation}
\eta_R = \frac{1 - |\sigma|^{-\frac{1}{2}}}{R},
\label{eq:xyBound}
\end{equation}
and notice $\eta_R > 0$.

%.................................................................................
\begin{lemma}
Suppose an infinite family of ${\rm SR}_k$-solutions with $\eta = \eta_R$
satisfies $\left( x^{(k)}_0, y^{(k)}_0 \right) \to (x^*,0)$ as $k \to \infty$.
Then
\begin{equation}
\left| y^{(k)}_0 \right| \le 2 y^* |\sigma|^{-\frac{k}{2}},
\label{eq:y0Bound}
\end{equation}
for all sufficiently large values of $k$.
\label{le:y0Bound}
\end{lemma}

%.................................................................................
\begin{proof}
By assumption $\left( x^{(k)}_j, y^{(k)}_j \right) \in \cN_{\eta_R}$
for each $j \in \{ 0,\ldots,k-1 \}$,
thus from \eqref{eq:T01}, \eqref{eq:FGBound}, and \eqref{eq:xyBound} we obtain
\begin{align*}
\left| y^{(k)}_{j+1} \right| &= |\sigma| \left| y^{(k)}_j \right|
\left| 1 + x^{(k)}_j y^{(k)}_j G \left( x^{(k)}_j, y^{(k)}_j \right) \right| \\
&\ge |\sigma| \left| y^{(k)}_j \right| \left( 1 - \eta_R R \right) \\
&= |\sigma|^{\frac{1}{2}} \left| y^{(k)}_j \right|.
\end{align*}
By applying this bound $k$ times we obtain
\begin{equation}
\left| y^{(k)}_0 \right| \le |\sigma|^{-\frac{k}{2}} \left| y^{(k)}_k \right|.
\nonumber
\end{equation}
But $y^{(k)}_k \to y^*$ as $k \to \infty$ \blue{because $y_{0}^{(k)} \rightarrow 0$ and therefore by mapping backwards under map $T_{1}$ we have $y^{(k)}_k \to y^*$}. Thus
$\left| y^{(k)}_k \right| \le 2 y^*$, say, for all sufficiently large values of $k$.
This verifies \eqref{eq:y0Bound}.
\end{proof}

%.................................................................................
\begin{lemma}
Suppose $|\lambda \sigma| \le 1$ and suppose an infinite family of ${\rm SR}_k$-solutions with $\eta = \eta_R$
satisfies $\left( x^{(k)}_0, y^{(k)}_0 \right) \to (x^*,0)$ as $k \to \infty$.
Then there exists $\omega > 0$ such that
\begin{equation}
\begin{split}
\left| x^{(k)}_j \right| &\le \omega |\lambda|^j, \\
\left| y^{(k)}_j \right| &\le \omega |\sigma|^{j-k},
\label{eq:xjyjBound}
\end{split}
\end{equation}
for all $j \in \{ 0, \ldots, k \}$ and all sufficiently large values of $k$.
\label{le:xjyjBound}
\end{lemma}

%.................................................................................
\begin{proof}
By Lemma \ref{le:y0Bound} there exists $\omega_0 > 0$ such that
\begin{align}
\left| x^{(k)}_0 \right| &\le \omega_0 \,,
\label{eq:x0Bound2} \\
\left| y^{(k)}_0 \right| &\le \omega_0 |\sigma|^{-\frac{k}{2}},
\label{eq:y0Bound2}
\end{align}
for all sufficiently large values of $k$.
For the remainder of the proof we assume $k$ is sufficiently large that
\begin{align}
\left( 1 + 4 \omega_0^2 R |\sigma|^{-\frac{k}{2}} \right)^k &< 2,
\label{eq:powerkBoundUpper} \\
\left( 1 - 4 \omega_0^2 R |\sigma|^{-\frac{k}{2}} \right)^k &> \tfrac{1}{2}.
\label{eq:powerkBoundLower}
\end{align}
\blue{Such $k$ exists because the left hand-sides of \eqref{eq:powerkBoundLower} and \eqref{eq:powerkBoundUpper} both converge to $1$ as $k \to \infty$.} Below we use induction on $j$ to prove that
\begin{align}
\left| x^{(k)}_0 \right| |\lambda|^j \left( 1 - 4 \omega_0^2 R |\sigma|^{-\frac{k}{2}} \right)^j
\le \left| x^{(k)}_j \right| \le
\left| x^{(k)}_0 \right| |\lambda|^j \left( 1 + 4 \omega_0^2 R |\sigma|^{-\frac{k}{2}} \right)^j,
\label{eq:xjDoubleBound} \\
\left| y^{(k)}_0 \right| |\sigma|^j \left( 1 - 4 \omega_0^2 R |\sigma|^{-\frac{k}{2}} \right)^j
\le \left| y^{(k)}_j \right| \le
\left| y^{(k)}_0 \right| |\sigma|^j \left( 1 + 4 \omega_0^2 R |\sigma|^{-\frac{k}{2}} \right)^j,
\label{eq:yjDoubleBound}
\end{align}
for all $j \in \{ 0,\ldots,k \}$.
This will complete the proof because, first, \eqref{eq:x0Bound2},
\eqref{eq:powerkBoundUpper}, and \eqref{eq:xjDoubleBound} combine to produce
\begin{equation}
\left| x^{(k)}_j \right| \le 2 \omega_0 |\lambda|^j.
\nonumber
\end{equation}
Second, \eqref{eq:powerkBoundLower} and \eqref{eq:yjDoubleBound} evaluated at $j=k$ combine to produce
$\frac{1}{2} \left| y^{(k)}_0 \right| |\sigma|^k \le \left| y^{(k)}_k \right|$.
But $y^{(k)}_k \to y^*$ as $k \to \infty$, thus
$\left| y^{(k)}_k \right| \le 2 y^*$ for sufficiently large values of $k$, and so
$\left| y^{(k)}_0 \right| \le 4 y^* |\sigma|^{-k}$.
Hence by \eqref{eq:powerkBoundUpper} and \eqref{eq:yjDoubleBound} we have
\begin{equation}
\left| y^{(k)}_j \right| \le 8 y^* |\sigma|^{j-k}.
\nonumber
\end{equation}

It remains to verify \eqref{eq:xjDoubleBound}--\eqref{eq:yjDoubleBound}.
Trivially these hold for $j=0$.
Suppose \eqref{eq:xjDoubleBound}--\eqref{eq:yjDoubleBound} hold for some $j \in \{ 0,\ldots,k-1 \}$
(this is our induction hypothesis).
Then by \eqref{eq:T01},
\begin{align}
|\lambda| \left| x^{(k)}_j \right| \left( 1 - \left| x^{(k)}_j y^{(k)}_j \right| R \right)
&\le \left| x^{(k)}_{j+1} \right| \le
|\lambda| \left| x^{(k)}_j \right| \left( 1 + \left| x^{(k)}_j y^{(k)}_j \right| R \right),
\label{eq:xjDoubleBoundProof} \\
|\sigma| \left| y^{(k)}_j \right| \left( 1 - \left| x^{(k)}_j y^{(k)}_j \right| R \right)
&\le \left| y^{(k)}_{j+1} \right| \le
|\sigma| \left| y^{(k)}_j \right| \left( 1 + \left| x^{(k)}_j y^{(k)}_j \right| R \right).
\label{eq:yjDoubleBoundProof}
\end{align}
The induction hypothesis implies
\begin{align}
\left| x^{(k)}_j y^{(k)}_j \right|
&\le \left| x^{(k)}_0 \right| \left| y^{(k)}_0 \right| |\lambda \sigma|^j
\left( 1 + 4 \omega_0^2 R |\sigma|^{-\frac{k}{2}} \right)^{2 j} \nonumber \\
&\le 4 \omega_0^2 |\sigma|^{-\frac{k}{2}},
\label{eq:xjyjBoundProof}
\end{align}
where we have used \eqref{eq:x0Bound2}, \eqref{eq:y0Bound2}, \eqref{eq:powerkBoundUpper}
and $|\lambda \sigma| \le 1$ in the second line.
By applying \eqref{eq:xjyjBoundProof} and
the induction hypothesis to \eqref{eq:xjDoubleBoundProof}--\eqref{eq:yjDoubleBoundProof}
we obtain \eqref{eq:xjDoubleBound}--\eqref{eq:yjDoubleBound} for $j+1$,
and this completes the induction step.
\end{proof}
%.................................................................................
% Two new lemma's for <= 1 case
For brevity we write
\begin{equation}
M_j = {\rm D} T_0 \left( x^{(k)}_j, y^{(k)}_j \right)
= \begin{bmatrix} p_j & q_j \\ r_j & s_j \end{bmatrix},
\label{eq:Mjl}
\end{equation}
and
\begin{equation}
P_j = M_{j-1} \cdots M_1 M_0
= \begin{bmatrix} t_j & u_j \\ v_j & w_j \end{bmatrix}.
\label{eq:Pjl}
\end{equation}
%.................................................................................
\begin{lemma}
Suppose $|\lambda \sigma| \le 1$ and suppose an infinite family of ${\rm SR}_k$-solutions with $\eta = \eta_R$
satisfies $\left( x^{(k)}_0, y^{(k)}_0 \right) \to (x^*,0)$ as $k \to \infty$.
Then there exists $\alpha > 0$ such that
\begin{equation}
\begin{aligned}
|p_j - \lambda| &\le \alpha |\sigma|^{-k}, &
|q_j| &\le \alpha |\sigma|^{-2 j}, \\
|r_j| &\le \alpha |\sigma|^{2 (j-k)}, &
\hspace{10mm}
\left|s_j - \sigma \right| &\le \alpha |\sigma|^{-k},
\end{aligned}
\label{eq:MjBoundsl}
\end{equation}
for all $j \in \{ 0,\ldots,k-1 \}$ and all sufficiently large values of $k$.
\label{le:MjBoundsl}
\end{lemma}

%.................................................................................
\begin{proof}
By Lemma \ref{le:xjyjBound} and $|\lambda \sigma| \leq 1$ there exists $\omega > 0$ such that
\begin{align}
\left| x^{(k)}_j \right| &\le \omega |\sigma|^{-j}, &
\left| y^{(k)}_j \right| &\le \omega |\sigma|^{j-k},
\label{eq:xjyjBound2l}
\end{align}
for all sufficiently large values of $k$.
By differentiating \eqref{eq:T01} we obtain
\begin{equation}
{\rm D} T_0(x,y) = \begin{bmatrix}
\lambda \big( 1 + x y (2 F + x \frac{\partial{F}}{\partial{x}}) \big) & \lambda x^2 (F +  y \frac{\partial{F}}{\partial{y}})\\
\sigma y^2 (G + x \frac{\partial{G}}{\partial{x}}) & \sigma \big( 1 + x y (2 G + y \frac{\partial{G}}{\partial{y}}) \big)
\end{bmatrix}.
\label{eq:DT0l}
\end{equation}
Since $F$, $G$, and their derivatives are continuous and $\cN_{\eta_R}$ is bounded
there exists $\alpha > 0$ such that throughout $\cN_{\eta_R}$ we have that
$\alpha$ is greater than
$\omega^2 |\lambda| |2 F + x \frac{\partial{F}}{\partial{x}}|$,
$\omega^2 |\lambda| |F +  y \frac{\partial{F}}{\partial{y}}|$,
$|\sigma| \omega^2 |G + x \frac{\partial{G}}{\partial{x}}|$, and
$|\sigma| \omega^2 |2 G +  y \frac{\partial{G}}{\partial{y}}|$.
Then \eqref{eq:MjBoundsl} follows immediately from
\eqref{eq:xjyjBound2l} and \eqref{eq:DT0l}.
\end{proof}

%.................................................................................
\begin{lemma}
Suppose $|\lambda \sigma| \le 1$ and suppose an infinite family of ${\rm SR}_k$-solutions with $\eta = \eta_R$
satisfies $\left( x^{(k)}_0, y^{(k)}_0 \right) \to (x^*,0)$ as $k \to \infty$.
Then there exists $\beta > 0$ such that
\begin{equation}
\begin{aligned}
|t_j - \lambda^j| &\le \beta j |\sigma|^{-(j+k)}, &
|u_j| &\le \beta j |\sigma|^{-j}, \\
|v_j| &\le \beta j |\sigma|^{j - 2 k}, &
\hspace{10mm}
\left|w_j - \sigma^{j} \right| &\le \beta j |\sigma|^{j-k},
\end{aligned}
\label{eq:PjBoundsl}
\end{equation}
for all $j \in \{ 1,\ldots,k \}$ and all sufficiently large values of $k$.
\label{le:PjBoundsl}
\end{lemma}

%.................................................................................
\begin{proof}
By Lemma \ref{le:xjyjBound} and $|\lambda\sigma| \leq 1$ there exists $\omega > 0$ satisfying
\eqref{eq:xjyjBound2l} for all sufficiently large values of $k$.
Let $\alpha > 0$ satisfy \eqref{eq:MjBoundsl}
and let $\beta = \frac{2 \alpha}{|\lambda|}$.

We now verify \eqref{eq:PjBoundsl} by induction on $j$.
Observe $P_1 = M_0$ so \eqref{eq:PjBoundsl} holds with $j = 1$
because \eqref{eq:MjBoundsl} holds with $j = 0$ and $\beta > \alpha$.
Now suppose \eqref{eq:PjBoundsl} holds for some $j \in \{ 1,\ldots,k-1 \}$
(this is our induction hypothesis).
Since $P_{j+1} = M_j P_j$ we have
\begin{equation}
\begin{bmatrix} t_{j+1} & u_{j+1} \\ v_{j+1} & w_{j+1} \end{bmatrix} =
\begin{bmatrix}
p_j t_j + q_j v_j & p_j u_j + q_j w_j \\
r_j t_j + s_j v_j & r_j u_j + s_j w_j
\end{bmatrix}.
\label{eq:Pjp1l}
\end{equation}
We now verify the four inequalities \eqref{eq:PjBoundsl} for $j+1$ in order.
First
\begin{align}
\left| t_{j+1} - \lambda^{j+1} \right|
&= \left| p_j t_j + q_j v_j - \lambda^{j+1} \right| \nonumber \\
&= \left| (p_j - \lambda) \left( t_j - \lambda^j \right) + (p_j - \lambda) \lambda^j
+ \lambda \left( t_j - \lambda^j \right) + q_j v_j \right| \nonumber \\
&\le |p_j - \lambda| \left| t_j - \lambda^j \right| + |p_j - \lambda| |\lambda|^j
+ |\lambda| \left| t_j - \lambda^j \right| + |q_j| |v_j|, \nonumber
\end{align}
and so by \eqref{eq:MjBoundsl} and the induction hypothesis we have
\begin{equation}
\left| t_{j+1} - \lambda^{j+1} \right|
\le \beta j |\sigma|^{-(k+j+1)}
+ \left( 2 \alpha \beta j |\sigma|^{-k+1} + \alpha|\sigma| \right) |\sigma|^{-(k+j+1)}.
\nonumber
\end{equation}
Thus for sufficiently large values of $k$,
\begin{align}
\left| t_{j+1} - \lambda^{j+1} \right|
&\le \beta j |\sigma|^{-(k+j+1)}
+ \tfrac{2 \alpha}{|\lambda|} |\sigma|^{-(k+j+1)}
\nonumber \\
&= \beta (j+1) |\sigma|^{-(k+j+1)},
\nonumber
\end{align}
where we have substituted our formula for $\beta$ in the last line.
In a similar fashion we obtain
\begin{align}
\left| u_{j+1} \right|
&= \left| p_j u_j + q_j w_j \right| \nonumber \\
&\le |p_j| |u_j| + |q_j| |w_j| \nonumber \\
&\le \left( |\lambda| + \alpha |\sigma|^{-k} \right) \beta j |\sigma|^{-j}
+ \alpha |\sigma|^{-2 j} \left( |\sigma|^{j} + \beta j |\sigma|^{j-k} \right) \nonumber \\
&= \beta j |\sigma|^{-(j+1)} + \left(2\alpha\beta j  |\sigma|^{-k+1} + \alpha |\sigma|\right)
|\sigma|^{-(j+1)}\\
%&\le \beta j |\lambda|^{j+1}
%+ \left( 2 \alpha \beta j |\lambda|^{k-1} + \alpha |\sigma| \right) |\lambda|^{j+1} \nonumber \\
%&\le \beta j |\lambda|^{j+1}
%+ \left( 2 \alpha \beta j |\lambda|^{k-1} + \tfrac{\alpha}{|\lambda|} \right) |\lambda|^{j+1} \nonumber \\
&\le \beta j |\sigma|^{-(j+1)} + \tfrac{2 \alpha}{|\lambda|} |\sigma|^{-(j+1)} \nonumber \\
&= \beta (j+1) |\sigma|^{-(j+1)}, \nonumber
\end{align}
\begin{align}
\left| v_{j+1} \right|
&= \left| r_j t_j + s_j v_j \right| \nonumber \\
&\le |r_j| |t_j| + |s_j| |v_j| \nonumber \\
&\le \alpha |\sigma|^{2 j - 2 k} \left( |\lambda|^j + \beta j |\sigma|^{-(k+j)} \right)
+ \left( |\sigma| + \alpha |\sigma|^{-k} \right) \beta j |\sigma|^{j - 2k} \nonumber \\
&\le \beta j |\sigma|^{(j+1) - 2 k}
+ \left( 2 \alpha \beta j |\sigma|^{-(k+1)} + \frac{\alpha}{|\sigma|} \right) |\sigma|^{(j+1) - 2 k} \nonumber \\
&\le \beta j |\sigma|^{(j+1) - 2 k} + \frac{2 \alpha}{|\sigma|} |\sigma|^{(j+1) - 2 k} \nonumber \\
&< \beta (j+1) |\sigma|^{(j+1) - 2 k}, \nonumber
\end{align}
and
\begin{align}
\left| w_{j+1} - \sigma^{(j+1)} \right|
&= \left| r_j u_j + s_j w_j - \sigma^{(j+1)} \right| \nonumber \\
&\le |r_j| |u_j| + \left| s_j - \sigma \right| \left| w_j - \sigma^{j} \right|
+ \left| s_j - \sigma \right| |\sigma|^{j} + |\sigma| \left| w_j - \sigma^{j} \right| \nonumber \\
&\le \beta j |\sigma|^{(j+1)-k}
+ \left( 2 \alpha \beta j |\sigma|^{-(k+1)} + \alpha |\sigma|^{-1} \right) |\sigma|^{(j+1)-k} \nonumber \\
&\le \beta j |\sigma|^{(j+1)-k}
+ \frac{2\alpha}{|\sigma|} |\sigma|^{(j+1)-k} \nonumber \\
&< \beta (j+1) |\sigma|^{(j+1)-k},
\end{align}
for sufficiently large values of $k$.
\end{proof}

%.................................................................................
\begin{lemma}
If $|\lambda \sigma| \leq 1$ then in a neighbourhood of $(x^*,0)$,
\begin{equation}
T_0^k(x,y) = \begin{bmatrix}
\lambda^k x \left( 1 + \cO(kxy) \right) \\
\sigma^k y \left(1 + \cO(kxy) \right)
\end{bmatrix}.
\label{eq:T0kl1}
\end{equation}
\label{le:T0k}
\end{lemma}
%.................................................................................
\begin{proof}
Write
\begin{equation}
    T_{0}^{k}(x,y) = \begin{bmatrix}
    \lambda^{k}x\left(1 + xyF_{k}(x,y)\right)\\
    \sigma^{k}y\left( 1 + xyG_{k}(x,y)\right)
    \end{bmatrix}.
    \label{eq:T0k1}
\end{equation}
Let $\delta = \frac{2^{\frac{1}{3}}-1}{2R}$. We will show that if $k|xy| \leq \delta$, then
\begin{equation}
    |F_{k}(x,y)| \leq 2kR, \hspace{0.2cm} |G_{k}(x,y)| \leq 2kR,
    \label{eq:FkGkl1}
\end{equation}
and this will complete the proof. We prove \eqref{eq:FkGkl1} by induction on $k$. Equation \eqref{eq:FkGkl1} is certainly true for $k=1$ because $|F_{1}(x,y)| = |F(x,y)| \leq R < 2R$, and similarly for $G_{1}$. Suppose \eqref{eq:FkGkl1} is true for some $k \geq 1$ (this is our induction hypothesis). By matching terms in $T_{0}^{k+1} = T_{0} \circ T_{0}^{k}$ we obtain  
\begin{equation}
    F_{k+1}(x,y) = F_{k}(x,y) + \lambda^{k}\sigma^{k}F\left(T_{0}^{k}(x,y)\right) \left( 1 + xyF_{k}(x,y)\right)^{2} \left( 1 + xyG_{k}(x,y)\right).
\end{equation}
By applying the induction hypothesis we obtain
\begin{equation}
    |F_{k+1}(x,y)| \leq 2kR + |F(T_{0}^{k}(x,y))|\left( 1 + |xy|2kR\right)^{3}.
\end{equation}
Since $k|xy|$ is small and $(x,y) \approx (x^*,0)$, we can assume $T_{0}^{k}(x,y) \in \cN$ and so $|F\left(T_{0}^{k}(x,y)\right)|\leq R$. Also $1 + |xy|2kR \leq 2^{\frac{1}{3}}$, so
\begin{equation}
    |F_{k+1}(x,y)| \leq 2(k+1)R,
\end{equation}
as required (and the result for $G_{k+1}$ is obtained similarly).
\end{proof}
%....................................................................................
%.................................................................................
\begin{lemma}
If $|\lambda \sigma| = 1$ then in a neighbourhood of $(x^*,0)$,
\begin{equation}
T_0^k(x,y) = \begin{bmatrix}
\lambda^k x \left( 1 + ka_{1}xy + \cO(k^2x^2y^2) \right) \\
\sigma^k y \left(1 + kb_{1}xy + \cO(k^2x^2y^2) \right)
\end{bmatrix}.
\label{eq:T0koripres1}
\end{equation}
\label{le:T0koripv}
\end{lemma}
%.................................................................................
\begin{proof}
The map $T_{0}$ can be written as below 
\begin{equation}
    T_{0}(x,y) = \begin{bmatrix}
    \lambda x\left(1 + a_{1}xy + x^2y^2A_{1}(xy)\right)\\
    \frac{1}{\lambda}y\left( 1 + b_{1}xy + x^2y^2B_{1}(xy)\right)
    \end{bmatrix},
    \label{eq:T0oripres2}
\end{equation}
where $a_{1} = A(0), b_{1} = B(0)$.
Since the region $\mathcal{N}$ is bounded, there exists a constant $S$ such that throughout $\mathcal{N}$ we have $|a_{1}| \leq S, |b_{1}| \leq S, |A_{1}(xy)| \leq S$ and $|B_{1}(xy)| \leq S$. Write 
\begin{equation}
    T_{0}^{k}(x,y) = \begin{bmatrix}
    \lambda^k x\left(1 + ka_{1}xy + x^2y^2A_{k}(xy)\right)\\
    \frac{1}{\lambda^k}y\left( 1 + kb_{1}xy + x^2y^2B_{k}(xy)\right)
    \end{bmatrix},
    \label{eq:T0koripresk}
\end{equation}
Let $\delta = \frac{1}{56S} = \frac{M}{S}$. We will show that if $k|xy| \leq \delta$ then 
\begin{equation}
    |A_{k}(x,y)| \leq 2k^2S^2, \hspace{0.2cm} |B_{k}(x,y)| \leq 2k^2S^2,
    \label{eq:AkBk}
\end{equation}
and this will complete the proof due to the meaning of big-$\mathcal{O}$ notation. We prove \eqref{eq:AkBk} by induction on $k$. Equation \eqref{eq:AkBk} is certainly true for $k=1$ because $|A_{1}(xy)| \leq S < 2S^2$, and similarly for $B_{1}$. Suppose \eqref{eq:AkBk} is true for some $k \geq 1$ (this is our induction hypothesis). By matching terms in $T_{0}^{k+1} = T_{0} \circ T_{0}^{k}$ we obtain  
\begin{equation}
\begin{aligned}
    A_{k+1}x^2y^2 &= T + x^2y^2\left(A_{k} + ka_{1}(2a_{1}+b_{1})\right)\\
                  &+ x^2y^2\left(A_{1}(T_{0}^{k})(1+ka_{1}xy+x^2y^2A_{k})^3(1+kb_{1}xy+x^2y^2B_{k})^2 \right)
\end{aligned}    
    \label{eq:Akp1el}
\end{equation}
where \begin{equation}
\begin{aligned}
    T &= a_{1}xy\left( 2x^2y^2A_{k} + x^2y^2B_{k}\right)\\
        &+a_{1}xy\left( (ka_{1}xy + x^2y^2A_{k})^2\right)\\
        &+a_{1}xy\left( 2(ka_{1}xy + x^2y^2A_{k})(kb_{1}xy + x^2y^2B_{k})\right)\\
        &+a_{1}xy\left( (ka_{1}xy + x^2y^2A_{k})^2 (kb_{1}xy + x^2y^2B_{k})\right).
    \end{aligned}
    \label{eq:formT}
\end{equation} We have
\begin{equation*}
    \begin{aligned}
    |A_{k+1}| &\leq T + |A_{k}| + |ka_{1}(2a_{1}+b_{1})|\\
              &+ |A_{1}(T_{0}^{k})(1+ka_{1}xy+x^2y^2A_{k})^3(1+kb_{1}xy+x^2y^2B_{k})^2|\\
              &\leq T + 2k^2S^2 + 3kS^2 + S(1+S\delta+2S^2\delta^2).
    \end{aligned}
    \label{eq:Akp1}
\end{equation*}
Using $S(1+S\delta+2S^2\delta^2)<2S\leq2S^2$ we have
$$|A_{k+1}| \leq T + 2k^2S^2 + 3kS^2 + 2S^2.$$
The proof will be completed once we show that $|T| \leq kS^2$. So, next we need to extract the coefficients of $x^2y^2$ from $T$ as mentioned in equation \eqref{eq:Akp1}. We have $T = (T_{a1} + T_{a2} + T_{a3} + T_{a4})$ where 
\begin{equation}
    \label{}
    \begin{aligned}
    T_{a1} &= a_{1}xy(2x^2y^2A_{k} + x^2y^2B_{k}),\\
    T_{a2} &= a_{1}xy((ka_{1}xy + x^2y^2A_{k})^2),\\
    T_{a3} &= a_{1}xy(2(ka_{1}xy + x^2y^2A_{k})(kb_{1}xy + x^2y^2B_{k})),\\
    T_{a4} &= a_{1}xy(ka_{1}xy + x^2y^2A_{k})^2(kb_{1}xy + x^2y^2B_{k}),
    \end{aligned}
\end{equation}
and we also write $T_{a1} = x^2y^2T_{1}, T_{a2} = x^2y^2T_{2}, T_{a3} = x^2y^2T_{3}$, and $ T_{a4} = x^2y^2T_{4}$. We try to bound each of $T_{1}, T_{2}, T_{3}, T_{4}$ which are the coefficients of $x^2y^2$ in each of $T_{a1}, T_{a2}, T_{a3}, T_{a4}$.

For $T_{a1}$ we have,
$$T_{a1} = x^2y^2(2a_{1}xyA_{k} + a_{1}B_{k}xy) = x^2y^2T_{1}$$
where $|T_{1}| \leq 4S^3k^2xy + 2S^3k^2xy = 6S^3k^2xy$, so $|T_{1}| \leq 6S^3\delta k = 6S^2 Mk$.

For $T_{a2}$ we have,
$$T_{a2} = x^2y^2(k^2a_{1}^3xy + a_{1}A_{k}^2 x^3 y^3 + 4k^3s^4x^2y^2) = x^2y^2T_{2}$$
and $|T_{2}| \leq kS^3\delta + 4S^5\delta^3k + 4S^4k\delta^2$ so $|T_{2}| \leq kMS^2 + 4kS^2M^3 + 4S^2kM^2$.

For $T_{a3}$ we have,
$$T_{a3} = x^2y^2T_{3}$$
where $|T_{3}| \leq 2S^3k^2xy + 8k^3S^4x^2y^2 + 8k^4S^5x^3y^3$ and so $|T_{3}| \leq 2kMS^2 + 8kS^2M^2 + 8kS^2M^3$.

For $T_{a4}$ we have,
$$T_{a4} = x^2y^2 T_{4}$$
where $|T_{4}| \leq k^3S^4x^2y^2 + 12k^5S^6x^4y^4 + 6k^4S^5x^3y^3 + 4k^6S^7x^5y^5$ and $|T_{4}| \leq kS^2M^2 + 12kS^2M^4 + 6kS^2M^3 + 4kS^2M^5$. Since we have $T = x^2y^2(T_{1} + T_{2} + T_{3} + T_{4})$ and the bounds for $T_{1},T_{2},T_{3},T_{4}$, so $$|T_{1} + T_{2} + T_{3} + T_{4}| \leq |T_{1}| + |T_{2}| + |T_{3}| + |T_{4}| \leq kS^2(56M) = kS^2.$$ Thus we have $|A_{k+1}| \leq 2k^2S^2 + 4kS^2 + 2S^2$. Thus $|A_{k+1}| \leq 2(k+1)^2S^2$ as desired.
By applying the induction hypothesis we obtain
\begin{equation}
    |F_{k+1}(x,y)| \leq 2kR + |F(T_{0}^{k}(x,y))|\left( 1 + |xy|2kR\right)^{3}.
\end{equation}
Since $k|xy|$ is small and $(x,y) \approx (x^*,0)$, we can assume $T_{0}^{k}(x,y) \in \cN$ and so $|F\left(T_{0}^{k}(x,y)\right)| \leq R$. Also $1 + |xy|2kR \leq 2^{\frac{1}{3}}$, so
\begin{equation}
    |F_{k+1}(x,y)| \leq 2(k+1)R,
\end{equation}
as required. The result for $G_{k+1}$ is obtained similarly.
\end{proof}
%....................................................................................
%.......................
\begin{proof}[Proof of Theorem \eqref{th:necessaryConditions}]
For brevity we only provide details in the case $|\lambda \sigma| \le 1$.
Since $T_0$ and $T_1$ are both locally invertible the case $|\lambda \sigma| > 1$
can be accommodated by considering $f^{-1}$ in place of $f$.

%{\sc This proof needs to be completed.}
%..................................................
% A try
The trace of \eqref{eq:DfkPlusm} is
\begin{equation}
\tau_k = {\rm trace} \left( P_k {\rm D} T_1 \left( x^{(k)}_k, y^{(k)}_k \right) \right),
\nonumber
\end{equation}
where $P_k$ is given by \eqref{eq:Pjl} with $j = k$, and so

\begin{equation}
\tau_k = {\rm trace} \left( \begin{bmatrix} t_k & u_k \\ v_k & w_k \end{bmatrix}
\begin{bmatrix} c_1 + \cO \left( x^{(k)}_k, y^{(k)}_k - y^* \right)
& c_2 + \cO \left( x^{(k)}_k, y^{(k)}_k - y^* \right) \\
d_1 + \cO \left( x^{(k)}_k, y^{(k)}_k - y^* \right)
& d_2 + \cO \left( x^{(k)}_k, y^{(k)}_k - y^* \right)
\end{bmatrix} \right).
\label{eq:tauk4l}
\end{equation}
From Lemma \ref{le:PjBoundsl} we see that the leading order term in $\tau_k$ is $w_{k}d_{2}$. In particular, if $d_2 \ne 0$ then $|\tau_{k}| \rightarrow \infty$ as $k \rightarrow \infty$. But ${\rm SR}_k$-solutions are assumed to be stable for arbitrarily large values of $k$, so this is only possible if $d_2 = 0$. This establishes equation \eqref{eq:tangencyCondition}.
To verify \eqref{eq:det1Condition} and \eqref{eq:globalResonanceCondition} we
use \eqref{eq:T11} and \eqref{eq:T0kl1} to solve for fixed points of $T_0^k \circ T_1$ in order to derive the leading order component of $\left(x_{k}^{(k)},y_{k}^{(k)}\right)$. Since $\left(x_{k}^{(k)},y_{k}^{(k)}\right) \rightarrow \left( 0, y^{*}\right)$ as $k \rightarrow \infty$, it is convenient to write
\begin{align}
    \left( x_{k}^{(k)},y_{k}^{(k)}\right) = \left( \lambda^{k}(x^{*} + \Phi_{k}), y^{*}+\Psi_{k} \right).
    \label{eq:fpcomposed}
\end{align}
Notice $\Psi_k \to 0$ because $y^{(k)}_k \to y^*$. 
Since $y^{(k)}_0 = \cO \left( \sigma^{-k} \right)$ by Lemma \ref{le:xjyjBound},
the first component of \eqref{eq:T0kl1} gives
$x^{(k)}_k = \lambda^k x^{(k)}_0 + \cO \left( k \lambda^{k} \sigma^{-k} \right)$.
Since $x^{(k)}_0 \to 0$ as $k \to \infty$
we can conclude that in \eqref{eq:fpcomposed} we also have $\Phi_k \to 0$.

Next we use \eqref{eq:T11} and \eqref{eq:T0kl1}
to obtain formulas for $\Phi_k$ and $\Psi_k$
based on the knowledge that $\Phi_k, \Psi_k \to 0$ as $k \to \infty$.
By substituting \eqref{eq:fpcomposed} into \eqref{eq:T11} we obtain
\begin{align}
x^{(k)}_0 &= x^* + \cO \left( \lambda^k, \Psi_k \right),
\label{eq:xk0l} \\
y^{(k)}_0 &=  d_{1} x^{*} \lambda^k + d_{1} \Phi_k \lambda^k
+ d_3 x^{*^2} \lambda^{2 k} + d_4 x^* \Psi_k \lambda^k  + d_5 \Psi_k^2
\nonumber \\ &\quad+ \cO \left( \lambda^{3 k}, \Phi_k \lambda^{2 k}, \Phi_k \Psi_k \lambda^k, \Psi_k^3 \right).
\label{eq:yk0l}
\end{align}
We then substitute \eqref{eq:xk0l} and \eqref{eq:yk0l} into \eqref{eq:T0kl1},
noting that $k xy = \cO \left( k \sigma^{-k} \right)$ by \eqref{eq:xjyjBound} with $j = 0$, to obtain for the $y$-component
\begin{align}
y^{(k)}_k &= d_{1}x^*\lambda^{k}\sigma^{k} + \sigma^{k} \Bigl[ d_{1} \Phi_k \lambda^k 
+ d_{3}x^{*^2}\lambda^{2k} + d_4 x^* \Psi_k \lambda^k + d_5 \Psi_k^2
\nonumber \\ &\quad+ \cO \left(k \sigma^{-k}, \lambda^{3 k}, \Phi_k \lambda^{2 k}, \Phi_k \Psi_k \lambda^k, \Psi_k^3 \right) \Bigr].
\label{eq:ykkl}
\end{align}
We now match the two expressions for $y^{(k)}_k$, \eqref{eq:fpcomposed} and \eqref{eq:ykkl}.  If $y^{*} \neq \lambda^{k}\sigma^{k}d_{1}x^{*}$ then because $d_5 \ne 0$ the $\sigma^k d_5 \Psi_k^2$ term in \eqref{eq:ykkl} must balance the $y^*$ term in \eqref{eq:fpcomposed}, and so we must have $\Psi_{k} \sim t \sigma^{-\frac{k}{2}}$ for some $t \neq 0$. In this case the $(2,2)$-element of ${\rm D}T_{1}\left(x_{k}^{(k)},y_{k}^{(k)}\right)$ is asymptotic to $2d_{5}t \sigma^{-\frac{k}{2}}$.  Consequently from \eqref{eq:tauk4l} and Lemma \ref{le:PjBoundsl}, $\tau_k$ is asymptotic to $2 d_5 t \sigma^{\frac{k}{2}}$.  This diverges as $k \to \infty$ contradicting the stability assumption of the ${\rm SR}_k$-solutions.

Therefore we must have $y^{*} = d_{1} \lambda^{k}\sigma^{k} x^{*}$, which implies $|\lambda \sigma| = 1$ and $|d_{1}| = \frac{y^{*}}{x^{*}}$. Moreover, if $\lambda \sigma = 1$ then $d_1 = \frac{y^*}{x^*}$.
\end{proof}
%===============================================================================
\section{The orientation-preserving case}
\label{sec:positiveCase}
In this section we prove Theorems~\ref{th:a1b1} and \ref{th:sufficientConditionsPositiveCase}.

%.................................................................................

%.................................................................................
\begin{proof}[Proof of Theorem \ref{th:a1b1}]
We write
\begin{equation}
\left( x^{(k)}_k, y^{(k)}_k \right) = \left( \lambda^k (x^* + \Phi_k), y^* + \Psi_k \right),
\label{eq:xkyk4}
\end{equation}
where $\Phi_k, \Psi_k \to 0$ as $k \to \infty$ as shown in the proof of Theorem~\ref{th:necessaryConditions}.
By substituting \eqref{eq:xkyk4} into \eqref{eq:T11} with $d_1 = \frac{y^*}{x^*}$ and $d_2 = 0$ we obtain
\begin{align}
x^{(k)}_0 &= x^* + c_1 x^* \lambda^k + c_2 \Psi_k
+ \cO \left( \lambda^{2 k}, \lambda^k \Phi_k, \Psi_k^2 \right),
\label{eq:xk04} \\
y^{(k)}_0 &=  y^* \lambda^k + \tfrac{y^*}{x^*} \Phi_k  \lambda^k
+ d_3 x^{*^2} \lambda^{2 k} + d_4 x^* \lambda^k \Psi_k + d_5 \Psi_k^2
+ \cO \left( \lambda^{3 k}, \lambda^{2 k} \Phi_k, \lambda^k \Phi_k \Psi_k,\Psi_k^3 \right).
\label{eq:yk04}
\end{align}
From Lemma \ref{le:T0koripv} we have
\begin{equation}
T_0^k(x,y) = \begin{bmatrix}
\lambda^k x \left( 1 + k a_1 x y + \cO \left( k^2 x^2 y^2 \right) \right) \\
\frac{1}{\lambda^k} y \left( 1 + k b_1 x y + \cO \left( k^2 x^2 y^2 \right) \right)
\end{bmatrix}.
\label{eq:T0k4}
\end{equation}
We then substitute \eqref{eq:xk04} and \eqref{eq:yk04} into \eqref{eq:T0k4},
noting that $k^2 x^2 y^2 = \cO \left( k^2 \lambda^{2 k} \right)$, to obtain
\begin{align}
x^{(k)}_k &=  x^* \lambda^k + a_1 x^{*^2} y^* k \lambda^{2 k}
+ \cO \left( \lambda^{2 k}, k \lambda^{2 k} \Phi_k, \lambda^k \Psi_k, k \lambda^k \Psi_k^2 \right),
\label{eq:xkk4} \\
y^{(k)}_k &= y^* + \lambda^{-k} \left[ \tfrac{y^*}{x^*} \lambda^k \Phi_k
+ b_1 x^* y^{*^2} k \lambda^{2 k} + d_5 \Psi_k^2
+ \cO \left( \lambda^{2 k}, k \lambda^{2 k} \Phi_k, \lambda^k \Psi_k, k \lambda^k \Psi_k^2, \Psi_k^3 \right) \right],
\label{eq:ykk4}
\end{align}
where we have only explicitly written the terms that will be important below.
By matching \eqref{eq:xkyk4} to \eqref{eq:ykk4} we obtain
\begin{equation}
\lambda^k \Psi_k = \tfrac{y^*}{x^*} \lambda^k \Phi_k
+ b_1 x^* y^{*^2} k \lambda^{2 k} + d_5 \Psi_k^2
+ \cO \left( \lambda^{2 k}, k \lambda^{2 k} \Phi_k, \lambda^k \Psi_k, k \lambda^k \Psi_k^2, \Psi_k^3 \right).
\label{eq:Psik4}
\end{equation}
By matching \eqref{eq:xkyk4} to \eqref{eq:xkk4} we obtain
a similar expression for $\lambda^k \Phi_k$ which we substitute into \eqref{eq:Psik4} to obtain
\begin{equation}
\lambda^k \Psi_k = (a_1 + b_1) x^* y^{*^2} k \lambda^{2 k} + d_5 \Psi_k^2
+ \cO \left( \lambda^{2 k}, k \lambda^{2 k} \Phi_k, \lambda^k \Psi_k, k \lambda^k \Psi_k^2, \Psi_k^3 \right).
\label{eq:Psik5}
\end{equation}
Notice that $\lambda^k \Psi_k \to 0$ faster than either $k \lambda^{2 k}$ or $\Psi_k^2$
(but possibly not both)
and by inspection the same is true for every error term in \eqref{eq:Psik5}.

We now show that $a_{1} + b_{1} = 0$. Suppose for a contradiction that $a_1 + b_1 \ne 0$.
Then since $d_5 \ne 0$ the
$k \lambda^{2 k}$ and $\Psi_k^2$ terms in \eqref{eq:Psik5} must balance.
Thus $\Psi_k \sim \zeta \sqrt{k} \lambda^k$, for some $\zeta \ne 0$, and~\eqref{eq:tauk4l} becomes
\begin{equation}
\tau_k = {\rm trace} \left( \begin{bmatrix} t_j & u_j \\ v_j & w_j \end{bmatrix}
\begin{bmatrix} c_1 + \cO \left( \sqrt{k} \lambda^k \right)
& c_2 + \cO \left( \sqrt{k} \lambda^k \right) \\
d_1 + \cO \left( \sqrt{k} \lambda^k \right)
& 2 d_5 \Psi_k + \cO \left( \lambda^k \right)
\end{bmatrix} \right).
\label{eq:tauk4}
\end{equation}
By Lemma \ref{le:PjBoundsl} with $j = k$
the term involving $\Psi_k$ provides the leading order contribution to \eqref{eq:tauk4},
specifically
\begin{equation}
\tau_k \sim 2 d_5 \zeta \sqrt{k}.
\label{eq:tauk}
\end{equation}
Thus $\tau_k \to \infty$ as $k \to \infty$
and so the ${\rm SR}_k$-solutions are unstable for sufficiently values of $k$.
This contradicts the stability assumption in the theorem statement,
therefore $a_1 + b_1 = 0$.
\end{proof}

%.................................................................................
%.................................................................................
\begin{proof}[Proof of Theorem \ref{th:sufficientConditionsPositiveCase}]
We look for ${\rm SR}_k$-solutions for which the $k^{\rm th}$ point has the form
\begin{equation}
\left( x^{(k)}_k, y^{(k)}_k \right) = \Big( \lambda^k \left( x^* + \phi_k \lambda^{k} + \cO \left( \lambda^{2 k} \right) \right),
y^* + \psi_k \lambda^k + \cO \left( \lambda^{2 k} \right) \Big),
\label{eq:xkyk3}
\end{equation}
where $\phi_k, \psi_k \in \mathbb{R}$. Recall $\left( x^{(k)}_k, y^{(k)}_k \right)$ is a fixed point of $T_0^k \circ T_1$, so is equal to its image under \eqref{eq:T11} and \eqref{eq:T0k4}. Through matching \eqref{eq:xkyk3} to this image we obtain
\begin{equation}
\begin{split}
\phi_k &= c_1 x^* + c_2 \psi_k + a_1 x^{*^2} y^* k, \\
\psi_k &= \frac{y^* \phi_k}{x^*} + d_3 x^{*^2} + d_4 x^* \psi_k + d_5 \psi_k^2 + b_1 x^* y^{*^2} k.
\end{split}
\label{eq:pkqk}
\end{equation}
By solving \eqref{eq:pkqk} simultaneously for $\phi_k$ and $\psi_k$ we find that the terms involving $k$ cancel because $a_1 + b_1 = 0$ and there are two solutions. The values of $\psi_k$ for these are given by
\begin{equation}
\psi_k^\pm = \frac{1}{2 d_5} \left( 1 - \frac{c_2 y^*}{x^*} - d_4 x^* \pm \sqrt{\Delta} \right).
\label{eq:qk}
\end{equation}
It is readily seen that these correspond to ${\rm SR}_k$-solutions (for some fixed $\eta > 0$)
for sufficiently large values of $k$ by Lemma \ref{le:xjyjBound}.

We now investigate the stability of the two solutions.
With \eqref{eq:xkyk3} equation~\eqref{eq:tauk4l} becomes
\begin{equation}
\tau_k = {\rm trace} \left(
\begin{bmatrix} t_k & u_k \\ v_k & w_k \end{bmatrix}
\begin{bmatrix} c_1 + \cO \left( \lambda^k \right) & c_2 + \cO \left( \lambda^k \right) \\
d_1 + \cO \left( \lambda^k \right) & d_4 x^* \lambda^k + 2 d_5 \psi_k \lambda^k + \cO \left( \lambda^{2 k} \right) \end{bmatrix} \right).
\label{eq:taukAgain}
\end{equation}
Thus by Lemma \ref{le:PjBoundsl} we have
\begin{align}
\lim_{k \to \infty} \tau_k = d_4 x^* + 2 d_5 \psi_k^\pm
= 1 - \frac{c_2 y^*}{x^*} \pm \sqrt{\Delta} \,.
\label{eq:tauk2}
\end{align}
Also the determinant of ${\rm D} \left( T_0^k \circ T_1 \right)(x_k,y_k)$ converges to
\begin{align}
\lim_{k \to \infty} \delta_k = c_1 d_2 - c_2 d_1
= -\frac{c_2 y^*}{x^*}.
\label{eq:deltak2}
\end{align}
To show that $\psi_k^-$ generates an asymptotically stable ${\rm SR}_k$-solution we verify
(i) $\delta_k - \tau_k + 1 > 0$,
(ii) $\delta_k + \tau_k + 1 > 0$, and 
(iii) $\delta_k < 1$,
for sufficiently large values of $k$ (see Fig.~\ref{fig:stability_triangle}).
From \eqref{eq:tauk2} and \eqref{eq:deltak2} with $\psi_k^-$ we have
\begin{align}
\lim_{k \to \infty} (\delta_k - \tau_k + 1) &= \sqrt{\Delta}, \label{eq:sqrtDelta} \\
\lim_{k \to \infty} (\delta_k + \tau_k + 1) &= 2 - \frac{2 c_2 y^*}{x^*} - \sqrt{\Delta}, \\
\lim_{k \to \infty} (1 - \delta_k) &= 1 + \frac{c_2 y^*}{x^*}.
\end{align}
These limits are all positive by \eqref{eq:stabilityCondition},
hence conditions (i)--(iii) (given just above Equation \eqref{eq:sqrtDelta})
are satisfied for sufficiently large values of $k$.
Finally observe that with instead $\psi^+_k$ we have
$\lim_{k \to \infty} (\delta_k - \tau_k + 1) = -\sqrt{\Delta} < 0$
and $\lim_{k \to \infty} |\delta_k| = \frac{|c_2| y^*}{x^*} < 1$, hence $\psi^+_k$ generates a saddle ${\rm SR}_k$-solution.
\end{proof}

%===============================================================================
\section{The orientation-reversing case}
\label{sec:negativeCase}
%\setcounter{equation}{0}

%.................................................................................

%.................................................................................

\begin{proof}[Proof of Theorem \ref{th:sufficientConditionsNegativeCase}]
As in the proof of Theorem \ref{th:sufficientConditionsPositiveCase} we assume $\left( x^{(k)}_k, y^{(k)}_k \right)$ has the form \eqref{eq:xkyk3}.  This point is a fixed point of $T_0^k \circ T_1$ where
$T_0^k$ again has the form \eqref{eq:T0k4} except now $a_1 = b_1 = 0$.
By composing this with \eqref{eq:T11} we obtain
\begin{equation}
\left( T_0^k \circ T_1 \right)\left( x^{(k)}_k, y^{(k)}_k \right) =
\begin{bmatrix}
x^* \lambda^k + \left( c_1 x^* + c_2 \psi_k \right) \lambda^{2 k} + \cO \left( \lambda^{3 k} \right) \\
(-1)^k d_1 x^* + (-1)^k \left( d_1 \phi_k + d_3 x^{*^2} + d_4 x^* \psi_k + d_5 \psi_k^2 \right) \lambda^k + \cO \left( \lambda^{2 k} \right)
\end{bmatrix},
\label{eq:T0kT1xkyk}
\end{equation}
where we have substituted $\lambda \sigma = -1$ and $d_2 = 0$.
For the remainder of the proof we assume $k$ is even in the case $d_1 = \frac{y^*}{x^*}$
and $k$ is odd in the case $d_1 = -\frac{y^*}{x^*}$.
In either case
\begin{equation}
(-1)^k d_1 = \frac{y^*}{x^*},
\label{eq:d1Parity}
\end{equation}
and so the leading-order terms of \eqref{eq:xkyk3} and \eqref{eq:T0kT1xkyk} are the same.
By matching the next order terms we obtain
\begin{equation}
\begin{split}
\phi_k &= c_1 x^* + c_2 \psi_k \,, \\
\psi_k &= \frac{y^* \phi_k}{x^*} + (-1)^k \left( d_3 x^{*^2} + d_4 x^* \psi_k + d_5 \psi_k^2 \right).
\end{split}
\label{eq:pkqkNegativeCase}
\end{equation}
These produce the following two solutions for the value of $\psi_k$
\begin{equation}
\psi_k^\pm = \frac{(-1)^k}{2 d_5} \left( 1 - \frac{c_2 y^*}{x^*} - \frac{d_4 y^*}{d_1} \pm \sqrt{\Delta} \right),
\label{eq:qkNegativeCase}
\end{equation}
where we have further used \eqref{eq:d1Parity}.
%(and with $\phi_k$ given by the first equation in \eqref{eq:pkqkNegativeCase}).
Analogous to the proof of Theorem \ref{th:sufficientConditionsPositiveCase} we obtain
\begin{align}
\lim_{k \to \infty} \tau_k &= (-1)^k \left( d_4 x^* + 2 d_5 \psi_k^\pm \right) 
\label{eq:limtau}\\
&= 1 - \frac{c_2 y^*}{x^*} \pm \sqrt{\Delta},
\nonumber \\
\lim_{k \to \infty} \delta_k &= -\frac{c_2 y^*}{x^*},
\nonumber
\end{align}
and the proof is completed via the same stability arguments.
\end{proof}
%===============================================================================
\section{Conclusions}
\label{sec:conc}
In this chapter we have considered single-round periodic solutions
associated with homoclinic tangencies of two-dimensional $C^\infty$ maps.
We have formalised these as ${\rm SR}_k$-solutions via Definition \ref{df:SRkSolution}.
The key arguments leading to our results are centred around calculations of $\tau_k$ ---
the sum of the eigenvalues associated with an ${\rm SR}_k$-solution.
Immediately we see from Fig.~\ref{fig:stability_triangle} that if $|\tau_k| > 2$ then the ${\rm SR}_k$-solution is unstable.

We first showed that if conditions \eqref{eq:tangencyCondition}--\eqref{eq:globalResonanceCondition} do not all hold then $|\tau_k| \to \infty$ as $k \to \infty$,
thus at most finitely many ${\rm SR}_k$-solutions can be stable (Theorem \ref{th:necessaryConditions}).
Equation \eqref{eq:det1Condition}, namely $|\lambda \sigma| = 1$, splits into two fundamentally distinct cases.
If $\lambda \sigma = -1$
then the resonant terms in $T_0$ that cannot be eliminated by a coordinate change
are of sufficiently high order that they have no bearing on the results.
In this case $\tau_k$ converges to the finite value \eqref{eq:limtau} and
infinitely many ${\rm SR}_k$-solutions can indeed be stable
as a codimension-three phenomenon (Theorem \ref{th:sufficientConditionsNegativeCase}).
If instead $\lambda \sigma = 1$ then
$|\tau_k|$ is asymptotically proportional to $\sqrt{k}$ unless $a_1 + b_1 = 0$,
that is unless the coefficients of the leading-order resonance terms cancel (Theorem \ref{th:a1b1}).
This is the only additional condition needed to have infinitely many ${\rm SR}_k$-solutions,
aside from the inequalities $\Delta > 0$ and \eqref{eq:stabilityCondition},
thus in this case the infinite coexistence is codimension-four (Theorem \ref{th:sufficientConditionsPositiveCase}).

We will compute and illustrate the infinite coexistence phenomenon in a specific  family of maps in \S\ref{sec:example}.

\chapter{Unfolding globally resonant homoclinic tangencies\label{cha:UnfGRHT}}
\graphicspath{ {./images/} }
\newcommand{\bn}{{\bf n}}
\newcommand{\bO}{{\bf 0}}
\newcommand{\bv}{{\bf v}}
\newcommand{\rD}{{\rm D}}
% To do:
%		- add figures
%		- Sishu to draft intro and concluding section
%		- edit throughout, obtain complete first draft
%		- get Robert to read it
%		- where to send it: DCDS, IJBC, J Diff Eq Appl?

% MSC codes:
% 	37G25 -- Bifurcations connected with nontransversal intersection
%		37G15 -- Bifurcations of limit cycles and periodic orbits
%		39A23 -- Periodic solutions
%===============================================================================
In Chapter \ref{cha:GRHT} we determined what degeneracies are needed \blue{for the coexistence of infinitely many stable single-round periodic solutions.} In this chapter we unfold about such a homoclinic tangency. The results in this chapter are in preparation to be published \citep{MuMcSi21b}. If a family of maps $f_\mu$, with $\mu \in \mathbb{R}^n$, has this phenomenon at some point $\mu^*$ in parameter space, then for any positive $j \in \mathbb{Z}$ there exists an open set containing $\mu^*$ in which $j$ asymptotically stable, single-round periodic solutions coexist.  Due to the high codimension, a precise description of the shape of these sets (for large $j$) is beyond the scope of this thesis.  The approach we take here is to consider one-parameter families that perturb from a globally resonant homoclinic tangency.  Some information about the size and shapes of the sets can then be inferred from our results. Globally resonant homoclinic tangencies are hubs for extreme multi-stability.  They should occur generically in some families of maps with three or more parameters, such as the
 generalised H\'enon map \citep{GoKu05,KuMe19}, but to our knowledge they \blue{are yet to be} identified.

We find that as the value of $\mu$ is varied from $\mu^*$, there occurs an infinite sequence of either saddle-node or period-doubling bifurcations that destroy the periodic solutions or make them unstable. Generically these sequences converge exponentially to $\mu^*$ with the distance (in parameter space) to the bifurcation asymptotically proportional to $|\lambda|^{2 k}$, where the periodic solutions have period $k + m$, for some fixed $m \ge 1$. \blue{In the parameter space, the codimension-one surface of homoclinic tangencies can be thought of as where a scalar function $h$ is zero. As parameters are varied, the value of $h$ changes and we cross the surface when $h=0$. If the derivative $h'$ at the crossing value is non-zero, we say that we have a ``linear change'' to the codimension-one condition and this corresponds to a generic transverse intersection}. If we move away from $\mu^*$ without a linear change to the codimension-one condition of a homoclinic tangency, the bifurcation values instead generically scale like $\frac{|\lambda|^k}{k}$. If the perturbation suffers further degeneracies, the scaling can be slower. Specifically we observe $|\lambda|^k$ and $\frac{1}{k}$ for an abstract example that we believe is representative of how the bifurcations scale in general.

Similar results have been obtained for more restrictive classes of maps. For area-preserving families the phenomenon is codimension-two and there exist infinitely many elliptic, single-round periodic solutions \citep{GoSh05}. As shown in \cite{DeGo15,GoGo09} the periodic solutions are destroyed or lose stability in bifurcations that scale like $|\lambda|^{2 k}$, matching our result. For piecewise-linear families the phenomenon is codimension-three \citep{Si14,DSH17}.  In this setting the bifurcation values instead scale like $|\lambda|^k$ \citep{Si14b}, see also \citep{DoLa08}.
%===============================================================================
\section{A quantitative description for the dynamics near a homoclinic connection}
\label{sec:localCoordinates}

\begin{figure}
\begin{center}
\includegraphics[width=6cm]{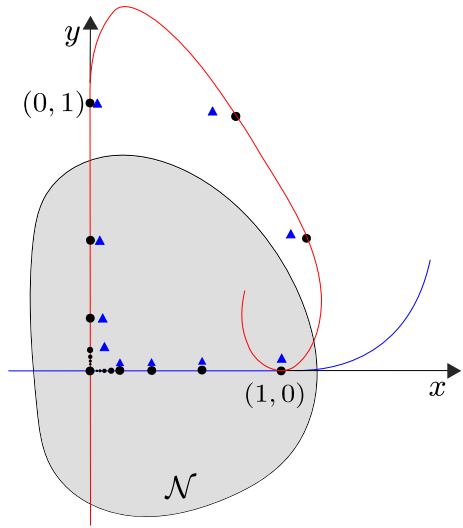}
\caption[A homoclinic tangency in a two-dimensional map with homoclinic orbit.]{A homoclinic tangency for a saddle fixed point of a two-dimensional map.
In this illustration the eigenvalues associated with the fixed points are positive, i.e.
$0 < \lambda < 1$ and $\sigma > 1$. A coordinate change has been applied so that in the region
$\cN$ (shaded) the coordinate axes coincide with the stable and unstable manifolds. The
homoclinic orbit $\Gamma_{\rm HC}$ is shown with black dots. A typical single-round periodic solution is shown with blue triangles.
\label{fig:TangUnf2}
} 
\end{center}
\end{figure}

Let $f$ be a $C^\infty$ map on $\mathbb{R}^2$.
Suppose the origin $(x,y) = (0,0)$ is a \blue{hyperbolic} saddle fixed point of $f$.
\blue{Further suppose $D f(0,0)$ does not have a zero eigenvalue, thus its} eigenvalues $\lambda, \sigma \in \mathbb{R}$ \blue{satisfy}
\begin{equation}
0 < |\lambda| < 1 < |\sigma|.
\label{eq:eigAssumption}
\end{equation}
By Theorem \ref{Th:StLinear}
there exists a $C^\infty$ coordinate change that transforms $f$ to
\begin{equation}
T_0(x,y) = \begin{bmatrix}
\lambda x \left( 1 + \cO \left( x y \right) \right) \\
\sigma y \left( 1 + \cO \left( x y \right) \right)
\end{bmatrix}.
\label{eq:T0general}
\end{equation}
In these new coordinates let $\cN$ be a convex neighbourhood of the origin for which
\begin{equation}
f(x,y) = T_0(x,y), \qquad \text{for all $(x,y) \in \cN$},
\label{eq:fEqualsT0}
\end{equation}
see Fig.~\ref{fig:TangUnf2}.  If $\lambda^p \sigma^q \ne 1$ for all integers \blue{$p, q \ge -1$}, then the eigenvalues are said to be {\em non-resonant} and the coordinate change can be chosen so that $T_0$ is linear.  If not, then $T_0$ must contain resonant terms that cannot be eliminated by the coordinate change. As explained in \S\ref{sec:parameters}, if $\lambda \sigma = 1$ we can reach the form
\begin{equation}
T_0(x,y) = \begin{bmatrix}
\lambda x \left( 1 + a_1 x y + \cO \left( x^2 y^2 \right) \right) \\
\sigma y \left( 1 + b_1 x y + \cO \left( x^2 y^2 \right) \right)
\end{bmatrix},
\label{eq:T0}
\end{equation}
where $a_1, b_1 \in \mathbb{R}$.
If $\lambda \sigma = -1$ we can obtain \eqref{eq:T0} with $a_1 = b_1 = 0$.

Now suppose there exists an orbit homoclinic to the origin, $\Gamma_{\rm HC}$.
By scaling $x$ and $y$ we may assume that $(1,0)$ and $(0,1)$ are points on $\Gamma_{\rm HC}$ and
\begin{equation}
(1,0), (\lambda,0), \left( 0, \tfrac{1}{\sigma} \right), \left( 0, \tfrac{1}{\sigma^2} \right) \in \cN,
\label{eq:pointsInN}
\end{equation}
By assumption there exists $m \ge 1$ such that $f^m(0,1) = (1,0)$.
We let $T_1 = f^m$ and expand $T_1$ in a Taylor series centred at $(x,y) = (0,1)$:
\begin{equation}
T_1(x,y) = \begin{bmatrix}
c_0 + c_1 x + c_2 (y-1) + \cO \left( (x,y-1)^2 \right) \\
d_0 + d_1 x + d_2 (y-1) + d_3 x^2 + d_4 x(y-1) + d_5 (y-1)^2 + \cO \left( (x,y-1)^3 \right)
\end{bmatrix},
\label{eq:T1}
\end{equation}
where $c_0 = 1$ and $d_0 = 0$. In \eqref{eq:T1} we have written explicitly the terms that will be important below.

%===============================================================================

%---------------------------------------------------------------------------------

%===============================================================================
\section{Smooth parameter dependence}
\label{sec:parameters}

Now suppose $f_\mu$ is a $C^\infty$ map on $\mathbb{R}^2$
with a $C^\infty$ dependence on a parameter $\mu \in \mathbb{R}^n$.
Let $\bO \in \mathbb{R}^n$ denote the origin in parameter space.
Suppose that for all $\mu$ in some region containing $\bO$,
the origin in phase space $(x,y) = (0,0)$ is a fixed point of $f_\mu$.
Let $\lambda = \lambda(\mu)$ and $\sigma = \sigma(\mu)$ be its associated eigenvalues
(these are $C^\infty$ functions of $\mu$) and suppose
\begin{align}
\lambda(\bO) &= \alpha, \label{eq:lambda} \\
\sigma(\bO) &= \frac{\chi_{\rm eig}}{\alpha}, \label{eq:sigma}
\end{align}
with $0 < \alpha < 1$ and $\chi_{\rm eig} \in \{ -1, 1 \}$.
With $\mu = \bO$ we have $|\lambda \sigma| = 1$, so as described above $T_0$ can be assumed to have the form \eqref{eq:T0}.
We now show we can assume $T_0$ has this form when the value of $\mu$ is sufficiently small.

%---------------------------------------------------------------------------------
\begin{lemma}
There exists a neighbourhood $\cN_{\rm param} \subset \mathbb{R}^n$ of $\bO$
and a $C^\infty$ coordinate change that puts $f_\mu$ in the form \eqref{eq:T0} for all $\mu \in \cN_{\rm param}$.
\label{le:T0}
\end{lemma}

%---------------------------------------------------------------------------------
\begin{proof}
Via a linear transformation $f_\mu$ can be transformed to
\begin{equation}
\begin{bmatrix} x \\ y \end{bmatrix} \mapsto
\begin{bmatrix}
\lambda x + \sum\limits_{i \ge 0, \,j \ge 0, \,i+j \ge 2} a_{ij} x^i y^j \\
\sigma y + \sum\limits_{i \ge 0, \,j \ge 0, \,i+j \ge 2} b_{ij} x^i y^j
\nonumber
\end{bmatrix},
\end{equation}
for some $a_{ij}, b_{ij} \in \mathbb{R}$.
It is a standard asymptotic matching exercise
to show that via an additional $C^\infty$ coordinate change we can achieve $a_{ij} = 0$ if $\lambda^{i-1} \sigma^j \ne 1$,
and $b_{ij} = 0$ if $\lambda^i \sigma^{j-1} \ne 1$. The remainder of the proof is based on this fact.

Assume the value of $\mu$ is small enough that \eqref{eq:eigAssumption} is satisfied.
Then $\lambda^p \sigma^q = 1$ is only possible with $p, q \ge 1$,
so a $C^\infty$ coordinate change can be performed to reduce the map to
\begin{equation}
\begin{bmatrix} x \\ y \end{bmatrix} \mapsto
\begin{bmatrix}
\lambda x + \sum\limits_{i \ge 2, \,j \ge 1} a_{ij} x^i y^j \\
\sigma y + \sum\limits_{i \ge 1, \,j \ge 2} b_{ij} x^i y^j
\nonumber
\end{bmatrix}.
\end{equation}
Since $|\lambda \sigma| = 1$ when $\mu = \bO$
we can assume $\mu$ is small enough that $\lambda^{p-1} \sigma \ne 1$ for all $p \ge 3$
and $\lambda \sigma^{q-1} \ne 1$ for all $q \ge 3$.
Consequently the map can further be reduced to
\begin{equation}
\begin{bmatrix} x \\ y \end{bmatrix} \mapsto
\begin{bmatrix}
\lambda x + a_{21} x^2 y + \sum\limits_{i \ge 3, \,j \ge 2} a_{ij} x^i y^j \\
\sigma y + b_{12} x y^2 + \sum\limits_{i \ge 2, \,j \ge 3} b_{ij} x^i y^j
\nonumber
\end{bmatrix},
\end{equation}
which can be rewritten as \eqref{eq:T0}.
\end{proof}

The product of the eigenvalues is
\begin{equation}
\lambda(\mu) \sigma(\mu) = \chi_{\rm eig} + \bn_{\rm eig}^{\sf T} \mu + \cO \left( \| \mu \|^2 \right),
\label{eq:nEig}
\end{equation}
where $\bn_{\rm eig}$ is the gradient of $\lambda \sigma$ evaluated at $\mu = \bO$ \blue{and the norm $ \|. \|$ introduced here is the Euclidean norm on $\mathbb{R}^{n}$. The same norm follows throughout the thesis}.
%\begin{equation}
%\bn_{\rm eig} = \nabla \big( \lambda(\mu) \sigma(\mu) \big) \big|_{\mu = \bO},
%\label{eq:nEig}
%\end{equation}
The following result is an elementary application of the implicit function theorem.

\begin{figure}
\begin{center}
\includegraphics[width=8cm]{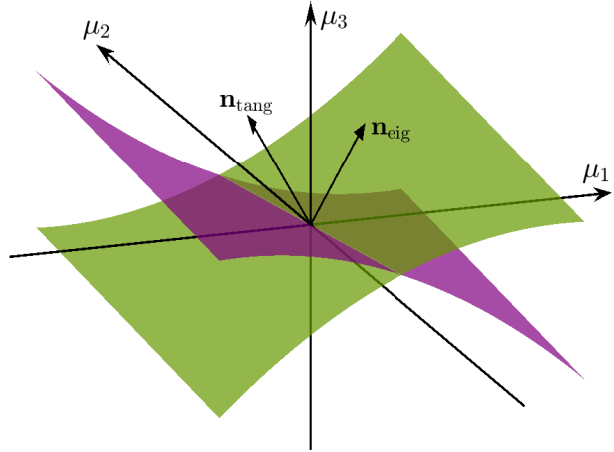}
\caption[A sketch of the codimension-one surfaces.]{A sketch of codimension-one surfaces of homoclinic tangencies (green) and where $\lambda(\mu) \sigma(\mu) = 1$ (purple).  The vectors ${\bf n}_{\rm tang}$ and ${\bf n}_{\rm eig}$, respectively, are normal to these surfaces at the origin $\mu = \bO$. 
\label{fig:Surf_Tang_Normal}
} 
\end{center}
\end{figure}

%---------------------------------------------------------------------------------
\begin{lemma}
Suppose $\bn_{\rm eig} \ne \bO$.
%In a neighbourhood of $\mu = \bO$
Then $|\lambda(\mu) \sigma(\mu)| = 1$ on a $C^\infty$ codimension-one surface
intersecting $\mu = \bO$ and with normal vector $\bn_{\rm eig}$ at $\mu = \bO$ (as illustrated in Fig.~\ref{fig:Surf_Tang_Normal}).
\label{le:areaPreservingSurface}
\end{lemma}

%%%%%%%%%%%%%%%%%%%%%%%%%%%%%%%%%%%%%%%%%%%%%%%%%%%%%%%%%%%%%%%%%%%%%%%%%%%%%%%%
%%%%%%%%%%%%%%%%%%%%%%%%%%%%%%%%%%%%%%%%%%%%%%%%%%%%%%%%%%%%%%%%%%%%%%%%%%%%%%%%

%{\color{red}
%FIG: sketch of three-dimensional parameter space:
%show axes $\mu_1$, $\mu_2$, and $\mu_3$;
%show a surface intersecting the origin and indicate and label its normal vector at the origin $\bn_{\rm eig}$;
%show another surface intersecting the origin and indicate and label its normal vector at the origin $\bn_{\rm tang}$.
%}

%===============================================================================
\section{The codimension-one surface of homoclinic tangencies}
\label{sec:tangency}
In this section we describe the codimension-one surface of homoclinic tangencies that intersects $\mu = \bO$ where we will be assuming that a globally resonant homoclinic tangency occurs. 
%For all $\mu$ in a sufficiently small neighbourhood of $\bO$ we can assume that $T_0$ and $T_1$
%take the form \eqref{eq:T0} and \eqref{eq:T1}
%but now the coefficients in these expressions are functions of $\mu$.
%Now suppose that when $\mu = \bO$ there exists an orbit homoclinic to the origin.
%Suppose $(1,0)$ and $(0,1)$ are points on this orbits and \eqref{eq:pointsInN} is satisfied.

Suppose \eqref{eq:pointsInN} is satisfied when $\mu = \bO$.
Write $f_\mu^m$ as \eqref{eq:T1} and suppose
\begin{align}
c_0(\bO) &= 1, \label{eq:c0} \\
d_0(\bO) &= 0, \label{eq:d0}
\end{align}
so that $f_\bO$ has an orbit homoclinic to the origin through $(x,y) = (1,0)$ and $(0,1)$.
Also suppose
\begin{align}
d_2(\bO) &= 0, \label{eq:d2} \\
d_5(\bO) &= d_{5,0} \ne 0, \label{eq:d5}
\end{align}
for a quadratic tangency.
Also write
\begin{equation}
d_0(\mu) = \bn_{\rm tang}^{\sf T} \mu + \cO \left( \| \mu \|^2 \right).
\label{eq:nTang}
\end{equation}

%---------------------------------------------------------------------------------
\begin{lemma}
Suppose \eqref{eq:pointsInN} is satisfied when $\mu = \bO$, \eqref{eq:c0}--\eqref{eq:d5} are satisfied, and  $\bn_{\rm tang} \ne \bO$.
Then $f_\mu$ has a quadratic homoclinic tangency to $(x,y) = (0,0)$ on a $C^\infty$ codimension-one surface intersecting $\mu = \bO$ and with normal vector $\bn_{\rm tang}$ at $\mu = \bO$.
Moreover, this tangency occurs at $(x,y) = (0,Y(\mu))$ where $Y$ is a $C^\infty$ function with $Y(\bO) = 1$.
\label{le:tangencySurface}
\end{lemma}

%---------------------------------------------------------------------------------
\begin{proof}
Let $T_{1,2}$ denote the second component of $T_1$ \eqref{eq:T1}.
The function
\begin{equation}
h(u;\mu) = \frac{\partial T_{1,2}}{\partial y}(0,1+u) = d_2(\mu) + 2 d_5(\mu) u + \cO(2),
\nonumber
\end{equation}
is a $C^\infty$ function of $u \in \mathbb{R}$ and $\mu$.
Since $h(0; \mathbf{0}) = 0$ by \eqref{eq:d2} and $\frac{\partial h}{\partial u}(0;\mathbf{0}) \ne 0$ by \eqref{eq:d5},
the implicit function theorem implies there exists a $C^\infty$ function $u_{\rm tang}(\mu)$ such that
$h(u_{\rm tang}(\mu);\mu) = 0$ locally.

By construction, the unstable manifold of $(x,y) = (0,0)$ is tangent to the $x$-axis at $T_1(0,1+u_{\rm tang})$.
Moreover this tangency is quadratic by \eqref{eq:d5}.
Thus a homoclinic tangency occurs if
$T_{1,2}(0,1+u_{\rm tang}) = d_0 + d_2 u_{\rm tang} + d_5 u_{\rm tang}^2 + \cO(3) = 0$.
This function is $C^\infty$ and
\begin{equation}
T_{1,2}(0,1+u_{\rm tang}) = \bn_{\rm tang}^{\sf T} \mu + \cO \left( \| \mu \|^2 \right).
\nonumber
\end{equation}
Since $\bn_{\rm tang} \ne \bO$ the result follows from another application of the implicit function theorem.
\end{proof}

%===============================================================================
\section{Sequences of saddle-node and period-doubling bifurcations}
\label{sec:mainResult}

Suppose
\begin{equation}
d_1(\bO) = \chi,
\label{eq:d1}
\end{equation}
where $\chi \in \{ -1, 1 \}$.
Write
\begin{equation}
\begin{aligned}
a_1(\bO) &= a_{1,0} \,,
&\qquad c_1(\bO) &= c_{1,0} \,,
&\qquad d_3(\bO) &= d_{3,0} \,, \\
b_1(\bO) &= b_{1,0} \,,
&\qquad c_2(\bO) &= c_{2,0} \,,
&\qquad d_4(\bO) &= d_{4,0} \,,
\end{aligned}
\label{eq:other}
\end{equation}
and, recalling \eqref{eq:discriminant}, let
\begin{equation}
\Delta_0 = \left( 1 - c_{2,0} - \chi d_{4,0} \right)^2 - 4 d_{5,0} (d_{3,0} + \chi c_{1,0}).
\label{eq:Delta0}
\end{equation}

%---------------------------------------------------------------------------------
Given $k_{\rm min} \in \mathbb{Z}$, let
\begin{equation}
K(k_{\rm min}) = \left\{ k \in \mathbb{Z} \,\big|\, k \ge k_{\rm min}, \left( \chi_{\rm eig} \right)^k = \chi \right\}.
\label{eq:K}
\end{equation}
If $\chi_{\rm eig} = \chi = 1$, then $K$ is the set of all integers greater than or equal to $k_{\rm min}$.
If $\chi_{\rm eig} = -1$ and $\chi = 1$ [resp.~$\chi = -1$],
then $K$ is the set of all even [resp.~odd] integers greater than or equal to $k_{\rm min}$.
\begin{theorem}
Suppose $f_\mu$ satisfies \eqref{eq:lambda}, \eqref{eq:sigma}, \eqref{eq:c0}--\eqref{eq:d5}, $a_{1,0} + b_{1,0} = 0$, $\Delta_0 > 0$, and $-1 < c_{2,0} < 1 - \frac{\sqrt{\Delta_0}}{2}$.
Let $\bv \in \mathbb{R}^n$.
If $\bn_{\rm tang}^{\sf T} \bv \ne 0$
then there exists $k_{\rm min} \in \mathbb{Z}$ such that for all $k \in K(k_{\rm min})$
there exist $\ee_k^- <0$ and $\ee_k^+ > 0$ with $\ee_k^\pm = \cO \left( \alpha^{2 k} \right)$
such that $f$ has an asymptotically stable period-$(k+m)$ solution for all $\mu = \ee \bv$ with $\ee_k^- < \ee < \ee_k^+$.
Moreover, one sequence, $\{ \ee_k^- \}$ or $\{ \ee_k^+ \}$, corresponds to saddle-node bifurcations of the periodic solutions, the other to period-doubling bifurcations.
If instead $\bn_{\rm tang}^{\sf T} \bv = 0$ and $\bn_{\rm eig}^{\sf T} \bv \ne 0$ the same results hold \blue{except} now $\ee_k^\pm = \cO \left( \frac{\alpha^k}{k} \right)$.
\label{th:main}
\end{theorem}

%===============================================================================
\section{Proof of the main result}
\label{sec:mainProof}

To prove Theorem \ref{th:main} we use the following lemma by carefully estimating the error terms in \blue{\eqref{eq:T0kun}}.
If $\chi_{\rm eig} = -1$ then \eqref{eq:T0kun} is true if $a_1 = b_1 = 0$
(and this can be proved in the same fashion).

%---------------------------------------------------------------------------------
\begin{lemma}
Suppose $\chi_{\rm eig} = 1$ and $\mu = \cO \left( \alpha^k \right)$.
If $|x-1|$ and \blue{$\left| \alpha^{-k} y  \right|$} are sufficiently small
for all sufficiently large values of $k$, then
\begin{equation}
T_0^k(x,y) = \begin{bmatrix}
\lambda^k x \left( 1 + k a_1 x y + \cO \left( k^2 \alpha^{2 k} \right) \right) \\
\sigma^k y \left( 1 + k b_1 x y + \cO \left( k^2 \alpha^{2 k} \right) \right)
\end{bmatrix}.
\label{eq:T0kun}
\end{equation}
\label{le:T0kun}
\end{lemma}
\begin{proof}
Write
\begin{equation}
T_0(x,y) = \begin{bmatrix}
\lambda x \left( 1 + a_1 x y + x^2 y^2 \tilde{F}(x,y) \right) \\
\sigma y \left( 1 + b_1 x y + x^2 y^2 \tilde{G}(x,y) \right)
\end{bmatrix}.
\label{eq:T02un}
\end{equation}
Let $R > 0$ be such that
\begin{equation}
|a_1|, |b_1|, |\tilde{F}(x,y)|, |\tilde{G}(x,y)| \le R
\label{eq:boundR}
\end{equation}
for all $(x,y) \in \cN$ and all sufficiently small values of $\mu$. For simplicity we assume $\alpha > 0$; if instead $\alpha < 0$ the proof can be completed in the same fashion.

We have $\lambda = \alpha + \cO \left( |\alpha|^k \right)$
and $\sigma = \frac{1}{\alpha} + \cO \left(|\alpha|^k \right)$.
Thus there exists $M \ge 2 R$ such that
$\lambda \le \alpha \left( 1 + M \alpha^k \right)$ and
$\sigma \le \frac{1}{\alpha} \left( 1 + M \alpha^k \right)$
for sufficiently large values of $k$.
%By using $M \alpha^k \le \frac{1}{2 k}$ (true for sufficiently large $k$)
%it is straight-forward to show (by induction on $j$) that
It follows (by induction on $j$) that
\begin{equation}
\sigma^j \le \alpha^{-j} \left( 1 + 2 M j \alpha^k \right),
\label{eq:T0kproof11}
\end{equation}
and
\begin{equation}
\left( \lambda \sigma \right)^j \le 1 + 4 M j \alpha^k,
\label{eq:T0kproof12}
\end{equation}
for all $j = 1,2,\ldots,k$
again assuming $k$ is sufficiently large.

Let $\ee > 0$.
Assume\blue{
\begin{align}
|x-1| &\le \ee, &
\left| \alpha^{-k} y  \right| &\le \ee,
\label{eq:T0kproof20}
\end{align}}
for sufficiently large values of $k$.
We can assume $\ee < 1 - \frac{1}{\sqrt{2}}$ so then
\begin{equation}
\frac{\alpha^k}{2} \le |x y| \le 2 \alpha^k.
\label{eq:T0kproof21}
\end{equation}

Write
\begin{equation}
T_0^j(x,y) = \begin{bmatrix}
\lambda^j x \left( 1 + j a_1 x y + x^2 y^2 \tilde{F}_j(x,y) \right) \\
\sigma^j y \left( 1 + j b_1 x y + x^2 y^2 \tilde{G}_j(x,y) \right)
\end{bmatrix}.
\label{eq:T0j}
\end{equation}
Below we will use induction on $j$ to show that
\begin{equation}
\left| \tilde{F}_j(x,y) \right|, \left| \tilde{G}_j(x,y) \right| \le 12 M j^2,
\label{eq:T0kproof30}
\end{equation}
for all $j = 1,2,\ldots,k$, assuming $k$ is sufficiently large.
This will complete the proof because with $j=k$, \eqref{eq:T0kproof30} implies \eqref{eq:T0kun}.

Clearly \eqref{eq:T0kproof30} is true for $j = 1$:
$\left| \tilde{F}_1(x,y) \right| = \left| \tilde{F}(x,y) \right| \le R \le \frac{M}{2} < 12 M$
and similarly for $\tilde{G}_1$.
Suppose \eqref{eq:T0kproof30} is true for some $j < k$.
It remains for us to verify \eqref{eq:T0kproof30} for $j+1$.
First observe that by using $|a_1| \le R$, \eqref{eq:T0kproof21}, and the induction hypothesis,
\begin{equation}
\left| 1 + j a_1 x y + x^2 y^2 \tilde{F}_j(x,y) \right|
\le 1 + 2 R j \alpha^k + 48 M j^2 \alpha^{2 k}.
\nonumber
\end{equation}
For sufficiently large $k$ this implies
\begin{equation}
\left| 1 + j a_1 x y + x^2 y^2 \tilde{F}_j(x,y) \right|
\le 1 + 2 M j \alpha^k,
\label{eq:T0kproof40}
\end{equation}
where we have also used $M \ge 2 R$.
Similarly
\begin{equation}
\left| 1 + j b_1 x y + x^2 y^2 \tilde{G}_j(x,y) \right|
\le 1 + 2 M j \alpha^k.
\label{eq:T0kproof41}
\end{equation}

Write $T_0^j(x,y) = (x_j,y_j)$.
%From \eqref{eq:T0j} we obtain
By using \eqref{eq:T0kproof11}, \eqref{eq:T0kproof20}, and \eqref{eq:T0kproof41} we obtain
\begin{equation}
|y_j| \le \alpha^{k-j} (1+\ee) \left( 1 + 2 M j \alpha^k \right)^2.
\nonumber
\end{equation}
Thus $|y_j| \le \alpha^{k-j} (1 + 2 \ee)$, say, for sufficiently large values of $k$.
Also $|x_j y_j|$ is clearly small, so we can conclude that $(x_j,y_j) \in \cN$
(in particular we have shown that $(x_{k-1},y_{k-1})$
can be made as close to $\left( 0, \frac{1}{\alpha} \right)$ as we like).

By matching the first components of
$T_0^{j+1}(x,y) = \left( T_0 \circ T_0^j \right)(x,y)$ we obtain
\begin{align}
\lambda^{j+1} x \left( 1 + (j+1) a_1 x y + x^2 y^2 \tilde{F}_{j+1}(x,y) \right) &=
\lambda^{j+1} x \left( 1 + j a_1 x y + x^2 y^2 \tilde{F}_j(x,y) \right) \nonumber \\
&\quad+ \lambda^{j+1} a_1 x^2 y (1 + P)
+ \lambda^{j+1} x^3 y^2 Q,
\label{eq:T0kproof50}
\end{align}
where
\begin{align}
1 + P &= \left( \lambda \sigma \right)^j 
\left( 1 + j a_1 x y + x^2 y^2 \tilde{F}_j(x,y) \right)^2
\left( 1 + j b_1 x y + x^2 y^2 \tilde{G}_j(x,y) \right),
\label{eq:T0kproof51} \\
Q &= \left( \lambda \sigma \right)^{2 j} 
\left( 1 + j a_1 x y + x^2 y^2 \tilde{F}_j(x,y) \right)^3
\left( 1 + j b_1 x y + x^2 y^2 \tilde{G}_j(x,y) \right)^2
\tilde{F}(x_j,y_j).
\label{eq:T0kproof52}
\end{align}
By \eqref{eq:T0kproof12}, \eqref{eq:T0kproof40}, and \eqref{eq:T0kproof41}, we obtain
\begin{align}
1 + P &\le \left( 1 + 4 M j \alpha^k \right) \left( 1 + 2 M j \alpha^k \right)^3 \nonumber \\
&\le 1 + 11 M j \alpha^k, \nonumber \\
Q &\le \left( 1 + 4 M j \alpha^k \right)^2 \left( 1 + 2 M j \alpha^k \right)^5 R \nonumber \\
&\le 2 R, \nonumber
\end{align}
assuming $k$ is sufficiently large and
where we have also used $|\tilde{F}(x_j,y_j)| \le R$ (valid because $(x_j,y_j) \in \cN$).
From \eqref{eq:T0kproof50},
\begin{equation}
\tilde{F}_{j+1}(x,y) = \tilde{F}_j(x,y) + \frac{P}{x y} + Q.
\nonumber
\end{equation}
Then by using the induction hypothesis,
the lower bound on $|x y|$ \eqref{eq:T0kproof21},
and the above bounds on $P$ and $Q$, we arrive at
\begin{align}
\left| \tilde{F}_{j+1}(x,y) \right| &\le 12 M j^2 + 22 M j + 2 R \nonumber \\
&\le 12 M j^2 + 24 M j \nonumber \\
&< 12 M (j+1)^2. \nonumber
\end{align}
In a similar fashion by matching the second components of
$T_0^{j+1}(x,y) = \left( T_0 \circ T_0^j \right)(x,y)$
we obtain $\left| \tilde{G}_{j+1}(x,y) \right| < 12 M (j+1)^2$.
This verifies \eqref{eq:T0kproof30} for $j+1$
and so completes the proof.
\end{proof}
%---------------------------------------------------------------------------------
\begin{proof}[Proof of Theorem \ref{th:main}]
\myStep{1}{Coordinate changes to distinguish the surface of homoclinic tangencies.}
First we perform two coordinate changes on parameter space.
There exists an $n \times n$ orthogonal matrix $A$ such that after $\mu$ is replaced with $A \mu$,
and $\mu$ is scaled,
we have $\bn_{\rm tang}^{\sf T} = [1,0,\ldots,0]$ --- the first coordinate vector of $\mathbb{R}^n$.
Then $d_0(\mu) = \mu_1 + \cO \left( \| \mu \|^2 \right)$.
Second, for convenience, we apply a near-identity transformation to remove the higher order terms, resulting in
\begin{equation}
d_0(\mu) = \mu_1 \,.
\label{eq:d0Exp}
\end{equation}
These coordinate changes do not alter the signs of the dot products ${\bf n}_{\rm tang}^{\sf T} {\bf v}$ and ${\bf n}_{\rm eig}^{\sf T} {\bf v}$.  Now write
\begin{align}
c_0(\mu) &= 1 + \sum_{i=1}^n p_i \mu_i + \cO \left( \| \mu \|^2 \right), \label{eq:c0Exp} \\
d_1(\mu) &= \chi + \sum_{i=1}^n q_i \mu_i + \cO \left( \| \mu \|^2 \right), \label{eq:d1Exp} \\
d_2(\mu) &= \sum_{i=1}^n r_i \mu_i + \cO \left( \| \mu \|^2 \right), \label{eq:d2Exp}
\end{align}
\begin{align}
\lambda(\mu) &= \alpha + \sum_{i=1}^n s_i \mu_i + \cO \left( \| \mu \|^2 \right), \label{eq:lambdaExp} \\
\sigma(\mu) &= \frac{\chi_{\rm eig}}{\alpha} + \sum_{i=1}^n t_i \mu_i + \cO \left( \| \mu \|^2 \right), \label{eq:sigmaExp}
\end{align}
where $p_i,\ldots,t_i \in \mathbb{R}$ are constants.

\myStep{2}{Apply a $k$-dependent scaling to $\mu$.}  In view of the coordinate changes applied in the previous step, the surface of homoclinic tangencies of Lemma \blue{\ref{le:tangencySurface}} is tangent to the $\mu_1 = 0$ coordinate hyperplane.  In order to show that bifurcation values scale like $|\alpha|^{2 k}$ if we adjust the value of $\mu$ in a direction transverse to this surface, and, generically, scale like $\frac{|\alpha|^k}{k}$ otherwise, we scale the components of $\mu$ as follows:
\begin{equation}
\mu_i = \begin{cases}
\alpha^{2 k} \tilde{\mu}_i \,, & i = 1, \\
\alpha^k \tilde{\mu}_i \,, & i \ne 1.
\end{cases}
\label{eq:muScaling}
\end{equation}
Below we will see that the resulting asymptotic expansions are consistent and this will justify \eqref{eq:muScaling}.  Notice that $\tilde{\mu}_1$-terms are higher order than $\tilde{\mu}_i$-terms, for $i \ne 1$.
%For example $c_0 = 1 + \sum_{i=2}^n r_i \tilde{\mu}_i \alpha^k + \cO \left( \alpha^{2 k} \right)$.
For example \eqref{eq:lambdaExp} now becomes $\lambda = \alpha + \sum_{i=2}^n s_i \tilde{\mu}_i \alpha^k + \cO \left( |\alpha|^{2 k} \right)$.
Further, let $k$ be such that $\lambda(\bO)^k \sigma(\bO)^k = d_1(\bO)$, that is $\chi_{\rm eig}^k = \chi$.
Then from \eqref{eq:lambdaExp}--\eqref{eq:muScaling} we obtain
\begin{align}
\lambda^k &= \alpha^k \left( 1 + \frac{k}{\alpha} \sum_{i=2}^n s_i \tilde{\mu}_i \alpha^k
+ \cO \left( k^2 |\alpha|^{2 k} \right) \right), \label{eq:lambdak} \\
\sigma^k &= \frac{\chi}{\alpha^k} \left( 1 + k \alpha \chi_{\rm eig} \sum_{i=2}^n t_i \tilde{\mu}_i \alpha^k
+ \cO \left( k^2 |\alpha|^{2 k} \right) \right). \label{eq:sigmak} \\
\end{align}

\myStep{3}{Calculate one point of each periodic solution.}  One point of a single-round periodic solution is a fixed point of $T_0^k \circ T_1$. We look for fixed points $\left( x^{(k)}, y^{(k)} \right)$ of $T_0^k \circ T_1$ of the form
%\begin{align}
\begin{equation}
\begin{split}
x^{(k)} &= \alpha^k \left( 1 + \phi_k \alpha^k + \cO \left( k^2 |\alpha|^{2 k} \right) \right), \\ %\label{eq:xk} \\
y^{(k)} &= 1 + \psi_k \alpha^k + \cO \left( k^2 |\alpha|^{2 k} \right), %\label{eq:yk}
\label{eq:xkyk}
\end{split}
\end{equation}
\blue{where $\phi_{k},\psi_{k} \rightarrow 0$ as $k \rightarrow \infty$.} 
%\end{align}
By substituting \eqref{eq:xkyk} into \eqref{eq:T1}
and the above various asymptotic expressions for the coefficients in $T_1$, we obtain
\begin{equation}
T_1 \left( x^{(k)}, y^{(k)} \right) =
\begin{bmatrix}
1 + \left( c_{1,0} + c_{2,0} \psi_k + \sum_{i=2}^n p_i \tilde{\mu}_i \right) \alpha^k + \cO \left( k^2 |\alpha|^{2 k} \right)\\

\begin{aligned}
\chi \alpha^k + \left( \tilde{\mu}_1 + d_{3,0} + d_{4,0} \psi_k + d_{5,0} \psi_k^2 + \chi \phi_k
+ \sum_{i=2}^n \left( q_i + r_i \psi_k \right) \tilde{\mu}_i \right) \alpha^{2 k}&\\[2mm]
+ \cO \left( k^2 |\alpha|^{3 k} \right)&
\end{aligned}
\end{bmatrix}.
%\label{eq:T1xkyk}
\nonumber
\end{equation}
Then by \eqref{eq:T0kun},
{\footnotesize
\begin{equation}
\left( T_0^k \circ T_1 \right) \left( x^{(k)}, y^{(k)} \right) =
\begin{bmatrix}
\begin{aligned}
\alpha^k + \bigg( a_{1,0} \chi + \frac{1}{\alpha} \sum_{i=2}^n s_i \tilde{\mu}_i \bigg) k \alpha^{2 k}
+ \bigg( c_{1,0} + c_{2,0} \psi_k + \sum_{i=2}^n p_i \tilde{\mu}_i \bigg) \alpha^{2 k} &\\[2mm]
+ \cO \left( k^2 |\alpha|^{3 k} \right) &
\end{aligned}\\
\begin{aligned}
1 + \bigg( b_{1,0} \chi + \alpha \chi_{\rm eig} \sum_{i=2}^n t_i \tilde{\mu}_i \bigg) k \alpha^k
+ \bigg( \chi \tilde{\mu}_1 + \chi d_{3,0} + \chi d_{4,0} \psi_k &\\
+ \chi d_{5,0} \psi_k^2 + \phi_k + \chi \sum_{i=2}^n \left( q_i + r_i \psi_k \right) \tilde{\mu}_i \bigg) \alpha^k + \cO \left( k^2 |\alpha|^{2 k} \right) &
\end{aligned}
\end{bmatrix}.
\label{eq:T0kT1xkykun}
\end{equation}
}
By matching \eqref{eq:xkyk} and \eqref{eq:T0kT1xkykun} and eliminating $\phi_k$
we obtain the following expression that determines the possible values of $\psi_k$:
\begin{equation}
\chi d_{5,0} \psi_k^2 - P \psi_k + Q=0,
\label{eq:psikQuadratic}
\end{equation}
where
\begin{align}
P &= 1 - c_{2,0} - \chi d_{4,0} - \chi \sum_{i=2}^n r_i \tilde{\mu}_i \,, \label{eq:P} \\
Q &= \chi \tilde{\mu}_1 + c_{1,0} + \chi d_{3,0} + \sum_{i=2}^n \left( p_i + \chi q_i \right) \tilde{\mu}_i
+ \sum_{i=2}^n \left( \frac{s_i}{\alpha} + \alpha \chi_{\rm eig} t_i \right) \tilde{\mu}_i k, \label{eq:Q}
\end{align}
and we have also used $a_{1,0} + b_{1,0} = 0$.
Of the two solutions to \eqref{eq:psikQuadratic}, the one that yields an asymptotically stable solution when $\mu = \bO$ is \blue{as shown below}
\begin{equation}
\psi_k = \frac{1}{2 \chi d_{5,0}} \left( P - \sqrt{P^2 - 4 \chi d_{5,0} Q} \right).
\label{eq:psik}
\end{equation}
Note that this solution exists and is real-valued for sufficiently small values of $\mu$ because when $\mu = \bO$
the discriminant is $P^2 - 4 \chi d_{5,0} Q = \Delta_0 > 0$.

\myStep{4}{Stability of the periodic solution.} By using \eqref{eq:T1}, \eqref{eq:T0kun}, \eqref{eq:lambdak}, \eqref{eq:sigmak}, and \eqref{eq:xkyk},
\begin{equation}
\rD \left( T_0^k \circ T_1 \right) \left( x^{(k)}, y^{(k)} \right) =
\begin{bmatrix}
\cO \left( |\alpha|^k \right) &
c_{2,0} \alpha^k + \cO \left( k |\alpha|^{2 k} \right) \\
\frac{1}{\alpha^k} \left( 1 + \cO \left( k |\alpha|^k \right) \right) &
\begin{aligned}
\\[2mm]\chi d_{4,0} + \chi \sum_{i=2}^n r_i \tilde{\mu}_i&\\[2mm] + 2 \chi d_{5,0} \psi_k + \cO \left( k |\alpha|^k \right)&
\end{aligned}
\end{bmatrix}.
\label{eq:Jacobian}
\end{equation}
Let $\tau_k$ and $\delta_k$ denote the trace and determinant of this matrix, respectively.
By \eqref{eq:P}, \eqref{eq:psik}, and \eqref{eq:Jacobian} we obtain
\begin{align}
\tau_k &= 1 - c_{2,0} - \sqrt{P^2 - 4 \chi d_{5,0} Q} + \cO \left( k |\alpha|^k \right), \label{eq:taukun} \\
\delta_k &= -c_{2,0} + \cO \left( k |\alpha|^k \right). \label{eq:deltak}
\end{align}
By substituting $\mu = \bO$ into \eqref{eq:taukun}
we obtain $\tau_k = 1 - c_{2,0} - \sqrt{\Delta_0}$.
It immediately follows from the assumption $-1 < c_{2,0} < 1 - \frac{\sqrt{\Delta_0}}{2}$
that the periodic solution is asymptotically stable for sufficiently large values of $k$.

\myStep{5}{The generic case ${\bf n}_{\rm tang}^{\sf T} {\bf v} \ne 0$.} Now suppose $\bn_{\rm tang}^{\sf T} \bv \ne 0$, that is, $v_1 \ne 0$ (in view of the earlier coordinate change).
Write $\mu = \ee \bv$ and $\ee = \tilde{\ee} \alpha^{2 k}$.
By \eqref{eq:muScaling}, $\tilde{\mu}_1 = \tilde{\ee} v_1$ and $\tilde{\mu}_i = \tilde{\ee} v_i \alpha^k$ for $i \ne 1$.
Then by \eqref{eq:P} and \eqref{eq:Q}, $P = 1 - c_{2,0} - \chi d_{4,0} + \cO \left( |\alpha|^k \right)$
and $Q = c_{1,0} + \chi d_{3,0} + \chi \tilde{\ee} v_1 + \cO \left( |\alpha|^k \right)$, so
\begin{equation}
P^2 - 4 \chi d_{5,0} Q = \left( 1 - c_{2,0} - \chi d_{4,0} \right)^2
- 4 \chi d_{5,0} \left( c_{1,0} + \chi d_{3,0} + \chi \tilde{\ee} v_1 \right) + \cO \left( |\alpha|^k \right).
\label{eq:DeltaCase1}
\end{equation}
Since $v_1 \ne 0$ and $d_{5,0} \ne 0$ we can solve $\delta_k - \tau_k + 1 = 0$ for $\tilde{\ee}$
(formally this is achieved via the implicit function theorem) and the solution is
\begin{equation}
\tilde{\ee}_{\rm SN} = \frac{\Delta_0}{4 d_{5,0} v_1} + \cO \left( k |\alpha|^k \right).
\label{eq:tildeeeSNCase1}
\end{equation}
Also we can use \eqref{eq:DeltaCase1} to solve $\delta_k + \tau_k + 1 = 0$ for $\tilde{\ee}$:
\begin{equation}
\tilde{\ee}_{\rm PD} = \frac{\Delta_0 - 4 \left( 1 - c_{2,0} \right)^2}{4 d_{5,0} v_1} + \cO \left( k |\alpha|^k \right).
\label{eq:tildeeePDCase1}
\end{equation}
Since $\tilde{\ee}_{\rm SN}$ and $\tilde{\ee}_{\rm PD}$ evidently have opposite signs, this completes the proof in the first case.

\myStep{6}{The degenerate case ${\bf n}_{\rm tang}^{\sf T} {\bf v} = 0$.}  Now suppose $v_1 = 0$ and $\bn_{\rm eig}^{\sf T} \bv \ne 0$.
Again write $\mu = \ee \bv$ but now write $\ee = \frac{\tilde{\ee} \alpha^k}{k}$.
By \eqref{eq:muScaling}, $\tilde{\mu}_1 = 0$ and $\tilde{\mu}_i = \frac{\tilde{\ee} v_i}{k}$ for $i \ne 1$.

We first evaluate $\bn_{\rm eig}^{\sf T} \bv$.
By multiplying \eqref{eq:lambdaExp} and \eqref{eq:sigmaExp} and comparing the result to \eqref{eq:nEig}
we see that the $i^{\rm th}$ element of $\bn_{\rm eig}$
is $\frac{\chi_{\rm eig} s_i}{\alpha} + \alpha t_i$.
Then since $v_1 = 0$ we have
$\bn_{\rm eig}^{\sf T} \bv = \sum_{i=2}^n \left( \frac{\chi_{\rm eig} s_i}{\alpha} + \alpha t_i \right) v_i$.
Further, $\mu = \ee \bv$ and $\ee = \frac{\tilde{\ee} \alpha^k}{k}$,
\begin{equation}
\chi_{\rm eig} \tilde{\ee} \bn_{\rm eig}^{\sf T} \bv
= \sum_{i=2}^n \left( \frac{s_i}{\alpha} + \alpha \chi_{\rm eig} t_i \right) \tilde{\mu}_i k,
\nonumber
\end{equation}
which is a term appearing in \eqref{eq:Q}.
So by \eqref{eq:P} and \eqref{eq:Q}, $P = 1 - c_{2,0} - \chi d_{4,0} + \cO \left( \frac{1}{k} \right)$ and $Q = c_{1,0} + \chi d_{3,0} + \frac{\tilde{\ee} \chi_{\rm eig} \bn_{\rm eig}^{\sf T} \bv}{k} + \cO \left( \frac{1}{k} \right)$.
By solving $\delta_k - \tau_k + 1 = 0$,
\begin{equation}
\tilde{\ee}_{\rm SN} = \frac{\Delta_0}{4 d_{5,0} \chi \chi_{\rm eig} \bn_{\rm eig}^{\sf T} \bv} + \cO \left( \frac{1}{k} \right),
\label{eq:tildeeeSNCase2}
\end{equation}
and by solving $\delta_k - \tau_k + 1 = 0$,
\begin{equation}
\tilde{\ee}_{\rm PD} = \frac{\Delta_0 - 4 \left( 1 - c_{2,0} \right)^2}{4 d_{5,0} \chi \chi_{\rm eig} \bn_{\rm eig}^{\sf T} \bv} + \cO \left( \frac{1}{k} \right).
\label{eq:tildeeePDCase2}
\end{equation}
\end{proof}
As in the previous case, $\tilde{\ee}_{\rm SN}$ and $\tilde{\ee}_{\rm PD}$ have opposite signs, and this completes the proof in the second case.  Notice how the assumptions we have made ensure the denominators of \eqref{eq:tildeeeSNCase2} and \eqref{eq:tildeeePDCase2} are nonzero.
%===============================================================================
\section{Conclusion}
\label{sec:concch4}
\blue{In this chapter we have studied what happens when parameters are varied to move away from a globally resonant homoclinic tangency at which infinitely many asymptotically stable single-round periodic solutions coexist.  We have shown that these periodic solutions are either sequentially destroyed in saddle-node bifurcations, or sequentially lose stability in period-doubling bifurcations.}

\blue{If the parameter change is performed in a generic fashion, then the amount by which the parameter needs to vary for the bifurcation of the ${\rm SR}_k$-solution to occur is asymptotically proportional to $|\alpha|^{2 k}$, where $\alpha$ is the stable eigenvalue associated with the fixed point at the globally resonant homoclinic tangency.  This scaling law forms the first part of Theorem \ref{th:main}.  Equations \eqref{eq:tildeeeSNCase1} and \eqref{eq:tildeeePDCase1} provide formulas for the values of the saddle-node and period-doubling bifurcations, to leading order.}

\blue{When we say that the parameter change is generic, we mean it produces a variation that is transverse to the surface of codimension-one homoclinic tangencies.  This surface passes through the point at which the globally resonant homoclinic tangency occurs and above we introduced coordinates so that this surface is simply given by $\mu_1 = 0$.  If the parameter change is not generic, so is tangent to the homoclinic tangency surface, then the bifurcation values will instead scale like $\frac{|\alpha|^k}{k}$, assuming other genericity conditions are satisfied.  This is the second part of Theorem \ref{th:main}.  Equations \eqref{eq:tildeeeSNCase2} and \eqref{eq:tildeeePDCase2} provide the leading-order contribution to the bifurcation values in this case.  These scaling laws are illustrated in the next chapter for a particular family of maps.}

\blue{As a final comment, observe that for any positive $j \in \mathbb{Z}$ there exists an open region $\Omega_j$ of parameter space in which the family of maps has $j$ asymptotically stable single-round periodic solutions.  Let $B_r$ be a ball (sphere) of radius $r$ centred at the globally resonant homoclinic tangency, and suppose $r$ is as big as possible subject to the constraint $B_r \subset \Omega_j$.  Then, since the fastest scaling law is $|\alpha|^{2 k}$, we can infer that the radius $r$ is proportional to $|\alpha|^{2 j}$.}

\chapter[A planar map with infinite coexistence]{An explicit example of a planar map with a globally resonant homoclinic tangency \label{cha:ExampleInfinte} }
In this chapter we study an explicit example of a planar map that exhibits a globally resonant homoclinic tangency. The map is constructed so that it is linear in a neighbourhood of the origin, has a homoclinic tangency, and all conditions in the above theorems for globally resonant homoclinic tangencies can be explicitly checked. We compute phase portraits, basins of attraction, and invariant manifolds that display the features studied in the earlier chapters. We then unfold to illustrate different scaling laws in different directions of parameter space. \blue{The numerical simulations below were all computed using {\sc matlab}.}
\section{Explicit examples of infinite coexistence}
\label{sec:example}

\begin{figure}[!htbp]
\begin{center}
\includegraphics{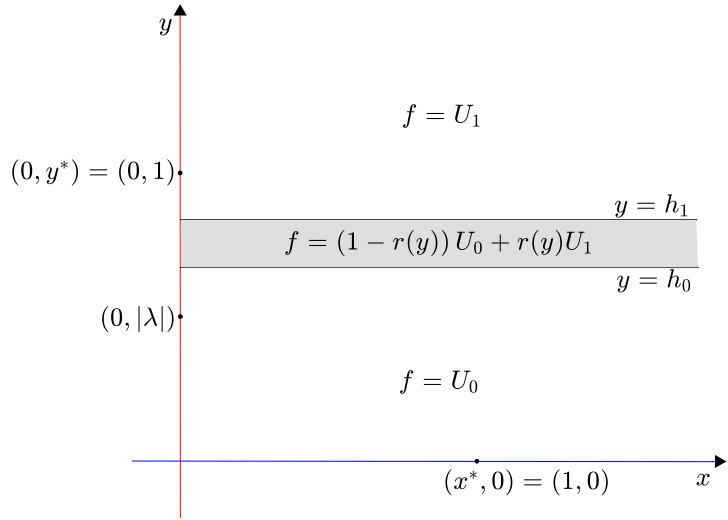}
\caption[Phase space of map exhibiting GRHT phenomenon.]{ The basic structure of the phase space of the map \eqref{eq:fEx}.}
\label{fig:smooth_example}
\end{center}
\end{figure}

Here we demonstrate Theorems \ref{th:sufficientConditionsPositiveCase} and \ref{th:sufficientConditionsNegativeCase}
with a piecewise-smooth $C^{1}$ \blue{family of} maps of the form
\begin{equation}
f(x,y) = \begin{cases}
U_0(x,y), & y \le h_0 \,, \\
(1 - r(y)) U_0(x,y) + r(y) U_1(x,y), & h_0 \le y \le h_1 \,, \\
U_1(x,y), & y \ge h_1 \,,
\end{cases}
\label{eq:fEx}
\end{equation}
where \blue{$U_{0}, U_{1}, r(y)$ are as defined in \eqref{eq:U0}, \eqref{eq:U1}, \eqref{eq:r} respectively, and}
\begin{align}
h_0 &= \frac{2 |\lambda| + 1}{3}, &
h_1 &= \frac{|\lambda| + 2}{3},
\label{eq:xi0xi1}
\end{align}
\blue{with $|\lambda|<1$.}
This family of maps was motivated by similar families exhibiting infinite coexistence introduced in \cite{GaTr83}. The piecewise nature of \eqref{eq:fEx} is illustrated in Fig.~\ref{fig:smooth_example}.
 
We recall that for a globally resonant homoclinic tangency with a local map $T_{0}$ of the form \eqref{eq:T02} for the orientation-preserving case and of the form \eqref{eq:T0rev} for the orientation-reversing case, and a global map $T_{1}$ of the form \eqref{eq:T11}, the necessary conditions for the \blue{infinite} coexistence  in the\\
(i) Orientation-reversing case \\
a) Eigenvalues multiplying to $-1$
\begin{equation}
    \lambda \sigma = -1
\end{equation}
b) Homoclinic tangency
\begin{equation}
    d_{2}  = 0
\end{equation}
c) Global resonance condition
\begin{equation}
    d_{1} = \begin{cases}
    \hspace*{3mm}1, {\rm even}\,k\\
    -1, {\rm odd}\,k
    \end{cases}
\end{equation}
%a) homoclinic tangency ($d_{2}  =0$), b) $\lambda \sigma = -1$, c) $d_{1} = 1$ if $k$ is even, $-1$ if $k$ is odd.\\
(ii) Orientation-preserving case \\
a) Eigenvalues multiplying to $1$
\begin{equation}
    \lambda \sigma = 1
\end{equation}
b) Homoclinic tangency
\begin{equation}
    d_{2} = 0
\end{equation}
c) Global resonance condition
\begin{equation}
    d_{1}  = 1
\end{equation}
d) Condition on coefficients of local map
\begin{equation}
    a_{1} + b_{1} = 0
\end{equation}

with the sufficient condition
\begin{equation}
-1 < \frac{c_2 y^*}{x^*} < 1 - \frac{\sqrt{\Delta}}{2}.
\end{equation}

We let
\begin{equation}
U_0(x,y) = \begin{bmatrix}
\lambda x \\
\sigma y
\end{bmatrix},
\label{eq:U0}
\end{equation}
so that $T_0 = U_0$, and
\begin{equation}
U_1(x,y) = \begin{bmatrix}
1 + c_2 (y - 1) \\
d_1 x + d_5 (y - 1)^2
\end{bmatrix},
\label{eq:U1}
\end{equation}
so that $T_1 = f^m = U_1$ using $m = 1$ and $x^* = y^* = 1$.
Since \eqref{eq:U1} neglects some terms in \eqref{eq:T11},
equation \eqref{eq:discriminant} reduces to $\Delta = (1 - c_2)^2$
and \eqref{eq:stabilityCondition} reduces to $|c_2| < 1$.
Notice that with \eqref{eq:xi0xi1}, equation \eqref{eq:fourPoints} is satisfied
and $(0,y^*) = (0,1)$ lies in the region $y \ge h_1$.

\begin{figure}[h!]
\begin{center}
\includegraphics{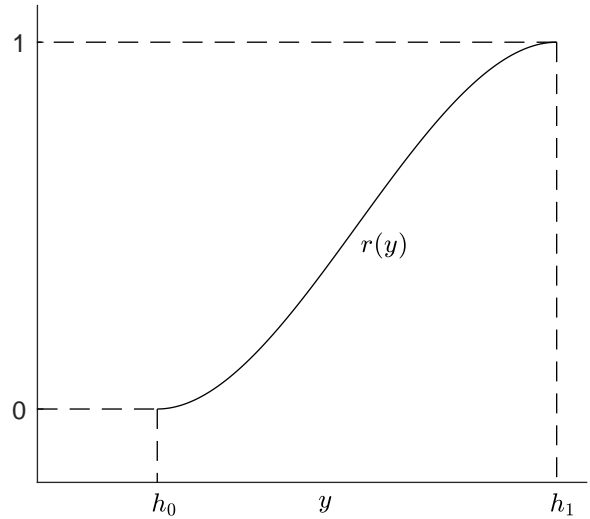}
\caption[The function used to smoothen the toy map.]{The function \eqref{eq:r} (with \eqref{eq:s}) that we use as a convex combination parameter in \eqref{eq:fEx}.}
\label{fig:smooth_ry}
\end{center}
\end{figure}

One could choose the function $r$ in \eqref{eq:fEx} such that $f$ is $C^\infty$ on the switching manifolds $y=h_0$ and $y=h_1$, and thus $C^\infty$ on $\mathbb{R}^2$. However, the nature of $f$ outside a neighbourhood of the homoclinic orbit is not important to Theorems \ref{th:sufficientConditionsPositiveCase} and \ref{th:sufficientConditionsNegativeCase} because the ${\rm SR}_k$-solutions converge to the homoclinic orbit, so for simplicity we choose $r$ such that $f$ is $C^1$. This requires
$r(h_0) = 0$,
$r'(h_0) = 0$,
$r(h_1) = 1$, and
$r'(h_1) = 0$.
The unique cubic polynomial satisfying these requirements is
\begin{equation}
r(y) = s \left( \frac{y - h_0}{h_1 - h_0} \right),
\label{eq:r}
\end{equation}
where
\begin{equation}
s(z) = 3 z^2 - 2 z^3,
\label{eq:s}
\end{equation}
see Fig.~\ref{fig:smooth_ry}.

\subsection{Non-invertibility of the map}

\begin{figure}[t!]
\begin{center}
\includegraphics[width=0.8\textwidth]{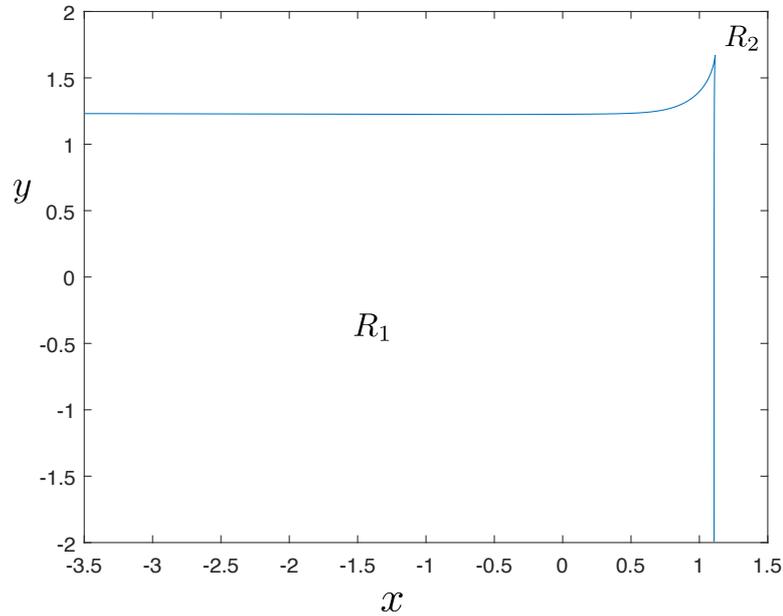}
\caption[Critical curve of the toy map.]{ The critical curve $LC$ of the map $f$. The non-invertible map $f$ is of type $Z_{1}-Z_{3}$. The parameters are set as $\lambda = \frac{3}{5}, \sigma = \frac{5}{3}, d_{1}=1, c_{2}=-0.8, d_{5} = 0.8$.}
\label{fig:LCToymap}
\end{center}
\end{figure}

It is important to note that the piecewise-smooth map $f$ is non-invertible. This is evident from the fact that the points $(0,1)$ and $(\frac{1}{\blue{\lambda}},0)$ both map to $(1,0)$ under $f$. Different regions of $\mathbb{R}^2$ have different numbers of preimages. The boundary which separates regions with distinct numbers of preimages is determined by the curve det$({\rm D}f) =0$, shown in Fig. \ref{fig:LCToymap}. This curve divides $\mathbb{R}^2$ into two regions, a region \blue{$R_{1}$} with three pre-images and a region \blue{$R_{2}$} with one pre-image, hence the map is of type $Z_{1}-Z_{3}$ using the terminology of \cite{CM96} discussed above in \S \ref{twodimnoninv}.  We also observe the curve $\det({\rm D} f) = 0$ has a cusp at $(x,y) \approx (1.6706,1.1156)$.

\subsection{Computation of stable and unstable manifolds}
\begin{figure}[h!]
\begin{center}
\includegraphics[width=0.75\textwidth]{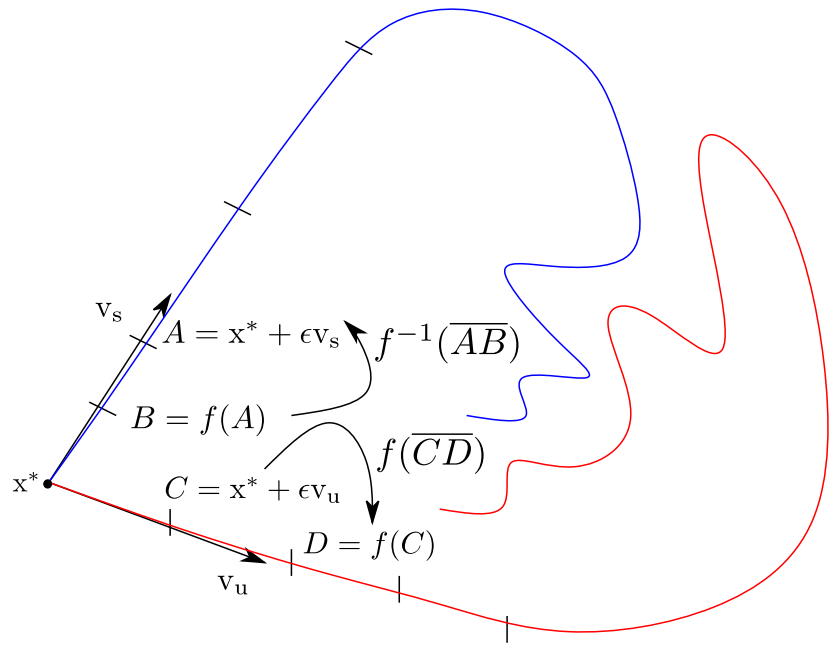}
\caption[Method of construction of the stable and unstable manifold of a fixed point.]{ A schematic illustrating our method for the numerical computation of the stable and unstable manifolds of a fixed point.}
\label{fig:ManifoldConstr}
\end{center}
\end{figure}

We next discuss the numerical computation of the stable and unstable manifolds of a saddle fixed point $\mathbf{x}^{*}$. To compute the stable manifold, we consider a nearby point $A$ to the saddle fixed point in the stable eigen-direction, see Fig. \ref{fig:ManifoldConstr}. We next consider the forward image of $A$ denoted as $B = f(A)$. In the case of a positive eigenvalue take a set of points on the line segment $\overline{AB}$. This is because $\overline{AB}$ approximates a {\em fundamental domain} for the stable manifold.  For a fundamental domain, every point in the domain corresponds to a different orbit in the stable manifold, and conversely every orbit in the stable manifold corresponds to one point in the domain.  With instead a negative eigenvalue we use $B = f^2(A)$. Using the line segment $\overline{AB}$ we determine a sequence of pre-images of these and our numerical approximation to the stable manifold is obtained. This method is simple and generates stable manifolds sufficiently accurately for our purposes.  More sophisticated methods for computing stable and invariant involve formulating and solving boundary value problems \citep{BK07}.

Similarly, for the unstable manifold, we consider a nearby point $C$ to the saddle fixed point in the unstable eigen-direction. In the case of a positive eigenvalue we next consider the forward image of $C$ denoted as $D = f(C)$. We take forward images of a set of points on the line segment $\overline{CD}$ and our numerical approximation to the unstable manifold is obtained.

\begin{figure}[h!]
\begin{center}
\includegraphics[width=0.75\textwidth]{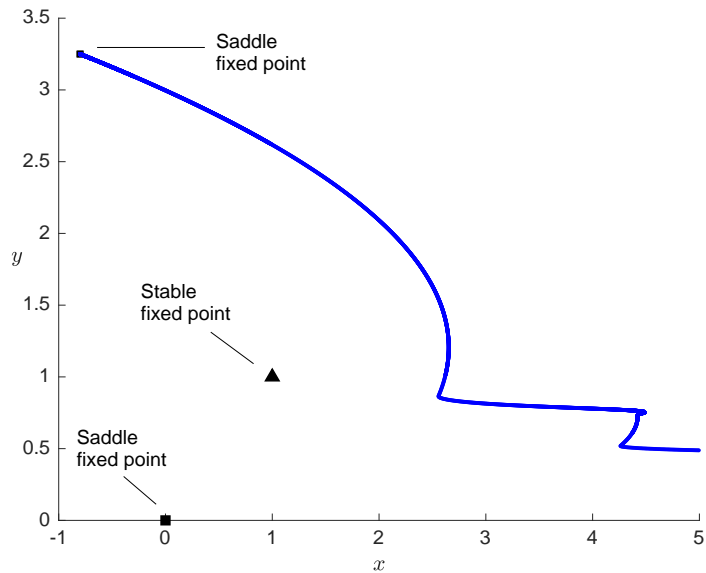}
\caption[Stable manifold of a saddle fixed point of map $U_{1}$. ]{The stable manifold in blue of the saddle fixed point $(-0.8,3.25)$ of map $U_{1}$. The stable fixed point $(1,1)$ of map $U_{1}$ is marked in black. The saddle fixed point $(0,0)$ of the map $U_{0}$ is also shown. The parameters are set as $\lambda = \frac{3}{5}, \sigma = \frac{5}{3}, d_{1}=1, c_{2}=-0.8, d_{5} = 0.8$.}
\label{fig:StableManifoldU1}
\end{center}
\end{figure}

Next we apply this to the piecewise-smooth $C^{1}$ map \eqref{eq:fEx}. We see that the fixed point of $U_{0}$ is $(0,0)$ and is a saddle fixed point. The fixed points of the map $U_{1}$ are given by solving the simultaneous equations
\begin{equation}
    \begin{aligned}
    1 + c_{2}(y-1) &= x,\\
    d_{1} x + d_{5}(y-1)^2 - 1 &= y-1,
    \end{aligned}
    \label{eq:U1fpseqn}
\end{equation}
which gives us two fixed points as 
\begin{equation}
    \begin{aligned}
    y^{*} &= 1 + \frac{(1-c_{2}d_{1} \pm \sqrt{(c_{2}d_{1})^2 - 4d_{5}(d_{1}-1)})}{2d_{5}},\\
    x^{*} &= 1 + c_{2}(y^* - 1).
    \end{aligned}
    \label{eq:U1fps}
\end{equation}
With $d_{5} = 0.8, c_{2} = -0.8, \lambda = \frac{3}{5}, \sigma = \frac{5}{3},$ and  $d_{1} = 1$, these fixed points are located at $(1,1)$ and $(-0.8,3.25)$. By computing the Jacobian matrix and analysing the eigenvalues of the map $U_{1}$ at its fixed points, we find that $(1,1)$ is asymptotically stable whereas $(-0.8,3.25)$ is a saddle. The stable manifold of this saddle is shown in Fig.~\ref{fig:StableManifoldU1}.

\begin{figure}[!htbp]
 %\hspace{0.8cm}
	\begin{tabular}{c c}

		\hspace{-0.5cm}\includegraphics[width=0.55\textwidth]{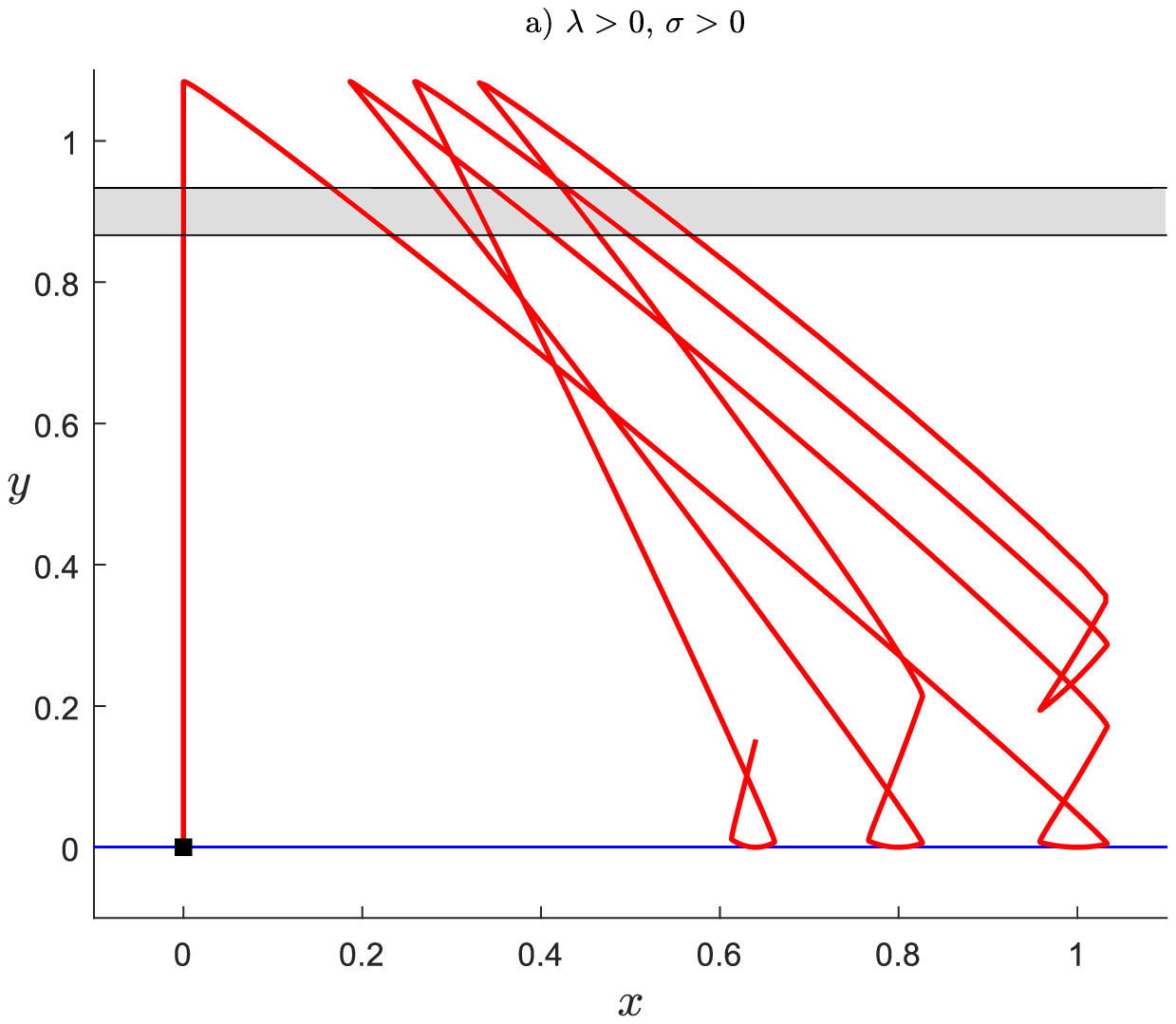}&
		\hspace{-1.3cm}\includegraphics[width=0.55\textwidth]{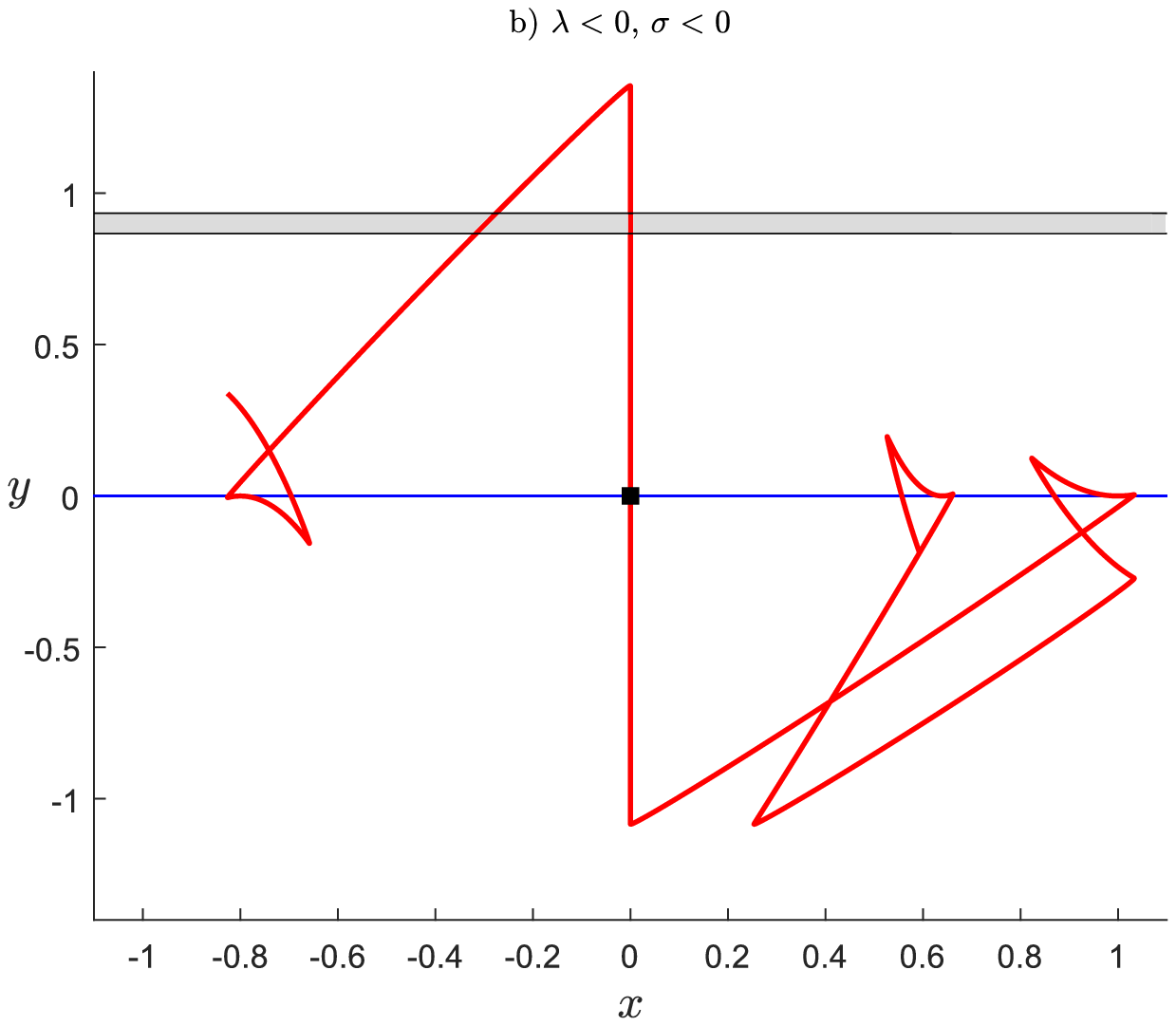}\\
%$a$ & $b$ \\
		%\includegraphics[width=0.55\textwidth]{case3.eps}&
		\hspace{-0.5cm}\includegraphics[width=0.55\textwidth]{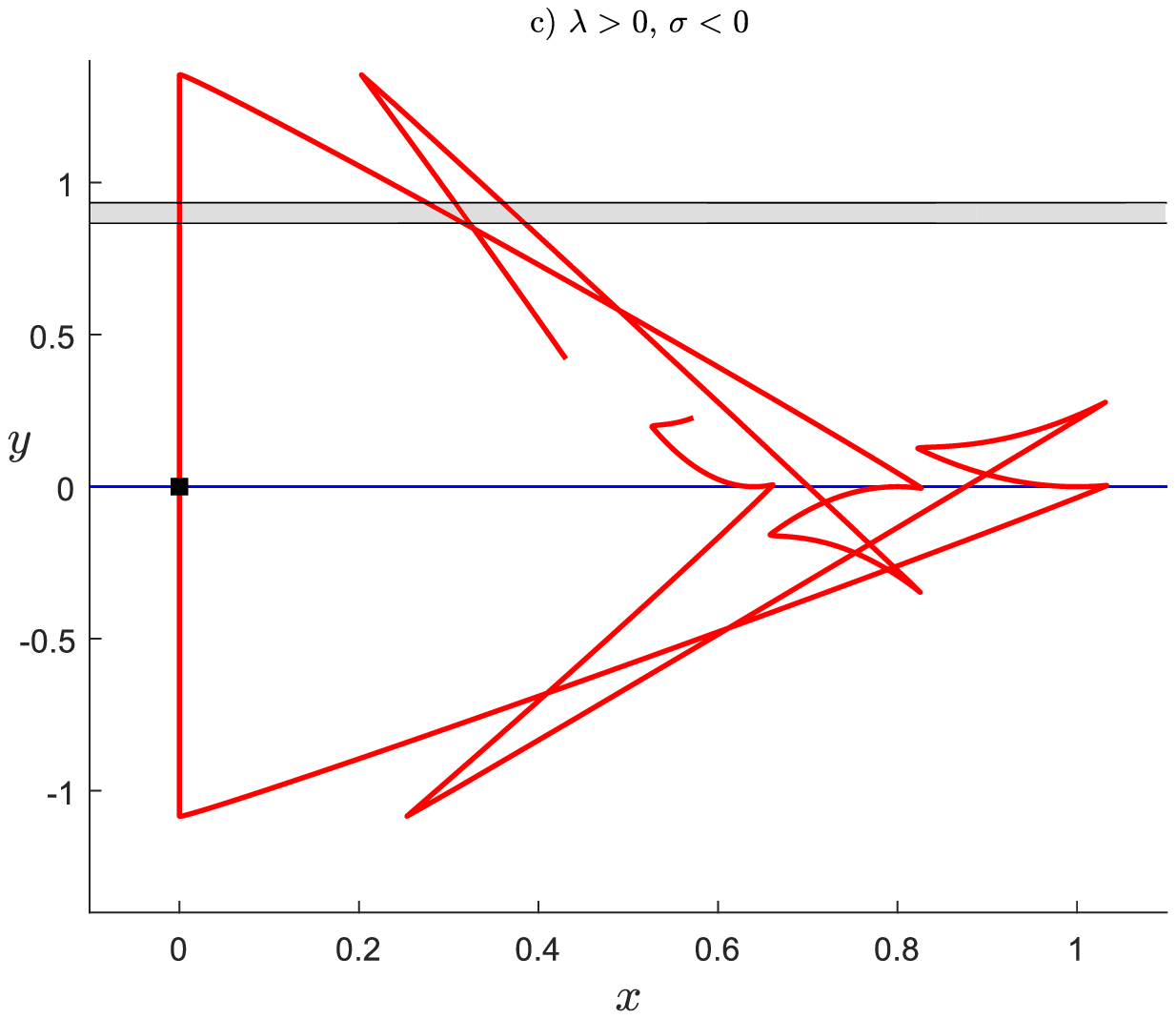}&
		\hspace{-1.3cm}\includegraphics[width=0.55\textwidth]{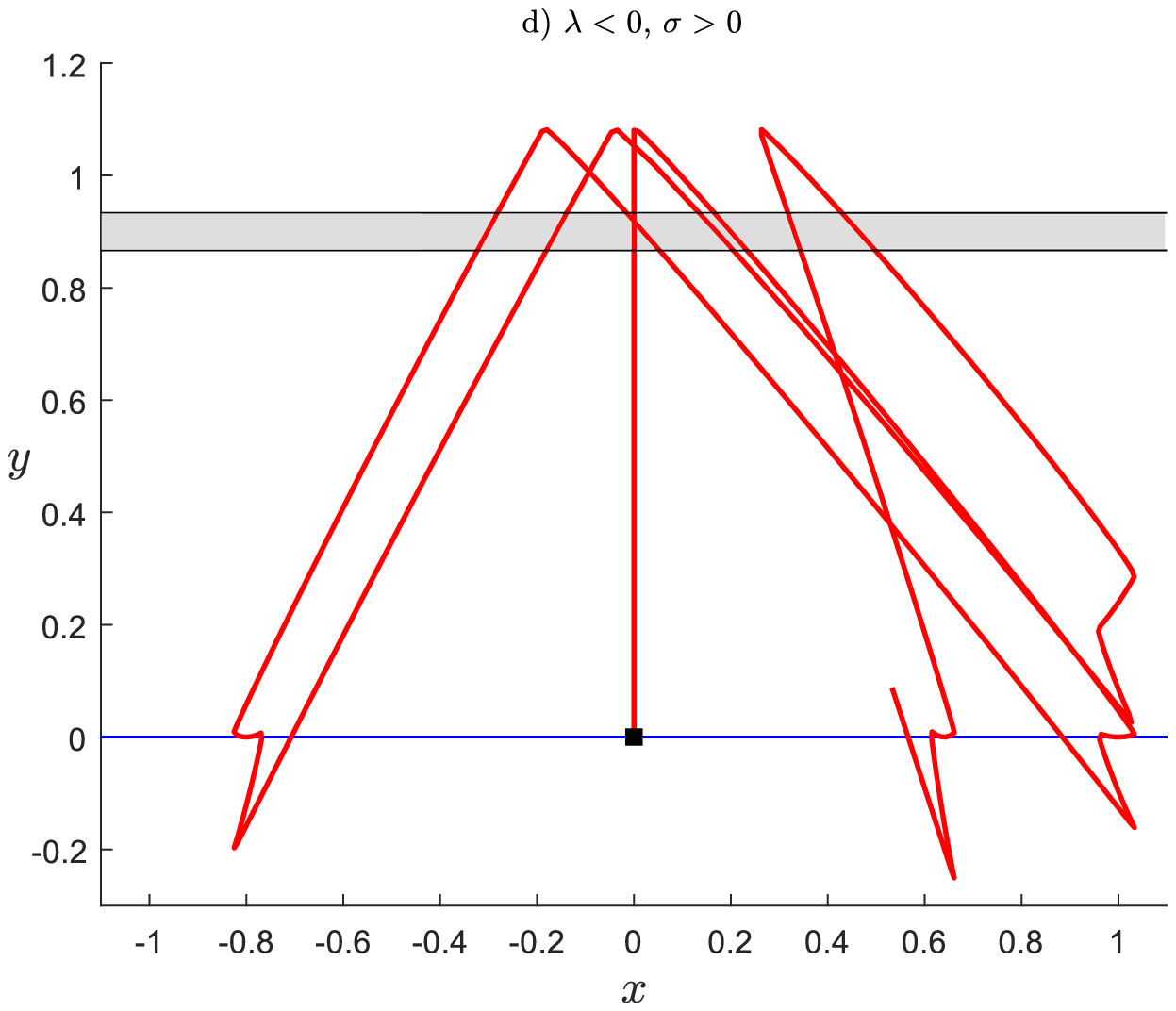}\\
%$c$ & $d$\\
	\end{tabular}
	\caption[Homoclinic tangency in piecewise-smooth $C^{1}$ map]{Parts of the stable [blue] and unstable [red] manifolds of the origin for the map \eqref{eq:fEx} with \eqref{eq:xi0xi1},\eqref{eq:c2d5}.
Panels (a)--(d) correspond to \eqref{eq:param1}--\eqref{eq:param4} respectively.
In each panel the region $h_0 < y < h_1$ is shaded. A tangential intersection of the stable and unstable manifolds can be seen in each of the cases. The homoclinic tangency in each case is quadratic.
}
	\label{fig:HT}
\end{figure}

Next we illustrate globally resonant homoclinic tangency in both orientation-preserving and orientation-reversing cases for \eqref{eq:fEx}.
The values of the parameters are chosen so that the necessary and sufficient conditions of Theorems \ref{th:sufficientConditionsPositiveCase} and  \ref{th:sufficientConditionsNegativeCase} for orientation-reversing and orientation-preserving cases respectively are satisfied. As examples we fix
\begin{align}
c_2 &= -\tfrac{1}{2}, & d_5 &= 1,
\label{eq:c2d5}
\end{align}
and consider the following four combinations for the remaining parameter values,
\begin{align}
\lambda &= \tfrac{4}{5}, & \sigma &= \tfrac{5}{4}, & d_1 &= 1, \label{eq:param1} \\
\lambda &= -\tfrac{4}{5}, & \sigma &= -\tfrac{5}{4}, & d_1 &= 1, \label{eq:param2} \\
\lambda &= \tfrac{4}{5}, & \sigma &= -\tfrac{5}{4}, & d_1 &= 1, \label{eq:param3} \\
\lambda &= -\tfrac{4}{5}, & \sigma &= \tfrac{5}{4}, & d_1 &= -1. \label{eq:param4}
\end{align}
The parameter combinations \eqref{eq:param1} and \eqref{eq:param2} correspond to the orientation-preserving case, while \eqref{eq:param3} and \eqref{eq:param4} correspond to the orientation-reversing case.
Fig.~\ref{fig:HT} shows the stable and unstable manifolds of the origin for all four combinations \eqref{eq:param1}--\eqref{eq:param4}. Specifically we have followed the unstable manifold upwards from the origin until observing three
tangential intersections with the stable manifold on the $x$-axis
(some transversal intersections are also visible).
In panels (b) and (c) the unstable manifold evolves outwards along both the positive and negative $y$-axes
because here the unstable eigenvalue $\sigma$ is negative.
In some places the unstable manifolds have extremely high curvature due to the high degree of
nonlinearity of~\eqref{eq:fEx} in the region $h_0 \le y \le h_1$.
Also notice the unstable manifolds have self-intersections because \eqref{eq:fEx} is non-invertible.

\subsection{Computation of single-round periodic solutions}
The single-round periodic solutions of \eqref{eq:fEx} converge to the homoclinic orbit \blue{of the saddle fixed point at the origin of the map $f$}.  The homoclinic orbit is bounded away from the horizontal strip in Fig.~\ref{fig:smooth_example} where the middle piece of \eqref{eq:fEx} applies.  Therefore, for sufficiently large values of $k$, one point of a single-round period-$k+1$ solution is given by a fixed point of $U_0^k \circ U_1$.  Fixed points of this composed map satisfy  
\begin{equation}
    \begin{aligned}
    \lambda^k(1 + c_{2}(y-1)) &= x,\\
    \sigma^k(d_{1} x + d_{5}(y-1)^2 - 1) &= y-1.
    \end{aligned}
    \label{eq:composedquad}
\end{equation}
This is a quadratic equation, so is trivial to solve. \blue{One of the solutions to \eqref{eq:composedquad} corresponds to a stable fixed point and the other to a saddle fixed point. When the stable fixed point of $U_{0}^{k} \circ U_{1}$ is iterated under $f$, we get an asymptotically stable single-round periodic solution of period $k+1$ of the map $f$. When the saddle fixed point of $U_{0}^{k} \circ U_{1}$ is iterated under the map $f$, we get a saddle single-round periodic solution of period $k+1$ of the map $f$. Observe that this is in accordance with Theorem \ref{th:sufficientConditionsPositiveCase} for the orientation-preserving case and Theorem \ref{th:sufficientConditionsNegativeCase} for the orientation-reversing case}.  The remaining points of the periodic solution are then given by iterating fixed points of $U_{0}^{k} \circ U_{1}$ under the map \blue{$f$}.

Fig.~\ref{fig:BasinManifoldU1} shows several of the periodic solutions for the parameter combination of $\lambda = \frac{3}{5}, \sigma = \frac{5}{3}, d_{1}=1, c_{2}=-0.8, d_{5} = 0.8$. By Theorem \ref{th:a1b1}, the map \blue{$f$} has an asymptotically stable, single-round periodic solution \blue{of period $k+1$} for all sufficiently large values of $k$.  In fact the map has these for all $k \ge 1$.  Moreover, the fixed point at $(x,y) = (1,1)$ can be viewed as corresponding to $k=0$.

\begin{figure}[!htbp]
	%\hspace{-1.3cm}	
	\begin{tabular}{c c }

		\hspace{-0.3cm}\includegraphics[width=0.5\textwidth]{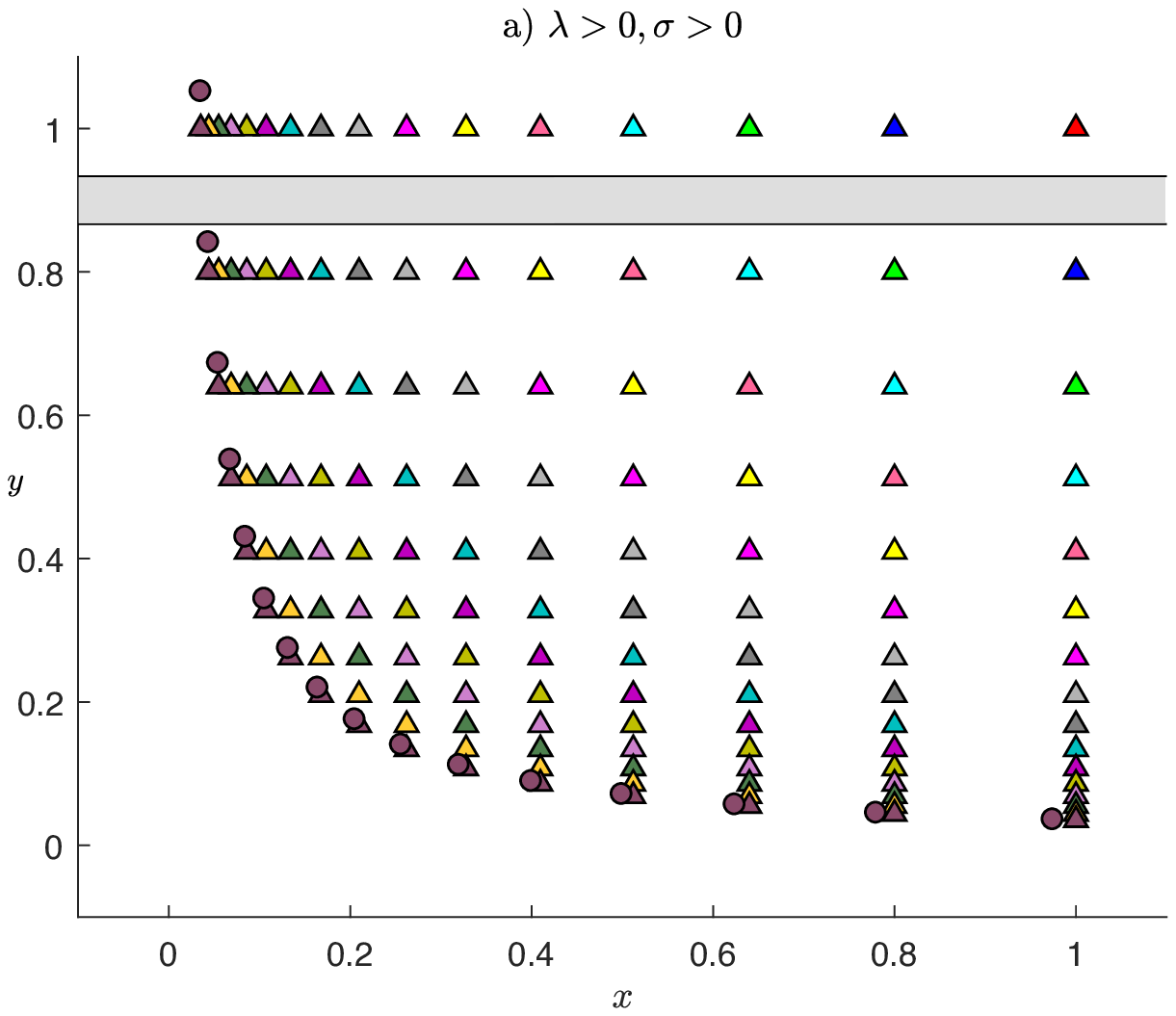}&
		
		\hspace{-0.1cm}\includegraphics[width=0.5\textwidth]{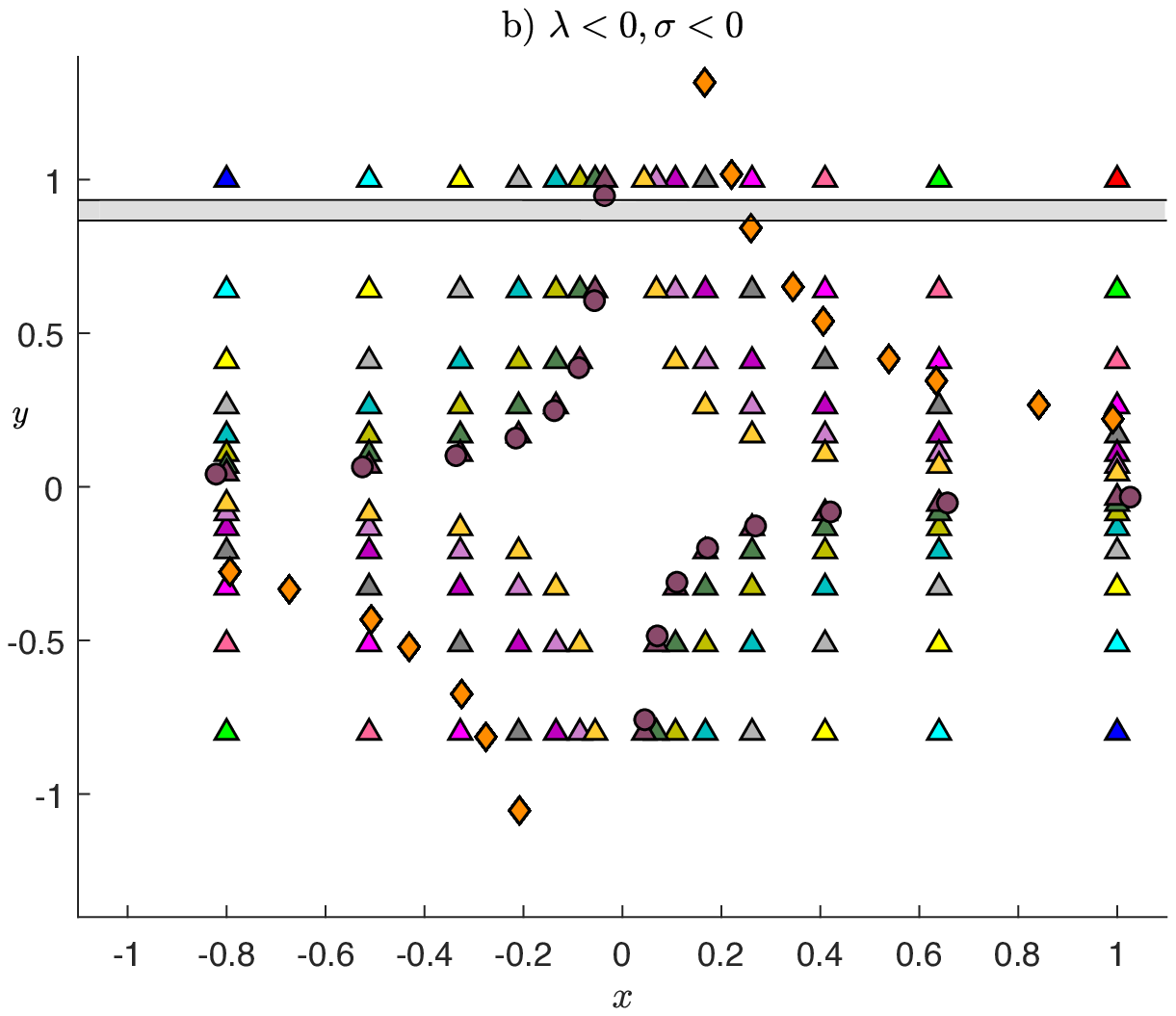}\\

		\hspace{-0.3cm}\includegraphics[width=0.5\textwidth]{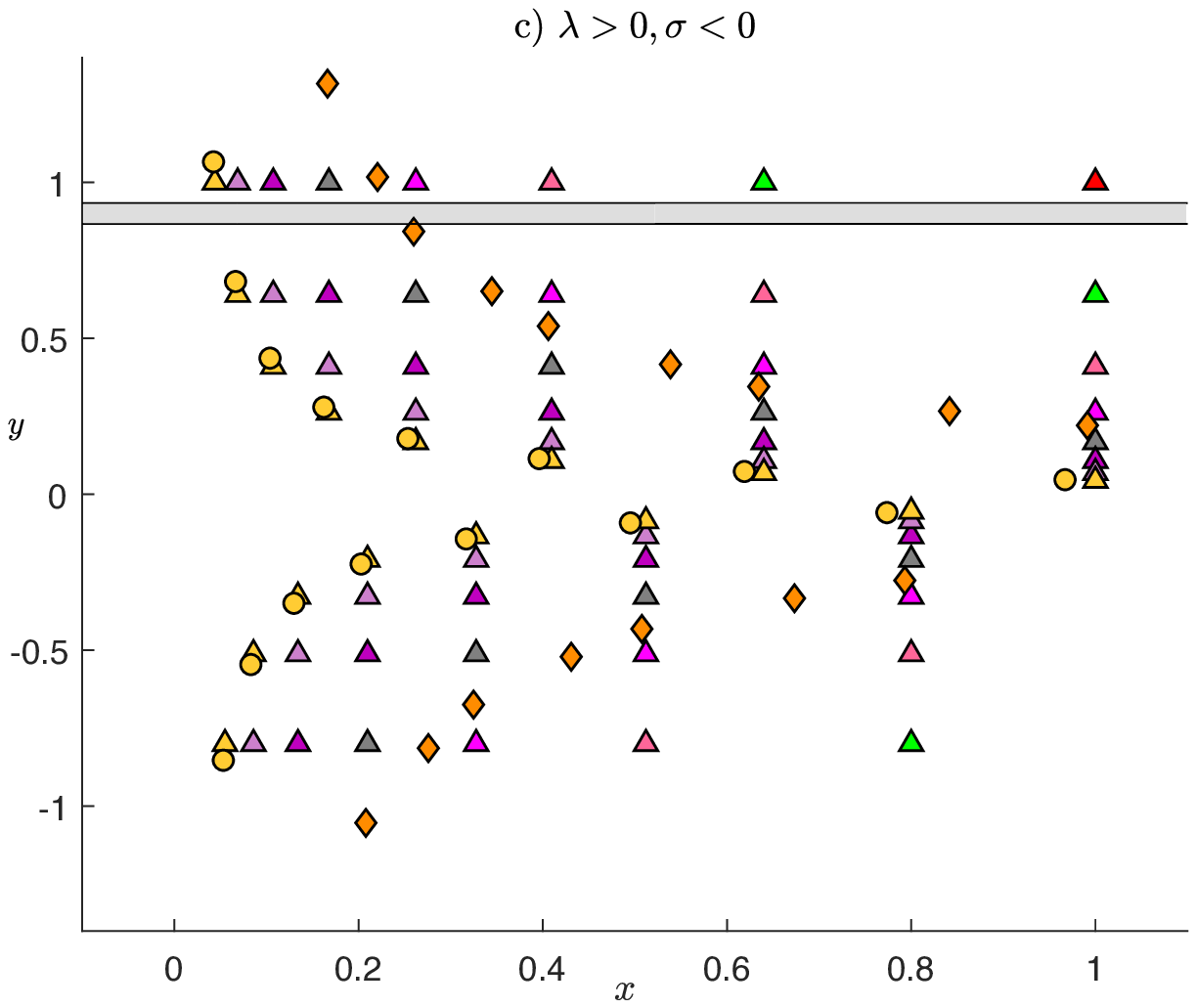}&
	
		\hspace{-0.1cm}\includegraphics[width=0.5\textwidth]{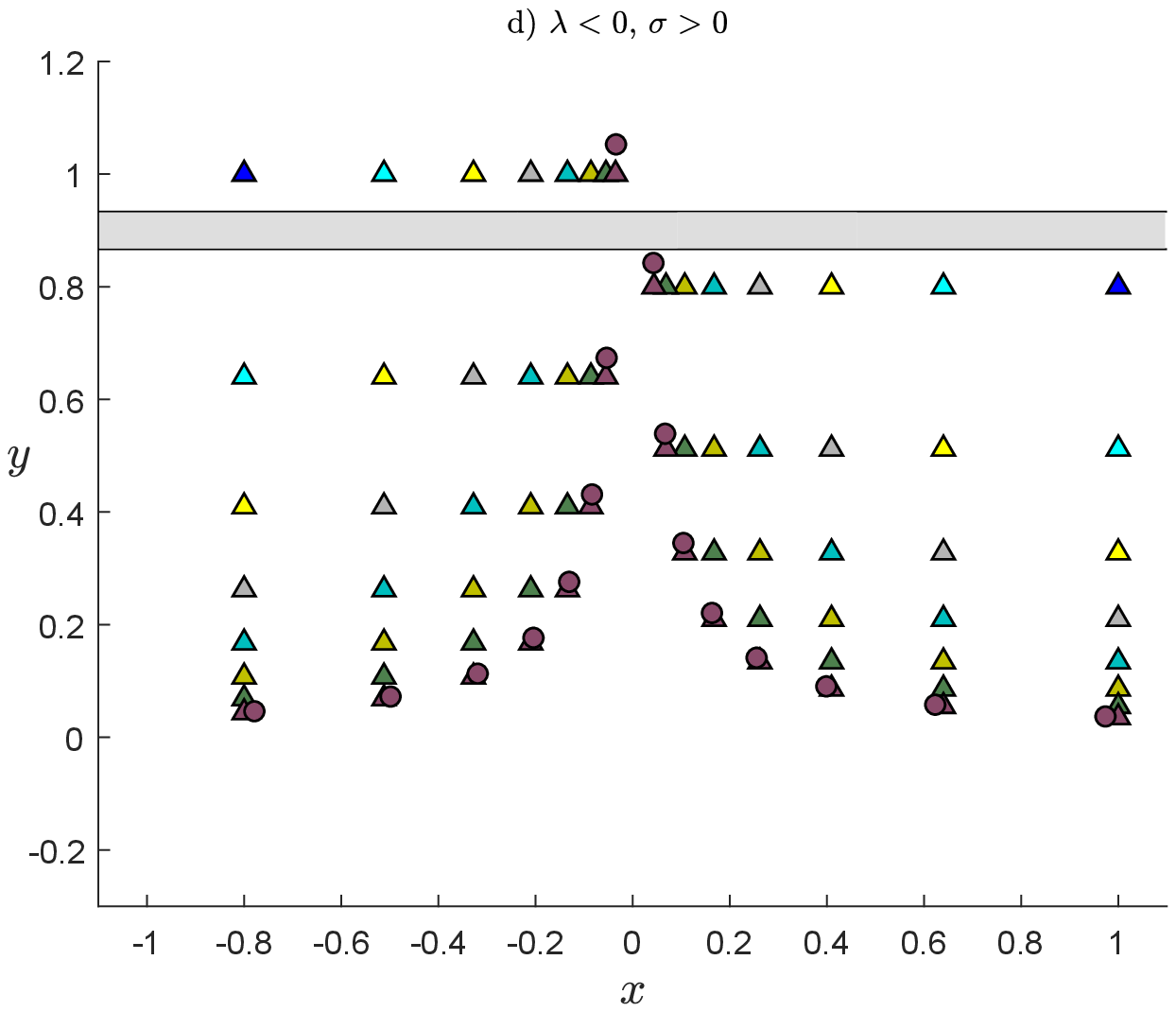}\\
		
%$c$ & $d$\\
	\end{tabular}
	\hspace*{4.5cm}\includegraphics{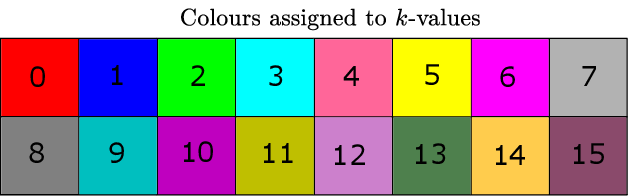}
	\vspace{0.3cm}\caption[Asymptotically stable single-round periodic solutions of map $U_{0}^{k} \circ U_{1}$.]{ Asymptotically stable ${\rm SR}_k$-solutions (period $k+1$ solution) of \eqref{eq:fEx} with \eqref{eq:xi0xi1},\eqref{eq:c2d5}.
Panels (a)--(d) correspond to \eqref{eq:param1}--\eqref{eq:param4} respectively.
Points of the stable ${\rm SR}_k$-solutions are indicated by triangles
and coloured by the value of $k$ (as indicated in the key).
In panels (a) and (b) the solutions are shown for $k=0$ (a fixed point in $y > h_1$) up to $k = 15$.
In panel (c) the solutions are shown for $k = 0,2,4,\ldots,14$
and in panel (d) the solutions are shown for $k = 1,3,5,\ldots,15$.
In each panel one saddle ${\rm SR}_k$-solution is shown with circles
(with $k = 14$ in panel (c) and $k = 15$ in the other panels).
In panels (b) and (c) asymptotically stable double-round periodic solutions are shown with diamonds.
} 
	\label{fig:SRP}
\end{figure}

With \eqref{eq:param1} and \eqref{eq:param2} the map satisfies the conditions of Theorem \ref{th:sufficientConditionsPositiveCase}, while with \eqref{eq:param3} and \eqref{eq:param4} the map satisfies the conditions of Theorem \ref{th:sufficientConditionsNegativeCase}. Thus in each case \eqref{eq:fEx} has infinitely many asymptotically stable ${\rm SR}_k$-solutions
and these are shown in Fig.~\ref{fig:SRP} up to $k = 15$ (each has period $n = k+1$).
The ${\rm SR}_k$-solutions converge to the homoclinic orbit that includes the points $(0,y^{*})$ and $(x^*,0)$
and occupy only certain quadrants of the $(x,y)$-plane as determined by the signs of $\lambda$ and $\sigma$.
In each case the ${\rm SR}_k$-solutions exist and are stable for relatively low values of $k$
(in fact for all possible $k \ge 0$).
This is due to the simplicity of~\eqref{eq:fEx}; to be clear Theorems \ref{th:sufficientConditionsPositiveCase} and \ref{th:sufficientConditionsNegativeCase} only
tell us about ${\rm SR}_k$-solutions for sufficiently large values of $k$.
In panel (c) the ${\rm SR}_k$-solutions exist for even values of $k$
while in panel (d) they exist for odd values of $k$,
in accordance with Theorem \ref{th:sufficientConditionsNegativeCase}.

Theorems \ref{th:sufficientConditionsPositiveCase} and \ref{th:sufficientConditionsNegativeCase} also guarantee the existence of saddle ${\rm SR}_k$-solutions.
For each of our four cases we show one such solution in Fig.~\ref{fig:SRP}
for the largest value of $k \le 15$. These saddle solutions were straight-forward to compute numerically because their computation reduces to solving the same quadratic equation as that mentioned above.
For smaller periods some saddle solutions involve points in $h_0 < y < h_1$.
Numerical root-finding methods are needed to compute them and this is beyond the scope of the thesis.
Also panels (b) and (c) each contain an asymptotically stable
double-round periodic solution of period $16$.

\subsection{Computing basins of attraction}
To numerically compute basins of attraction, we consider a  $1000 \times 1000$ grid of $x$ and $y$ values. Each point of the grid is iterated forward under $f$ a total of $1000$ times. This produces an orbit $(x_i,y_i)$ for $i = 0,1,\ldots,1000$.  We then look at the norm of the difference between $(x_{1000},y_{1000})$ and $(x_{1000-p},y_{1000-p})$ for successively increasing values of $p = 1,2,\ldots$ (up to a maximum period we wish to consider).  If this norm is smaller than a tolerance (we used $10^{-13}$), we say that the initial point $(x_0,y_0)$ belongs to the basin of attraction of a periodic orbit with period $p$.  \blue{If they converge to some periodic orbit it is colour coded according to its period.} 

\begin{figure}[!htbp]
\begin{center}
\includegraphics[width=0.95\textwidth]{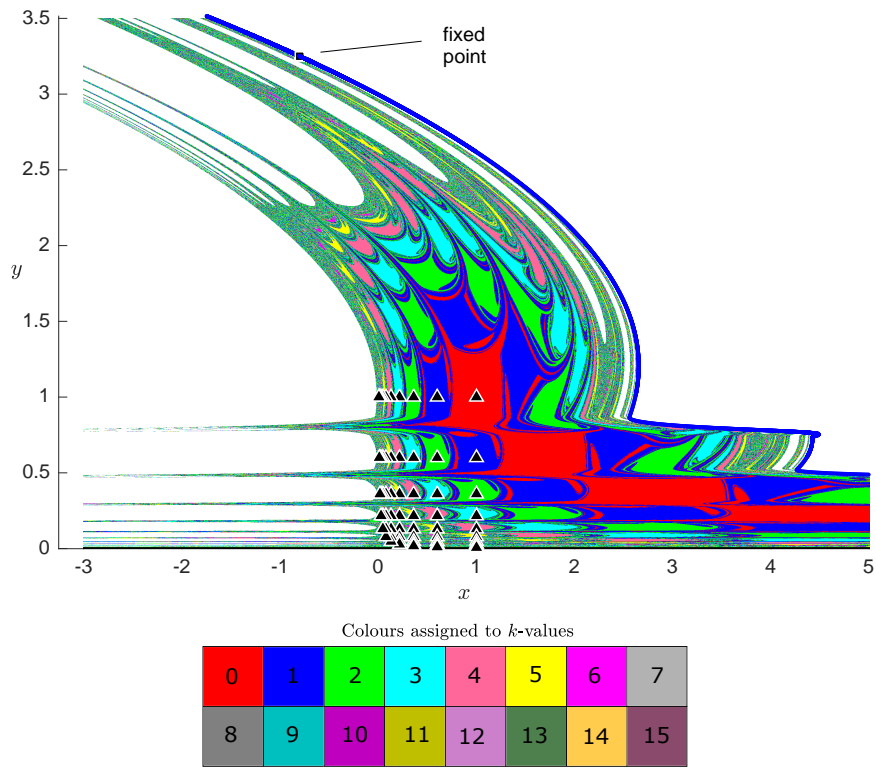}\\
\includegraphics{colour-legend-periodic.eps}
\caption[The stable manifold bounding the basins of attraction]{Basin of attraction of single-round periodic solutions from period $1$ to $16$ color coded according to the legend shown. The stable manifold of the saddle fixed point of the map $U_{1}$ appears to bound the union of the basins of attraction. The parameters are set as $\lambda = \frac{3}{5}, \sigma = \frac{5}{3}, d_{1}=1, c_{2}=-0.8, d_{5} = 0.8$.}
\label{fig:BasinManifoldU1}
\end{center}
\end{figure}

If the norm is greater than a threshold (we used $10^2$), we say that the initial point $(x_{0},y_{0})$ is diverging \blue{and we colour it white}. If it is not converging to the set of periodic orbits considered nor is diverging, it is coloured black. Fig. \ref{fig:BasinManifoldU1} shows the plot of basins of attraction of single-round periodic solutions from $k=0$ to $k=15$. We also observe that the stable manifold of the saddle fixed point at $(-0.8,3.25)$ appears to bound the basins of attraction of all the single-round periodic solutions. We believe that if we grow this stable manifold outwards further, it would fill in all of the boundary of the white region, particularly in the left half-plane $x < 0$.

\begin{figure}[!htbp]
%		\hspace{-1.5cm}
	\begin{tabular}{c c}
					
		\hspace{-0.5cm}\includegraphics[width=0.55\textwidth]{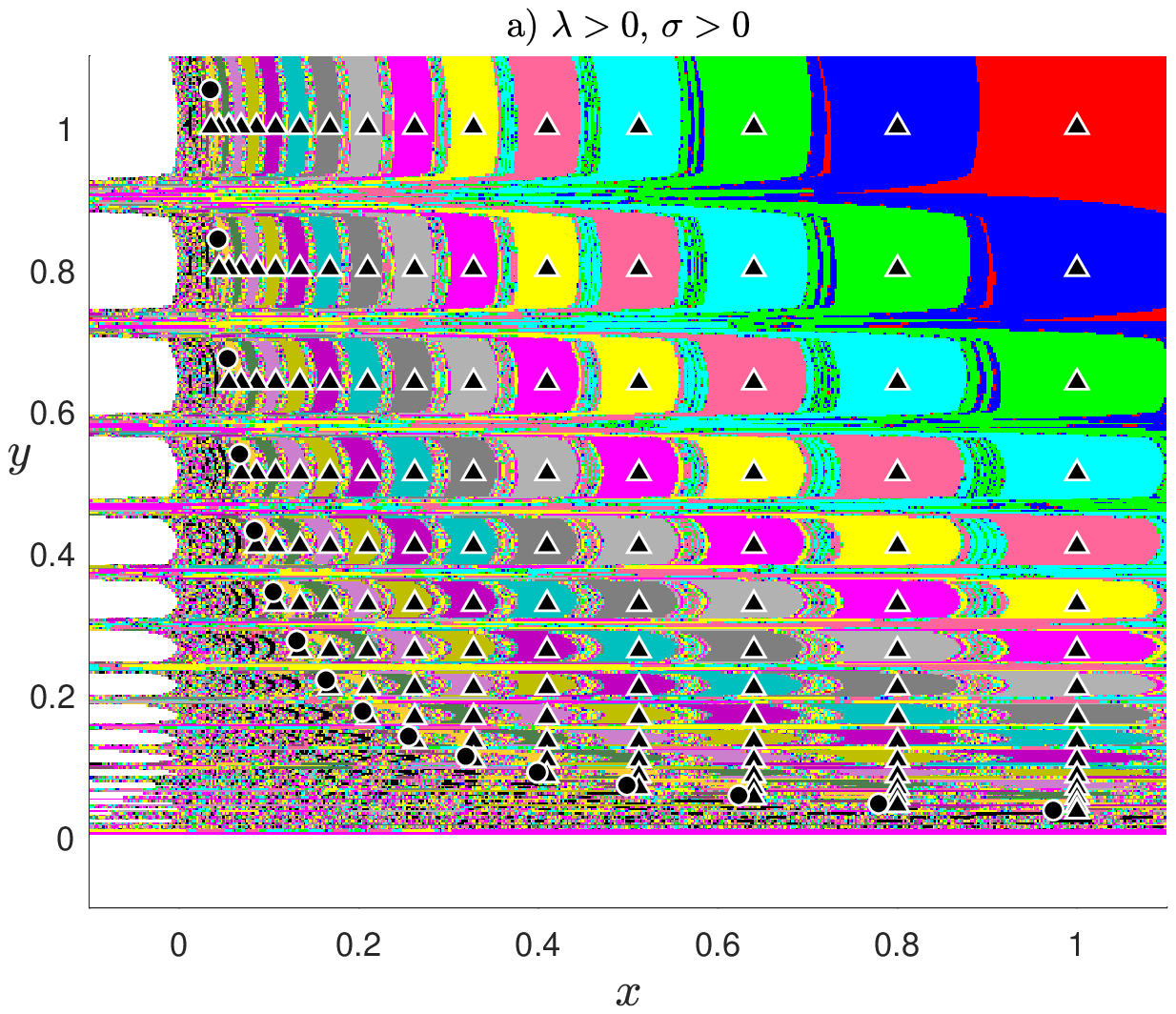}&
		
		\hspace{-1.3cm}\includegraphics[width=0.55\textwidth]{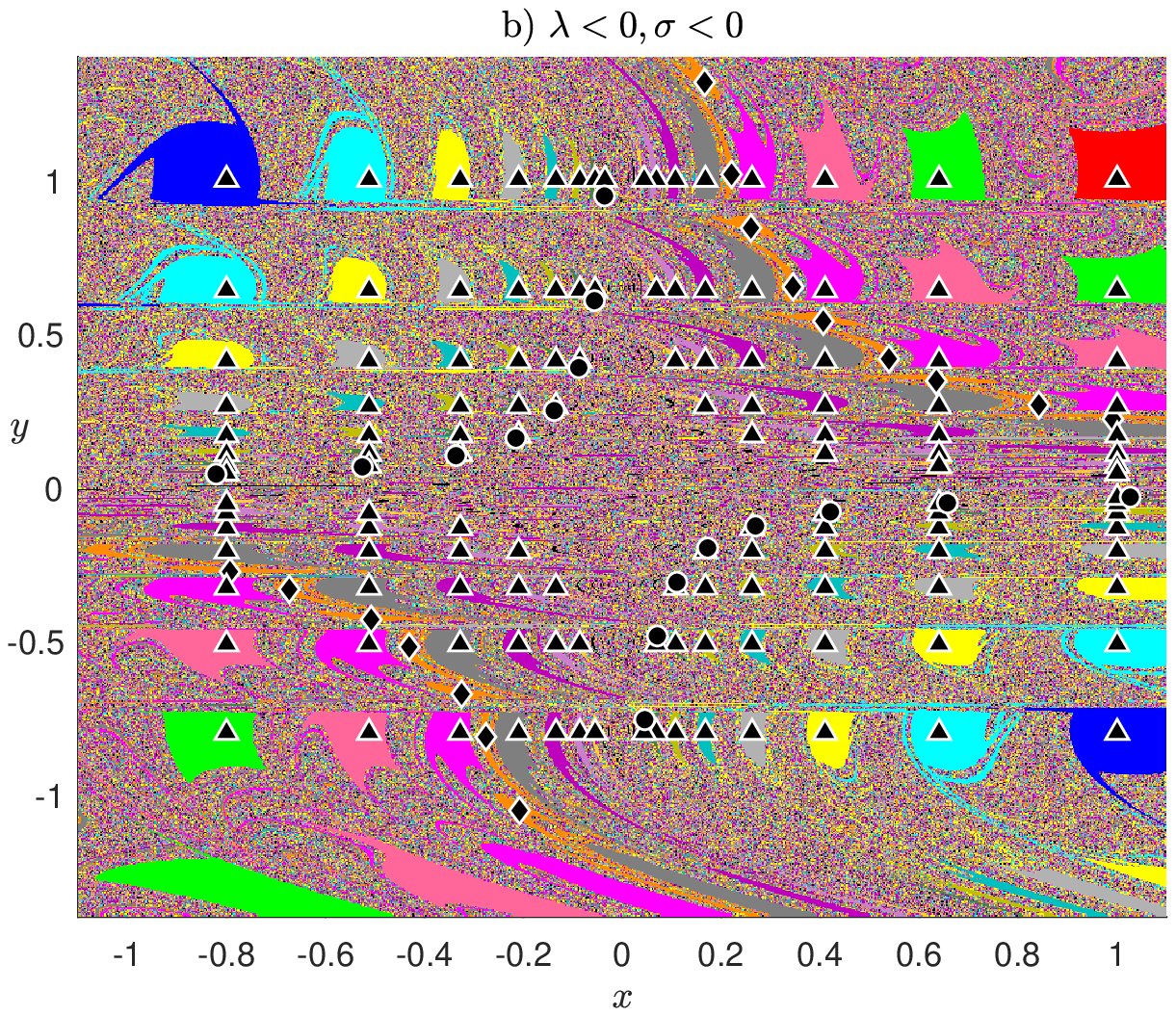}\\
%$a$ & $b$ \\
	
		\hspace{-0.5cm}\includegraphics[width=0.55\textwidth]{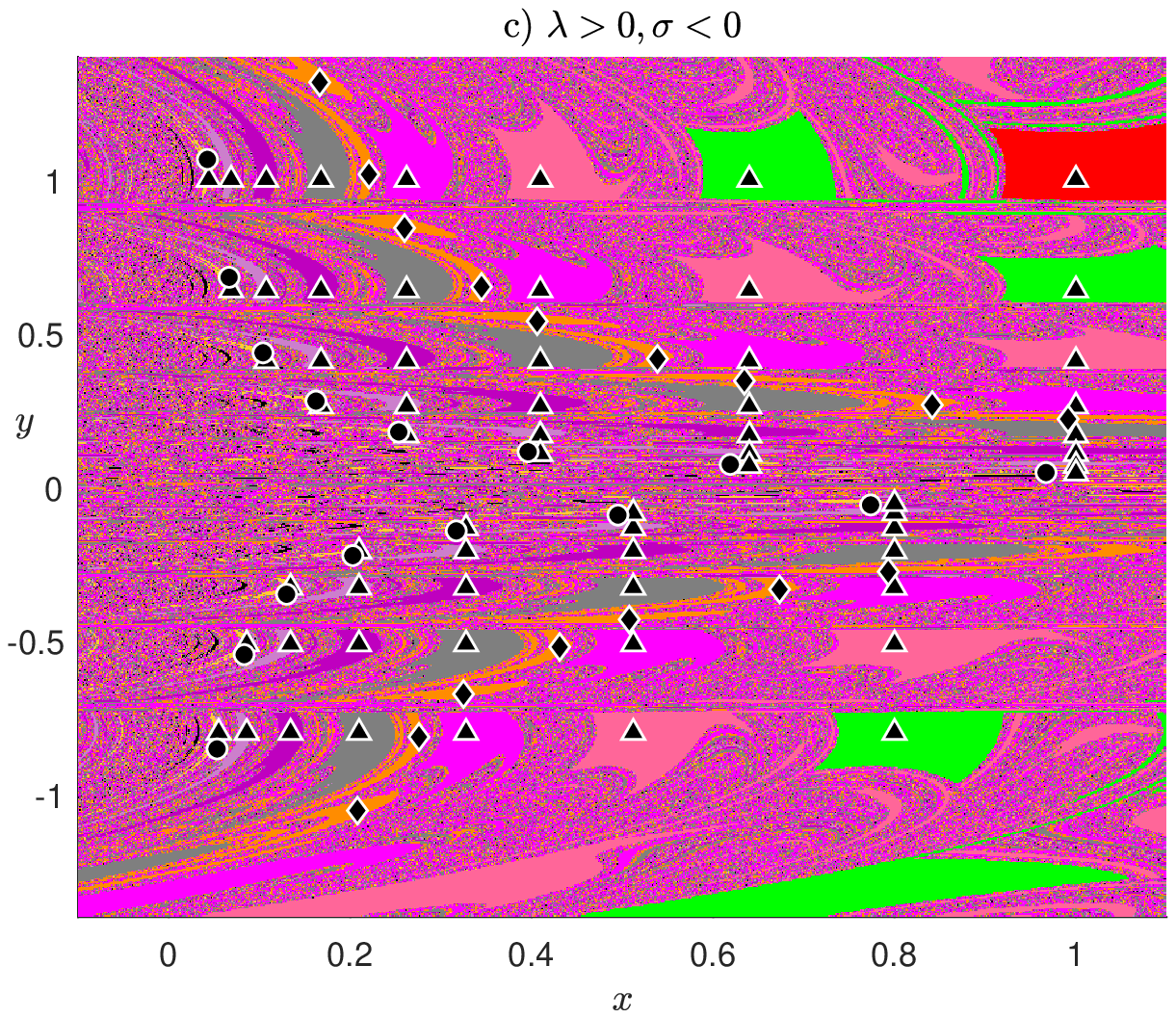}&
		\hspace{-1.3cm}\includegraphics[width=0.55\textwidth]{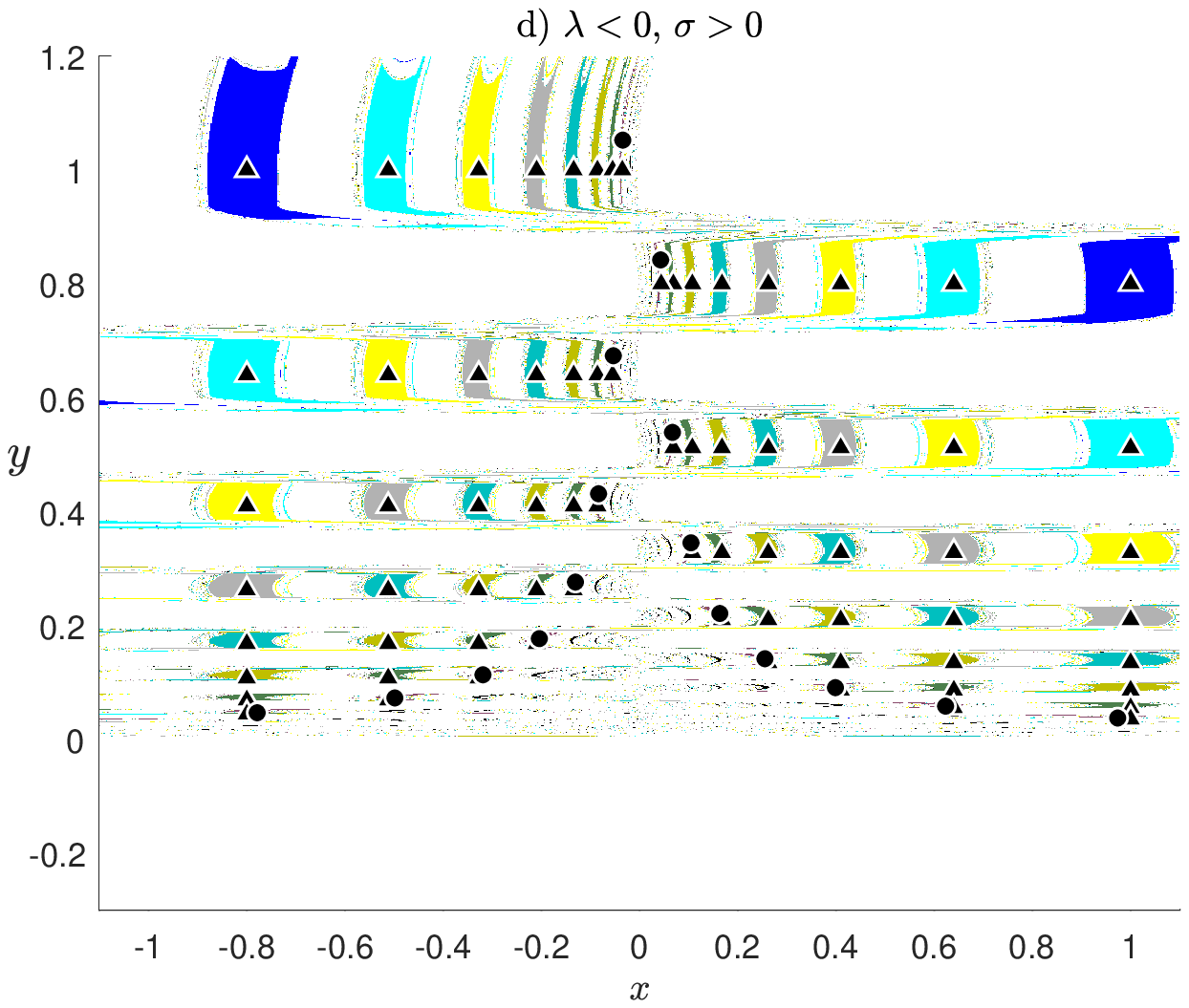}\\
		
%$c$ & $d$\\
	\end{tabular}
    \hspace*{4.5cm}\includegraphics{colour-legend-periodic.eps}
	\vspace*{0.3cm}\caption[Basins of attraction of single round periodic solutions.]{ Basins of attraction for the asymptotically stable ${\rm SR}_k$-solutions (period $k+1$ solutions) shown in Fig \ref{fig:SRP}. Specifically each point in a $1000 \times 1000$ grid is coloured by that of the ${\rm SR}_k$-solution to which its forward orbit under $f$ converges to.} 
	\label{fig:Basin}
\end{figure}

We next compute the basins of attraction for the parameter combinations in \eqref{eq:c2d5}--\eqref{eq:param4}. These are shown in Fig.~\ref{fig:Basin}.

%figures of tangency, SRP, basins goes below:
In the next section we illustrate the unfolding features near a globally resonant homoclinic tangency discussed in Chapter \ref{cha:UnfGRHT}.
\newpage
\section{Perturbations from a globally resonant homoclinic tangency}
\label{sec:numerics}

Here we introduce parameters to the example given in \S \ref{sec:example}. Specifically we consider \eqref{eq:fEx} with now

\begin{align}
U_0(x,y) &= \begin{bmatrix}
(\blue{\lambda} + \mu_2) x \big( 1 + (a_{1,0} + \mu_4) x y \big) \\
\blue{\frac{1}{\lambda}} y (1 - a_{1,0} x y)
\end{bmatrix},\\
%\label{eq:U0}
U_1(x,y) &= \begin{bmatrix}
1 + c_{2,0} (y-1) \\
\mu_1 + (1 + \mu_3) x + d_{5,0} (y-1)^2
\end{bmatrix},\\
h_0 &= \frac{2 |\lambda| + 1}{3},\\
h_1 &= \frac{|\lambda| + 2}{3}.
\end{align}
%\label{eq:r}
Below we fix
\begin{equation}
\begin{split}
\blue{\lambda} &= 0.8, \\
a_{1,0} &= 0.2, \\
c_{2,0} &= -0.5, \\
d_{5,0} &= 1,
\end{split}
\label{eq:toyParameters}
\end{equation}
and vary $\mu = (\mu_1,\mu_2,\mu_3,\mu_4) \in \mathbb{R}^4$.

%%%%%%%%%%%%%%%%%%%%%%%%%%%%%%%%%%%%%%%%%%%%%%%%%%%%%%%%%%%%%%%%%%%%%%%%%%%%%%%%
\begin{figure}
\begin{center}
\includegraphics[width=12cm]{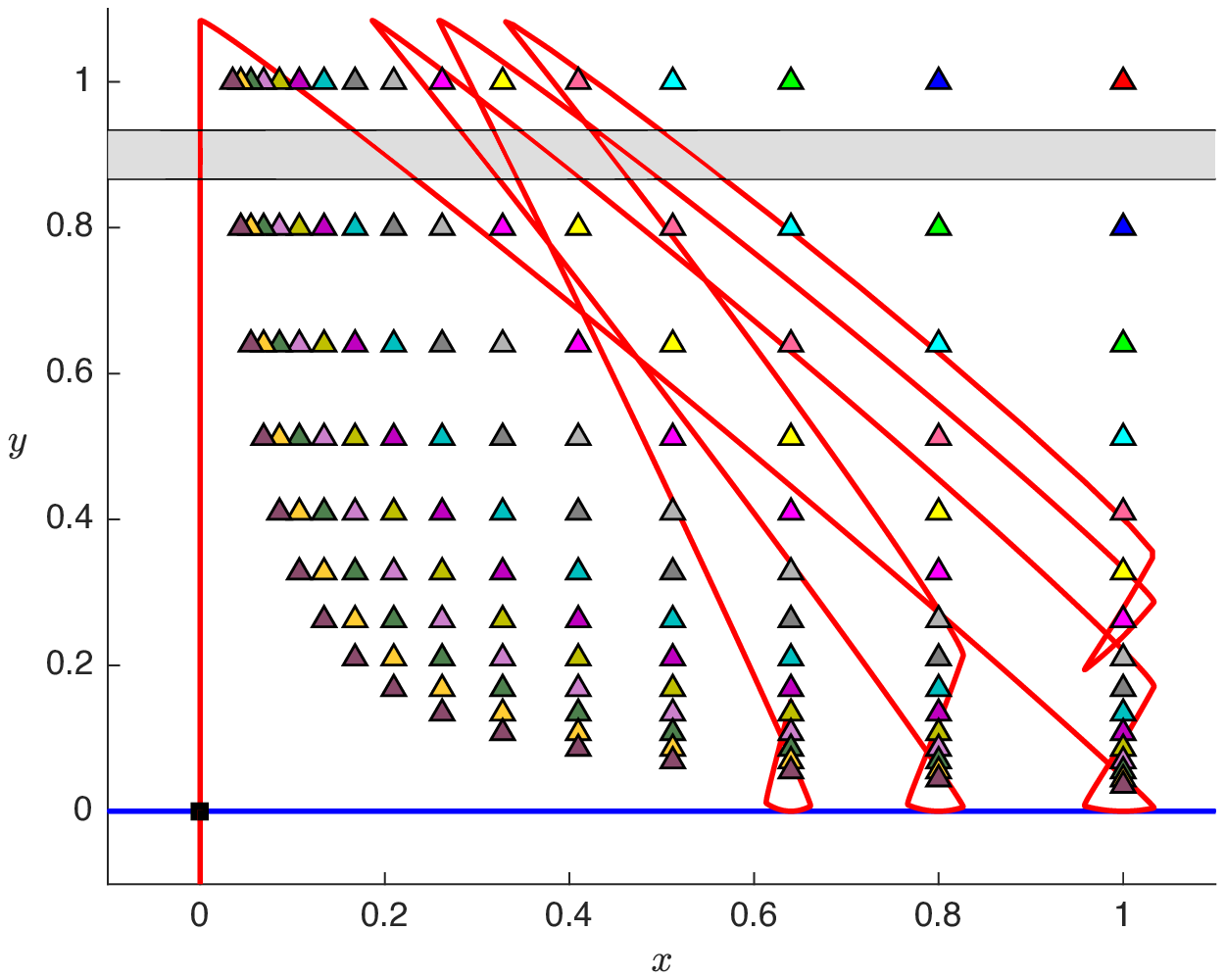}\\
\includegraphics{colour-legend-periodic.eps}
\caption[Phase portrait of the map at a homoclinic tangency]{A phase portrait of \eqref{eq:fEx} \blue{with} $\mu = \bO$.  The shaded horizontal strip is where the middle component of \eqref{eq:fEx} applies.  We show parts of the stable and unstable manifolds of $(x,y) = (0,0)$.  Note the unstable manifold has very high curvature at $(x,y) \approx (0,1.1)$ because \eqref{eq:fEx} is highly nonlinear in the horizontal strip.  For the given parameter values \eqref{eq:fEx} has an asymptotically stable, single-round periodic solutions of period $k+1$ for all $k \ge 1$.  These are shown for $k = 1,2,\ldots,15$; different colours correspond to different values of $k$.  The map also has an asymptotically stable fixed point at $(x,y) = (1,1)$. \blue{The parameters are fixed as $a_{1,0} = b_{1,0} = 0, \alpha = 0.8, c_{2,0} = -0.5, d_{5,0} = 1, \mu_{3} = 0, \mu_{4} = 0.$}
\label{fig:HTSRPUNF}
} 
\end{center}
\end{figure}
%%%%%%%%%%%%%%%%%%%%%%%%%%%%%%%%%%%%%%%%%%%%%%%%%%%%%%%%%%%%%%%%%%%%%%%%%%%%%%%%

With $\mu = \bO$, \eqref{eq:fEx} satisfies the conditions of Theorem \ref{th:sufficientConditionsPositiveCase}. In particular $\Delta_0 = 2.25$ so $\Delta_0 > 0$ and $-1 < c_{2,0} < 1 - \frac{\Delta_0}{2}$.  Therefore \eqref{eq:fEx} has an asymptotically stable single-round periodic solution for all sufficiently large values of $k$.
In fact these exist for all $k \ge 1$, see Fig.~\ref{fig:HTSRPUNF}, plus there exists an asymptotically stable fixed point at $(x,y) = (1,1)$ that can be interpreted as corresponding to $k=0$.\\
\subsection{Bifurcation of periodic solutions}
\label{sec:BifPeriodic}
\begin{figure}[t!]
\begin{center}
\includegraphics[width=0.6\textwidth]{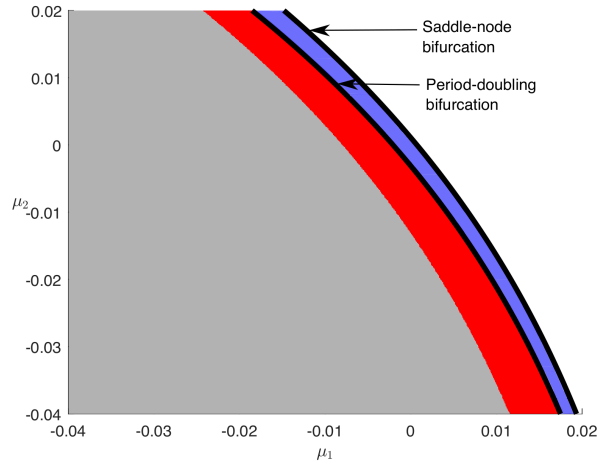}
\caption[Stability region of period-$16$ solution in the $\mu_{1}-\mu_{2}$ parameter plane.]{In the $\mu_{1}-\mu_{2}$ parameter plane the region of stability of a period-$(k+1)$ solution \blue{for $k=9$} is investigated.  The region of stability is shown in blue and is bounded between the saddle-node and period-doubling bifurcation curves. \blue{The region of instability is shown in red.} In the grey region the periodic solution has at least one point in the horizontal strip that corresponds to the middle piece of \eqref{eq:fEx}. The parameters are fixed as $a_{1,0} = b_{1,0} = 0, \alpha = 0.8, c_{2,0} = -0.5, d_{5,0} = 1, \mu_{3} = 0, \mu_{4} = 0.$}
\label{fig:mu1mu2SinglePeriodic}
\end{center}
\end{figure}

\begin{figure}[t!]
\begin{center}
\includegraphics[width=0.6\textwidth]{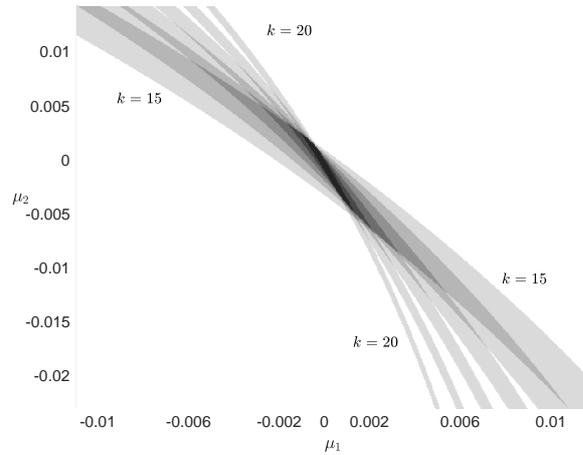}
\caption[Overlapping of stability regions from period-$16$ to period-$21$ solution in the $\mu_{1}-\mu_{2}$ parameter plane.]{In the $\mu_{1}-\mu_{2}$ parameter plane overlapping regions of stability of period-$(k+1)$ solutions for $k=15$ to $k=20$. The parameters are fixed as $a_{1,0} = b_{1,0} = 0, \alpha = 0.8, c_{2,0} = -0.5, d_{5,0} = 1, \mu_{3} = 0, \mu_{4} = 0.$}
\label{fig:mu1mu2Overlapping}
\end{center}
\end{figure}

\begin{figure}[t!]
\begin{center}
\includegraphics[width=0.6\textwidth]{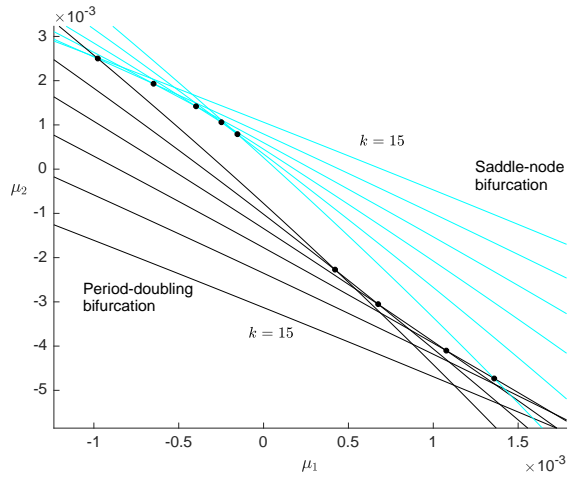}
\caption[Stability region of the periodic solutions from period $16$ to $21$ (nine-sided polygon).]{Stability region of the period-$(k+1)$ solutions for $k=15$ to $k=20$. Curves of saddle-node bifurcations are cyan; curves of period-doubling bifurcations are black. The values of the parameters are $\alpha = 0.8, a_{1,0} = b_{1,0} = 0, \mu_{3} = \mu_{4}=0, c_{2,0} = -0.5,$  and $d_{5,0}=1$.}
\label{fig:ninesided}
\end{center}
\end{figure}

\begin{figure}[t!]
\begin{center}
\includegraphics[width=0.6\textwidth]{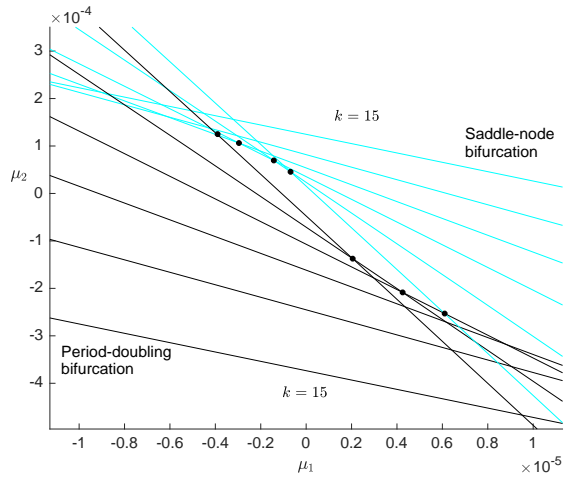}
\caption[Stability region of the periodic solutions from period $16$ to $21$ (seven-sided polygon).]{This figure shows stability regions as in Fig.~\ref{fig:ninesided} but with $\alpha = 0.7$.}
\label{fig:Sevensided}
\end{center}
\end{figure}

We compute a two-parameter bifurcation diagram showing the regions of stability and that the stable region is bounded by saddle-node and period-doubling bifurcations. For simplicity, we take the following set of parameter values as in \eqref{eq:toyParameters} but with $\mu_{3} = \mu_{4} = 0$ \blue{and neglecting resonance terms, that is $a_{1,0} = b_{1,0} = 0$. We have neglected the resonance terms for the computation of Figs.~\ref{fig:HTSRPUNF}--\ref{fig:Sevensided} just because the absence of these terms facilitates analytical calculations, computations of the stable, unstable manifolds, single-round periodic solutions, and saddle-node and period-doubling bifurcation curves. The resonance terms are considered in the next section as they are necessary during the process of unfolding.} 

Fig. \ref{fig:mu1mu2SinglePeriodic} shows the region of stability in the case of map \eqref{eq:fEx}. The blue region denotes the region of stability. The red region denotes where the periodic solution is unstable and the grey region denotes where the periodic solution has at least one point inside the horizontal strip corresponding to the middle piece of \eqref{eq:fEx}.  In this case the periodic solution does not correspond to a fixed point of $U_1 \circ U_0^k$. 
Observe that the region of stability (in blue) is bounded by the saddle-node bifurcations (in cyan) and period-doubling bifurcations (in black). The saddle-node and period-doubling bifurcations can be computed analytically due to the simplicity of the map with the chosen parameters. The saddle-node bifurcation is given by 
\begin{equation}
\mu_{1} = -\left( \frac{(a-c_{2}(\alpha + \mu_{2})^k c^k (\mu_{3}+1))^2 - (a-1)^2 + 4d_{5}((\mu_{3}+1)b-1)c^k}{4d_{5}c^{2k}} \right),
 \label{eq:SNtoyoverlap}
 \end{equation}
and the period-doubling bifurcation is given by
\begin{equation}
\mu_{1} = -\left(\frac{(a+c_{2}(\alpha + \mu_{2})^k c^k (\mu_{3} + 1) - 2)^2 - (a-1)^2 + 4d_{5}((\mu_{3}+1)b  -1)c^k}{4d_{5}c^{2k}} \right),
\label{eq:PDtoyoverlap}
\end{equation}
where 
\begin{align*}
a &= c_{2}b,\\
b &= ((\alpha + \mu_{2}) c)^k,\\
c &= \frac{1}{\alpha}.   
\end{align*}

Fig.~\ref{fig:mu1mu2Overlapping} shows where period-$(k+1)$ solutions are stable for $k = 15,16,\ldots,20$. Observe that the regions intersect with each other so in the black region all six periodic solutions exist and are stable.

The complicated shape of the overlapped black region can be seen below. Overlapping saddle-node (in cyan) and period-doubling bifurcations (in black) are computed from $k=15$ to $k=20$. For example in Fig.~\ref{fig:ninesided} the overlapped region is a nine-sided polygon (with black dots as its vertices) when $\alpha = 0.8$. With instead $\alpha = 0.7$ this region is a seven-sided polygon, Fig.~\ref{fig:Sevensided}.  These examples serve to show that regions where a certain number of stable periodic solutions coexist often have a complicated shape and with different shapes for different values of the parameters.

\subsection{Computation of bifurcation points}

In Chapter \ref{cha:UnfGRHT}, we showed that there exist sequences of saddle-node and period-doubling bifurcations near a globally resonant homoclinic tangency. Here we numerically compute these bifurcations in the map \eqref{eq:fEx}. While this could be achieved by using software, e.g.~{\sc auto} \citep{DoCh07}, \blue{~{\sc matcontm} \citep{KuMe19},} to numerically continue the periodic solutions, we found it sufficient to \blue{perform what could be called a manual bisection method} and carefully study how the phase portrait changes.  This approach is straight-forward but limited; certainly to obtain a more detailed bifurcation analysis some use of numerical continuation software would be almost essential. 

\begin{figure}[t!]
\begin{center}
\includegraphics[width=0.6\textwidth]{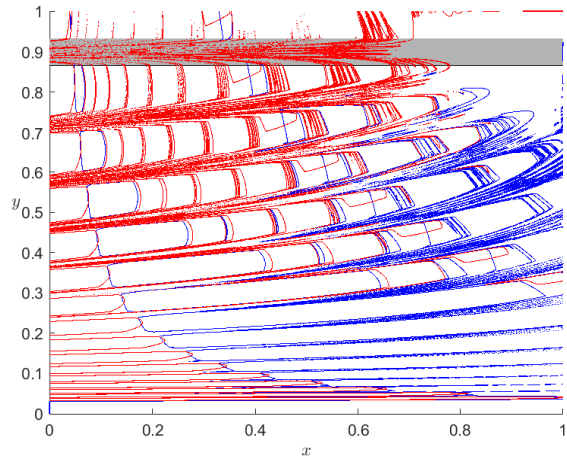}
\caption[Contour plot of the $(k+1)$th iterate of the map $f$]{Contour plot of $f^{k+1}(x,y) = (x,y)$ for $k=15$. The first component is shown in blue and the second component is shown in red. The parameters are set as $\mu_{1} = \mu_{2} = \mu_{3} = \mu_{4} = 0$ \blue{and as in \eqref{eq:toyParameters}.}}
\label{fig:ContourPlot}
\end{center}
\end{figure}

\begin{figure}[!htbp]
\begin{center}
\setlength{\unitlength}{1cm}
\begin{picture}(14,5.1)
\put(0,0){\includegraphics[width=6.8cm]{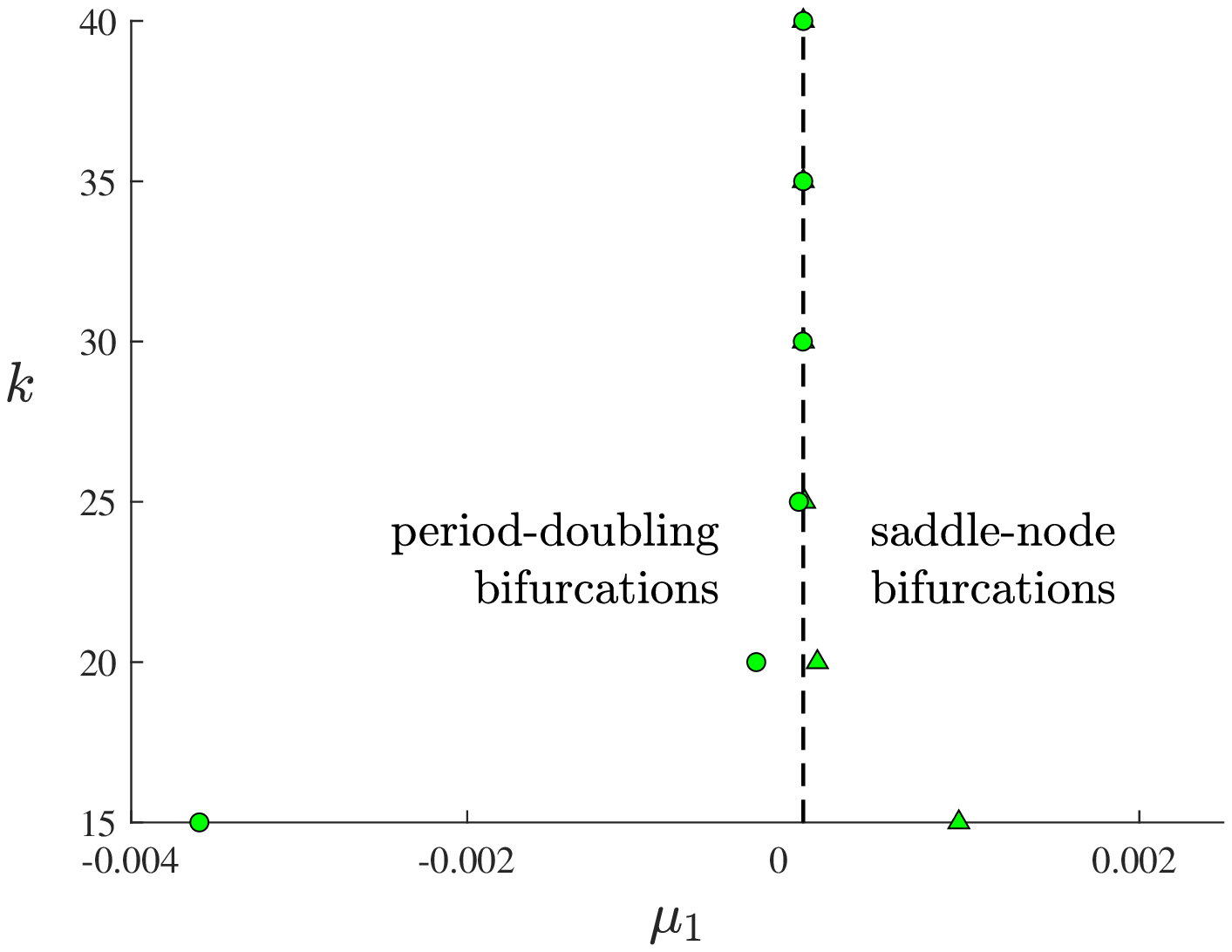}}
\put(7.2,0){\includegraphics[width=6.8cm]{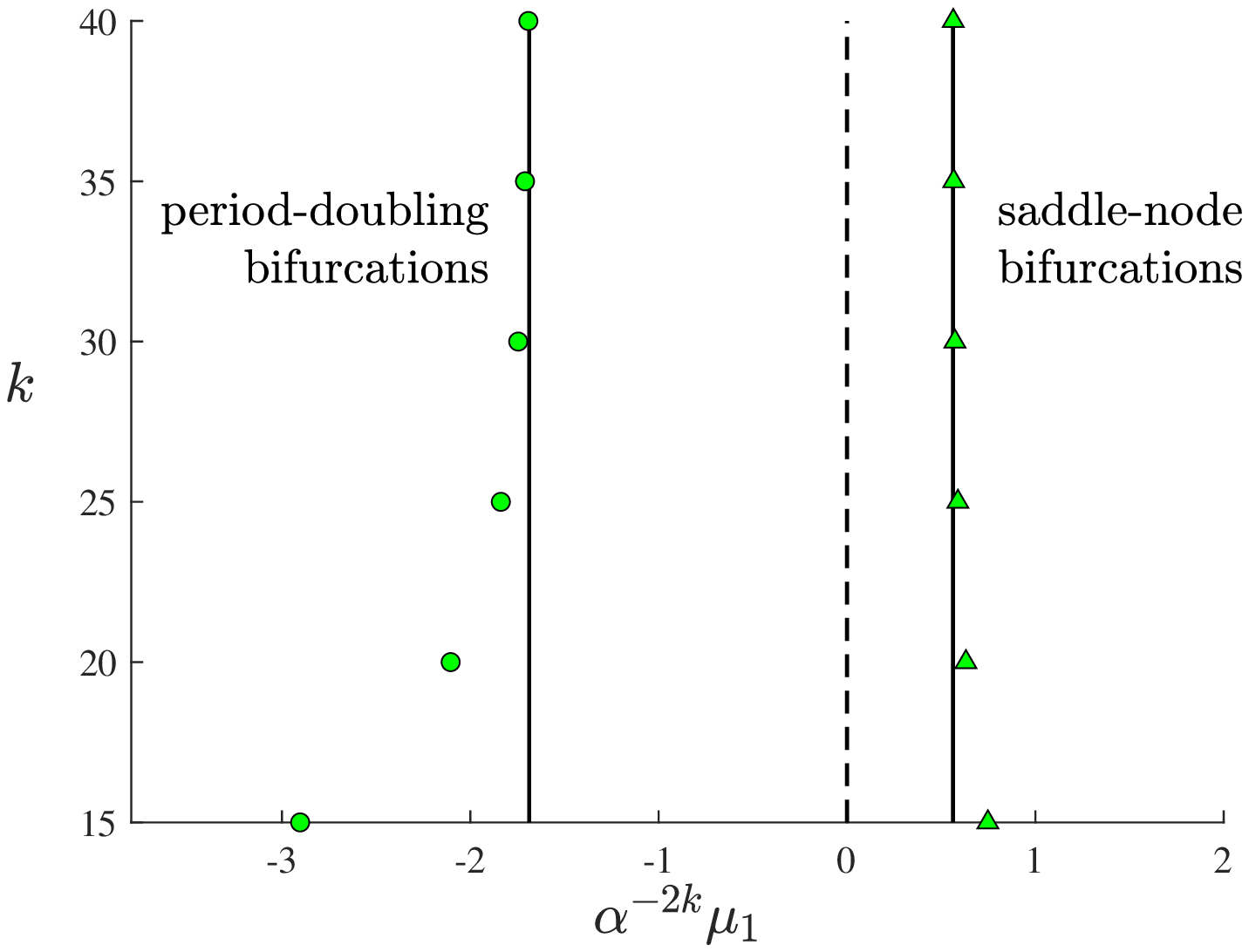}}
\put(0,5.1){\small \bf a)}
\put(7.2,5.1){\small \bf b)}
\end{picture}
\caption[Bifurcation diagram with the variation of parameter $\mu_{1}$]{Panel (a) is a numerically computed bifurcation diagram of \eqref{eq:fEx} with \eqref{eq:toyParameters} and $\mu_2 = \mu_3 = \mu_4 = 0$.  The triangles [circles] are saddle-node [period-doubling] bifurcations of single-round periodic solutions of period $k+1$.  Panel (b) shows the same points but with the horizontal axis scaled in such a way that the asymptotic approximations to these bifurcations, given by the leading-order terms in \eqref{eq:tildeeeSNCase1} and \eqref{eq:tildeeePDCase1}, appear as vertical lines.
\label{fig:A}
} 
\end{center}
\end{figure}

To compute saddle-node bifurcations we consider the contour plot of $f^{k+1}(x,y) = (x,y)$. For example let us consider $k=15$ which corresponds to a period-$16$ solution. For brevity write $(x',y') = f^{k+1}(x,y)$. Fig.~\ref{fig:ContourPlot} shows where $x' - x = 0$ in blue and where $y' - y = 0$ in red. \blue{This was} computed by using {\sc matlab}'s `contour' command \blue{
by evaluating $f^{k+1}(x,y) - (x,y)$ throughout a fine grid of $x$ and $y$ values and having the contour command use the resulting matrix of data points to fit curves to where each component is zero.} \blue{This determines the zero contours of the two components of $f^{k+1}(x,y) - (x,y)$ and can be viewed as the discrete-time analogue of the nullclines of a two-dimensional system of first-order ODEs.} Any point of intersection of the blue and red curves is a fixed point of $f^{k+1}$ and one point of a period-$(k+1)$ solution of $f$.  Fig.~\ref{fig:SNUnfold} shows the behaviour of the contours as the value of $\mu_1$ is varied and illustrates how we have detected saddle-node bifurcations. In Fig.~\ref{fig:SNUnfold}, at $\mu_{1} = 0.0005$, we observe that the contours intersect (two fixed points) and as we vary $\mu_{1}$, the two fixed points come closer and collide when a saddle-node bifurcation takes place. This happens at approximately $\mu_{1} = 0.0009257$. This point is marked in Fig.~\ref{fig:A}. Other saddle-node bifurcation points are computed similarly.

\begin{figure}[htbp!]
\begin{center}
\setlength{\unitlength}{1cm}
\hspace*{-3cm}
\begin{picture}(12,5.1)
\put(0,0){\includegraphics[width=0.45\textwidth]{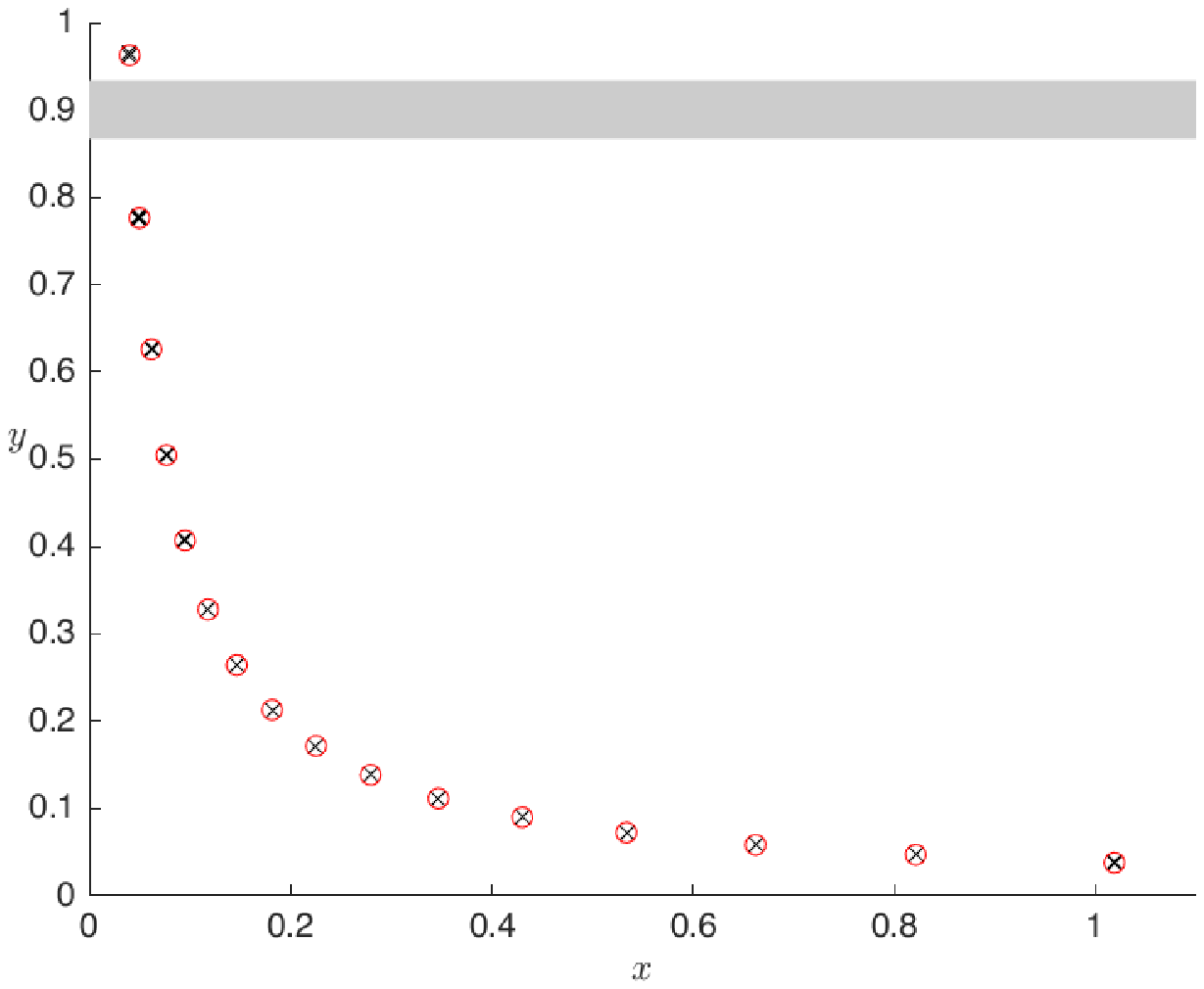}}
\put(8.2,0){\includegraphics[width=0.45\textwidth]{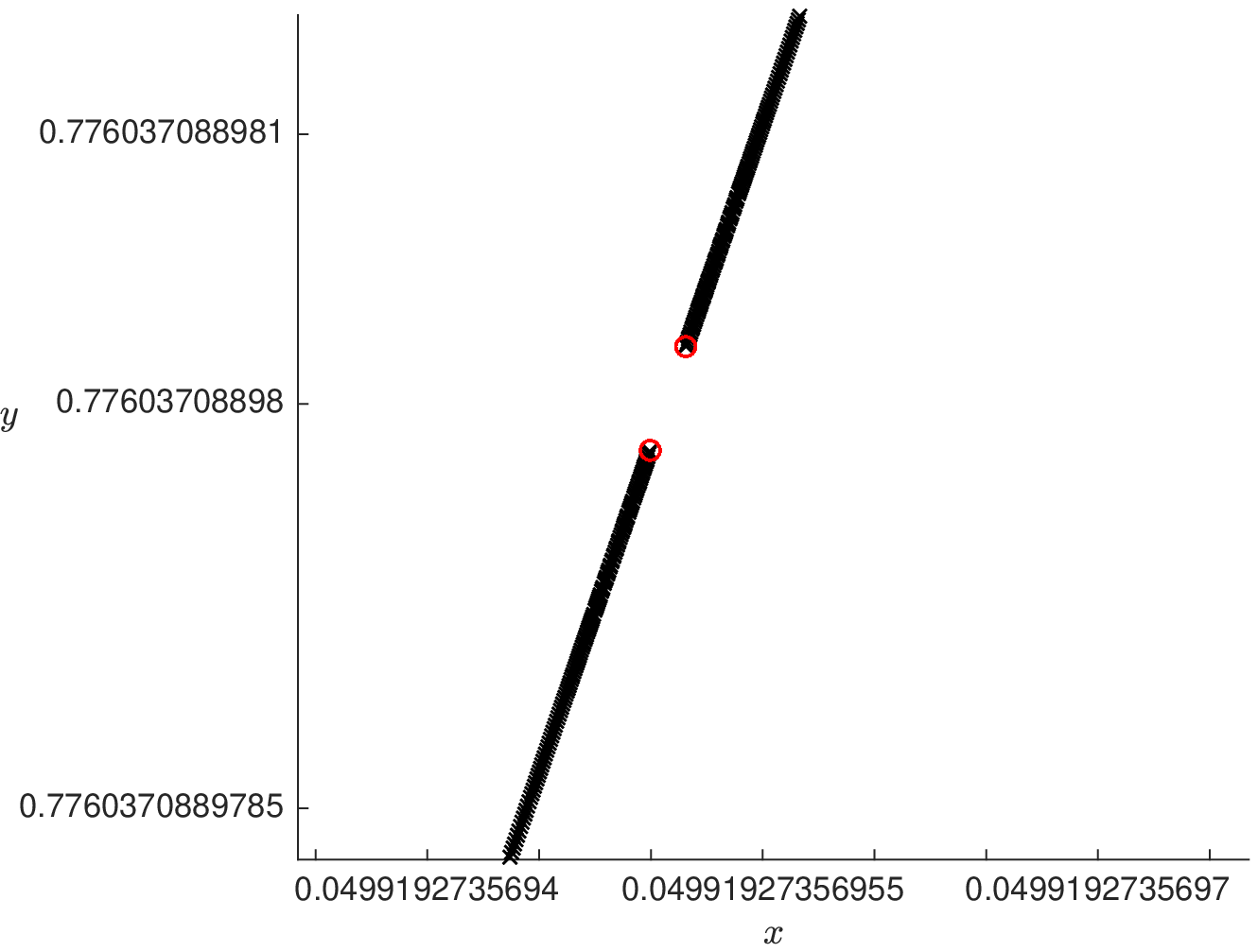}}
%\put(0,0){\small \bf a)}
%\put(8.2,0){\small \bf b)}
\put(0,5.1){\small \bf a)}
\put(7.2,5.1){\small \bf b)}
\end{picture}
\caption[Computation of period doubling bifurcation numerically.]{ In (a), an orbit is shown in \blue{red circles} with iterates shown in black \blue{crosses} converging to a stable period-$32$ orbit. In (b), a zoomed version \blue{near a periodic point is shown} which clarifies that two points are very close to each other, just after a period-$16$ orbit has undergone a period-doubling bifurcation. The red circles denote the final few iterates of the computed orbit. \blue{ In panel (b), the sequences of black crosses are so close that they appear to be straight lines but are discrete set of close points marked by black crosses, converging to the two periodic points approximated by the red circles.} \blue{The parameters are set as $\mu_{2} =\mu_{3} =\mu_{4} = 0$  and as in \eqref{eq:toyParameters}.}}
	\label{fig:PDUnfold}
\end{center}
\end{figure}

\iffalse
\begin{figure}[htbp!]
\begin{center}
\setlength{\unitlength}{1cm}
\hspace*{-3cm}
\begin{picture}(12,4.1)
\put(0,0){\includegraphics[width=0.35\textwidth]{IntnContours.eps}}
\put(5.2,0){\includegraphics[width=0.35\textwidth]{BifPointContours.eps}}
\put(10.1,0){\includegraphics[width=0.35\textwidth]{NoIntnContours.eps}}
\put(0,0){\small \bf a)}
\put(5.2,0){\small \bf b)}
\put(10.1,0){\small \bf c)}
\end{picture}
\caption[Computation of saddle-node bifurcation numerically.]{ Contour plot of $f^{k+1}(x,y) = (x,y)$ for $k=15$. The first component is marked in blue and the second component is marked in red. In (a), the two contours intersect suggesting two fixed points at parameter $\mu_{1} = 0.0005$. In (b), there is a tangential intersection between the two contours at $\mu_{1} = 0.0009257$ and hence a saddle node bifurcation is detected. In (c), there is no intersection between the two contours at $\mu_{1} = 0.002$.}
	\label{fig:SNUnfold}
\end{center}
\end{figure}
\fi
\begin{figure}[htbp!]
\begin{center}
\setlength{\unitlength}{1cm}
\hspace*{-7cm}
\begin{picture}(5.1,22)
\put(0,0){\includegraphics[width=10cm]{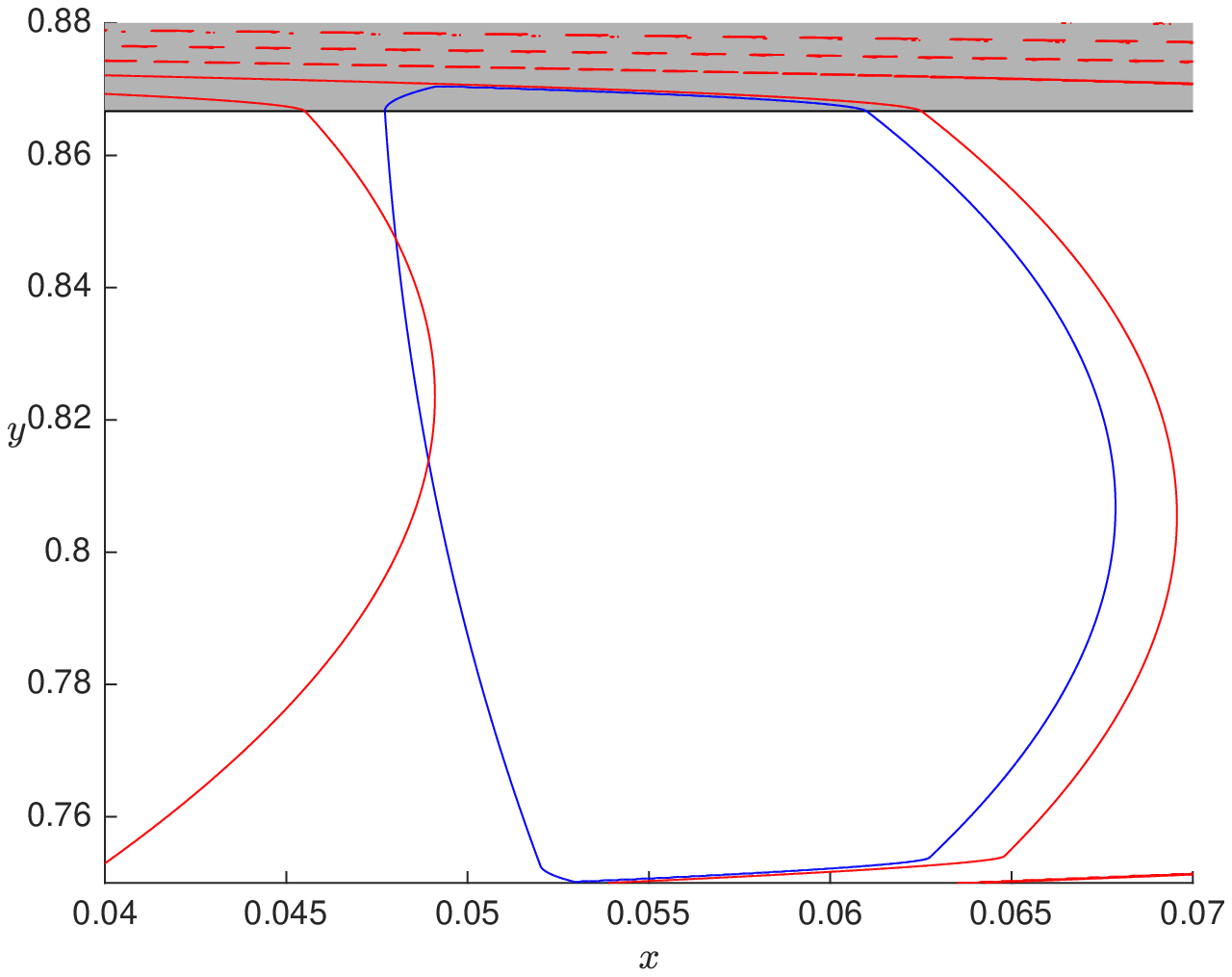}}
\put(0,7.2){\includegraphics[width=10cm]{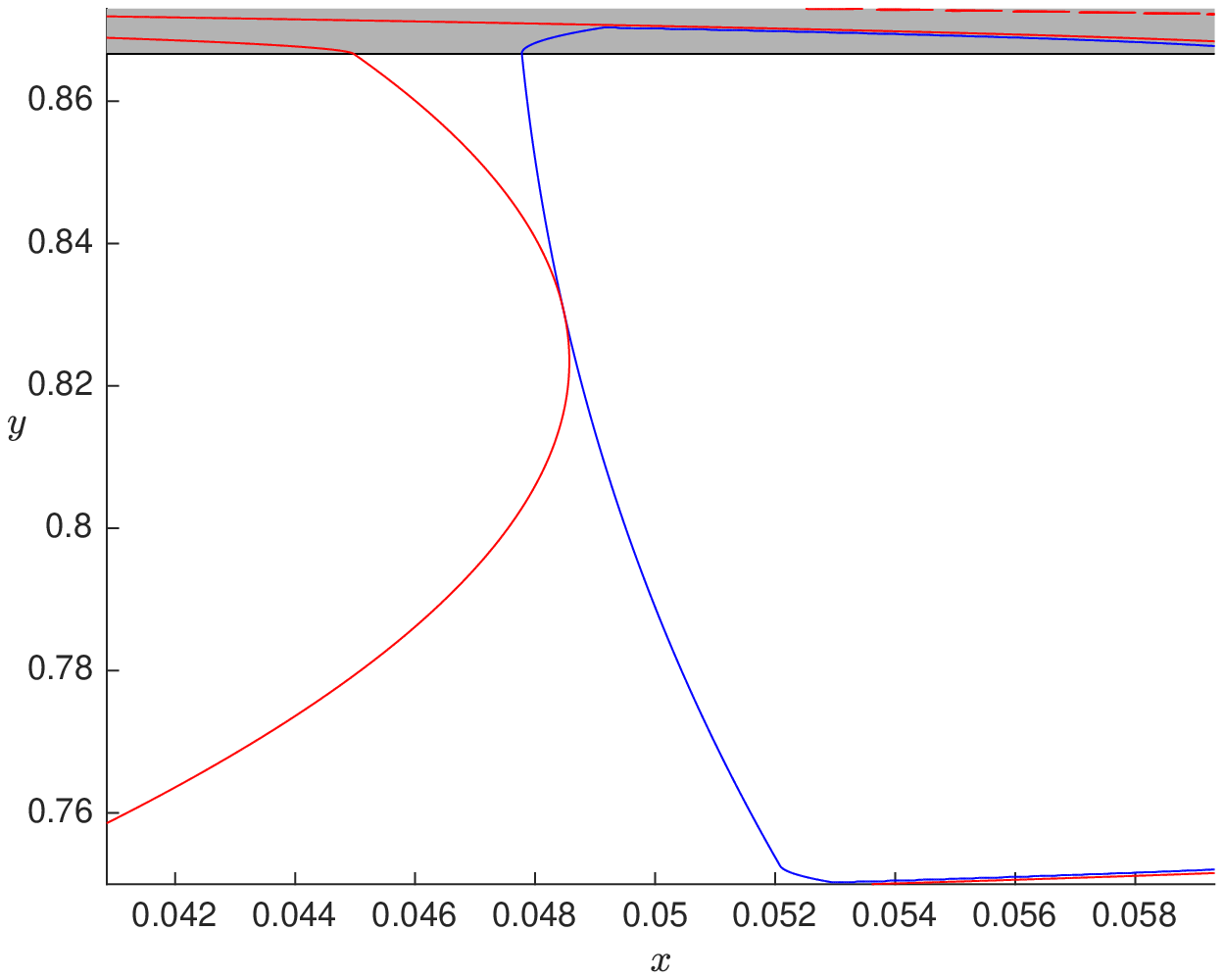}}
\put(0,14.4){\includegraphics[width=10cm]{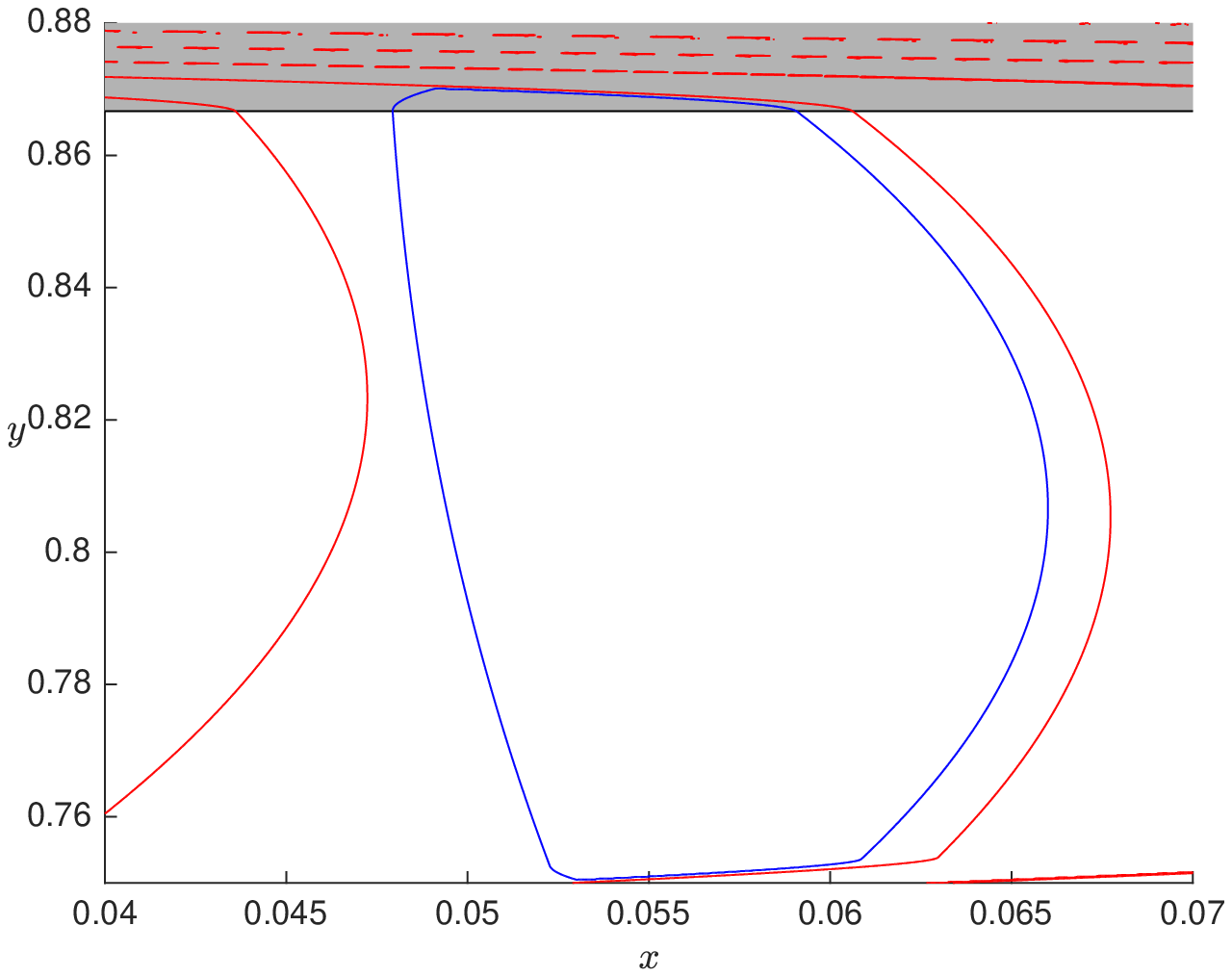}}
\put(0,7){\small \bf c)}
\put(0,14.2){\small \bf b)}
\put(0,21.4){\small \bf a)}
\end{picture}
\caption[Computation of saddle-node bifurcation numerically.]{ Contour plot of $f^{k+1}(x,y) = (x,y)$ for $k=15$. The first component is marked in blue and the second component is marked in red. In (a), the two contours intersect suggesting two fixed points at parameter $\mu_{1} = 0.0005$. In (b), there is a tangential intersection between the two contours at $\mu_{1} = 0.0009257$ and hence a saddle node bifurcation is detected. In (c), there is no intersection between the two contours at $\mu_{1} = 0.002$. \blue{The parameters are set as $\mu_{2} =\mu_{3} =\mu_{4} = 0$ and as in \eqref{eq:toyParameters}}.}
	\label{fig:SNUnfold}
\end{center}
\end{figure}
To compute period-doubling bifurcations, we start from an initial guess close to the periodic solution and as we vary parameters check that the numerically computed forward orbit is converging to the periodic solution. If it does, we compute the eigenvalues of the Jacobian matrix of the map $f^{k+1}$. When we are approaching the period doubling bifurcation, the eigenvalues approach $-1$. For example, consider $k=15$ and variations to the value of $\mu_1$. Fig.~\ref{fig:PDUnfold} shows that at $\mu_{1} = -0.003594$, we have just passed a supercritical period-doubling bifurcation as the period has doubled to $32$. It is difficult to see the doubled orbit, so we present a zoomed version near one of the periodic points in Fig.~\ref{fig:PDUnfold} (b).

%%%%%%%%%%%%%%%%%%%%%%%%%%%%%%%%%%%%%%%%%%%%%%%%%%%%%%%%%%%%%%%%%%%%%%%%%%%%%%%%
%%%%%%%%%%%%%%%%%%%%%%%%%%%%%%%%%%%%%%%%%%%%%%%%%%%%%%%%%%%%%%%%%%%%%%%%%%%%%%%%
In this remainder of this section we study how the infinite coexistence is destroyed by varying each the components of $\mu$ from zero in turn.  We identify saddle-node and period-doubling bifurcations numerically and compare these to our above asymptotic results. First, by varying the value of $\mu_1$ from zero we destroy the homoclinic tangency.
Indeed in \eqref{eq:nTang} we have ${\bf n}_{\rm tang}^{\sf T} = [1,0,0,0]$.
Thus if we fix $\mu_2 = \mu_3 = \mu_4 = 0$ and vary the value of $\mu_1$,
by Theorem \ref{th:main} there must exist sequences of saddle-node and period-doubling bifurcations
occurring at values of $\mu_1$ that are asymptotically proportional to $\alpha^{2 k}$.
Fig.~\ref{fig:A} (a) shows the bifurcation values (obtained numerically) for six different values of $k$.
We have designed \eqref{eq:fEx} so that it satisfies \eqref{eq:d0Exp}.  Consequently the formulas \eqref{eq:tildeeeSNCase1} and \eqref{eq:tildeeePDCase1} for the bifurcation values can be applied directly.  In panel (b) we observe that the numerically computed bifurcation values indeed  converge to their leading-order approximations.

%%%%%%%%%%%%%%%%%%%%%%%%%%%%%%%%%%%%%%%%%%%%%%%%%%%%%%%%%%%%%%%%%%%%%%%%%%%%%%%%
\begin{figure}[!htbp]
\begin{center}
\setlength{\unitlength}{1cm}
\begin{picture}(14,5.1)
\put(0,0){\includegraphics[width=6.8cm]{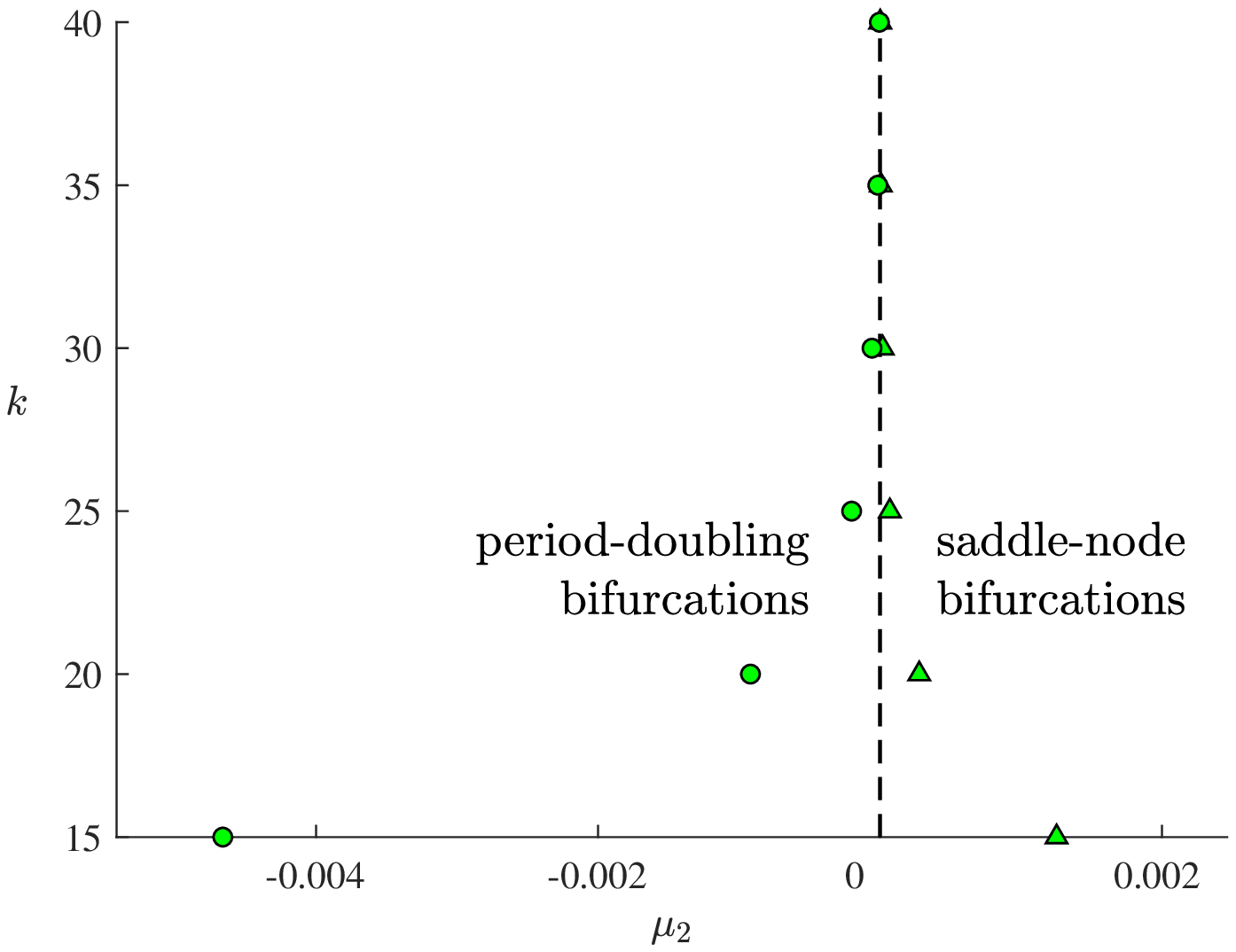}}
\put(7.2,0){\includegraphics[width=6.8cm]{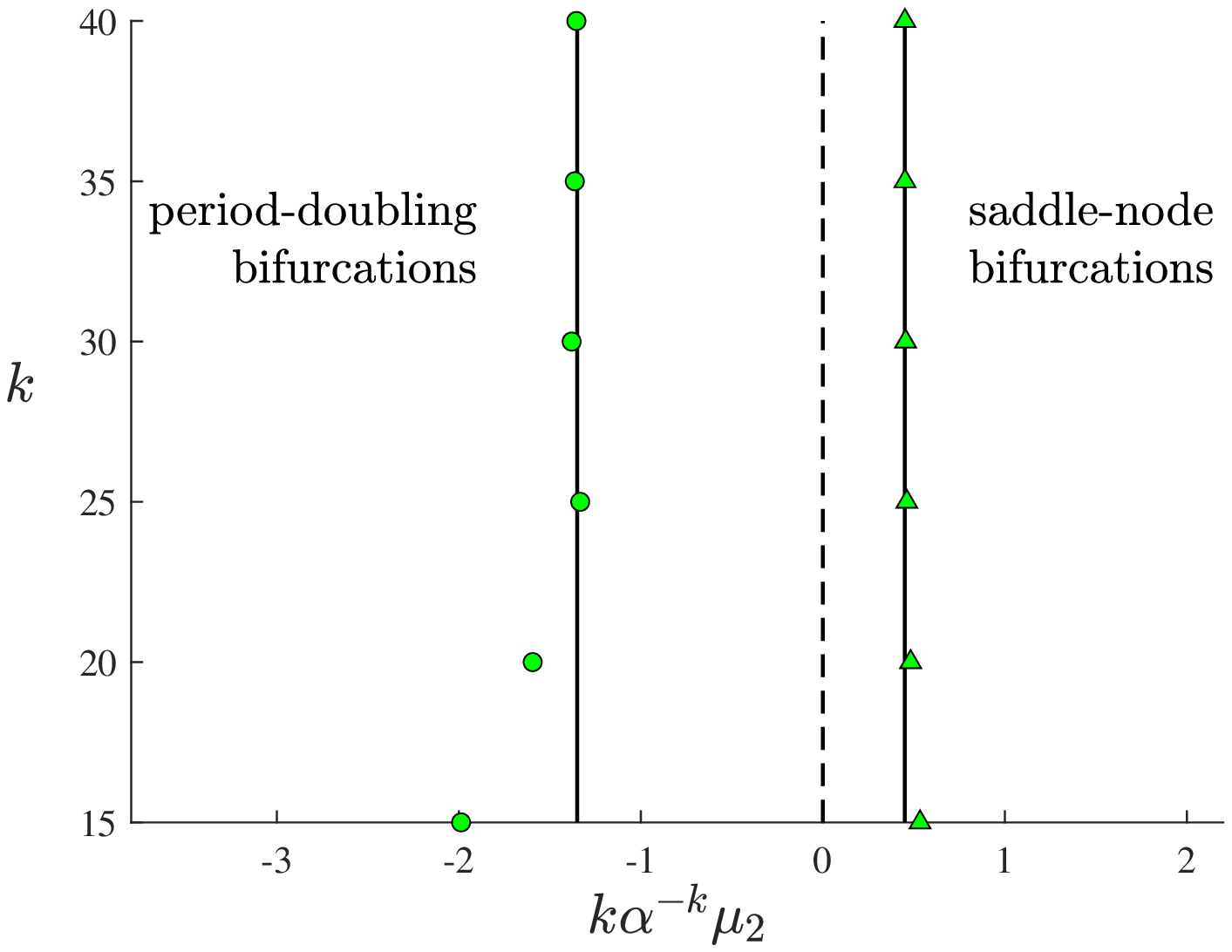}}
\put(0,5.1){\small \bf a)}
\put(7.2,5.1){\small \bf b)}
\end{picture}
\caption[Bifurcation diagram with variation of the parameter $\mu_{2}$]{Panel (a) is a numerically computed bifurcation diagram of \eqref{eq:fEx} with \eqref{eq:toyParameters} and $\mu_1 = \mu_3 = \mu_4 = 0$.  Panel (b) shows convergence to the leading-order terms of \eqref{eq:tildeeeSNCase2} and \eqref{eq:tildeeePDCase2}.
\label{fig:B}
} 
\end{center}
\end{figure}
%%%%%%%%%%%%%%%%%%%%%%%%%%%%%%%%%%%%%%%%%%%%%%%%%%%%%%%%%%%%%%%%%%%%%%%%%%%%%%%%

We now fix $\mu_1 = \mu_3 = \mu_4=0$ and vary the value of $\mu_2$.
This parameter variation alters the product of the eigenvalues associated with the origin.
Specifically ${\bf n}_{\rm eig}^{\sf T} = \left[ 0,\frac{1}{\alpha},0,0 \right]$ in \eqref{eq:nEig}
so by Theorem \ref{th:main} the bifurcation values are asymptotically proportional to $\frac{\alpha^k}{k}$.
In Fig.~\ref{fig:B} we see the numerically computed bifurcation values converging to their leading order approximations \eqref{eq:tildeeeSNCase2} and \eqref{eq:tildeeePDCase2}.

%%%%%%%%%%%%%%%%%%%%%%%%%%%%%%%%%%%%%%%%%%%%%%%%%%%%%%%%%%%%%%%%%%%%%%%%%%%%%%%%
\begin{figure}[!htbp]
\begin{center}
\setlength{\unitlength}{1cm}
\begin{picture}(14,5.1)
\put(0,0){\includegraphics[width=6.8cm]{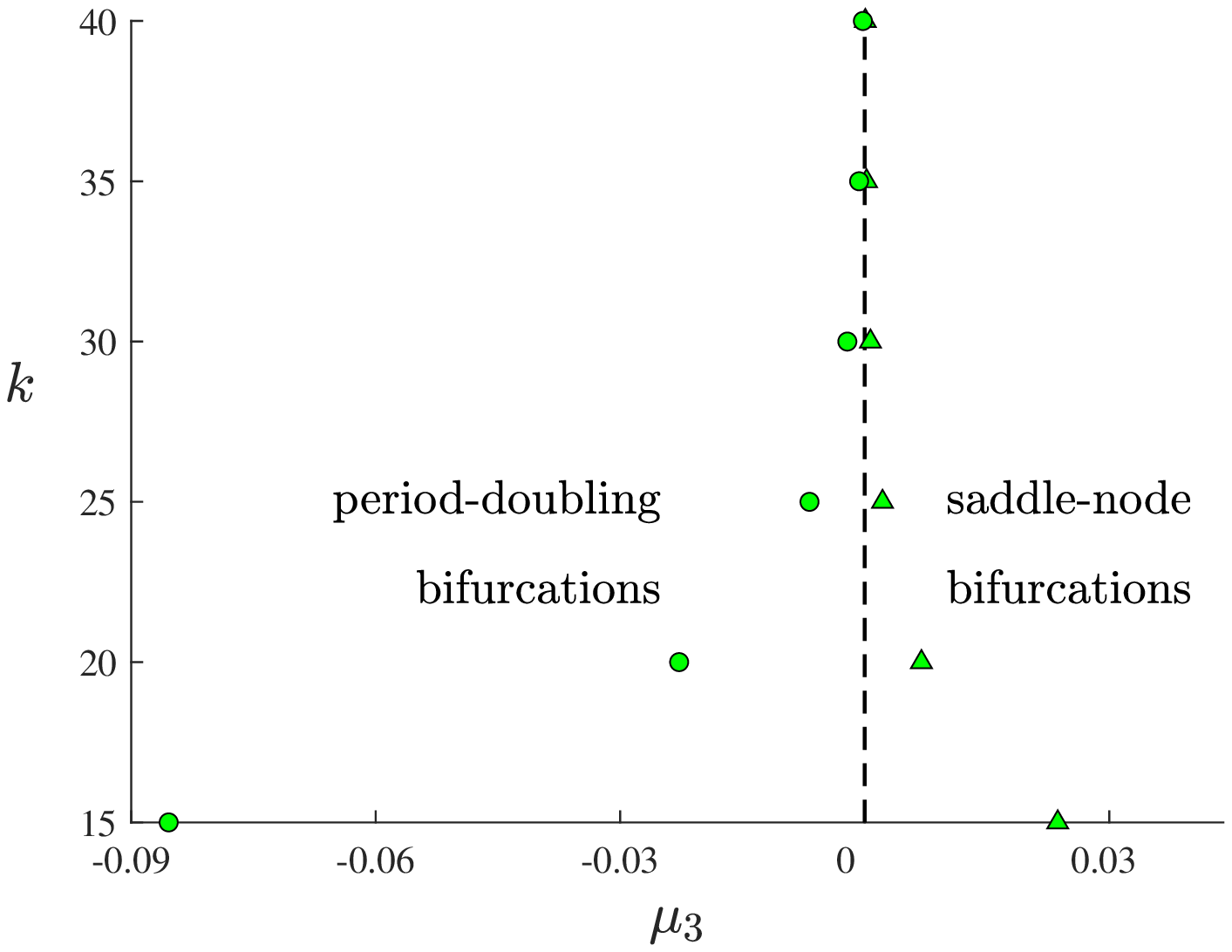}}
\put(7.2,0){\includegraphics[width=6.8cm]{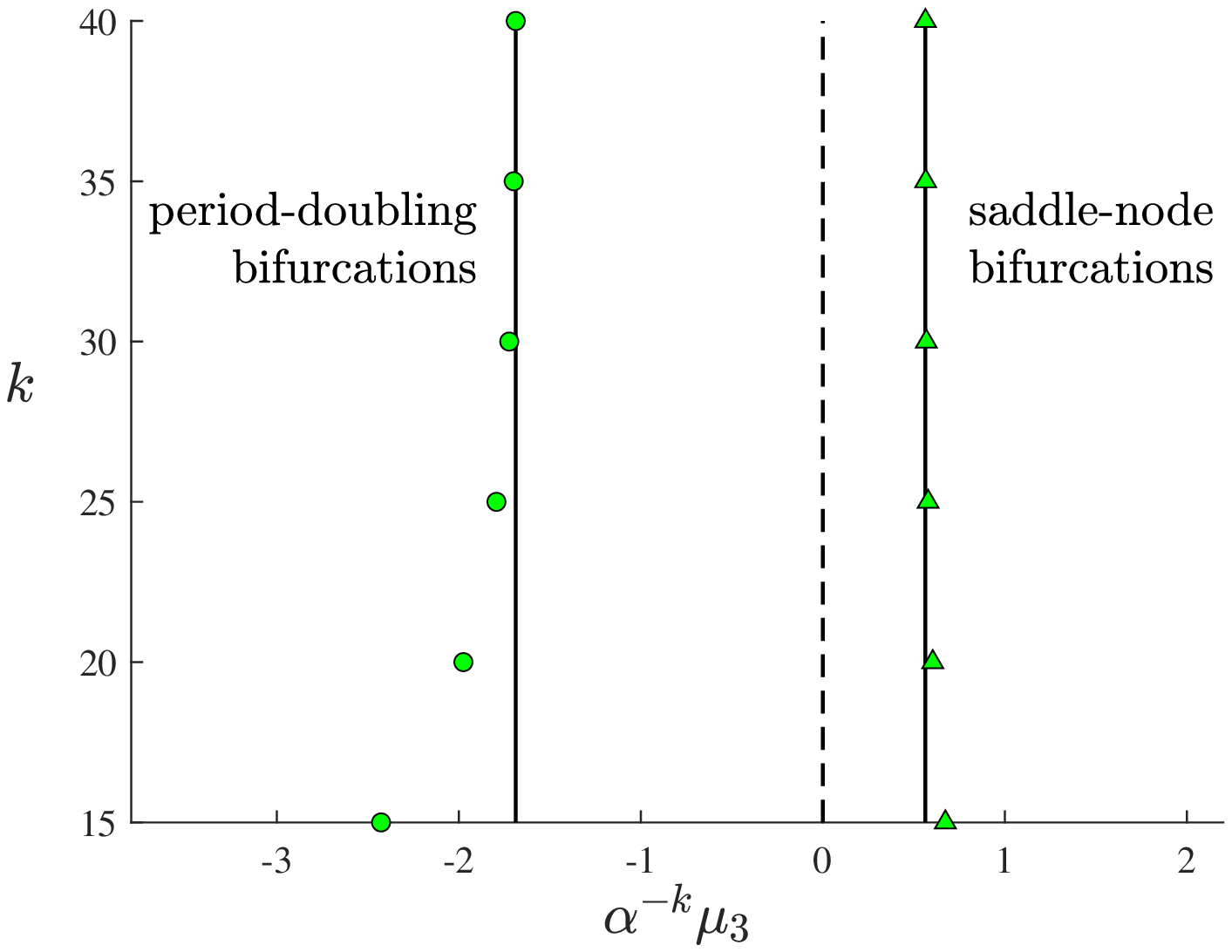}}
\put(0,5.1){\small \bf a)}
\put(7.2,5.1){\small \bf b)}
\end{picture}
\caption[Bifurcation diagram with variation of the parameter $\mu_{3}$]{Panel (a) is a numerically computed bifurcation diagram of \eqref{eq:fEx} with \eqref{eq:toyParameters} and $\mu_1 = \mu_2 = \mu_4 = 0$.  Panel (b) shows convergence to the leading-order terms of \eqref{eq:CO3SN} and \eqref{eq:CO3PD}.
\label{fig:C}
} 
\end{center}
\end{figure}
%%%%%%%%%%%%%%%%%%%%%%%%%%%%%%%%%%%%%%%%%%%%%%%%%%%%%%%%%%%%%%%%%%%%%%%%%%%%%%%%

Next we fix $\mu_1 = \mu_2 = \mu_4 = 0$ and vary the value of $\mu_3$ which breaks the global resonance condition.
Here ${\bf n}_{\rm tang}^{\sf T} {\bf v} = 0$ and ${\bf n}_{\rm eig}^{\sf T} {\bf v} = 0$ so Theorem \eqref{th:main} does not apply.  But by performing calculations analogous to those given above in the proof of Theorem \ref{th:main} directly to the map \eqref{eq:fEx}, 
we obtain the following expressions for the saddle-node and period-doubling bifurcation values
\begin{equation}
    \epsilon_{\rm SN} = \frac{(1-c_{2,0})^2}{4d_{5,0}} \alpha^{k} + \mathcal{O}(|\alpha|^{2k}),
    \label{eq:CO3SN}
\end{equation}

\begin{equation}
    \epsilon_{\rm PD} = -\frac{3(1-c_{2,0})^2}{4d_{5,0}} \alpha^{k}+ \mathcal{O}(|\alpha|^{2k}).
    \label{eq:CO3PD}
\end{equation}
Fig.~\ref{fig:C} shows that the numerically computed bifurcations do indeed appear to be converging to the leading-order components of \eqref{eq:CO3SN} and \eqref{eq:CO3PD}.
Notice the bifurcation values are asymptotically proportional to $\alpha^k$
(a slightly slower rate than that in Fig.~\ref{fig:B}).

%%%%%%%%%%%%%%%%%%%%%%%%%%%%%%%%%%%%%%%%%%%%%%%%%%%%%%%%%%%%%%%%%%%%%%%%%%%%%%%%
\begin{figure}[!htbp]
\begin{center}
\setlength{\unitlength}{1cm}
\begin{picture}(14,5.1)
\put(0,0){\includegraphics[width=6.8cm]{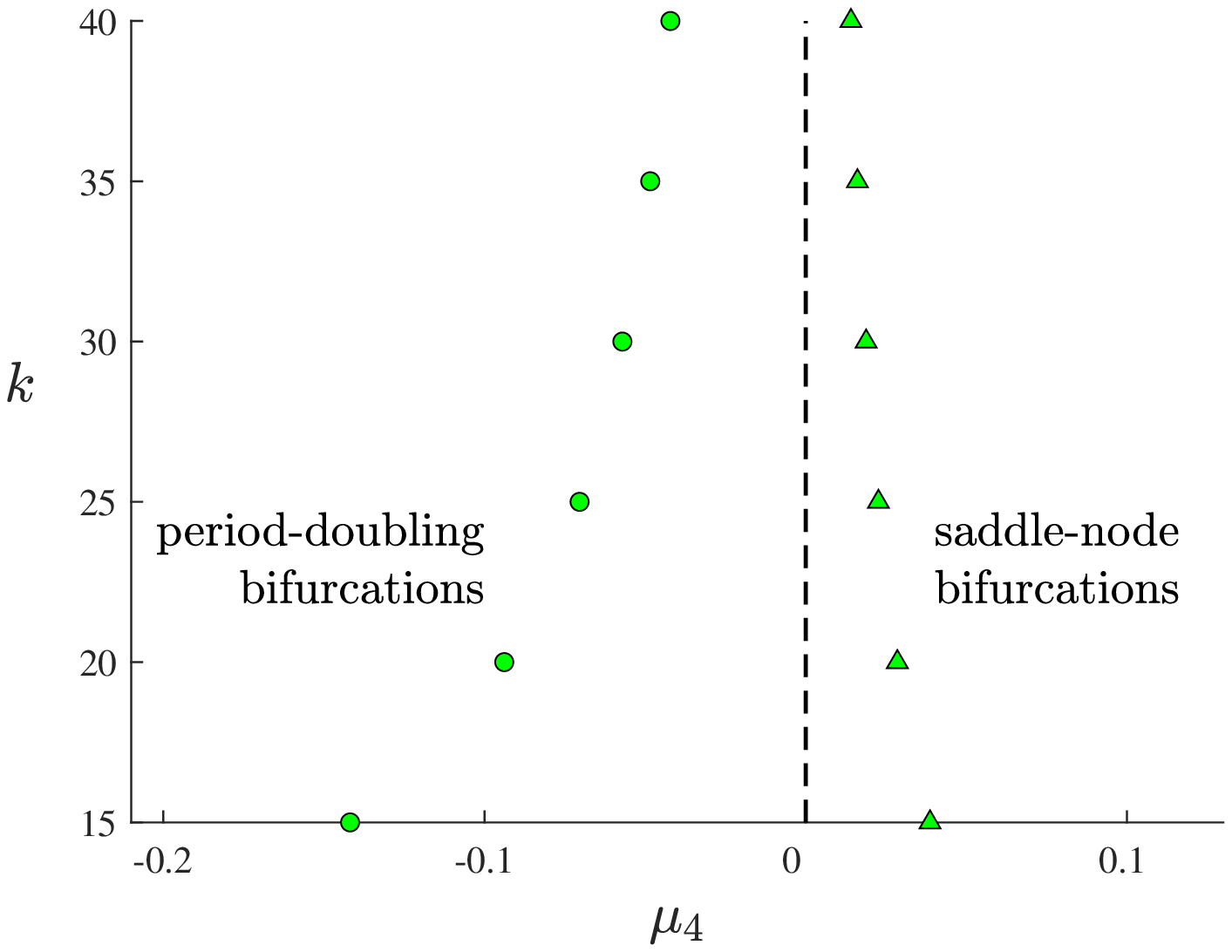}}
\put(7.2,0){\includegraphics[width=6.8cm]{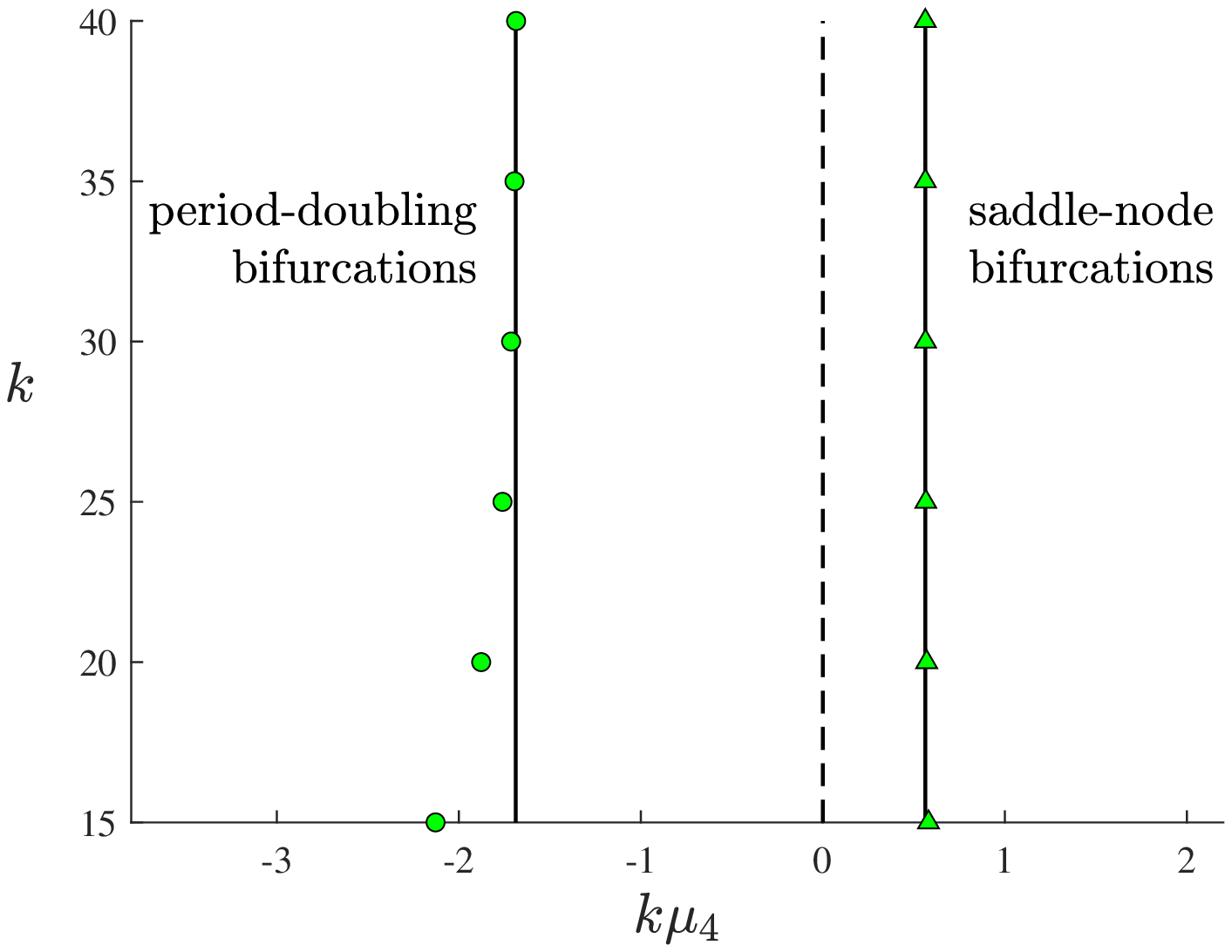}}
\put(0,5.1){\small \bf a)}
\put(7.2,5.1){\small \bf b)}
\end{picture}
\caption[Bifurcation diagram with variation of the parameter $\mu_{4}$]{Panel (a) is a numerically computed bifurcation diagram of \eqref{eq:fEx} with \eqref{eq:toyParameters} and $\mu_1 = \mu_2 = \mu_3 = 0$.  Panel (b) shows convergence to the leading-order terms of \eqref{eq:CO4SN} and \eqref{eq:CO4PD}.
\label{fig:D}
} 
\end{center}
\end{figure}
%%%%%%%%%%%%%%%%%%%%%%%%%%%%%%%%%%%%%%%%%%%%%%%%%%%%%%%%%%%%%%%%%%%%%%%%%%%%%%%%

Finally we fix $\mu_1 = \mu_2 = \mu_3 = 0$ and vary the value of $\mu_4$
which breaks the condition on the resonance terms in $T_0$.
As in the previous case Theorem \ref{th:main} does not apply.  By again calculating the bifurcations as above we obtain
\begin{equation}
    \epsilon_{{\rm SN}} = \frac{(1-c_{2,0})^2}{4d_{5,0}k} + \mathcal{O}\left(\frac{1}{k^2}\right),
    \label{eq:CO4SN}
\end{equation}
\begin{equation}
    \epsilon_{{\rm PD}} = -\frac{3(1-c_{2,0})^2}{4d_{5,0}k}+ \mathcal{O}\left(\frac{1}{k^2}\right),
    \label{eq:CO4PD}
\end{equation}
and these agree with the numerically computed bifurcation values as shown in Fig.~\ref{fig:D}. The bifurcation values are asymptotically proportional to $\frac{1}{k}$ which is substantially slower than in the previous three cases.
\chapter[Conclusions and open problems]{Conclusions and open problems \label{cha:OpenProblems}}
In Chapter \ref{cha:ExampleInfinte} we illustrated the phenomenon of globally resonant homoclinic tangencies (GRHT) in an abstract family of $C^{1}$ smooth maps. Next we explore the open problem of finding instances of GRHT in previously studied prototypical smooth maps that have enough parameters to fulfill the codimension-three and codimension-four criteria. We sketch a process which could be carried out to detect the high-codimension points. It is difficult to carry out the process as we need to be extremely close to the points to observe the phenomenon as we have seen in Chapter \ref{cha:UnfGRHT}. We showcase our trials in the \blue{GHM} to observe the GRHT phenomenon. We also discuss the challenges we faced while unfolding the codimension-four scenario in general.  

\section{GRHT in prototypical planar maps}
\label{sec:GHManalyse}
After exploring the phenomenon of GRHT in the abstract piecewise-smooth $C^{1}$ map above, we try to find more instances in prototypical maps. The \blue{GHM} is such an example of a prototypical map and is a generalisation of the well known H{\'e}non map. \blue{We have observed homoclinic tangles in the \blue{GHM} in Chapter \ref{cha:intro}.} The \blue{GHM} contains four parameters. So we expect a curve of solutions in the four-dimensional space of parameters for the codimension-three scenario and a point for the codimension-four scenario. 

\begin{figure}[t!]
\begin{center}
\includegraphics[width=0.8\textwidth]{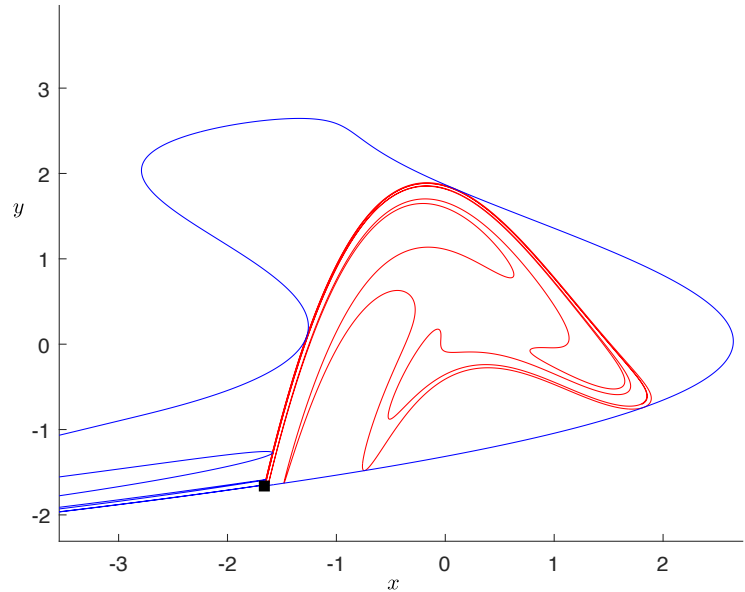}
\caption[Homoclinic tangency in the generalised H\'enon map.]{ Homoclinic tangency in the \blue{GHM}. The unstable manifold (red) and the stable manifold (blue) of a neutral saddle fixed point (marked by a black square dot) develop a tangential intersection at the parameter values $R = 0.1170465574,\,S = 0.3,\,\beta = 0.8055,\,\alpha = 0.8145$.}
\label{GHM_tangency_fig}
\end{center}
\end{figure}

\begin{figure}[t!]
\begin{center}
\includegraphics[width=0.8\textwidth]{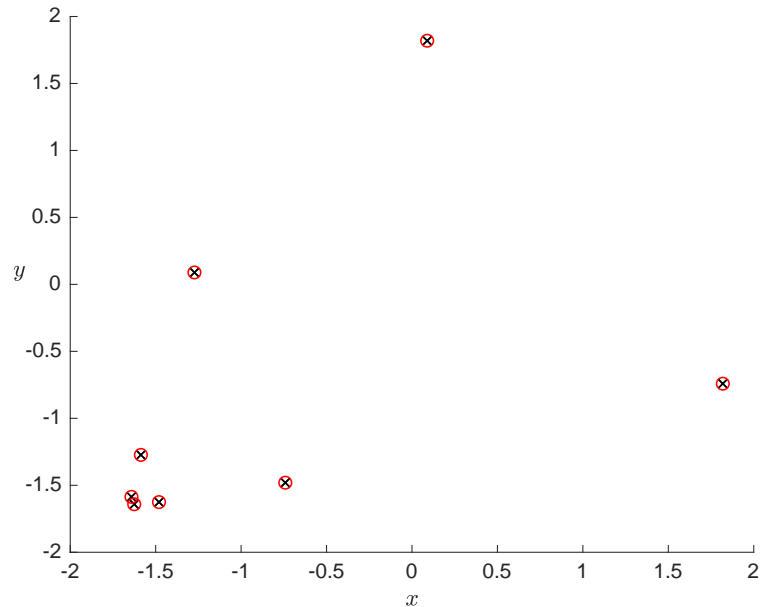}
\caption[A stable single-round period-eight solution for the generalised H\'enon map.]{ A stable period-eight solution for the \blue{GHM} at the parameter values $R = 0.1170465574, \, S = 0.3,\, \beta = 0.8055,\, \alpha = 0.8145$.}
\label{period8stab_fig}
\end{center}
\end{figure}

The explicit expression of the \blue{GHM} is below:
\begin{align*}
x' &= y,\\
y' &= \alpha - \beta x - y^2 + Rxy + Sy^3,   
\end{align*}
where $S, R, \alpha, \beta$ are parameters of the map. The fixed points of the GHM are given by solving the simultaneous equations 
\begin{equation}
    \begin{aligned}
    x &= y,\\
    y &= \alpha - \beta x - y^2 + Rxy + Sy^3.
    \end{aligned}
\end{equation}
Thus fixed points correspond to solutions of the cubic
\begin{equation}
Sy^3 + y^2 (R-1) + y(-\beta - 1) + \alpha = 0.
\label{eq:CubicGHM}
\end{equation}
The inverse of the \blue{GHM} is 
\begin{equation*}
f^{-1}\begin{bmatrix} x\\ y\\ \end{bmatrix} = \begin{bmatrix} \frac{y' - \alpha + x'^2 - Sx'^3}{Rx' - \beta}\\ x'\\ \end{bmatrix} .
\end{equation*}
We observe that the inverse of the \blue{GHM} exists if $x' \neq \frac{\beta}{R}$ which is equivalent to $y \neq \frac{\beta}{R}$. 
The Jacobian matrix of the \blue{GHM} is 
\begin{equation*}
J =
\begin{bmatrix}
0 & 1\\
(-\beta + Ry) & (-2y + Rx + 3Sy^2)\\
\end{bmatrix}.
\end{equation*}
In order to have a neutral saddle (i.e.~a fixed point whose associated eigenvalues multiply to $\pm 1$), we need the determinant of the Jacobian matrix to be $\pm 1$, that is $\beta - Ry = \pm 1$.  For this reason we first substitute $y = \frac{(\beta - 1)}{R}$ into the cubic equation \eqref{eq:CubicGHM} satisfied by the fixed points, so we have
$$S\left( \frac{\beta - 1}{R} \right)^3 + (R-1) \left(\frac{(\beta - 1)}{R}\right)^2 + (-\beta - 1)\left(\frac{\beta - 1}{R}\right) + \alpha = 0.$$
Solving for $\alpha$ gives
\begin{align}
\alpha &= \frac{-S(\beta - 1)^3 + (\beta - 1)R(\beta - 1 + 2R)}{R^3}.
\label{eq:alphaNeutral}
\end{align}
So in order to obtain a neutral saddle  with $\lambda \sigma  =1$, we need to use this value of $\alpha$.
Similarly, for the orientation-reversing case $\lambda\sigma=-1$, we need 
\begin{align*}
\alpha &= \frac{-S(\beta + 1)^3 + R(\beta + 1)^2}{R^3}.
\end{align*}
The \blue{GHM} exhibits homoclinic tangencies. Fig. \ref{GHM_tangency_fig} shows a tangential intersection between the stable and unstable manifolds of a neutral saddle fixed point. We vary the parameters until we observe a transversal intersection. The bisection method is then performed to obtain the tangency accurately. In an example, choosing initial conditions very close to the tangential points reveals a stable single-round periodic solution of period-$8$ as shown in Fig. \ref{period8stab_fig}.  
\subsection{Continuation of a curve of homoclinic tangencies}
\begin{figure}[!htbp]
\begin{center}
\includegraphics[width=0.8\textwidth]{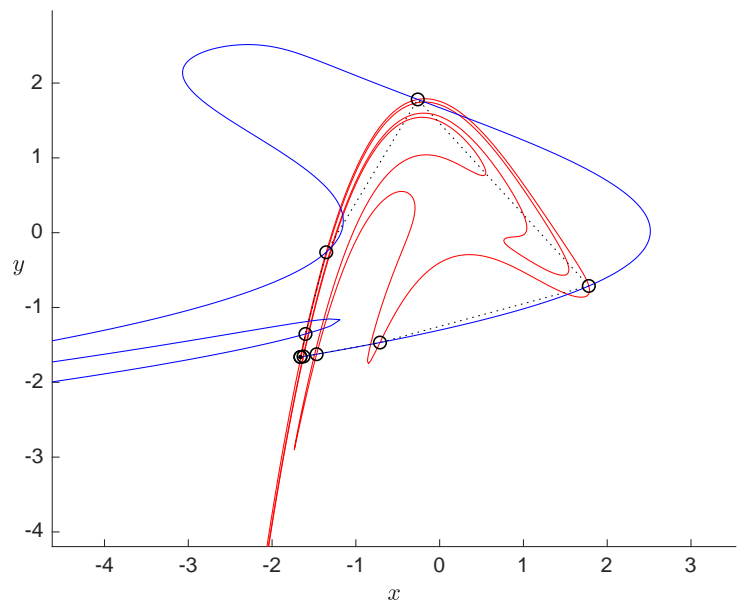}
\caption[Transversal intersection points in generalised H\'enon map using {\sc Matcontm}]{Computing Transversal intersections of the stable and unstable manifolds of the \blue{GHM} at the parameter values $R =0.1170465574, S =0.28, \beta = 0.8055$ and with the value of $\alpha$ given by \eqref{eq:alphaNeutral}.  For the purposes of numerical continuation some transversal intersections are marked with black circles.}
\label{fig:IntnTransvIntnGHMmm}
\end{center}
\end{figure}

We proceed to numerically continue the homoclinic tangency as described in \cite{KuMe19} using software {\sc matcontm}.
We begin with the parameters set as $R =0.1170465574, S =0.28, \beta = 0.8055$, at which the stable and unstable manifolds of saddle fixed point $(-1.6617,-1.6617)$ develop transversal intersections as shown in Fig.~\ref{fig:IntnTransvIntnGHMmm}.

We next compute some intersection points of the stable and unstable manifolds. In Fig.~\ref{fig:IntnTransvIntnGHMmm}, these are marked with black circles. We continue this homoclinic connection with respect to a single parameter (here, $S$) to obtain a tangential intersection between the manifolds.  Once they are detected they are marked by so-called limit points ${\rm LP_{HO}}$. \blue{We note that ${\rm LP_{HO}}$ are just the starting tangential points of the homoclinic tangency continuation in two parameters}. We were able to use {\sc matcontm} to find such limit points at $S =0.30007909$ and $S = 0.12398895$.

\begin{figure}[!htbp]
 %\hspace{0.8cm}
	\begin{tabular}{c c c}
		\includegraphics[width=0.65\textwidth]{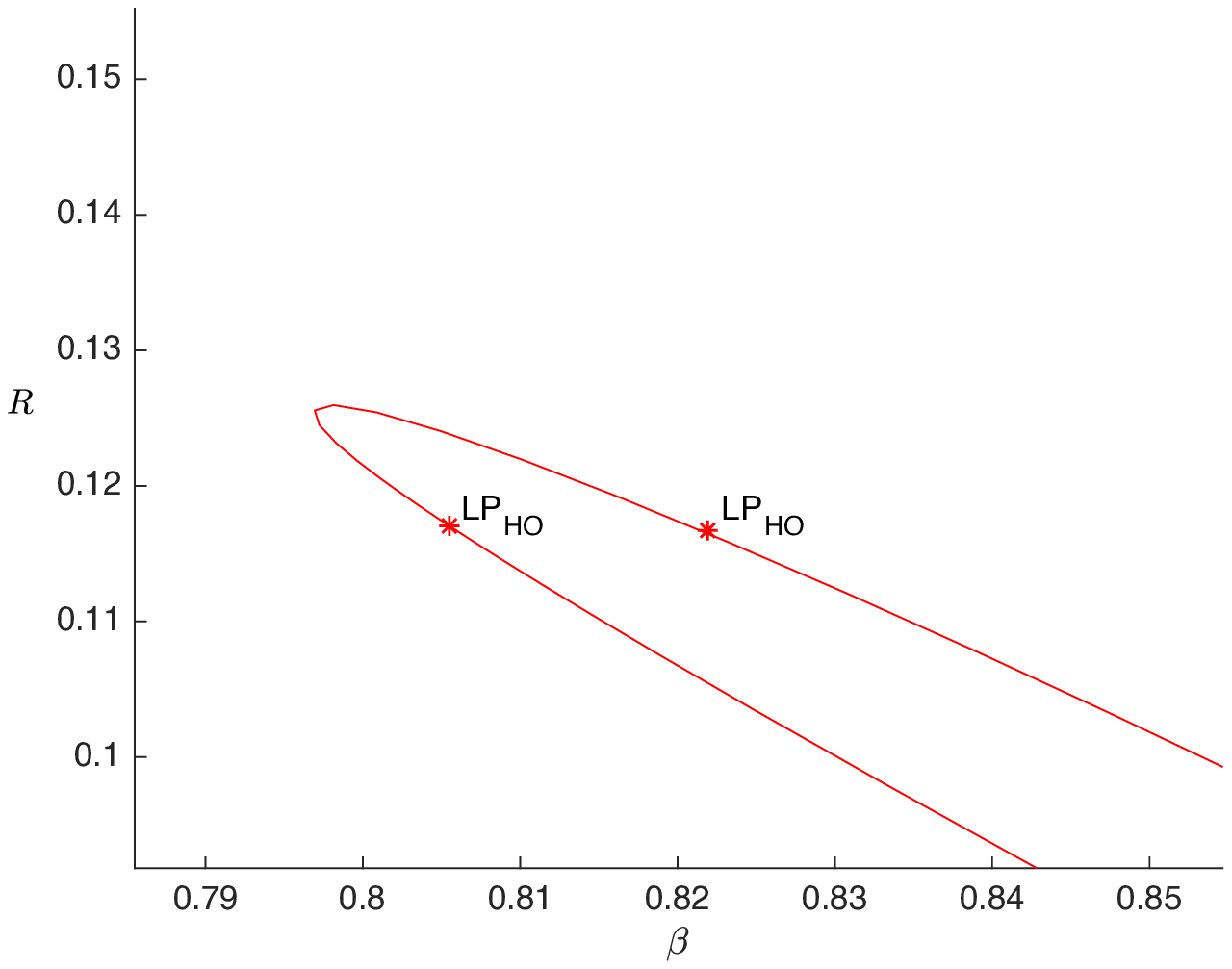}&\\
		\includegraphics[width=0.65\textwidth]{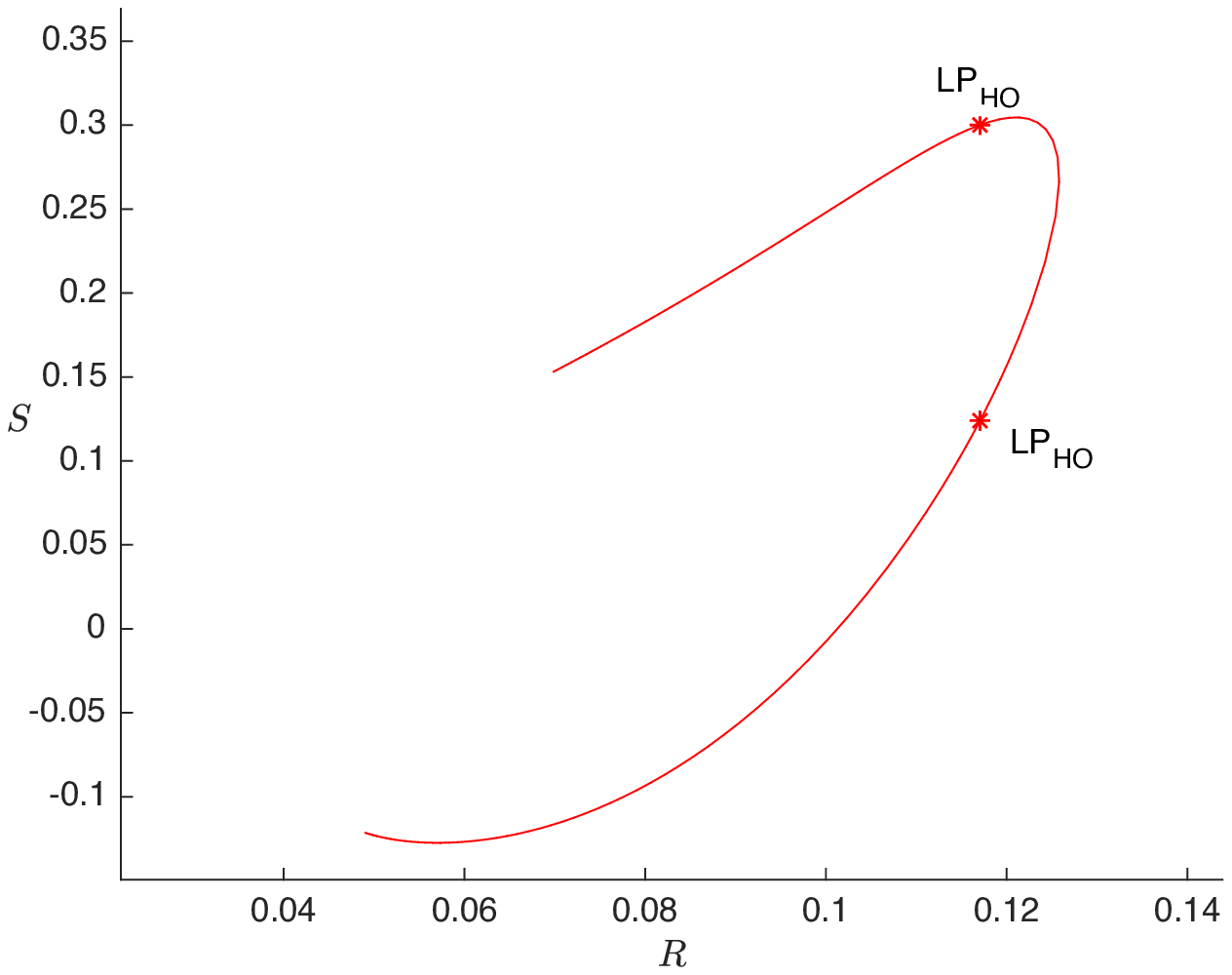}&\\
	\includegraphics[width=0.65\textwidth]{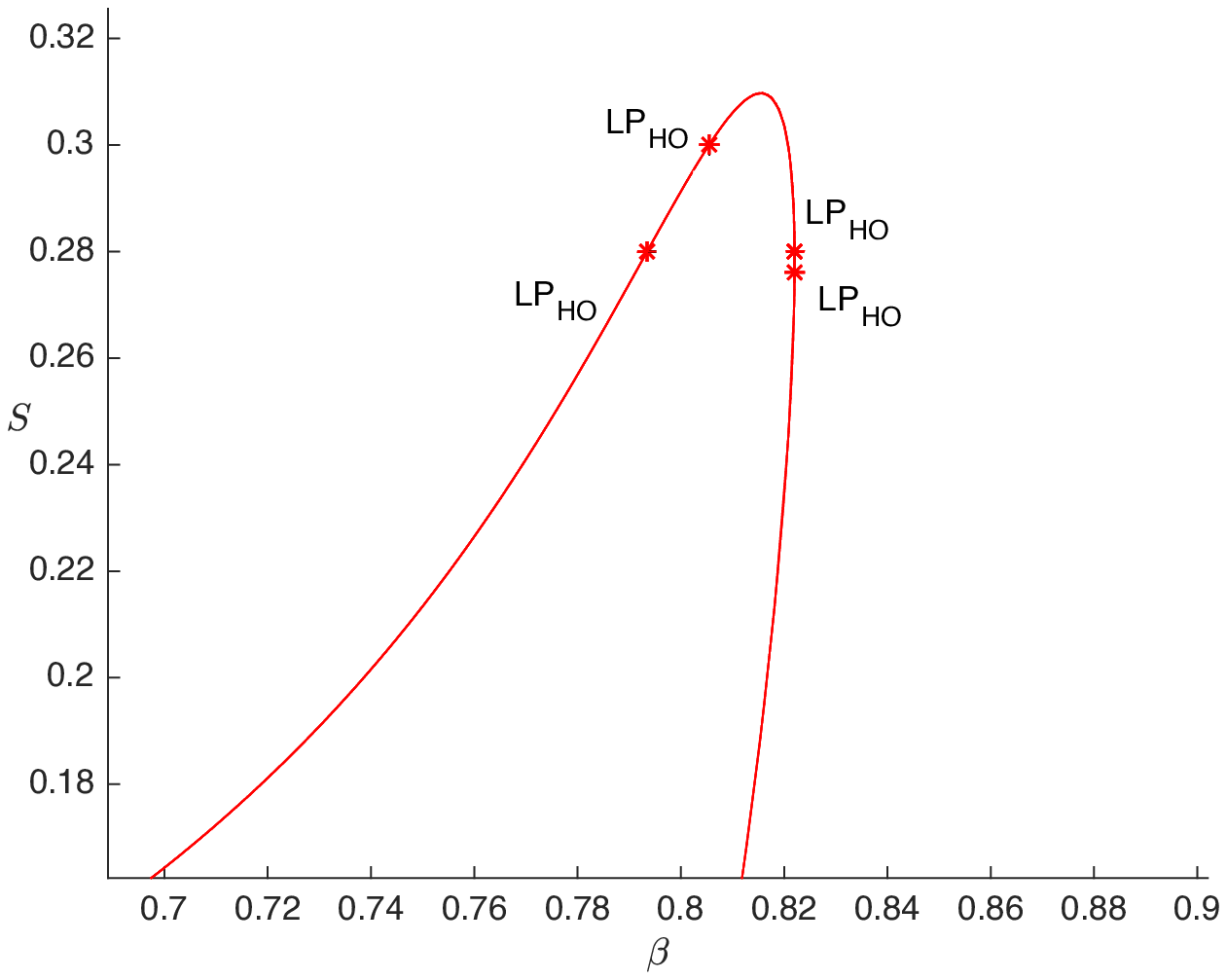}\\
	%	\put(8.1,9.1){\small \bf a)}
%\put(2.1,3.1){\small \bf b)}
	\end{tabular}
	\caption[Curves of homoclinic tangencies in two parameter space]{Curve of homoclinic tangencies continued in two-dimensional slices of parameter space.}
	\label{fig:TangencyContinuation}
\end{figure}

By numerically computing the stable and unstable manifolds at these values we were able to independently verify that these are indeed points of homoclinic tangencies. We then used {\sc matcontm} to numerically continue the homoclinic tangencies. Two parameters are needed to continue homoclinic tangency. 
Figure \ref{fig:TangencyContinuation} shows the continuation of homoclinic tangencies when it is continued with parameters a) $b$ and $S$, b) $R$ and $S$, and c) $b$ and $R$. In each plot we can see the curves of homoclinic tangencies at which each point on the curve gives a tangential intersection between the stable and unstable manifolds.  

We need to evaluate the global resonance condition along the curve of tangencies and ascertain when we have reached the codimension-three or four points. Checking the global resonance conditions and the inequalities involved in the sufficient conditions is a challenging task that needs to be addressed. 

\subsection{Sketch of a method to compute the conditions for GRHT}
\begin{figure}[!htbp]
\begin{center}
\ifx\du\undefined
  \newlength{\du}
\fi
\setlength{\du}{15\unitlength}
\begin{tikzpicture}[thick,scale=0.6, every node/.style={scale=0.7}]
\pgftransformxscale{1.000000}
\pgftransformyscale{-1.000000}
\definecolor{dialinecolor}{rgb}{0.000000, 0.000000, 0.000000}
\pgfsetstrokecolor{dialinecolor}
\definecolor{dialinecolor}{rgb}{1.000000, 1.000000, 1.000000}
\pgfsetfillcolor{dialinecolor}
\pgfsetlinewidth{0.0500000\du}
\pgfsetdash{}{0pt}
\pgfsetdash{}{0pt}
\pgfsetbuttcap
{
\definecolor{dialinecolor}{rgb}{0.000000, 0.000000, 0.000000}
\pgfsetfillcolor{dialinecolor}
% was here!!!
\definecolor{dialinecolor}{rgb}{0.000000, 0.000000, 0.000000}
\pgfsetstrokecolor{dialinecolor}
\draw (49.750000\du,28.100000\du)--(49.900000\du,28.200000\du);
}
\definecolor{dialinecolor}{rgb}{0.000000, 0.000000, 0.000000}
\pgfsetstrokecolor{dialinecolor}
\draw (49.750000\du,28.100000\du)--(49.900000\du,28.200000\du);
\pgfsetlinewidth{0.0500000\du}
\pgfsetdash{}{0pt}
\pgfsetmiterjoin
\pgfsetbuttcap
\definecolor{dialinecolor}{rgb}{0.000000, 0.000000, 0.000000}
\pgfsetfillcolor{dialinecolor}
\pgfpathmoveto{\pgfpoint{49.900000\du}{28.200000\du}}
\pgfpathcurveto{\pgfpoint{49.830662\du}{28.304006\du}}{\pgfpoint{49.657319\du}{28.338675\du}}{\pgfpoint{49.553312\du}{28.269338\du}}
\pgfpathcurveto{\pgfpoint{49.449306\du}{28.200000\du}}{\pgfpoint{49.414637\du}{28.026656\du}}{\pgfpoint{49.483975\du}{27.922650\du}}
\pgfpathcurveto{\pgfpoint{49.553312\du}{27.818644\du}}{\pgfpoint{49.726656\du}{27.783975\du}}{\pgfpoint{49.830662\du}{27.853312\du}}
\pgfpathcurveto{\pgfpoint{49.934669\du}{27.922650\du}}{\pgfpoint{49.969338\du}{28.095994\du}}{\pgfpoint{49.900000\du}{28.200000\du}}
\pgfusepath{fill}
\definecolor{dialinecolor}{rgb}{0.000000, 0.000000, 0.000000}
\pgfsetstrokecolor{dialinecolor}
\pgfpathmoveto{\pgfpoint{49.900000\du}{28.200000\du}}
\pgfpathcurveto{\pgfpoint{49.830662\du}{28.304006\du}}{\pgfpoint{49.657319\du}{28.338675\du}}{\pgfpoint{49.553312\du}{28.269338\du}}
\pgfpathcurveto{\pgfpoint{49.449306\du}{28.200000\du}}{\pgfpoint{49.414637\du}{28.026656\du}}{\pgfpoint{49.483975\du}{27.922650\du}}
\pgfpathcurveto{\pgfpoint{49.553312\du}{27.818644\du}}{\pgfpoint{49.726656\du}{27.783975\du}}{\pgfpoint{49.830662\du}{27.853312\du}}
\pgfpathcurveto{\pgfpoint{49.934669\du}{27.922650\du}}{\pgfpoint{49.969338\du}{28.095994\du}}{\pgfpoint{49.900000\du}{28.200000\du}}
\pgfusepath{stroke}
\pgfsetlinewidth{0.0500000\du}
\pgfsetdash{}{0pt}
\pgfsetdash{}{0pt}
\pgfsetbuttcap
{
\definecolor{dialinecolor}{rgb}{0.498039, 0.498039, 0.498039}
\pgfsetfillcolor{dialinecolor}
% was here!!!
\pgfsetarrowsend{stealth}
\definecolor{dialinecolor}{rgb}{0.498039, 0.498039, 0.498039}
\pgfsetstrokecolor{dialinecolor}
\draw (49.950000\du,27.450000\du)--(51.200000\du,24.200000\du);
}
\pgfsetlinewidth{0.0500000\du}
\pgfsetdash{}{0pt}
\pgfsetdash{}{0pt}
\pgfsetbuttcap
{
\definecolor{dialinecolor}{rgb}{0.000000, 0.000000, 0.000000}
\pgfsetfillcolor{dialinecolor}
% was here!!!
\definecolor{dialinecolor}{rgb}{0.000000, 0.000000, 0.000000}
\pgfsetstrokecolor{dialinecolor}
\draw (51.400000\du,23.550000\du)--(51.500000\du,23.850000\du);
}
\definecolor{dialinecolor}{rgb}{0.000000, 0.000000, 0.000000}
\pgfsetstrokecolor{dialinecolor}
\draw (51.400000\du,23.550000\du)--(51.386607\du,23.509822\du);
\pgfsetdash{}{0pt}
\pgfsetmiterjoin
\pgfsetbuttcap
\definecolor{dialinecolor}{rgb}{0.000000, 0.000000, 0.000000}
\pgfsetfillcolor{dialinecolor}
\fill (51.500000\du,23.850000\du)--(51.183772\du,23.691886\du)--(51.341886\du,23.375658\du)--(51.658114\du,23.533772\du)--cycle;
\pgfsetlinewidth{0.0500000\du}
\pgfsetdash{}{0pt}
\pgfsetmiterjoin
\pgfsetbuttcap
\definecolor{dialinecolor}{rgb}{0.000000, 0.000000, 0.000000}
\pgfsetstrokecolor{dialinecolor}
\draw (51.500000\du,23.850000\du)--(51.183772\du,23.691886\du)--(51.341886\du,23.375658\du)--(51.658114\du,23.533772\du)--cycle;
\pgfsetlinewidth{0.0500000\du}
\pgfsetdash{}{0pt}
\pgfsetdash{}{0pt}
\pgfsetbuttcap
{
\definecolor{dialinecolor}{rgb}{0.000000, 0.000000, 0.000000}
\pgfsetfillcolor{dialinecolor}
% was here!!!
\definecolor{dialinecolor}{rgb}{0.000000, 0.000000, 0.000000}
\pgfsetstrokecolor{dialinecolor}
\draw (41.750000\du,28.050000\du)--(42.050000\du,28.250000\du);
}
\definecolor{dialinecolor}{rgb}{0.000000, 0.000000, 0.000000}
\pgfsetstrokecolor{dialinecolor}
\draw (41.750000\du,28.050000\du)--(42.050000\du,28.250000\du);
\pgfsetlinewidth{0.0500000\du}
\pgfsetdash{}{0pt}
\pgfsetmiterjoin
\pgfsetbuttcap
\definecolor{dialinecolor}{rgb}{0.000000, 0.000000, 0.000000}
\pgfsetfillcolor{dialinecolor}
\pgfpathmoveto{\pgfpoint{42.050000\du}{28.250000\du}}
\pgfpathcurveto{\pgfpoint{41.980662\du}{28.354006\du}}{\pgfpoint{41.807319\du}{28.388675\du}}{\pgfpoint{41.703312\du}{28.319338\du}}
\pgfpathcurveto{\pgfpoint{41.599306\du}{28.250000\du}}{\pgfpoint{41.564637\du}{28.076656\du}}{\pgfpoint{41.633975\du}{27.972650\du}}
\pgfpathcurveto{\pgfpoint{41.703312\du}{27.868644\du}}{\pgfpoint{41.876656\du}{27.833975\du}}{\pgfpoint{41.980662\du}{27.903312\du}}
\pgfpathcurveto{\pgfpoint{42.084669\du}{27.972650\du}}{\pgfpoint{42.119338\du}{28.145994\du}}{\pgfpoint{42.050000\du}{28.250000\du}}
\pgfusepath{fill}
\definecolor{dialinecolor}{rgb}{0.000000, 0.000000, 0.000000}
\pgfsetstrokecolor{dialinecolor}
\pgfpathmoveto{\pgfpoint{42.050000\du}{28.250000\du}}
\pgfpathcurveto{\pgfpoint{41.980662\du}{28.354006\du}}{\pgfpoint{41.807319\du}{28.388675\du}}{\pgfpoint{41.703312\du}{28.319338\du}}
\pgfpathcurveto{\pgfpoint{41.599306\du}{28.250000\du}}{\pgfpoint{41.564637\du}{28.076656\du}}{\pgfpoint{41.633975\du}{27.972650\du}}
\pgfpathcurveto{\pgfpoint{41.703312\du}{27.868644\du}}{\pgfpoint{41.876656\du}{27.833975\du}}{\pgfpoint{41.980662\du}{27.903312\du}}
\pgfpathcurveto{\pgfpoint{42.084669\du}{27.972650\du}}{\pgfpoint{42.119338\du}{28.145994\du}}{\pgfpoint{42.050000\du}{28.250000\du}}
\pgfusepath{stroke}
\pgfsetlinewidth{0.0500000\du}
\pgfsetdash{}{0pt}
\pgfsetdash{}{0pt}
\pgfsetbuttcap
{
\definecolor{dialinecolor}{rgb}{0.498039, 0.498039, 0.498039}
\pgfsetfillcolor{dialinecolor}
% was here!!!
\pgfsetarrowsend{stealth}
\definecolor{dialinecolor}{rgb}{0.498039, 0.498039, 0.498039}
\pgfsetstrokecolor{dialinecolor}
\draw (42.050000\du,27.250000\du)--(43.233693\du,24.201825\du);
}
\pgfsetlinewidth{0.0500000\du}
\pgfsetdash{}{0pt}
\pgfsetdash{}{0pt}
\pgfsetbuttcap
{
\definecolor{dialinecolor}{rgb}{0.000000, 0.000000, 0.000000}
\pgfsetfillcolor{dialinecolor}
% was here!!!
\definecolor{dialinecolor}{rgb}{0.000000, 0.000000, 0.000000}
\pgfsetstrokecolor{dialinecolor}
\draw (43.500000\du,23.700000\du)--(43.150000\du,23.600000\du);
}
\definecolor{dialinecolor}{rgb}{0.000000, 0.000000, 0.000000}
\pgfsetstrokecolor{dialinecolor}
\draw (43.500000\du,23.700000\du)--(43.494782\du,23.698509\du);
\pgfsetdash{}{0pt}
\pgfsetmiterjoin
\pgfsetbuttcap
\definecolor{dialinecolor}{rgb}{0.000000, 0.000000, 0.000000}
\pgfsetfillcolor{dialinecolor}
\fill (43.150000\du,23.600000\du)--(43.459061\du,23.428299\du)--(43.630762\du,23.737361\du)--(43.321701\du,23.909061\du)--cycle;
\pgfsetlinewidth{0.0500000\du}
\pgfsetdash{}{0pt}
\pgfsetmiterjoin
\pgfsetbuttcap
\definecolor{dialinecolor}{rgb}{0.000000, 0.000000, 0.000000}
\pgfsetstrokecolor{dialinecolor}
\draw (43.150000\du,23.600000\du)--(43.459061\du,23.428299\du)--(43.630762\du,23.737361\du)--(43.321701\du,23.909061\du)--cycle;
\pgfsetlinewidth{0.0500000\du}
\pgfsetdash{{\pgflinewidth}{0.200000\du}}{0cm}
\pgfsetdash{{\pgflinewidth}{0.200000\du}}{0cm}
\pgfsetmiterjoin
\pgfsetbuttcap
{
\definecolor{dialinecolor}{rgb}{0.000000, 0.000000, 0.000000}
\pgfsetfillcolor{dialinecolor}
% was here!!!
\pgfsetarrowsend{stealth}
\definecolor{dialinecolor}{rgb}{0.000000, 0.000000, 0.000000}
\pgfsetstrokecolor{dialinecolor}
\pgfpathmoveto{\pgfpoint{48.850000\du}{26.450000\du}}
\pgfpathcurveto{\pgfpoint{47.350000\du}{25.100000\du}}{\pgfpoint{46.350000\du}{24.650000\du}}{\pgfpoint{43.850000\du}{26.250000\du}}
\pgfusepath{stroke}
}
\pgfsetlinewidth{0.0500000\du}
\pgfsetdash{}{0pt}
\pgfsetdash{}{0pt}
\pgfsetbuttcap
{
\definecolor{dialinecolor}{rgb}{0.000000, 0.000000, 0.000000}
\pgfsetfillcolor{dialinecolor}
% was here!!!
\definecolor{dialinecolor}{rgb}{0.000000, 0.000000, 0.000000}
\pgfsetstrokecolor{dialinecolor}
\draw (33.800000\du,27.950000\du)--(34.050000\du,28.050000\du);
}
\definecolor{dialinecolor}{rgb}{0.000000, 0.000000, 0.000000}
\pgfsetstrokecolor{dialinecolor}
\draw (33.800000\du,27.950000\du)--(34.050000\du,28.050000\du);
\pgfsetlinewidth{0.0500000\du}
\pgfsetdash{}{0pt}
\pgfsetmiterjoin
\pgfsetbuttcap
\definecolor{dialinecolor}{rgb}{0.000000, 0.000000, 0.000000}
\pgfsetfillcolor{dialinecolor}
\pgfpathmoveto{\pgfpoint{34.050000\du}{28.050000\du}}
\pgfpathcurveto{\pgfpoint{34.003576\du}{28.166060\du}}{\pgfpoint{33.841093\du}{28.235695\du}}{\pgfpoint{33.725033\du}{28.189272\du}}
\pgfpathcurveto{\pgfpoint{33.608974\du}{28.142848\du}}{\pgfpoint{33.539338\du}{27.980364\du}}{\pgfpoint{33.585762\du}{27.864305\du}}
\pgfpathcurveto{\pgfpoint{33.632185\du}{27.748245\du}}{\pgfpoint{33.794669\du}{27.678609\du}}{\pgfpoint{33.910728\du}{27.725033\du}}
\pgfpathcurveto{\pgfpoint{34.026788\du}{27.771457\du}}{\pgfpoint{34.096424\du}{27.933940\du}}{\pgfpoint{34.050000\du}{28.050000\du}}
\pgfusepath{fill}
\definecolor{dialinecolor}{rgb}{0.000000, 0.000000, 0.000000}
\pgfsetstrokecolor{dialinecolor}
\pgfpathmoveto{\pgfpoint{34.050000\du}{28.050000\du}}
\pgfpathcurveto{\pgfpoint{34.003576\du}{28.166060\du}}{\pgfpoint{33.841093\du}{28.235695\du}}{\pgfpoint{33.725033\du}{28.189272\du}}
\pgfpathcurveto{\pgfpoint{33.608974\du}{28.142848\du}}{\pgfpoint{33.539338\du}{27.980364\du}}{\pgfpoint{33.585762\du}{27.864305\du}}
\pgfpathcurveto{\pgfpoint{33.632185\du}{27.748245\du}}{\pgfpoint{33.794669\du}{27.678609\du}}{\pgfpoint{33.910728\du}{27.725033\du}}
\pgfpathcurveto{\pgfpoint{34.026788\du}{27.771457\du}}{\pgfpoint{34.096424\du}{27.933940\du}}{\pgfpoint{34.050000\du}{28.050000\du}}
\pgfusepath{stroke}
\pgfsetlinewidth{0.0500000\du}
\pgfsetdash{}{0pt}
\pgfsetdash{}{0pt}
\pgfsetbuttcap
{
\definecolor{dialinecolor}{rgb}{0.498039, 0.498039, 0.498039}
\pgfsetfillcolor{dialinecolor}
% was here!!!
\pgfsetarrowsend{stealth}
\definecolor{dialinecolor}{rgb}{0.498039, 0.498039, 0.498039}
\pgfsetstrokecolor{dialinecolor}
\draw (34.050000\du,27.350000\du)--(35.000000\du,24.250000\du);
}
\pgfsetlinewidth{0.0500000\du}
\pgfsetdash{}{0pt}
\pgfsetdash{}{0pt}
\pgfsetbuttcap
{
\definecolor{dialinecolor}{rgb}{0.000000, 0.000000, 0.000000}
\pgfsetfillcolor{dialinecolor}
% was here!!!
\definecolor{dialinecolor}{rgb}{0.000000, 0.000000, 0.000000}
\pgfsetstrokecolor{dialinecolor}
\draw (35.300000\du,23.650000\du)--(35.300000\du,23.550000\du);
}
\definecolor{dialinecolor}{rgb}{0.000000, 0.000000, 0.000000}
\pgfsetstrokecolor{dialinecolor}
\draw (35.300000\du,23.650000\du)--(35.300000\du,23.908579\du);
\pgfsetdash{}{0pt}
\pgfsetmiterjoin
\pgfsetbuttcap
\definecolor{dialinecolor}{rgb}{0.000000, 0.000000, 0.000000}
\pgfsetfillcolor{dialinecolor}
\fill (35.300000\du,23.550000\du)--(35.550000\du,23.800000\du)--(35.300000\du,24.050000\du)--(35.050000\du,23.800000\du)--cycle;
\pgfsetlinewidth{0.0500000\du}
\pgfsetdash{}{0pt}
\pgfsetmiterjoin
\pgfsetbuttcap
\definecolor{dialinecolor}{rgb}{0.000000, 0.000000, 0.000000}
\pgfsetstrokecolor{dialinecolor}
\draw (35.300000\du,23.550000\du)--(35.550000\du,23.800000\du)--(35.300000\du,24.050000\du)--(35.050000\du,23.800000\du)--cycle;
\pgfsetlinewidth{0.0500000\du}
\pgfsetdash{{\pgflinewidth}{0.200000\du}}{0cm}
\pgfsetdash{{\pgflinewidth}{0.200000\du}}{0cm}
\pgfsetmiterjoin
\pgfsetbuttcap
{
\definecolor{dialinecolor}{rgb}{0.000000, 0.000000, 0.000000}
\pgfsetfillcolor{dialinecolor}
% was here!!!
\pgfsetarrowsend{stealth}
\definecolor{dialinecolor}{rgb}{0.000000, 0.000000, 0.000000}
\pgfsetstrokecolor{dialinecolor}
\pgfpathmoveto{\pgfpoint{40.700000\du}{26.300000\du}}
\pgfpathcurveto{\pgfpoint{39.500000\du}{24.500000\du}}{\pgfpoint{37.150000\du}{25.050000\du}}{\pgfpoint{35.800000\du}{26.250000\du}}
\pgfusepath{stroke}
}
\pgfsetlinewidth{0.0500000\du}
\pgfsetdash{{\pgflinewidth}{0.200000\du}}{0cm}
\pgfsetdash{{\pgflinewidth}{0.200000\du}}{0cm}
\pgfsetmiterjoin
\pgfsetbuttcap
{
\definecolor{dialinecolor}{rgb}{0.000000, 0.000000, 0.000000}
\pgfsetfillcolor{dialinecolor}
% was here!!!
\definecolor{dialinecolor}{rgb}{0.000000, 0.000000, 0.000000}
\pgfsetstrokecolor{dialinecolor}
\pgfpathmoveto{\pgfpoint{61.100000\du}{19.650000\du}}
\pgfpathcurveto{\pgfpoint{64.600000\du}{13.850000\du}}{\pgfpoint{24.850000\du}{4.150000\du}}{\pgfpoint{33.600000\du}{27.200000\du}}
\pgfusepath{stroke}
}
\pgfsetlinewidth{0.0500000\du}
\pgfsetdash{}{0pt}
\pgfsetdash{}{0pt}
\pgfsetbuttcap
{
\definecolor{dialinecolor}{rgb}{0.000000, 0.000000, 0.000000}
\pgfsetfillcolor{dialinecolor}
% was here!!!
\pgfsetarrowsend{stealth}
\definecolor{dialinecolor}{rgb}{0.000000, 0.000000, 0.000000}
\pgfsetstrokecolor{dialinecolor}
\draw (44.077292\du,12.601890\du)--(44.709498\du,12.601890\du);
}
\pgfsetlinewidth{0.0500000\du}
\pgfsetdash{}{0pt}
\pgfsetdash{}{0pt}
\pgfsetbuttcap
{
\definecolor{dialinecolor}{rgb}{0.000000, 0.000000, 0.000000}
\pgfsetfillcolor{dialinecolor}
% was here!!!
\definecolor{dialinecolor}{rgb}{0.000000, 0.000000, 0.000000}
\pgfsetstrokecolor{dialinecolor}
\draw (63.950000\du,20.550000\du)--(64.150000\du,20.500000\du);
}
\definecolor{dialinecolor}{rgb}{0.000000, 0.000000, 0.000000}
\pgfsetstrokecolor{dialinecolor}
\draw (63.950000\du,20.550000\du)--(64.150000\du,20.500000\du);
\pgfsetlinewidth{0.0500000\du}
\pgfsetdash{}{0pt}
\pgfsetmiterjoin
\pgfsetbuttcap
\definecolor{dialinecolor}{rgb}{0.000000, 0.000000, 0.000000}
\pgfsetfillcolor{dialinecolor}
\pgfpathmoveto{\pgfpoint{64.150000\du}{20.500000\du}}
\pgfpathcurveto{\pgfpoint{64.180317\du}{20.621268\du}}{\pgfpoint{64.089366\du}{20.772853\du}}{\pgfpoint{63.968098\du}{20.803170\du}}
\pgfpathcurveto{\pgfpoint{63.846830\du}{20.833486\du}}{\pgfpoint{63.695246\du}{20.742536\du}}{\pgfpoint{63.664929\du}{20.621268\du}}
\pgfpathcurveto{\pgfpoint{63.634612\du}{20.500000\du}}{\pgfpoint{63.725563\du}{20.348415\du}}{\pgfpoint{63.846830\du}{20.318098\du}}
\pgfpathcurveto{\pgfpoint{63.968098\du}{20.287781\du}}{\pgfpoint{64.119683\du}{20.378732\du}}{\pgfpoint{64.150000\du}{20.500000\du}}
\pgfusepath{fill}
\definecolor{dialinecolor}{rgb}{0.000000, 0.000000, 0.000000}
\pgfsetstrokecolor{dialinecolor}
\pgfpathmoveto{\pgfpoint{64.150000\du}{20.500000\du}}
\pgfpathcurveto{\pgfpoint{64.180317\du}{20.621268\du}}{\pgfpoint{64.089366\du}{20.772853\du}}{\pgfpoint{63.968098\du}{20.803170\du}}
\pgfpathcurveto{\pgfpoint{63.846830\du}{20.833486\du}}{\pgfpoint{63.695246\du}{20.742536\du}}{\pgfpoint{63.664929\du}{20.621268\du}}
\pgfpathcurveto{\pgfpoint{63.634612\du}{20.500000\du}}{\pgfpoint{63.725563\du}{20.348415\du}}{\pgfpoint{63.846830\du}{20.318098\du}}
\pgfpathcurveto{\pgfpoint{63.968098\du}{20.287781\du}}{\pgfpoint{64.119683\du}{20.378732\du}}{\pgfpoint{64.150000\du}{20.500000\du}}
\pgfusepath{stroke}
\pgfsetlinewidth{0.0500000\du}
\pgfsetdash{}{0pt}
\pgfsetdash{}{0pt}
\pgfsetbuttcap
{
\definecolor{dialinecolor}{rgb}{0.498039, 0.498039, 0.498039}
\pgfsetfillcolor{dialinecolor}
% was here!!!
\pgfsetarrowsend{stealth}
\definecolor{dialinecolor}{rgb}{0.498039, 0.498039, 0.498039}
\pgfsetstrokecolor{dialinecolor}
\draw (63.250000\du,20.850000\du)--(60.387929\du,22.061261\du);
}
\pgfsetlinewidth{0.0500000\du}
\pgfsetdash{}{0pt}
\pgfsetdash{}{0pt}
\pgfsetbuttcap
{
\definecolor{dialinecolor}{rgb}{0.000000, 0.000000, 0.000000}
\pgfsetfillcolor{dialinecolor}
% was here!!!
\definecolor{dialinecolor}{rgb}{0.000000, 0.000000, 0.000000}
\pgfsetstrokecolor{dialinecolor}
\draw (59.800000\du,22.250000\du)--(59.950000\du,22.400000\du);
}
\definecolor{dialinecolor}{rgb}{0.000000, 0.000000, 0.000000}
\pgfsetstrokecolor{dialinecolor}
\draw (59.800000\du,22.250000\du)--(59.696447\du,22.146447\du);
\pgfsetdash{}{0pt}
\pgfsetmiterjoin
\pgfsetbuttcap
\definecolor{dialinecolor}{rgb}{0.000000, 0.000000, 0.000000}
\pgfsetfillcolor{dialinecolor}
\fill (59.950000\du,22.400000\du)--(59.596447\du,22.400000\du)--(59.596447\du,22.046447\du)--(59.950000\du,22.046447\du)--cycle;
\pgfsetlinewidth{0.0500000\du}
\pgfsetdash{}{0pt}
\pgfsetmiterjoin
\pgfsetbuttcap
\definecolor{dialinecolor}{rgb}{0.000000, 0.000000, 0.000000}
\pgfsetstrokecolor{dialinecolor}
\draw (59.950000\du,22.400000\du)--(59.596447\du,22.400000\du)--(59.596447\du,22.046447\du)--(59.950000\du,22.046447\du)--cycle;
% setfont left to latex
\definecolor{dialinecolor}{rgb}{0.000000, 0.000000, 0.000000}
\pgfsetstrokecolor{dialinecolor}
\node[anchor=west] at (47.700000\du,29.350000\du){$\mathbf{W_{2}}=(u_{2},v_{2})$};
% setfont left to latex
\definecolor{dialinecolor}{rgb}{0.000000, 0.000000, 0.000000}
\pgfsetstrokecolor{dialinecolor}
\node[anchor=west] at (51.177030\du,22.606760\du){$\mathbf{W_{2}} + \delta v_{s}$};
% setfont left to latex
\definecolor{dialinecolor}{rgb}{0.000000, 0.000000, 0.000000}
\pgfsetstrokecolor{dialinecolor}
\node[anchor=west] at (41.400000\du,29.200000\du){$f(\mathbf{W_{2}})$};
% setfont left to latex
\definecolor{dialinecolor}{rgb}{0.000000, 0.000000, 0.000000}
\pgfsetstrokecolor{dialinecolor}
\node[anchor=west] at (46.591890\du,26.054069\du){$f$};
% setfont left to latex
\definecolor{dialinecolor}{rgb}{0.000000, 0.000000, 0.000000}
\pgfsetstrokecolor{dialinecolor}
\node[anchor=west] at (42.850000\du,22.550000\du){$f(\mathbf{W_{2}}+\delta v_{s})$};
% setfont left to latex
\definecolor{dialinecolor}{rgb}{0.000000, 0.000000, 0.000000}
\pgfsetstrokecolor{dialinecolor}
\node[anchor=west] at (38.467570\du,25.952719\du){$f$};
% setfont left to latex
\definecolor{dialinecolor}{rgb}{0.000000, 0.000000, 0.000000}
\pgfsetstrokecolor{dialinecolor}
\node[anchor=west] at (33.350000\du,28.900000\du){$f^{2}(\mathbf{W_{2}})$};
% setfont left to latex
\definecolor{dialinecolor}{rgb}{0.000000, 0.000000, 0.000000}
\pgfsetstrokecolor{dialinecolor}
\node[anchor=west] at (34.650000\du,22.700000\du){$f^{2}(\mathbf{W_{2}}+\delta v_{s})$};
% setfont left to latex
\definecolor{dialinecolor}{rgb}{0.000000, 0.000000, 0.000000}
\pgfsetstrokecolor{dialinecolor}
\node[anchor=west] at (42.050000\du,11.650000\du){applying $f$ for $(n-2)$ times};
% setfont left to latex
\definecolor{dialinecolor}{rgb}{0.000000, 0.000000, 0.000000}
\pgfsetstrokecolor{dialinecolor}
\node[anchor=west] at (64.712160\du,20.494600\du){$\mathbf{W_{1}} = f^{n}(\mathbf{W_{2}})$};
% setfont left to latex
\definecolor{dialinecolor}{rgb}{0.000000, 0.000000, 0.000000}
\pgfsetstrokecolor{dialinecolor}
\node[anchor=west] at (59.627020\du,23.398650\du){$\mathbf{W} = \mathbf{W_{1}} + \delta v_{\text{diff}}$};
\end{tikzpicture}
 \caption{A sketch illustrating a method to compute  $d_{1} = \frac{\partial{T_{1}}}{\partial{x}}\Bigr|_{(0,y^*)}$.}
    \label{fig:f_method2}
\end{center}
\end{figure}
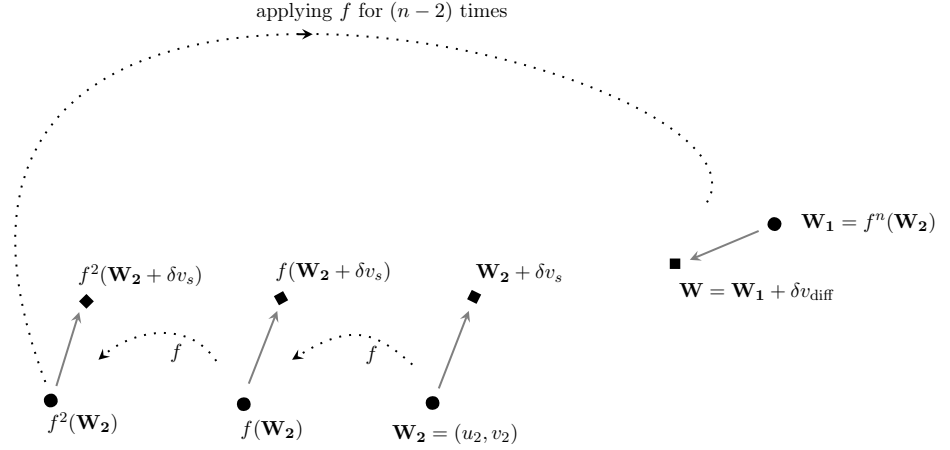

We sketch approaches taken to compute the global resonance condition $|d_{1}| = \frac{y^*}{x^*}$. Observe that  $ \left|\frac{\partial{T_{1}}}{\partial{x}}|_{(0,y^*)}\right| = \frac{y^*}{x^*}$. Here we explain how we calculate the value of $d_1$ at a tangency, see Fig.~\ref{fig:f_method2}.
Suppose $\mathbf{W_{1}} \, {\rm and} \, \mathbf{W_{2}}$ represent two homoclinic points of the homoclinic orbit. Let $\mathbf{W_{2}} = (u_{2},v_{2})$ map to the homoclinic point $\mathbf{W_{1}} = (u_{1},v_{1})$ under a number of iterations (say $n$) of the map $f$. We have $\mathbf{W_{2}} = (u_{2},v_{2})$ and we perturb it by $\delta$ in the direction of the stable eigenvector $v_{s}$. Then iterate both $\mathbf{W_{2}}$ and the perturbed point under the map $f$. So the perturbed point is at $$f(\mathbf{W_{2}} + \delta v_{s}) = f(\mathbf{W_{2}}) + \delta {\rm D}f(\mathbf{W_{2}}) v_{s} + \mathcal{O}(\delta^2).$$ Next we iterate them under $f$ again, then the new perturbed point becomes $$f(f(\mathbf{W_{2}} + \delta v_{s})) = f^{2}(\mathbf{W_{2}}) + \delta {\rm D}f(f(\mathbf{W_{2}})) {\rm D}f(\mathbf{W_{2}}) v_{s} + \mathcal{O}(\delta ^ 2).$$
We apply $f$ another $(n-2)$ times and $\mathbf{W_{1}} = f^n(\mathbf{W_{2}})$ and the perturbed point has the form $\mathbf{W} = \mathbf{W_{1}} + v_{\text{diff}}$ where $$v_{\text{diff}} = {\rm D}f(f^{(n-1)}(\mathbf{W_{2}})) \ldots {\rm D}f(f(\mathbf{W_{2}})) {\rm D}f(\mathbf{W_{2}}) v_{s}\, ,$$ to leading order. Following the geometric interpretation of the global resonance in Fig.~\ref{fig:third_condn} of Chapter \ref{cha:GRHT}, $d_{1}$ can then be written as
$$d_{1} = [0,1]M^{-1}v_{{\rm diff}},$$ where $M = \begin{bmatrix} v_{u} \quad v_{s}\end{bmatrix}$. 

The above method requires the manifolds to be aligned with the orthogonal coordinate axis. This can be achieved via  a series of nonlinear coordinate changes. Performing nonlinear coordinate changes numerically and automating it for efficient computation is a challenging task which remains for future work. 

%\subsection{General unfolding of the codimension-four scenario}
%Another open problem is to determine the general unfolding of the codimension-four phenomenon. This requires calculating three-dimensional surfaces of saddle-node and period-doubling bifurcation in a way that accommodates all the various scaling laws that we have observed. Doing this in general is a challenging task as we have to handle four parameters simultaneously.

\section{Conclusion}
\label{sec:conclusion}
In this thesis we have addressed a new phenomenon of the coexistence of infinitely many asymptotically stable single-round periodic solutions and the mechanism associated to it in discrete two-dimensional maps. The phenomenon is termed a globally resonant homoclinic tangency due to the behaviour required for the reinjection mechanism. We have identified both necessary and sufficient conditions behind the occurrence of the globally resonant homoclinic tangencies.  These were discussed in detail in Chapter \ref{cha:GRHT}. We have proved that the phenomenon is codimension-three in the case of orientation-reversing maps while it is codimension-four in the case of orientation-preserving case. The idea behind the proof was to bound the trace and determinant so that they lie inside the stability triangle for higher iterates of the map. 

To understand the stability region in parameter space, we then studied the unfolding scenario about such a tangency. In a one-parameter variation of the tangency, \blue{ it was found that there is a sequence of saddle-node and period-doubling bifurcations nearby the tangency.}  It was interesting to observe that the scaling laws of bifurcation values were different in different directions of parameter space. Generically the sequences scale like $|\lambda|^{2k}$, where $-1<\lambda<1$ is the stable eigenvalue of the fixed point. If we vary parameters without a linear change to the codimension-one condition for a homoclinic tangency the bifurcations scale as $\frac{|\lambda|^k}{k}$. If the perturbation is further degenerate, the scaling laws are found to be slower. Specifically they scale like $|\lambda|^k$ and $\frac{1}{k}$ for the example \eqref{eq:fEx}.  Preliminary studies found that the shape of the stability region is extremely complicated with a precise nature that heavily depends on the parameters. A detailed structure is yet to be determined.

We also note that the generic scaling law of $|\lambda|^{2 k}$ differs from that for families of piecewise-linear maps which instead admit a $|\lambda|^k$ scaling law. 

\subsection{Open problems}
%Refer subsection of chapter 1, \ref{subsec:PWL}.
It remains to identify infinite coexistence in previously studied prototypical maps.
The relatively high codimension possibly explains why this does not appear to have been done already. Yet the codimension is not so high that the phenomenon cannot be expected to have an important role in some applications.

\blue{We have encountered double-round periodic solutions in the orientation-preserving case. It would be interesting to understand the one or more bifurcations that lead to its formation. A more general study on multi-round periodic solutions in the context of the results of this thesis remains to be carried out.}

\blue{Above we have provided some description of how we have searched for globally resonant homoclinic tangencies in the generalised H\'{e}non map.  It remains to identify an effective automated numerical procedure that one could use to do this for any appropriate family of maps.}

Performing coordinate changes numerically to straighten the manifolds and computing the global resonance condition and conditions involving the coefficients of resonance terms of the local map remains a challenging task to address. Also the typical bifurcation structure that surrounds these codimension-three and four situations has yet to be determined in general. Unfolding the tangency by varying four parameters simultaneously in general is a challenging task.

Numerically computing the saddle-node and period doubling bifurcation points near the unfolding scenario using bifurcation software is a possible future direction that may help to understand the complicated stability region in the codimension-three or codimension-four scenarios better.

Finally it also may be of interest to generalise the results to incorporate cubic and higher order tangencies and extend the results to higher dimensional maps.
\blue{
\section{Broader significance}
\label{sec:significance}
This thesis reveals the mechanism behind extreme multistability or coexistence of infinitely many asymptotically stable single-round periodic solutions arising from homoclinic tangencies in two-dimensional maps. Multistability has been found in many real world systems and applications \citep{Fe08}. Multistability is undesirable in engineering systems because the system can switch to undesirable attractors in an unpredictale way depending on the initial state of the system \citep{HDC90}, but it is advantageous in the case of neuron systems, where coexistence of many attractors is responsible for multi component information processing \citep{TM11} and has far reaching implications for motor control, and decision making. The mechanism behind such infinite coexistence has been answered for piecewise linear maps \citep{DSS14}. Lots of attractors have been found in impacting systems governed by the Nordmark map in \cite{Nord20}. A single asymptotically stable periodic orbit could bifurcate into infinitely many asymptotically stable periodic orbits via grazing- sliding bifurcations \citep{Si17d} commonly found in impacting systems. 
}
% Appendices start here
%\appendix
%\include{appendix1}

% Bibliography starts here
% First: name of .bib file (e.g., refs)

\addcontentsline{toc}{chapter}{Bibliography}
\bibliography{refs}             
% Second: style in which printed. There is no official style, but I like plainnat
\bibliographystyle{plainnat}  
\end{document}